\documentclass[11pt,reqno]{article}
\usepackage{amsfonts}
\usepackage{graphicx}
\usepackage{imakeidx}
\makeindex

\begin{document}

\begin{center}
{\large\bf Spectral properties of Ruelle transfer operators for regular Gibbs measures  and decay of correlations for contact Anosov flows}
\end{center}

\begin{center}
{Luchezar Stoyanov\footnote{Department of Mathematics and Statistics,
University of Western Australia;
luchezar.stoyanov@uwa.edu.au}}
\end{center}


\begin{abstract}
In this work we study strong spectral properties of Ruelle  transfer operators related to 
a large family of Gibbs measures for contact Anosov flows. The ultimate aim is to establish exponential decay of
correlations for H\"older observables with respect to a very general class of Gibbs measures.
The approach invented in 1997 by Dolgopyat  \cite{D1} and further developed in \cite{St2} is substantially refined here, 
allowing to deal with much more general situations than before, although we still restrict ourselves to the uniformly
hyperbolic case. A rather general procedure is established which produces the desired 
estimates whenever the Gibbs measure admits a Pesin set with exponentially small tails, that is a Pesin set whose
preimages along the flow have measures decaying exponentially fast. We call such Gibbs measures regular.
Recent results in \cite{GSt} prove
existence of such Pesin sets for hyperbolic diffeomorphisms and flows for a large variety of Gibbs measures
determined by H\"older continuous potentials. The strong spectral estimates for
Ruelle operators and well-established techniques lead to exponential decay of correlations for H\"older continuous 
observables, as well as to some other consequences such as: (a) existence of a non-zero analytic continuation of the 
Ruelle zeta function  with a pole at the entropy in a vertical strip containing  the entropy in its interior; 
(b) a Prime Orbit Theorem with an exponentially small error.

\end{abstract}









\def\R{{\mathbb R}}
\def\T{{\mathbb T}}
\def\S{{\mathbb S}}
\def\C{{\mathbb C}}
\def\Z{{\mathbb Z}}
\def\N{{\mathbb N}}
\def\H{{\mathbb H}}
\def\B{{\mathbb B}}
\def\diam{\mbox{\rm diam}}
\def\sn{\S^{n-1}}
\def\rr{{\cal R}}
\def\mt{{\Lambda}}
\def\e{\emptyset}
\def\dQ{\partial Q}
\def\dk{\partial K}
\def\endofproof{{\rule{6pt}{6pt}}}
\def\di{\displaystyle}
\def\dist{\mbox{\rm dist}}
\def\sa+{\Sigma_A^+}
\def\du{\frac{\partial}{\partial u}}
\def\dv{\frac{\partial}{\partial v}}
\def\dt{\frac{d}{d t}}
\def\dx{\frac{\partial}{\partial x}}
\def\con{\mbox{\rm const }}
\def\mat{\mm_{at}}
\def\ma{\mm_{a}}
\def\lab{L_{ab}}
\def\labt{L_{abt}}
\def\mabn{\mm_{a}^N}
\def\man{\mm_a^N}
\def\labn{L_{ab}^N}
\def\fa{f^{(a)}}
\def\i{{\bf i}}
\def\gge{{\mathcal G}_\epsilon}
\def\gej{\chi^{(j)}_\mu}
\def\ge{\chi_\epsilon}
\def\chio{\chi^{(1)}}
\def\chit{\chi^{(2)}}
\def\chii{\chi^{(i)}}
\def\chil{\chi^{(\ell)}}
\def\gett{\chi^{(2)}_{\mu}}
\def\geol{\chi^{(1)}_{\ell}}
\def\getl{\chi^{(2)}_{\ell}}
\def\geil{\chi^{(i)}_{\ell}}
\def\gee{\chi_{\ell}}
\def\wloc{W_{\epsilon}}
\def\Int{\mbox{\rm Int}}
\def\dist{\mbox{\rm dist}}
\def\pr{\mbox{\rm pr}}
\def\tpp{\widetilde{\pp}}
\def\supp{\mbox{\rm supp}}
\def\Arg{\mbox{\rm Arg}}
\def\In{\mbox{\rm Int}}
\def\con{\mbox{\rm const}\;}
\def\Re{\mbox{\rm Re}}
\def\li{\mbox{\rm li}} 
\def\Seo{S^*_\epsilon(\Omega)}
\def\sdk{S^*_{\dk}(\Omega)}
\def\lae{\Lambda_{\epsilon}}
\def\ep{\epsilon}
\def\be{\begin{equation}}
\def\ee{\end{equation}}
\def\beqn{\begin{eqnarray*}}
\def\eeqn{\end{eqnarray*}}
\def\Pr{\mbox{\rm Pr}}

\def\gi{\gamma^{(i)}}
\def\ii{{\imath }}
\def\jj{{\jmath }}
\def\II{{\mathcal I}}
\def\ccij{ \cc_{i'_0,j'_0}[\eta]}
\def\la{\langle}
\def\ra{\rangle}
\def\bs{\bigskip}
\def\xio{\xi^{(0)}}
\def\xo{x^{(0)}}
\def\zo{z^{(0)}}
\def\Con{\mbox{\rm Const}\;}
\def\do{\partial \Omega}
\def\dk{\partial K}
\def\dl{\partial L}
\def\pr{{\rm pr}}
\def\dist{{\rm dist}}
\def\dds{\frac{d}{ds}}
\def\con{{\rm const}\;}
\def\Con{{\rm Const}\;}
\def\di{\displaystyle}
\def\oo{\mbox{\rm O}}
\def\hess{\mbox{\rm Hess}}
\def\gi{\gamma^{(i)}}
\def\endofproof{{\rule{6pt}{6pt}}}
\def\xm{x^{(m)}}
\def\vm{\varphi^{(m)}}
\def\km{k^{(m)}}
\def\dm{d^{(m)}}
\def\kam{\kappa^{(m)}}
\def\dem{\delta^{(m)}}
\def\xim{\xi^{(m)}}
\def\ep{\epsilon}
\def\ms{\medskip}
\def\ex{\mbox{\rm extd}}

\def\clip{C^{\mbox{\footnotesize \rm Lip}}}
\def\wlocs{W^s_{\mbox{\footnote\rm loc}}}
\def\Lip{\mbox{\rm Lip}}

\def\Xr{X^{(r)}}
\def\lip{\mbox{{\footnotesize\rm Lip}}}
\def\Vol{\mbox{\rm Vol}}

\def\naf{\nabla f(z)}
\def\so{\sigma_0}
\def\Xo{X^{(0)}}
\def\z1{z^{(1)}}
\def\Vo{V^{(0)}}
\def\Yo{Y{(0)}}

\def\uo{u^{(0)}}
\def\vo{v^{(0)}}
\def\no{\nu^{(0)}}
\def\psa{\partial^{(s)}_a}
\def\hcd{\hc^{(\delta)}}
\def\Md{M^{(\delta)}}
\def\Uo{U^{(1)}}
\def\Ut{U^{(2)}}
\def\Uj{U^{(j)}}
\def\no{n^{(1)}}
\def\nt{n^{(2)}}
\def\nj{n^{(j)}}
\def\ccm{\cc^{(m)}}

\def\ooo{\oo^{(1)}}
\def\oot{\oo^{(2)}}
\def\ooj{\oo^{(j)}}
\def\fo{f^{(1)}}
\def\ft{f^{(2)}}
\def\fj{f^{(j)}}
\def\wo{w^{(1)}}
\def\wt{w^{(2)}}
\def\wj{w^{(j)}}
\def\Vo{V^{(1)}}
\def\Vt{V^{(2)}}
\def\Vj{V^{(j)}}

\def\Ul{U^{(\ell)}}
\def\Uj{U^{(j)}}
\def\wl{w^{(\ell)}}
\def\Vl{V^{(\ell)}}
\def\Ujj{U^{(j+1)}}
\def\wjj{w^{(j+1)}}
\def\Vjj{V^{(j+1)}}
\def\Ujo{U^{(j_0)}}
\def\wjo{w^{(j_0)}}
\def\Vjo{V^{(j_0)}}
\def\vj{v^{(j)}}
\def\vl{v^{(\ell)}}

\def\f0{f^{(0)}}

\def\gl{\gamma_\ell}
\def\id{\mbox{\rm id}}
\def\piU{\pi^{(U)}}

\def\cca{C^{(a)}}
\def\bba{B^{(a)}}
\def\saa{\Sigma^+_A}
\def\sa{\Sigma_A}

\def\Int{\mbox{\rm Int}}
\def\epo{\ep^{(0)}}
\def\pH{\partial \H^{n+1}}
\def\sh{S^*(\H^{n+1})}
\def\zoo{z^{(1)}}
\def\yoo{y^{(1)}}
\def\xoo{x^{(1)}}


\def\supp{\mbox{\rm supp}}
\def\Arg{\mbox{\rm Arg}}
\def\In{\mbox{\rm Int}}
\def\diam{\mbox{\rm diam}}
\def\e{\emptyset}
\def\endofproof{{\rule{6pt}{6pt}}}
\def\di{\displaystyle}
\def\dist{\mbox{\rm dist}}
\def\con{\mbox{\rm const }}
\def\Box{\spadesuit}
\def\Int{\mbox{\rm Int}}
\def\dist{\mbox{\rm dist}}
\def\pr{\mbox{\rm pr}}
\def\be{\begin{equation}}
\def\ee{\end{equation}}
\def\beqn{\begin{eqnarray*}}
\def\eeqn{\end{eqnarray*}}
\def\la{\langle}
\def\ra{\rangle}
\def\bs{\bigskip}
\def\Con{\mbox{\rm Const}\;}
\def\clip{C^{\alpha}}
\def\wlocs{W^s_{\mbox{\footnote\rm loc}}}
\def\Lip{\mbox{\rm Lip}}
\def\lip{\mbox{\footnotesize\rm Lip}}
\def\Re{\mbox{\rm Re}}
\def\li{\mbox{\rm li}} 
\def\ep{\epsilon}
\def\ms{\medskip}
\def\dds{\frac{d}{ds}}
\def\oo{\mbox{\rm O}}
\def\hess{\mbox{\rm Hess}}
\def\id{\mbox{\rm id}}
\def\ii{{\imath }}
\def\jj{{\jmath }}
\def\graph{\mbox{\rm graph}}
\def\span{\mbox{\rm span}}
\def\Intu{\Int^u}
\def\Ints{\Int^s}

\def\i{{\bf i}}

\def\G{{\mathcal G}}
\def\nn{{\mathcal N}}
\def\mm{{\mathcal M}}
\def\kk{{\mathcal K}}
\def\ll{{\mathcal L}}
\def\vv{{\mathcal V}}
\def\bb{{\mathcal B}}
\def\ff{{\mathcal F}}
\def\tt{{\mathcal T}}
\def\uu{{\mathcal U}}
\def\pp{{\mathcal P}}
\def\aa{{\mathcal A}}
\def\cc{{\mathcal C}}
\def\oo{{\mathcal O}}
\def\dd{{\mathcal D}}
\def\ss{{\mathcal S}}
\def\Ee{{\mathcal E}}
\def\tEe{{\widetilde{\cal E}}}
\def\rr{{\mathcal R}}
\def\hh{{\mathcal H}}
\def\II{{\mathcal I}}
\def\jj{{\mathcal J}}

\def\hs{\hat{s}}
\def\hz{\hat{z}}
\def\hL{\widehat{L}}
\def\hl{\hat{l}}
\def\hl{\hat{l}}
\def\hc{\hat{\cc}}
\def\hbb{\widehat{\cal B}}
\def\hu{\hat{u}}
\def\hX{\hat{X}}
\def\hx{\hat{x}}
\def\hu{\hat{u}}
\def\hv{\hat{v}}
\def\hQ{\hat{Q}}
\def\hC{\widehat{C}}
\def\hF{\hat{F}}
\def\hf{\hat{\varphi}}
\def\hvarphi{\hat{\varphi}}
\def\hii{\hat{\ii}}
\def\hr{\hat{r}}
\def\hq{\hat{q}}
\def\hy{\hat{y}}
\def\hZ{\widehat{Z}}
\def\hz{\hat{z}}
\def\hE{\widehat{E}}
\def\hR{\widehat{R}}
\def\hell{\hat{\ell}}
\def\hs{\hat{s}}
\def\hW{\widehat{W}}
\def\hS{\widehat{S}}
\def\hV{\widehat{V}}
\def\hB{\widehat{B}}
\def\hhh{\widehat{\cal H}}
\def\hK{\widehat{K}}
\def\hU{\widehat{U}}
\def\hhh{\widehat{\hh}}
\def\hdd{\widehat{\dd}}
\def\hZ{\widehat{Z}}
\def\hGa{\widehat{\Gamma}}

\def\hal{\hat{\alpha}}
\def\hbe{\hat{\beta}}
\def\hg{\hat{\gamma}}
\def\hrho{\hat{\rho}}
\def\hd{\hat{\delta}}
\def\hphi{\hat{\phi}}
\def\hmu{\hat{\mu}}
\def\hnu{\hat{\nu}}
\def\hsi{\hat{\sigma}}
\def\htau{\hat{\tau}}
\def\hpi{\hat{\pi}}
\def\hep{\hat{\epsilon}}
\def\hxi{\hat{\xi}}
\def\hLa{\widehat{\Lambda}^u}
\def\hPhi{\widehat{\Phi}}
\def\hPsi{\widehat{\Psi}}
\def\hPhii{\widehat{\Phi}^{(i)}}
\def\hath{\hat{h}}

\def\tc{\tilde{C}}
\def\tg{\tilde{\gamma}}  
\def\tV{\widetilde{V}}
\def\tC{\widetilde{\cc}}
\def\tb{\tilde{b}}
\def\tt{\tilde{t}}
\def\tx{\tilde{x}}
\def\tp{\tilde{p}}
\def\tz{\tilde{Z}}
\def\tZ{\tilde{Z}}
\def\tF{\tilde{F}}
\def\tf{\tilde{\varphi}}
\def\tvarphi{\tilde{\varphi}}
\def\tp{\tilde{p}}
\def\te{\tilde{e}}
\def\tv{\tilde{v}}
\def\tu{\tilde{u}}
\def\tw{\tilde{w}}
\def\ts{\tilde{\sigma}}
\def\talpha{\tilde{\alpha}}
\def\tr{\tilde{r}}
\def\tU{\widetilde{U}}
\def\tS{\tilde{S}}
\def\tP{\widetilde{P}}
\def\ttau{\tilde{\tau}}
\def\tLip{\widetilde{\Lip}}
\def\tz{\tilde{z}}
\def\tS{\tilde{S}}
\def\tts{\tilde{\sigma}}
\def\tVl{\widetilde{V}^{(\ell)}}
\def\tVj{\widetilde{V}^{(j)}}
\def\tVo{\widetilde{V}^{(1)}}
\def\tVj{\widetilde{V}^{(j)}}
\def\tPsi{\tilde{\Psi}}
 \def\tp{\tilde{p}}
 \def\tVjo{\widetilde{V}^{(j_0)}}
\def\tvj{\tilde{v}^{(j)}}
\def\tVjj{\widetilde{V}^{(j+1)}}
\def\tvl{\tilde{v}^{(\ell)}}
\def\tVt{\widetilde{V}^{(2)}}
\def\tR{\widetilde{R}}
\def\tQ{\widetilde{Q}}
\def\oL{\tilde{\Lambda}}
\def\tq{\tilde{q}}
\def\tk{\tilde{k}}
\def\tx{\tilde{x}}
\def\ty{\tilde{y}}
\def\tz{\tilde{z}}
\def\txo{\tilde{x}^{(0)}}
\def\tso{\tilde{\sigma}_0}
\def\tmt{\tilde{\Lambda}}
\def\tg{\tilde{g}}
\def\tsi{\tilde{\sigma}}
\def\tC{\tilde{C}}
\def\tell{\tilde{\ell}}
\def\trho{\tilde{\rho}}
\def\ts{\tilde{s}}
\def\tB{\widetilde{B}}
\def\thh{\widetilde{\cal H}}
\def\tV{\widetilde{V}}
\def\trr{\tilde{r}}
\def\te{\tilde{e}}
\def\tv{\tilde{v}}
\def\tu{\tilde{u}}
\def\tw{\tilde{w}}
\def\trho{\tilde{\rho}}
\def\tell{\tilde{\ell}}
\def\tz{\tilde{Z}}
\def\tF{\tilde{F}}
\def\tf{\tilde{f}}
\def\tp{\tilde{p}}
\def\ttau{\tilde{\tau}}
\def\tz{\tilde{z}}
\def\tg{\tilde{\gamma}}  
\def\tV{\widetilde{V}}
\def\tC{\widetilde{\cc}}
\def\tLa{\widetilde{\Lambda}^u}
\def\tR{\widetilde{R}}
\def\tr{\tilde{r}}
\def\tc{\widetilde{C}}
\def\tD{\widetilde{D}}
\def\tt{\tilde{t}}
\def\tp{\tilde{p}}
\def\tS{\tilde{S}}
\def\tts{\tilde{\sigma}}
\def\tZ{\widetilde{Z}}
\def\tdelta{\tilde{\delta}}
\def\th{\tilde{h}}
\def\tB{\widetilde{B}}
\def\thh{\widetilde{\hh}}
\def\tep{\tilde{\ep}}
\def\tE{\widetilde{E}}
\def\tu{\tilde{u}}
\def\txi{\tilde{\xi}}
\def\teta{\tilde{\eta}}
\def\tRR{\widetilde{\rr}}

\def\sr{{\sc r}}
\def\mt{{\Lambda}}
\def\do{\partial \Omega}
\def\dk{\partial K}
\def\dl{\partial L}
\def\wloc{W_{\epsilon}}
\def\piU{\pi^{(U)}}
\def\Rio{\R_{i_0}}
\def\Ri{\R_{i}}
\def\Rii{\R^{(i)}}
\def\Riii{\R^{(i-1)}}
\def\hRii{\widehat{\R}_i}
\def\hRiio{\widehat{\R}_{(i_0)}}
\def\Eii{E^{(i)}}
\def\Eio{E^{(i_0)}}
\def\Rj{\R_{j}}
\def\Vio{{\cal V}^{i_0}}
\def\Vi{{\cal V}^{i}}
\def\Wio{W^{i_0}}
\def\Wioo{W^{i_0-1}}
\def\hi{h^{(i)}}
\def\Psii{\Psi^{(i)}}
\def\pii{\pi^{(i)}}
\def\piii{\pi^{(i-1)}}
\def\gxyii{g_{x,y}^{i-1}}
\def\span{\mbox{\rm span}}
\def\Jac{\mbox{\rm Jac}}
\def\Vol{\mbox{\rm Vol}}
\def\limp{\lim_{p\to\infty}}

\def\xijl{X_{i,j}^{(\ell)}}
\def\xij{X_{i,j}}
\def\hyijl{\widehat{Y}_{i,j}^{(\ell)}}
\def\hxijl{\widehat{X}_{i,j}^{(\ell)}}
\def\hxij{\widehat{X}_{i,j}}
\def\eijl{\omega_{i,j}^{(\ell)}}
\def\eij{\omega_{i,j}}
\def\Gl{\Gamma_\ell}

\def\cB{\check{B}}
\def\tpi{\tilde{\pi}}
\def\J{{\sf J}}
\def\bJ{{\mathbb J}}

\def\hcc{\widehat{\cc}}
\def\hpp{\widehat{\pp}}
\def\ttP{\widetilde{\pp}}
\def\tP{\widetilde{P}}
\def\hP{\widehat{P}}
\def\hY{\widehat{Y}}

\def\diamtef{{\footnotesize \diam_\theta}}


\def\tc{\tilde{C}}
\def\tg{\tilde{\gamma}}  
\def\tV{\widetilde{V}}
\def\tW{\widetilde{W}}
\def\tC{\widetilde{\cc}}
\def\tKo{\widetilde{K_0}}
\def\tUKo{\widetilde{U\setminus K_0}}

\def\wo{w^{(1)}}
\def\vo{v^{(1)}}
\def\uo{u^{(1)}}
\def\wt{w^{(2)}}
\def\xio{\xi^{(1)}}
\def\xit{\xi^{(2)}}
\def\etao{\eta^{(1)}}
\def\etat{\eta^{(2)}}
\def\zetao{\zeta^{(1)}}
\def\zetat{\zeta^{(2)}}
\def\vt{v^{(2)}}
\def\ut{u^{(2)}}
\def\Wo{W^{(1)}}
\def\Vo{V^{(1)}}
\def\Uo{U^{(1)}}
\def\Wt{W^{(2)}}
\def\Vt{V^{(2)}}
\def\Ut{U^{(2)}}
\def\tmu{\tilde{\mu}}
\def\tla{\tilde{\lambda}}
\def\diamf{{\rm\footnotesize diam}}
\def\Intu{\mbox{\rm Int}^u}
\def\Ints{\mbox{\rm Int}^s}

\def\Bmt{\overline{B_{\ep_0}(\mt)}}
\def\Lye{L_{y,\eta}}
\def\Lyep{L^{(p)}_{y,\eta}}
\def\Fyp{F^{(p)}_y}
\def\Fxp{F^{(p)}_x}
\def\Lxx{L_{x,\xi}}
\def\Lxxp{L^{(p)}_{x,\xi}}

\def\Wuo{W^{u,1}}
\def\Wui{W^{u,i}}
\def\Wuj{W^{u,j}}
\def\Wut{\tW^{u,2}}
\def\Wuk{W^{u,k}}
\def\Wuh{\hW^{u}}
\def\tWuo{\tW^{u,1}}
\def\tWui{\tW^{u,i}}
\def\tWuj{\tW^{u,j}}
\def\tWuk{\tW^{u,k}}
\def\hWuo{\hW^{u,1}}
\def\hWui{\hW^{u,i}}
\def\hWuj{\hW^{u,j}}
\def\hWuk{\hW^{u,k}}
\def\dj{\delta^{(j)}}
\def\do{\delta^{(1)}}
\def\epj{\ep^{(j)}}
\def\epo{\ep^{(1)}}
\def\hSj{\widehat{S}^{(j)}}
\def\hSo{\widehat{S}^{(1)}}

\def\tmu{\tilde{\mu}}
\def\tla{\tilde{\lambda}}
\def\hE{\widehat{E}}
\def\uk{u^{(k)}}
\def\ui{u^{(i)}}
\def\uj{u^{(j)}}
\def\vk{v^{(k)}}
\def\vl{v^{(l)}}
\def\vi{v^{(i)}}
\def\vj{v^{(j)}}
\def\wk{w^{(k)}}
\def\wi{w^{(i)}}
\def\wj{w^{(j)}}
\def\etak{\eta^{(k)}}
\def\etai{\eta^{(i)}}
\def\etaj{\eta^{(j)}}
\def\zetak{\zeta^{(k)}}
\def\zetai{\zeta^{(i)}}
\def\zetaj{\zeta^{(j)}}


\def\yj{y^{(j)}}
\def\yi{y^{(i)}}
\def\tyi{\ty^{(i)}}
\def\yo{y^{(1)}}
\def\zj{z^{(j)}}
\def\zo{z^{(1)}}
\def\vj{v^{(j)}}
\def\vo{v^{(1)}}
\def\kaj{\kappa^{(j)}}
\def\kao{\kappa^{(1)}}

\def\tyj{\tilde{y}^{(j)}}
\def\yl{y^{(\ell)}}
\def\tyl{\tilde{y}^{(l)}}
\def\wo{w^{(1)}}
\def\vo{v^{(1)}}
\def\vi{v^{(i)}}
\def\vj{v^{(j)}}
\def\vk{v^{(k)}}
\def\uo{u^{(1)}}
\def\wt{w^{(2)}}
\def\xio{\xi^{(1)}}
\def\xit{\xi^{(2)}}
\def\xii{\xi^{(i)}}
\def\xij{\xi^{(j)}}
\def\hxio{\hxi^{(1)}}
\def\hxit{\hxi^{(2)}}
\def\hxii{\hxi^{(i)}}
\def\hxij{\hxi^{(j)}}

\def\cxi{\check{\xi}}
\def\cxio{\cxi^{(1)}}
\def\cxit{\cxi^{(2)}}
\def\cet{\check{\eta}}
\def\ceto{\cet^{(1)}}
\def\cett{\cet^{(2)}}
\def\cv{\check{v}}
\def\cvo{\cv^{(1)}}
\def\cvt{\cv^{(2)}}
\def\cu{\check{u}}
\def\cuo{\cu^{(1)}}
\def\cut{\cu^{(2)}}
\def\cj{c^{(j)}}
\def\fj{f^{(j)}}
\def\gji{g^{(j,i)}}
\def\tPsi{\widetilde{\Psi}}
\def\chU{\check{U}}


\def\Ulo{U^{(\ell_0)}}
\def\dte{D_\theta}
\def\diamte{\mbox{\rm diam}_{\theta}}
\def\Ial{I^{(\alpha)}}
\def\uml{u_m^{(\ell)}}
\def\yl{y^{(\ell)}}
\def\tyl{\tilde{y}^{(\ell)}}
\def\ool{\oo^{(\ell)}}
\def\fl{f^{(\ell)}}
\def\hep{\hat{\ep}}
\def\dl{d^{(\ell)}}
\def\dli{d_{\ell,i}}
\def\dlo{d_{\ell,1}}
\def\dlt{d_{\ell,2}}
\def\Lipt{{\Lip_\theta}}
\def\lipt{{\footnotesize \Lip_\theta}}
\def\tm{\tilde{m}}
\def\tj{\tilde{j}}
\def\lengthf{\mbox{\rm\footnotesize length}}
\def\length{\mbox{\rm length}}


\def\Xijl{X^{(\ell)}_{i,j}}
\def\hXijl{\widehat{X}^{(\ell)}_{i,j}}
\def\Wl{W^{(\ell)}}
\def\omijl{\omega^{(\ell)}_{i,j}}

\def\hXitl{\widehat{X}^{(\ell)}_{i,t}}
\def\Vl{V^{(\ell)}}
\def\omitl{\omega^{(\ell)}_{i,t}}
\def\Xisl{X^{(\ell)}_{i,s}}
\def\hXisl{\widehat{X}^{(\ell)}_{i,s}}
\def\omisl{\omega^{(\ell)}_{i,s}}
\def\hGa{\widehat{\Gamma}}
\def\hOm{\widehat{\Omega}}
\def\tGa{\widetilde{\Gamma}}
\def\hA{\widehat{A}}
\def\tnu{\tilde{\nu}}
\def\tX{\widetilde{X}}

\def\ww{{\mathcal W}}
\def\Zl{Z^{(\ell)}}
\def\hpp{\widehat{\pp}}
\def\tnn{\widetilde{\nn}}

\def\ftt{f^{(t)}}
\def\f0{f^{(0)}}
\def\fat{f^{(at)}}
\def\Fat{F^{(at)}}
\def\Fa{F^{(a)}}
\def\F0{F^{(0)}}
\def\tu{\tilde{u}}
\def\tD{\widetilde{D}}
\def\tchi{\tilde{\chi}}
\def\tC{\widetilde{C}}
\def\hC{\widehat{C}}
\def\hQ{\widehat{Q}}
\def\hF{\widehat{F}}
\def\hD{\widehat{D}}
\def\hr{\hat{r}}
\def\psid{\psi^\dag}
\def\taud{\tau^\dag}
\def\Omn{\Omega^{(n)}}
\def\Omm{\Omega^{(m)}}
\def\Omk{\Omega^{(k)}}
\def\Conf{{\mbox{\footnotesize\rm Const}}}
\def\hp{\hat{p}}

\def\tj{t^{(j)}}
\def\tyj{\tilde{y}^{(j)}}
\def\tyjo{\tilde{y}_{j,1}}
\def\tyjt{\tilde{y}_{j,2}}
\def\tyji{\tilde{y}_{j,i}}
\def\yjo{y_{j,1}}
\def\yjt{y_{j,2}}
\def\yji{y_{j,i}}
\def\tylo{\tilde{y}_{\ell,1}}
\def\tylt{\tilde{y}_{\ell,2}}
\def\tyli{\tilde{y}_{\ell,i}}
\def\ylo{y_{\ell,1}}
\def\ylt{y_{\ell,2}}
\def\yli{y_{\ell,i}}

\def\ulo{u_{\ell,1}}
\def\ult{u_{\ell,2}}
\def\uli{u_{\ell,i}}
\def\tulo{\tilde{u}_{\ell,1}}
\def\tult{\tilde{u}_{\ell,2}}
\def\tuli{\tilde{u}_{\ell,i}}

\def\tdlo{\tilde{d}_{\ell,1}}
\def\tdlt{\tilde{d}_{\ell,2}}
\def\tdli{\tilde{d}_{\ell,i}}
\def\dlo{d_{\ell,1}}
\def\dlt{d_{\ell,2}}
\def\dli{d_{\ell,i}}
\def\wjo{w_{j,1}}
\def\wjt{w_{j,2}}
\def\wji{w_{j,i}}
\def\sj{s^{(j)}}
\def\Yj{Y^{(j)}}
\def\Vj{V^{(j)}}
\def\Zj{Z^{(j)}}
\def\vj{v^{(j)}}
\def\wj{w^{(j)}}
\def\twj{\tilde{w}^{(j)}}
\def\gj{g^{(j)}}
\def\tgj{\tilde{g}^{(j)}}
\def\tg{\tilde{g}}
\def\hn{\hat{n}}

\def\hbeta{\hat{\beta}}
\def\hmu{\hat{\mu}}
\def\piS{\pi^{(S)}}
\def\hb{\hat{b}}
\def\shP{P^\sharp}
\def\tshP{\widetilde{P}^\sharp}
\def\T{\mathcal T}
\def\tut{\tu^{(2)}}
\def\twt{\tw^{(2)}}
\def\piS{\pi^{(S)}}
\def\tsigma{\tilde{\sigma}}
\def\td{\tilde{d}}
\def\m{{\sf m}}
\def\tXi{\widetilde{\Xi}}

\def\Omb{\Omega^{(\hb)}}
\def\Xib{\Xi^{(L\hb)}}

\def\tomega{\tilde{\omega}}
\def\vlim{v^{(\ell)}_{i,m}}
\def\tvlim{\tilde{v}^{(\ell)}_{i,m}}
\def\thetlim{\theta^{(\ell)}_{i,m}}
\def\lamblim{\lambda^{(\ell)}_{i,m}}
\def\vl{v^{(\ell)}}
\def\vlik{v^{(\ell)}_{i,k}}
\def\vlom{v^{(\ell)}_{1,m}}
\def\vltm{v^{(\ell)}_{2,m}}
\def\tvlik{\tilde{v}^{(\ell)}_{i,k}}
\def\tvlom{\tilde{v}^{(\ell)}_{1,m}}
\def\tvltm{\tilde{v}^{(\ell)}_{2,m}}
\def\Cum{\cc^{(u_m)}}
\def\Cu{\cc^{(u)}}
\def\gao{\gamma^{(1)}}
\def\gat{\gamma^{(2)}}
\def\diamtef{{\footnotesize\mbox{\rm  diam}_\theta}}
\def\tcc{\widetilde{\cc}}
\def\halpha{\hat{\alpha}}
\def\hgamma{\hat{\gamma}}
\def\tXijl{\widetilde{X}^{(\ell)}_{i.j}}
\def\tdd{\widetilde{\dd}}
\def\hcc{\widehat{\cc}}
\def\lambdam{\lambda^{(m)}}
\def\thetam{\theta^{(m)}}
\def\tomijl{\tilde{\omega}^{(\ell)}_{i,j}}
\def\oloij{\omega^{(\ell,0)}_{i,j}}
\def\olooj{\omega^{(\ell,0)}_{1,j}}
\def\olotj{\omega^{(\ell,0)}_{2,j}}
\def\olkij{\omega^{(\ell,k)}_{i,j}}
\def\olkoj{\omega^{(\ell,k)}_{1,j}}
\def\olktj{\omega^{(\ell,k)}_{2,j}}
\def\oloim{\omega^{(\ell,0)}_{i,m}}
\def\oloom{\omega^{(\ell,0)}_{1,m}}
\def\olotm{\omega^{(\ell,0)}_{2,m}}
\def\olkim{\omega^{(\ell,k)}_{i,m}}
\def\olkom{\omega^{(\ell,k)}_{1,m}}
\def\olktm{\omega^{(\ell,k)}_{2,m}}
\def\vlo{v^{(\ell)}_{1}}
\def\vlt{v^{(\ell)}_2}
\def\vli{v^{(\ell)}_{i}}
\def\vloij{v^{(\ell,0)}_{i,j}}
\def\vlooj{v^{(\ell,0)}_{1,j}}
\def\vlotj{v^{(\ell,0)}_{2,j}}
\def\vlkij{v^{(\ell,k)}_{i,j}}
\def\vlkoj{v^{(\ell,k)}_{1,j}}
\def\vlktj{v^{(\ell,k)}_{2,j}}

\def\vloim{v^{(\ell,0)}_{i,m}}
\def\vloom{v^{(\ell,0)}_{1,m}}
\def\vlotm{v^{(\ell,0)}_{2,m}}
\def\vlkim{v^{(\ell,k)}_{i,m}}
\def\vlkom{v^{(\ell,k)}_{1,m}}
\def\vlktm{v^{(\ell,k)}_{2,m}}

\def\vlrim{v^{(\ell,r)}_{i,m}}
\def\vlrom{v^{(\ell,r)}_{1,m}}
\def\vlrtm{v^{(\ell,r)}_{2,m}}

\def\Xloij{X^{(\ell,0)}_{i,j}}
\def\Xlooj{X^{(\ell,0)}_{1,j}}
\def\Xlotj{X^{(\ell,0)}_{2,j}}
\def\Xlkij{X^{(\ell,k)}_{i,j}}
\def\Xlkoj{X^{(\ell,k)}_{1,j}}
\def\Xlktj{X^{(\ell,k)}_{2,j}}

\def\Eekj{\Ee^{(k)}_{j}}
\def\Eerj{\Ee^{(r)}_{j}}

\def\Wr{W^{(r)}}
\def\Wk{W^{(k)}}

\def\tvlo{\tilde{v}^{(\ell)}_{1}}
\def\tvlt{\tilde{v}^{(\ell)}_{2}}
\def\tvli{\tilde{v}^{(\ell)}_{i}}

\def\omegak{\omega^{(k)}}
\def\nnk{\nn^{(k)}}
\def\omegao{\omega^{(o)}}
\def\nno{\nn^{(o)}}
\def\htheta{\hat{\theta}}

\def\ho{h^{(1)}}
\def\Ho{H^{(1)}}
\def\tcc{\widetilde{\cc}}
\def\utk{u^{(\tk)}}
\def\hXi{\widehat{\Xi}}
\def\chR{\check{R}}
\def\llm{\ll^{(m)}}
\def\llk{\ll^{(k)}}

\def\hB{\widehat{B}}
\def\Omb{\Omega^{(\hb)}}
\def\tnn{\widetilde{\nn}}

\def\Xim{\Xi^{(m)}}
\def\Omm{\Omega^{(m)}}
\def\Xin{\Xi^{(n)}}
\def\Omn{\Omega^{(n)}}
\def\hP{\widehat{P}}

\def\halpha{\hat{\alpha}}
\def\hgamma{\hat{\gamma}}

\def\E{\mathcal E}
\def\tp{\tilde{p}}
\def\tq{\tilde{q}}

\def\tyj{\tilde{y}^{(j)}}
\def\yl{y^{(\ell)}}
\def\wo{w^{(1)}}
\def\vo{v^{(1)}}
\def\vi{v^{(i)}}
\def\vj{v^{(j)}}
\def\vk{v^{(\tk)}}
\def\uo{u^{(1)}}
\def\wt{w^{(2)}}
\def\xio{\xi^{(1)}}
\def\xit{\xi^{(2)}}
\def\xii{\xi^{(i)}}
\def\xij{\xi^{(j)}}
\def\hxio{\hxi^{(1)}}
\def\hxit{\hxi^{(2)}}
\def\hxii{\hxi^{(i)}}
\def\hxij{\hxi^{(j)}}

\def\cxi{\check{\xi}}
\def\cxio{\cxi^{(1)}}
\def\cxit{\cxi^{(2)}}
\def\cet{\check{\eta}}
\def\ceto{\cet^{(1)}}
\def\cett{\cet^{(2)}}
\def\cv{\check{v}}
\def\cvo{\cv^{(1)}}
\def\cvt{\cv^{(2)}}
\def\cu{\check{u}}
\def\cuo{\cu^{(1)}}
\def\cut{\cu^{(2)}}
\def\cj{c^{(j)}}
\def\fj{f^{(j)}}
\def\gji{g^{(j,i)}}

\def\ao{a^{(1)}}
\def\bo{b^{(1)}}
\def\zi{z^{(i)}}
\def\hGa{\widehat{\Gamma}}

\def\tvo{\tv^{(1)}}
\def\tetao{\teta^{(1)}}

\def\Bmt{\overline{B_{\ep_0}(\mt)}}
\def\chBo{\check{B}^{u,1}}
\def\tBo{\tB^{u,1}}
\def\hBo{\hB^{u,1}}
\def\hpi{\hat{\pi}}
\def\tU{\widetilde{U}}
\def\tr{\tilde{r}}
\def\chBo{\check{B}^{u,1}}
\def\tBo{\tB^{u,1}}
\def\hBo{\hB^{u,1}}
\def\hpi{\hat{\pi}}
\def\tL{\widetilde{L}}

\def\pij{\pi^{(j)}}
\def\wo{w^{(1)}}
\def\vo{v^{(1)}}
\def\uo{u^{(1)}}
\def\wt{w^{(2)}}
\def\vt{v^{(2)}}
\def\twt{\tw^{(2)}}
\def\tvt{\tv^{(2)}}
\def\xio{\xi^{(1)}}
\def\xit{\xi^{(2)}}
\def\txio{\txi^{(1)}}
\def\txit{\txi^{(2)}}
\def\ut{u^{(2)}}
\def\tut{\tu^{(2)}}
\def\tuo{\tu^{(1)}}
\def\tut{\tu^{(2)}}
\def\Wo{W^{(1)}}
\def\Vo{V^{(1)}}
\def\Uo{U^{(1)}}
\def\tUo{\tU^{(1)}}
\def\Wt{W^{(2)}}
\def\Vt{V^{(2)}}
\def\Ut{U^{(2)}}
\def\tWt{\tW^{(2)}}
\def\Vt{\tV^{(2)}}
\def\tUt{\tU^{(2)}}

\def\tde{\tilde{\delta}}
\def\tSt{\widetilde{S}^{(2)}}
\def\So{S^{(1)}}
\def\hr{\hat{r}}
\def\li{\mbox{\rm li}}

\section{Introduction and Results}
\renewcommand{\theequation}{\arabic{section}.\arabic{equation}}

\subsection{Introduction}
The study of statistical properties of dynamical systems has a long history and has been the subject 
of a considerable interest due to their applications in statistical mechanics and thermodynamics. 
Many physical systems poses some kind of `strong hyperbolicity' and are known to have or expected 
to have strong mixing properties. For example in the 70's, due to works by Sinai, Bowen and Ruelle,
it was already known that  for Anosov diffeomorphisms exponential decay of correlations takes place 
for H\"older continuous observables (see e.g. the survey article \cite{ChY}). However the continuous 
case proved to be much more difficult and it took more than twenty years until the breakthrough work of 
Dolgopyat \cite{D1}, where he established exponential decay of correlations for H\"older continuous 
potentials in two major cases: (i) geodesic flows on compact surfaces of negative curvature (with respect to any 
Gibbs measure);  (ii) transitive Anosov flows  on compact Riemannian manifolds with $C^1$ jointly non-integrable local 
stable and unstable foliations (with respect to the Sinai-Bowen-Ruelle measure). The work of Dolgopyat
was preceded (and possibly partially inspired) by that of Chernov \cite{Ch1} who proved sub-exponential decay of 
correlations for Anosov flows on 3D Riemannian manifolds  (with respect to the Sinai-Bowen-Ruelle measure). 

Dolgopyat's work was followed by a considerable activity to establish exponential and other types of decay of
correlations for various  kinds of systems, most of the results dealing with the measure determined by the Riemannian volume.
Without trying to provide a comprehensive review of the literature, a sample of important works in this area is the following:


(a)  The so called functional-analytic approach initiated by the work of Blank, Keller and Liverani
\cite{BKL} which involves the so called Ruelle-Perron-Frobenius operators  
$\di \ll_tg = \frac{g\circ \phi_{-t}}{|(\det d\phi_t)| \circ \phi_{-t} }$, $t\in \R$, was further developed by various authors,
notably Liverani, Baladi, Tsujii, Gou\"ezel,  and many others (\cite{BaG}, \cite{BaT}, \cite{GL1}, \cite{GL2}, \cite{T})  
-  see e.g.  the lectures of Liverani \cite{L3} for a nice exposition 
of the main ideas. Using this method and also (generally speaking) "Dolgopyat's cancellation mechanism" from \cite{D1},
Liverani \cite{L2} proved exponential decay of correlations for $C^4$ contact Anosov flows with respect to the measure 
determined by the Riemannian volume. Some finer results were obtained later by  Tsujii \cite{T} (for $C^3$ contact Anosov flows).

A similar approach, however studying Ruelle-Perron-Frobenius operators acting on currents, 
was used by Giulietti, Liverani and Pollicott in \cite{GLP} where 
they proved some remarkable results. For example, they established that
for $C^\infty$ Anosov flows the Ruelle zeta function is meromorphic in the whole complex plane. 
In  \cite{GLP} the authors derived also (amongst other things) 
exponential decay of correlations for contact Anosov flows with respect to the measure of maximal 
entropy (generated by the potential $F = 0$) under a bunching condition (which implies that the 
stable/unstable foliations are $\frac{2}{3}$-H\"older).

Various other results on decay of correlations for uniformly hyperbolic systems  have been established using different methods as well, 
see e.g.  \cite{Ch1}, \cite{D2}, \cite{D3}, \cite{St1}, \cite{N}, \cite{ABV},  \cite{FMT},  \cite{OWi},  \cite{Wi},
and the references there.

\ms

(b) The ideas of Dolgopyat were used for example in \cite{BaV} and \cite{AGY} to prove exponential
decay of correlations for systems with infinite Markov partitions  (with respect to SRB measures).
For such systems a general approach was invented by L.-S. Young \cite{Y1}, \cite{Y2} who introduced the so called "Young towers". 
This approach was later used by many authors in
a variety of papers dealing with decay of correlations for diffeomorphisms and flows (uniformly and non-uniformly
hyperbolic) - see e.g. \cite{M1}, \cite{M2}, \cite{AM}, and the historical remarks and references there.
See also \cite{Sa} for a different approach.
 
\ms

(c) Decay of correlations for hyperbolic systems with singularities (e.g. billiards) have been studied for a very long time.
The first results in this area deal with the corresponding discrete dynamical system, generated by the 
billiard ball map from boundary to boundary. To my knowledge, these results were: the subexponential decay of correlations 
established for a very large class of dispersing billiards by Bunimovich, Sinai and Chernov \cite{BSC} and the exponential 
decay of correlations for some classes of dispersing billiards in the plane and on the two-dimensional torus established by 
Young \cite{Y1} and Chernov \cite{Ch2} as  consequences of their more general arguments. See also \cite{ChY} and
the comments and references there. Later on various other results 
were established both for the discrete and continuous dynamical systems, all of them dealing with SRB measures. Notably, 
Melbourne \cite{M1} proved super-polynomial decay of correlations (rapid mixing) for Lorentz  billiard flows with finite
horizon on a two-dimensional torus, while Chernov \cite{Ch4} established stretched exponential decay of correlations
for such flows. More recently Baladi and Liverani \cite{BaL} proved exponential decay of correlations for 
piecewise hyperbolic contact flows on three-dimensional  manifolds. Finally, using techniques, methods and ideas
from \cite{D1}, \cite{L2}, \cite{BaL}, \cite{ChM}, \cite{DZ1} and \cite{DZ2}, in a remarkable
recent paper Baladi, Demers and Liverani \cite{BaDL} established exponential decay of correlations for Sinai
billiards with finite horizon on a two-dimensional torus. See also \cite{BNST} and the historical remarks and references in \cite{BaDL}.

\ms

(d) During the last several years there has been a lot of activity in applying methods and tools that are usually seen 
in the analysis of PDE's and scattering theory to dynamical problems such as the study of decay of correlations, dynamical
zeta functions  and  the distribution of Ruelle-Pollicott resonances. For example: Nonnenmacher and Zworski \cite{NZ} 
established exponential decay of correlations for  a class of $C^\infty$ flows which includes the $C^\infty$ contact Anosov flows;
Dyatlov and Zworski \cite{DyZ2} gave a proof using microlocal analysis of the meromorphic continuation of the Ruelle
zeta function for $C^\infty$ Anosov flows. Various other interesting and deep results have been obtained in \cite{DDyZ}, \cite{DyFG},
\cite{DyG}, \cite{DyZ1}, \cite{DyZ3}, \cite{FaSj}, \cite{FaT1}, \cite{FaT2}, \cite{JZ}. See also the comments and references in those papers.
To my knowledge, {\bf all works in this area deal with the measure determined by the Riemannian volume}.

\ms

(e) What concerns extension and further development of the ideas of Dolgopyat \cite{D1}, we should mention here our 
papers \cite{St1} - \cite{St3}. Strictly speaking the first result on exponential decay of correlations for billiard flows was the one 
in \cite{St1} for open billiard flows in the plane, however the system in this case is uniformly hyperbolic and admits a finite 
Markov family\footnote{In \cite{St1} the so called {\it triple intersection property} of cylinders on unstable manifolds was introduced, 
and it was proved in the case of open billiard flows in the plane. As a consequence, it was established that for such flows all Gibbs measures 
generated by H\"older potentials had the so called {\it Federer property}, and then an appropriate modification of the approach 
in \cite{D1} could be used. Later, Naud \cite{N} used a similar procedure in the case of geodesic flows on 
 convex co-compact hyperbolic surfaces.}. Higher-dimensional open billiards were considered in \cite{St3}, 
 however under an additional condition.
The results in  \cite{St2} are much more general - they deal with Axiom A flows on basic sets having Lipschitz stable and
unstable foliation and satisfying a certain non-integrability condition. E.g. as shown in \cite{St5}, {\bf contact Axiom A flows on basic sets
satisfying a pinching condition} (similar to the $1/4$-pinching condition for geodesic flows on manifolds of negative curvature)
always satisfy the assumptions in \cite{St2}, and therefore they {\bf have exponential decay of correlations with respect to
any Gibbs measure generated by a H\"older potential}. In the present work we generalise the results and the
developments  in \cite{St2} by far. 

Ruelle transfer operators with two complex parameters were studied in \cite{PeS5} and
results in the spirit of these in \cite{St2} were established for the same kind of Axiom A flows on basic sets as the one considered
in \cite{St2}. It seems this was the first time transfer operators depending on two parameters have been considered. It should be
mentioned that the transition from one to two parameters is non-trivial.

\ms

In this work, as a consequence of the main result, we derive exponential decay of correlations 
for $C^5$ contact Anosov flows on Riemannian manifolds $M$ of {\bf any dimension and with respect to any regular
Gibbs measure on $M$, i.e. a Gibbs measure  admitting a Pesin set with exponentially small tails} (see Sect. 1.2 for 
the precise definition of the latter).  
We should stress again the fact that we deal with a large variety of Gibbs measures, not just SRB measures.
It appears that so far the only results  of this kind has been that of Dolgopyat \cite{D1} for 
geodesic flows on $C^5$ compact  surfaces and the one in \cite{St2} for Axiom A flows on basic sets (under some additional assumptions).

\ms


In \cite{D1} Dolgopyat developed  a certain technique   involving estimates of  
Laplace transforms of  correlations functions (following previous works of Pollicott \cite{Po1} 
and Ruelle \cite{R3}) that  leads more or less automatically to  exponential decay of correlations for 
H\"older continuous potentials, once certain strong spectral properties of Ruelle transfer operators 
have been established. Given an Anosov flow
 $\phi_t : M \longrightarrow M$ on a  Riemannian manifold $M$, consider a Markov family 
consisting of rectangles  $R_i = [U_i ,S_i ] = \{ [x,y] : x \in U_i, y \in S_i\}$, where $U_i$ and  $S_i$ are pieces of unstable/stable 
manifolds at some $z_i\in M$ , the first return time function 
$$\tau : R = \cup_{i=1}^{k_0} R_i \longrightarrow [0,\infty)$$ 
and the standard Poincar\'e map $\pp: R \longrightarrow R$ (see the first two pages in Ch. 2 for the detailed definitions of $R_i$, $\tau$ and $\pp$).
The {\it shift map} 
$$\sigma : U = \cup_{i=1}^{k_0} U_i\longrightarrow U ,$$
given by $\sigma = \piU\circ \pp$, where 
$\piU : R \longrightarrow U$ is the projection along the leaves of local stable manifolds, 
defines a dynamical system which is essentially isomorphic to an one-sided Markov shift.  Given a bounded function 
$f \in B (U)$, one defines the {\it Ruelle transfer operator} $L_{f} : B (U) \longrightarrow B (U)$ by 
$$(L_f h)(x) = \sum_{\sigma(y) = x} e^{f(y)} h(y) .$$
Assuming that $f$ is real-valued and H\"older continuous, 
let $P_f \in \R$ be such that the topological pressure of $f-P_f\tau$ with respect to $\sigma$ is zero 
(cf. e.g. \cite{PP}). Dolgopyat proved (for the type of flows he considered in \cite{D1}) that for small 
$|a|$ and large $|b|$ the spectral radius of the Ruelle operator 
$$L_{f-(P_f+ a+\i b)\tau} : C^\alpha (U) \longrightarrow C^\alpha (U)$$
acting on $\alpha$-H\"older continuous functions 
($0 < \alpha \leq 1$) is uniformly bounded by a constant $\rho < 1$. 

More general results of this kind were proved in \cite{St2} for mixing Axiom A flows on basic sets under some 
additional regularity assumptions,  amongst them -- Lipschitzness of the so local stable holonomy 
maps\footnote{In general these are only H\"older continuous -- see  \cite{Ha1}, \cite{Ha2}.} (see Sect. 2). 

{\bf Our main result in this paper is that for contact Anosov flows on compact  Riemannian manifolds $M$ correlations 
for H\"older continuous observables decay exponentially fast with respect to any  regular Gibbs measure on $M$. }

It was proved recently in \cite{GSt} that Pesin sets with exponentially small tails exist for Gibbs measures for Axiom A flows (and
diffeomorphism) satisfying a certain condition, called exponential large deviations for all Lyapunov exponents
(see Sect. 3 below). In fact, under such a condition,  Pesin sets with exponentially small tails exist for every
continuous linear cocycle over a transitive subshift of finite type (see Theorem 1.7 in \cite{GSt}). And it turns
out that in this generality, exponential large deviations for all exponents is a generic condition (see Theorem 1.5
in \cite{GSt}). 

The main results mentioned above  are in fact consequences of a more general result.
Given $\theta \in (0,1)$, the metric $\dte$ on $U$ is defined by
$\dte(x,y) = 0$ if $x = y$, $\dte(x,y) = 1$ if $x,y$ belong to different $U_i$'s and $\dte(x,y) = \theta^N$ if 
$\pp^j(x)$ and $\pp^j(y)$ belong to the same rectangle $R_{i_j}$ for all $j = 0,1, \ldots,N-1$, 
and $N$ is the largest integer with this property. Denote by $\ff_\theta(U)$ {\it the space of all functions 
$h : U \longrightarrow \C$ with Lipschitz constants }
$$|h|_\theta = \sup \left\{ \frac{|h(x) - h(y)|}{\dte(x,y)} : x\neq y  \, , \, x,y \in U\right\} < \infty .$$
The central Theorem 1.3 below says that
for sufficiently large $\theta \in (0,1)$ and any real-valued function $f \in \ff_\theta(U)$ the Ruelle transfer 
operators related to $f$  are eventually  contracting on $\ff_\theta (U)$ (see the definition of this in Sect. 1.2). 
A similar result holds for H\"older  continuous functions on $U$ --  see Corollary 1.4 below. 
\index{$\ff_\theta(U)$} 
\index{$\dte$} 
\index{$|h|_\theta$} 

In the proof of the central Theorem 1.3  we generalise significantly  the approach of Dolgopyat \cite{D1} and 
its development in \cite{St2}. The general framework described in Sections 6, 7 and 8 below is indeed rather 
general and is expected to work in a variety of other situations (and possibly for some non-uniformly 
hyperbolic systems, as well). The {\bf assumption of fundamental importance is the existence of a Pesin set
with exponentially small tails}. Certain technical troubles, such as the lack of regularity of the local stable/unstable 
manifolds and related local stable/unstable holonomy maps\footnote{E.g. the local stable holonomy maps are defined by 
sliding along local stable manifolds -- as we mentioned earlier, in general these are only H\"older continuous. In \cite{D1} 
and \cite{St2} these were assumed to be $C^1$ and  Lipschitz, respectively. Since the definition of Ruelle operators itself
involves sliding along local stable manifolds, it appears to be a significant technical problem to overcome the 
lack of regularity in general.} are overcome here by using the fact that the flow is contact and by using Lyapunov exponents
on an appropriately chosen Pesin set. The most significant part in 
overcoming these difficulties is Sect. 9 below dealing with non-integrability matters -  it is technical, lengthy and non-trivial.
However, as far as ideas are concerned, 
{\bf the most significant ideas in this work are those in Sections 6 and 7} below. See Sect 1.3 for some more details.

It has been well known since Dolgopyat's paper \cite{D1} that
strong spectral estimates for Ruelle transfer operators as the ones described  in Theorem 1.3 
lead to deep results concerning zeta functions and related topics which are difficult to obtain by 
other means. For example, such estimates were fundamental in \cite{PoS1}, where the statements in
Theorem 1.2, stated in Sect. 1.2 below, was proved for geodesic flows on compact surfaces of negative curvature. 
In particular, a sharp asymptotic was obtained for the number $\pi(\lambda)$ of closed orbits of length not exceeding 
$\lambda$. For the same kind  of flows, fine and very interesting asymptotic estimates for pairs of closed geodesics were established 
in \cite{PoS3}, again by using the strong spectral estimates in \cite{D1}. Counting functions similar 
to $\pi(\lambda)$ however with some homological constraints were studied in \cite{An} and \cite{PoS2}.
For Anosov flows  with $C^1$  jointly non-integrable horocycle foliations  full asymptotic expansions for such counting functions 
were obtained.
In \cite{PeS2} Theorem 1.3 above was used to obtain results similar to these in \cite{PoS3} about correlations 
for pairs of closed billiard trajectories  for billiard flows in $\R^n \setminus K$,
where $K$ is a finite disjoint union of strictly convex compact bodies with smooth boundaries satisfying the so 
called `no eclipse condition'  (and some additional conditions as well). For the same kind of models and using 
Theorem 1.3 again, a rather non-trivial result was established in \cite{PeS1}  about analytic  continuation of  
the cut-off resolvent of the Dirichlet Laplacian in $\R^n\setminus K$, which appears to be the first of its kind in 
the field of quantum chaotic scattering. In \cite{PeS3}, 
using the spectral estimate in \cite{St2} and under the  assumptions there, a fine asymptotic was obtained  
for the number of closed trajectories in $M$ with primitive periods 
lying in exponentially shrinking intervals  $(x - e^{-\delta x}, x + e^{-\delta x})$,  $\delta > 0$, $x \to + \infty.$ 
In \cite{PeS4} a sharp large deviation principle was established concerning intervals shrinking with 
sub-exponential speed  for the Poincar\'e map related to a Markov family for an Axiom A flow on a basic set 
satisfying the assumptions in \cite{St2}. Finally, the spectral estimates in \cite{PeS5}, which we mentioned above,
were used to derive some interesting applications: (i) combining them with some arguments from \cite{PoS5},
we proved the so called Hannay-Ozorio de Almeida sum formula for the distribution of periods of periodic orbits of the 
flow for a rather general class of hyperbolic attractors; (ii) for Axiom A flows on basic sets $\Lambda$
satisfying the assumptions in \cite{St2} and for any H\"older continuous function $F: \Lambda \longrightarrow \R$
there exists $\ep > 0$ such that the counting function
$$\pi_F(T) = \sum_{\lambda(\gamma) \leq T} e^{\lambda_F(\gamma)} \quad , \quad
\lambda_F(\gamma) = \int_0^{\lambda(\gamma)} F(\phi_t(x_{\gamma})) dt,$$
where $\gamma$ is a primitive period orbit of the flow, $\lambda(\gamma)$ is the least period of $\gamma$,
and $x_{\gamma} \in \gamma$, has the asymptotic 
$$\pi_F(T) = \li(e^{Pr(F) T})(1 + {\mathcal O} (e^{-\epsilon T})) \quad ,\quad T \to \infty ,$$
where 
$$\li(x):= \int_2^{x} \frac{1}{\log y} dy \sim \frac{x}{\log x} \quad , \quad \: x \to +\infty .$$ 
This generalises a result of Pollicott \cite{Po2} for geodesic flows on compact manifolds of  negative $\frac{1}{4}$-pinched curvature.

\newpage

\subsection{Statement of results}

Let $\phi_t : M \longrightarrow M$ be a $C^2$ contact Anosov flow on a $C^2$ compact Riemannian manifold $M$.

Let $\varphi = \phi_1$ be the time-one map of the flow, and let $\m$ be an $\varphi$-invariant probability measure on $M$.
We refer the reader to the beginning of Sect. 3 below for the definition of a Pesin set with respect to $\m$.

A compact subset $P$ of $M$ will be called a {\it Pesin set with exponentially small tails} 
\index{Pesin set with exponentially small tails} with respect to $\m$
if $P$ is a Pesin set with respect to $\m$ and for every $\delta > 0$ there exist $C > 0$ and $c > 0$ such that
$$\m \left( \left\{ x\in \ll : \:\: \sharp \, \left\{ j : 0 \leq j \leq n-1 \: 
\mbox{\rm and } \:  \varphi^j(x) \notin P \right\} \geq \delta n \right\}\right) \leq C e^{- c n} ,$$
for all $n\geq 1$. The measure $\m$ will be called {\it regular} 
\index{regular measure} 
if it admits a Pesin set with exponentially small tails. See Sect. 3  for a sufficient condition for the
existence of Pesin sets with exponentially small tails. As explained below this sufficient condition is `generic' in a certain sense. 

Given $\alpha > 0$ denote by $C^{\alpha} (M)$ is {\it the space of all $\alpha$-H\"older 
complex-valued functions} on $M$, i.e. functions $h : M \longrightarrow \C$ for which there exists
$L \geq 0$ with  $|h(x) - h(y)| \leq L\, (d(x,y))^\alpha$ for all $x,y \in M$. For such $h$,  let $|h|_\alpha$ be the smallest possible
choice for $L$. Set $\|h\|_0 = \sup_{x\in M} |h(x)|$, and $\|h\|_\alpha = \|h\|_0 + |h|_\alpha$.

The main result in this work is the following.

\bs

\noindent
{\bf Theorem 1.1.} {\it Let $\phi_t : M \longrightarrow M$ be a $C^5$ contact Anosov flow,
let $F_0$ be a H\"older continuous function on $M$ and let $\m$ be the Gibbs measure
determined by $F_0$ on $M$. Assume that $\m$ is regular.
Then for every $\alpha > 0$ there exist constants $C = C(\alpha) > 0$ and $c = c(\alpha) > 0$ such that 
$$\left| \int_{M} A(x) B(\phi_t(x)) d\m (x) -  \left( \int_{M} A(x)  d\m (x)\right)\left(\int_{M} B(x)  d\m (x)\right)\right|
\leq C e^{-ct} \|A\|_\alpha  \|B\|_\alpha \;$$
for any two functions $A, B\in C^\alpha(M)$.}

\bs

We obtain this as a consequence of Theorem 1.3 below  and the procedure described in \cite{D1}.
The assumption that the flow is $C^5$ is made so that one can apply the procedure in \cite{D1}.
In particular this is essential when estimating the Laplace transform of the correlation function
$$\rho(t) = \int_{M} A(x) B(\phi_t(x))\; d\m (x) - 
\left( \int_{M} A(x)\; d\m (x)\right)\left(\int_{M} B(x) \; d\m (x)\right) $$
(see part VI in Sect. 4 in \cite{D1}).


As we mentioned earlier, it appears that so far the only results concerning exponential decay of correlations for 
general Gibbs potentials have been that of Dolgopyat \cite{D1} for  geodesic flows on
compact surfaces and the one in \cite{St2} for Axiom A flows on basic sets (under some additional  assumptions).

\index{$\zeta(s)$} 
Next, consider the {\it Ruelle zeta function} 
\index{Ruelle zeta function}
$$\zeta(s) = \prod_{\gamma} (1- e^{-s\ell(\gamma)})^{-1} \quad, \quad s\in \C ,$$
where $\gamma$ runs over the set of primitive  closed orbits of $\phi_t: M \longrightarrow M$
and $\ell(\gamma)$ is the least period of $\gamma$.  Denote by $h_T$ the  {\it topological entropy} 
of $\phi_t$ on $M$.
\index{$h_T$} 

Using Theorem 1.3 below  and an argument of Pollicott and Sharp \cite{PoS1}, one derives the 
following\footnote{Instead of using the norm $\|\cdot \|_{1,b}$ as in \cite{PoS1}, in the present 
case one has to work with $\|\cdot \|_{\theta,b}$ for some  $\theta \in (0,1)$, 
and then one has to use the so called Ruelle's Lemma in the form proved in  \cite{W}. 
This is enough to prove the estimate (2.3) 
for $\zeta(s)$ in \cite{PoS1}, and from there the arguments are the same.}.

\bs

\noindent
{\bf Theorem 1.2.} {\it Let $\phi_t : M \longrightarrow M$ be a $C^2$ contact Anosov flow on a $C^2$ compact  
Riemannian manifold $M$. Assume\footnote{It is still not proven that every Gibbs measure related to a contact 
Anosov flow is regular.} that there exists a Pesin set with exponentially small tails with respect to the
Sinai-Bowen-Ruelle measure\footnote{This is known to be true under some standard pinching conditions -- see
e.g. the comments at the end of Sect. 1 in \cite{GSt}. However  we expect that this condition should be satisfied in
much more general circumstances.}. Then:} 

\ms

(a) {\it The Ruele zeta function $\zeta(s)$ of the flow  $\phi_t: M \longrightarrow M$  has an analytic  and non-vanishing 
continuation in a half-plane $\Re(s) > c_0$ for some $c_0 < h_T$ except for a  simple pole at $s = h_T$.  }

\ms

(b) {\it There exists $c \in (0, h_T)$ such that
$$\di \pi(\lambda) = \# \{ \gamma : \ell(\gamma) \leq \lambda\} = \li(e^{h_T \lambda}) + O(e^{c\lambda})$$
as $\lambda\to \infty$, where 
$\di \li(x) = \int_2^x \frac{du}{\log u} \sim \frac{x}{\log x}$ as  $x \to \infty$. }

\bigskip

Parts (a) and (b) were first established by Pollicott and Sharp \cite{PoS1} for geodesic flows on 
compact surfaces of negative curvature (using \cite{D1}), and then similar results were proved in 
\cite{St2} for mixing Axiom A flows on basic sets satisfying certain additional assumptions 
(as mentioned above). Recently, using different methods, it was proved in \cite{GLP} that: 
(i) for volume preserving three dimensional Anosov flows (a) holds, and moreover, in the case of 
$C^\infty$ flows,  the Ruelle zeta function $\zeta(s)$ is meromorphic in $\C$ and  $\zeta(s) \neq 0$ for $\Re(s) > 0$; 
(ii) (b) holds for geodesic flows on $\frac{1}{9}$-pinched compact Riemannian manifolds of negative 
curvature. These were obtained as consequences of more general results in \cite{GLP}.

Let $\rr = \{R_i\}_{i=1}^{k_0}$ be a (pseudo-) Markov partition for $\phi_t$  consisting of 
rectangles $R_i = [U_i ,S_i ]$, where $U_i$ (resp. $S_i$) are (admissible) subsets of  $W^u_{\ep}(z_i)$
(resp. $W^s_{\ep}(z_i)$) for some $\ep > 0$ and $z_i\in M$ (cf. Sect. 2 for details). 
The first return time function $\tau : R = \cup_{i=1}^{k_0} R_i  \longrightarrow [0,\infty)$ 
is essentially $\alpha_1$-H\"older continuous on $R$  for some $\alpha_1 > 0$,
i.e. there exists a constant $L > 0$ such that if $x,y \in R_i \cap \pp^{-1}(R_j)$  for some $i,j$, 
where $\pp: R \longrightarrow R$ is the 
standard Poincar\'e map,  then $|\tau(x) - \tau(y)| \leq L\, (d(x,y))^{\alpha_1}$.
The {\it shift map} $\sigma : U = \cup_{i=1}^{k_0} U_i \longrightarrow U$ 
is defined by $\sigma = \piU\circ \pp$, where 
$\piU : R \longrightarrow U$ is the projection along the leaves of local stable manifolds.
Let $\hU$ be the set of all $x \in U$ whose orbits do not have common points with the boundary of $R$.
Given $\theta \in (0,1)$, recall the metric $\dte$ on $\hU$ from Sect. 1.1.
Denote by $\ff_\theta(\hU)$ {\it the space of all bounded functions $h : \hU \longrightarrow \C$ 
with Lipschitz constants }
$$|h|_\theta = \sup \left\{ \frac{|h(x) - h(y)|}{\dte(x,y)} : x\neq y; \ ; x,y \in \hU\right\} < \infty .$$
Define the norm $\|.\|_{\theta,b}$ on $\ff_\theta (\hU)$ by  
$$\| h\|_{\theta,b} = \|h\|_0 + \frac{|h|_{\theta}}{|b|} ,$$
where $\|h\|_0 = \sup_{x\in \hU} |h(x)|$.
\index{$\|h\|_0$}  
\index{$\|h\|_{\theta,b}$} 

Given a  real-valued function $f \in \ff_\theta (\hU)$, set $g = g_f = f - P_f\tau$, where   $P_f\in \R$ 
\index{$P_f$}
is the unique  number such that the topological pressure $\Pr_\sigma(g)$ of $g$ with respect to $\sigma$  is zero (cf. \cite{PP}). 

We say  that  {\it Ruelle transfer operators related to $f$ are eventually contracting }
\index{eventually contracting Ruelle operator}
on $\ff_\theta(\hU)$ if there exist constants $0 < \rho < 1$, $a_0 > 0$, 
$b_0 \geq 1$, $T \geq 1$ and  $C > 0$ such that if $a,b\in \R$  satisfy $|a| \leq a_0$ and $|b| \geq b_0$, 
then 
$$\|L_{f -(P_f+a+ \i b)\tau}^m h \|_{\theta,b} \leq C \;\rho^m \; \| h\|_{\theta,b}$$
for any integer $m \geq T \log |b|$ and any  $h\in \ff_\theta (\hU)$.

This condition implies that the spectral radius  of $L_{f-(P_f+ a+\i b)\tau}$ on 
$\ff_\theta (\hU)$  does not exceed  $\rho$. It is also easy to see that it implies the following\footnote{Which is the way we defined 
eventual contraction of Ruelle transfer operators in \cite{St2}, and it agrees with the way 
the main result in \cite{D1} is stated.}: for every $\epsilon > 0$ there exist  constants $0 < \rho < 1$, $a_0 > 0$, 
$b_0 \geq 1$ and  $C > 0$ such that if $a,b\in \R$  satisfy $|a| \leq a_0$ and $|b| \geq b_0$,  then
$$\|L_{f -(P_f+a+ \i b)\tau}^m h \|_{\theta,b} \leq C \;\rho^m \; |b|^{\ep}\; \| h\|_{\theta,b}$$
for any integer $m \geq 0$ and any  $h\in \ff_\theta (\hU)$.

The central result in this work is the following.

\bs

\noindent
{\bf Theorem 1.3.} {\it  Let $\phi_t : M \longrightarrow M$ be a $C^2$ contact Anosov flow 
on a $C^2$ compact  Riemannian manifold $M$, let  $\rr = \{R_i\}_{i=1}^{k_0}$ be a (pseudo-) Markov 
partition for $\phi_t$ as above and let $\sigma : U \longrightarrow U$ be the corresponding shift map. 
There exists a constant  $0 <  \hat{\theta} < 1$ such that for any  $\theta \in [\hat{\theta}, 1)$ and any real-valued function 
$f\in \ff_{\theta}(\hU)$ which is the restriction of a H\"older continuous function $F_0$ on $M$ so that 
the Gibbs measure $\nu_{F_0}$ is regular, the Ruelle transfer operators related to $f$ are eventually contracting on $\ff_\theta(\hU)$}.

\bs

Here $\hat{\theta}$ is the minimal number in $(0,1)$ such that the first-return time function
$\tau \in \ff_{\hat{\theta}}(\hU)$.

A similar result for H\"older continuous functions (with respect to the Riemannian metric) looks a 
bit more complicated, since in general Ruelle transfer operators do not preserve any of the spaces 
$C^\alpha(\hU)$. However, they preserve a certain `filtration' $\cup_{0 < \alpha \leq \alpha_0} \clip (\hU)$.
Here  $\alpha > 0$ and $C^{\alpha} (\hU)$ is {\it the space of all $\alpha$-H\"older  complex-valued functions} on $\hU$. 
Then $|h|_\alpha$ is the smallest non-negative number so that 
$|h(x) - h(y)| \leq |h|_\alpha (d(x,y))^\alpha$ for all $x,y \in \hU$.
Define the norm $\|.\|_{\alpha,b}$ on $\clip (\hU)$ by  $\| h\|_{\alpha,b} = \|h\|_0 + \frac{|h|_{\alpha}}{|b|}$.

\bs

\noindent
{\bf Corollary 1.4.} {\it  Under the assumptions of Theorem 1.3,  there exists a constant $\alpha_0 > 0$ 
such that for any real-valued function  $f \in C^{\alpha_0}(\hU)$ the  Ruelle transfer operators related to  $f$  are 
eventually contracting on $\cup_{0 < \alpha \leq \alpha_0} \clip (\hU)$. More precisely, there 
exists a constant $\hbeta\in (0,1]$
and for each $\epsilon > 0$ there exist constants $0 < \rho < 1$, $a_0 > 0$, $b_0 \geq 1$,
$C > 0$ and $M > 0$ such that if $a,b\in \R$  
satisfy $|a| \leq a_0$ and $|b| \geq b_0$, then for every integer $m  \geq M\, \log |b|$ and every 
$\alpha \in (0,\alpha_0]$ the operator
$$L_{f -(P_f+a+ \i b)\tau}^m : C^\alpha (\hU) \longrightarrow C^{\alpha \hbeta}  (\hU)$$ 
is well-defined and
$$\|L_{f -(P_f+a+ \i b)\tau}^m h \|_{\alpha\hbeta,b} \leq C \;\rho^m \; \| h\|_{\alpha,b}$$ 
for every  $h\in C^\alpha  (\hU)$.}

\ms

The maximal constant $\alpha_0 \in (0,1]$ that one can choose above (which is determined by 
the minimal $\hat{\theta}$ one can choose in Theorem 1.3) is 
related to the regularity of the local stable/unstable foliations. Estimates for this constant 
can be derived from certain bunching 
condition concerning the rates of expansion/contraction of the  flow along local unstable/stable 
manifolds (see \cite{Ha1},  \cite{Ha2}, \cite{PSW}). In the proof of Corollary 1.4 
in Sect. 8 below we give some rough estimate for $\alpha_0$.

The above  was first proved by Dolgopyat (\cite{D1}) in the case of geodesic flows on compact 
surfaces of negative curvature with $\alpha_0 = 1$ (then one can choose $\hbeta = 1$ as well).
The second main result in \cite{D1} concerns transitive Anosov flows  on compact Riemannian
manifolds with $C^1$ jointly  non-integrable local stable and unstable foliations. For such 
flows Dolgopyat proved that the conclusion of Corollary 1.4 with 
$\alpha_0 = 1$ holds for the Sinai-Bowen-Ruelle  potential  $F_0 = \log \det (d\phi_\tau)_{|E^u}$. 
More general results were proved in \cite{St2} for mixing Axiom A flows on basic sets
(again for $\alpha_0 = 1$) under some additional regularity assumptions. As mentioned earlier,
the results apply e.g. to  $C^2$ mixing Axiom A flows on basic sets satisfying a certain pinching condition (similar to the
$1/4$-pinching condition for geodesic flows on manifolds of negative curvature).

\subsection{Plan of the work and comments on the proof of Theorem 1.3}

 Unlike previous works on Ruelle transfer operators, {\bf here we make a heavy use of Pesin's theory of Lyapunov exponents}.
 Sections 2 and 3 contain some basic definitions and facts from hyperbolic dynamics and  Pesin's theory, respectively.

Let $F_0 : M \longrightarrow \R$ be a H\"older continuous function and let $\m$ be the 
Gibbs measure determined by $F_0$. Given a pseudo-Markov family $\rr = \{ R_i\}_{i=1}^{k_0}$  
for $\phi_t$ (see Sect. 2 for details), let 
$$\tau: R = \cup_{i=1}^{k_0} R_i \longrightarrow [0,1/2] \quad \mbox{\rm and} \quad \pp : R\longrightarrow R$$ 
be the corresponding first return map and the Poincar\'e map. The measure $\m$ induces a  Gibbs measure $\mu$ on $R$ 
(with respect to the Poincar\'e map $\pp$) for the function
$$F(x) = \int_0^{\tau(x)} F_0 (\phi_s(x))\, ds \quad, \quad x\in R .$$
The function $F$ is H\"older and, by Sinai's Lemma, is cohomologous to a H\"older function 
$f: R \longrightarrow \R$  which is constant on stable leaves in rectangles $R_i$ in $R$. Thus, $\mu$ 
coincides with the Gibbs measure on $R$ determined by $f$.

In Sect. 4 we state several lemmas concerning the non-integrability of the flow
due to the preservation of a contact form. Lemma 4.1 is Liverani's Lemma B.7 from \cite{L2} -- see below.
The rest of Sect. 4 has to do with Lyapunov exponents and all
statements involve a certain fixed Pesin set $P_0$ with exponentially small tails and fixed constants $\hep_0, \hd_0 > 0$.
Given an integer $m$, let $\Xi_m$ be the set of those Lyapunov regular points $x \in R$ such that
$\pp^j(x) \notin P_0$ `relatively frequently' for $0 \leq j < m$, more precisely
$$\# \{ j : 0 \leq j \leq m-1 \:, \: \pp^j(x) \notin P_0\} \geq \hd_0\, m .$$
By the choice of $P_0$, the sets $\Xi_m$ have exponentially small measure: $\mu(\Xi_m) \leq C\, e^{-cm}$
for some constants $C, c > 0$. Lemmas 4.2 - 4.4 deal with unstable cylinders $\cc$ of length $m$  in $R$ that have common points
with $P_0 \setminus \Xi_m$. For such cylinders $\cc$ we can estimate their diameters, or rather the diameter of {\it their projections
$\check{\cc}$  to a true unstable manifold}, by using powers of the smallest `unstable' Lyapunov exponent $\lambda_1 > 1$. 
We also get an important estimate for the temporal distance function:
$$|\Delta(x,y') - \Delta(x,y'')| \leq C\; \diam(\check{\cc})\, (d(y',y''))^\beta$$
for some constants $C, \beta > 0$ independent of $\cc$ and $m$, whenever $x,z \in \cc$, $y',y'' \in W^s_\ep(z)$. 
This is the contents of Lemma 4.2. Its proof is given in Sect. 9. 

The non-integrability Lemma 4.3 plays a very important  role in the proof of the main result.
It has to do with an important feature in the general construction of Dolgopyat's contraction operators,
although here the situation is significantly more complicated than the one in \cite{D1} (and also the one in \cite{St2}).
Given a large integer $N > > 1$ for cylinders $\cc$ of length $m$ with $\cc \cap P_0\setminus \Xi_m \neq \e$ we construct
families of points $y_1, y_2$ in $\pp^N(V) \cap W^s_\ep(Z)$, where $Z \in \cc \cap P_0\setminus \Xi_m$ and $V$ is a small
neighbourhood of $Z$ in $W^u_\ep(Z)$ such that 
$$c\; \diam(\check{\cc}) \leq |\Delta(x,\pi_{y_1}(z)) - \Delta(x,\pi_{y_2}(z)) |$$
for some constant $c > 0$ independent of the cylinder $\cc$ and its length $m$, for all `appropriately positioned' points $x,z \in \cc$. 
To prove Lemma 4.3 we use Liverani's Lemma 4.1 which says that there exist constants 
$C_0  > 0$, $\vartheta > 0$ and $\ep_0 > 0$ such that
for any $z\in M$, any $u \in E^u(z)$ and $v \in E^s(z)$ with $\|u\|, \|v\| \leq \ep_0$ we have
\be
|\Delta(\exp^u_z(u) , \exp^s_z(v))  - d\omega_z(u,v)| \leq C_0\,\left[ \|u\|^2\, \|v\|^\vartheta 
+ \|u\|^\vartheta \|v\|^2 \right] ,
\ee
where $d\omega_z$ is the symplectic form defined by the contact form on $M$ and $\Delta$ is 
the temporal distance function (see Sect. 2). We want to use this when $z\in P_0$, $v\neq 0$ is fixed and $\|u\|$ is small. 
Then however  the right-hand-side of (1.1) is only $O( \|u\|^\vartheta)$ which is not good enough. As Liverani suggests in 
Remark B.8 in \cite{L2}, one might be  able to improve the estimate pushing the points $x = \exp^u_z(u)$ and $y= \exp^s_z(v)$ 
forwards or backwards along the  flow\footnote{Liverani says that the best one can hope for is to get $o(\|u\|)$ in the right-hand-side 
of (1.1) if $\vartheta > \sqrt{3} -1$,  and this may be so, however we are interested in particular directions $u$ and for these without 
any restrictions on  $\vartheta$ we succeed to get a bit more.}. We go forwards roughly until 
$\|d\phi_t(z) \cdot u\| \geq \|d\phi_t(z) \cdot v\|$ for some $t > 0$. Moreover we  are only interested in directions $u \in E^u_1(z)$, where
$E^u_1$ is the sub-bundle of $E^u$ corresponding to the smallest `unstable' Lyapunov  exponent $\lambda_1 > 1$.
The proof of Lemma 4.3 is given in  Sect. 9. As we mentioned earlier, it is rather technical, lengthy and non-trivial.

As in \cite{St2}, here we work a lot with cylinders defined by the Markov family.  In \cite{St2} we worked under a 
certain regularity  assumption --  the so called {\it regular distortion along unstable manifolds}.
A range of examples of flows having this property was described in \cite{St5}.  However it seems unlikely that it holds for 
every (contact) Anosov flow. Lemma 4.4 states that unstable cylinders $\cc$ of length $m$ with $\cc \cap P_0 \neq \e$
and $\pp^m(\cc) \cap P_0 \neq \e$ have properties similar to the properties considered in \cite{St2}.
Its proof is given in Sect. 10. Although it is technical and takes quite a bit of space and effort, it is using ideas similar to these 
in \cite{St5}, so, in some sense, it cannot be regarded as something  that requires a significant intellectual effort.

Sect. 5 contains the main application of Lemma 4.3, namely Lemma 5.4 which provides estimates from below of 
differences of temporal distances in a form convenient for the estimates of contraction operators in Sect. 6. 

As we have mentioned already, in the proof of Theorem 1.3 we use the general framework of 
Dolgopyat's method from \cite{D1} and its development in \cite{St2}. 
As in \cite{D1}, we deal with the normalized operators 
$$ \lab = L_{\fa - \i\,b\tau} ,$$
where
$$\fa(u) = f (u) - (P_f+a) \tau(u) + \ln h_{a}(u) -  \ln h_{a}(\sigma(u)) - \ln \lambda_{a} ,$$
$\lambda_a > 0$ being the largest eigenvalue of $L_{f-(P_f+a)\tau}$, and $h_a$ a particular 
corresponding positive eigenfunction (see Sect. 5.1 for details). 
Their real parts  $\ma = L_{\fa}$ satisfy $\ma  1 = 1$. Now instead of dealing with these 
operators on some $C^\alpha(U)$, we consider them on the space $\ff_\theta(\hU)$ of 
$\dte$-Lipschitz functions on $\hU$. Recall that for $x \neq y \in U$, $\dte (x,y) = \theta^{N}$ if $x,y$ belong to a 
cylinder of length $N-1$ and $N$ is maximal with this property. 
We choose $\theta \in (0,1)$ so that $\tau \in \ff_\theta(\hU)$.
The main benefit in working with $\dte$ is that the local stable holonomy maps are 
isometries in this metric, and every $h \in \ff_\theta(\hU)$ can be considered as a function 
on $R$ which is constant on stable leaves. Such $h$ has the same `trace' on each 
$W^u_R(x)$, so we have the freedom to choose whichever unstable leaf is more convenient to work on.

Another simple thing that helps to avoid the lack of regularity is to approximate the partition  
$\rr = \{ R_i\}_{i=1}^{k_0}$ (which we call a {\it pseudo-Markov partition} below) by a true 
(at least according to the standard definition, see \cite{B}) Markov
partition $\{ \tR_i\}_{i=1}^{k_0}$, where each $\tR_i$ is contained in a submanifold 
$D_i$ of $M$ of codimension one.  We can take $D_i$ so that $U_i \cup S_i \subset D_i$. 
The shift along the flow determines a bi-H\"older continuous bijection  
$$\tPsi : R \longrightarrow \tR =  \cup_{i=1}^{k_0} \tR_i ,$$
and whenever we need to
measure the `size' of a cylinder $\cc$ lying in some $W^u_R(x)$ we use 
$\diam(\tPsi(\cc))$, instead of $\diam(\cc)$. Since the Poincar\'e map 
$$\tpp: \tR =  \cup_{i=1}^{k_0} \tR_i \longrightarrow \tR$$
is essentially Lipschitz, estimates involving $\diam(\tPsi(\cc))$
are much nicer. Roughly speaking, whenever we deal with Ruelle operators, measures and 
Gibbsian properties of measures, we work on $R$, and whenever we have to estimate distances and 
diameters in the Riemannian metric we use projections\footnote{In fact, sometimes it is more 
convenient to use  the projections $\Psi_i: R_i \longrightarrow \check{R}_i = \cup_{z\in S_i} \chU_i(z)$, 
where $\chU_i(z)$ is the part of the true unstable manifold $W^u_{\ep_0}(z)$ corresponding to 
$W^u_{\tR_i}(z)$ via the shift along the flow. This is particularly convenient 
when using Liverani's Lemma 4.1.} in $\tR$.

As in \cite{D1}, the main result would follow if we show that, given $f\in \ff_\theta(\hU)$, 
there exist constants $C > 0$, $\rho \in (0,1)$, $a_0 > 0$,
$b_0 \geq 1$ and an integer $N \geq 1$ such that for $a, b\in \R$ with $|a|\leq a_0$ 
and $|b|\geq b_0$ and any $h\in \ff_\theta(\hU)$  we have
\be
\di\int_{U} |\lab^{mN} h|^2\; d\nu \leq C\, \rho^m \, \|h \|_{\theta,b} 
\ee
for every positive integer $m$. Here $\nu$ is the Gibbs measure on $U$  determined by $g = f - P_f \, \tau$, which
is naturally related to the Gibbs ($\pp$-invariant) measure $\mu$ on $R$.

In order to prove (1.2) in \cite{D1}, Dolgopyat constructs  for any choice of $a$ and $b$, a family of (what we call below) 
{\it contraction operators}  $\nn_J$, and the proof of (1.2) goes as follows. Given $h$ with (in his case) $\|h\|_{1,b} \leq 1$, 
define $h^{(m)} = \lab^{Nm} h$,  $H^{(0)} = 1$ and $H^{(m)} =  \nn_{J_m}(H^{(m-1)})$ for an appropriately chosen sequence 
of contraction operators $\nn_{J_m}$, so that  $|h^{(m)}| \leq H^{(m)}$ for all $m$. In \cite{D1} (and \cite{St2}) the contraction 
operators indeed contract in the $L^1$ norm so that
\be
\int_U (\nn_{J} H)^2 \, d\nu \leq \rho \, \int_U H^2 \, d\nu
\ee
for some constant $\rho\in (0,1)$ independent of $a$, $b$, $J$ and $H$. Thus,
\be
\int_{U} |\lab^{mN} h|^2\; d\nu  = \int_U |h^{(m)}|^2 \, d\nu \leq  \int_U (H^{(m)})^2 \, d\nu \leq \rho^m .
\ee

In the present work our contraction operators do not satisfy (1.3). Moreover we cannot deal 
immediately with functions  $f \in \ff_\theta(\hU)$; instead we fix a sufficiently small $\theta_1 \in (0,\theta)$ 
and assume initially  $f\in \ff_{\theta_1}(\hU)$. The general case is dealt with using approximations.  

Even with $\theta$ replaced by a very small $\theta_1$, we cannot prove an analogue of (1.2); instead we establish that 
for every $s > 0$ there  exist integers $N \geq 1$, $k = k(N) \geq 1$ and a constant $C = C (N)> 0$ such that\footnote{The role
of the integer $N$ is special here, since it is involved in the operators $\nn_J$; we need another constant $k \geq 1$ 
to get the action of multiple contraction operators.}
\be
\int_{U} |\lab^{k N \log |b|} h| \; d\nu \leq  \frac{C}{|b|^{s}} \,\|h\|_{\theta_1,b} 
\ee
for all $h \in \ff_{\theta_1}(U)$. 
From this it follows by a relatively standard procedure (see Ch. 8) that for any $h \in \ff_{\theta_1}(\hU)$ we have 
\be
\|\lab^{k N\, \log |b|} h\|_{\theta_1,b} \leq \frac{C}{|b|^{s}}\, \|h\|_{\theta_1,b} 
\ee
for some (possibly different) constants $C > 0$, $k \geq 1$, depending on $N$ but not on $a$ and $b$.
The estimate in Theorem 1.3 is derived from (1.6) using another `standard procedure'.

Here we come to the most central part of this work --  the construction of the contraction operators and the study 
of their main properties\footnote{We believe that the scope of applicability of the arguments
developed in Sections 6-7 is significantly wider than what is actually stated as results in this paper.} in Sections 6 and 7.
Having fixed a Pesin set $P_0$ with exponentially small tails, for any sufficiently large value of the parameter $b$
that appears in $ \lab = L_{\fa - \i\,b\tau}$ we construct several objects that depend on $b$. First, we construct in a special
way a `sufficiently representative' compact subset\footnote{A different notation is used in Sect. 6; in fact
the whole construction is more complicated than what we say here.} $P'_0$ of $P_0$ so that 
$\mu(P_0\setminus P'_0) \leq \frac{C}{|b|^c}$ for some constants $C, c > 0$ independent of $b$ and set $K_0 = \piU(P'_0)$, 
where $\piU : R \longrightarrow U$ is the projection along stable leaves in $R$. Then we choose a special essentially 
disjoint{\footnote{They can only have common points at their boundaries.} family of
unstable maximal cylinders $\cc_1, \ldots, \cc_{m_0}$ such that their projections $\cc'_m = \piU(\cc_m)$ cover $K_0$
and their {\it projections $\tcc_m = \tPsi(\cc_m)$ along the flow to $\tR$} satisfy $\diam(\tcc_m) \leq \frac{C_1}{|b|}$ 
for some constant $C_1 > 0$. 

Next, given a large integer $N$, independent of $b$, we  use Lemma 5.4
(consequence of Lemma 4.3) and construct a family of pairs of local inverses of the map $\sigma^N$: 
$$\vli = \vli(Z_m, \cdot) : \cc'_m = \piU(\cc_m) \longrightarrow U  \quad , \quad \ell = 1, \ldots, \ell_0 \: , \: i = 1,2 ,$$
where $Z_m$ is fixed point  in $\cc_m \cap P'_0$ and 
$$\sigma^N(\vli(Z_m, x)) = x \quad , \quad x \in \cc'_m .$$
Then define `contraction' operators $\nn_J(a,b)$ for a family of symbols $J$, small $|a|$ and large $|b|$, very much as in \cite{St2}.
However proving some kind of contraction properties of these operators is non-trivial. To achieve this in Sect. 6.2 we 
define a special metric $\dd$ on $\hU$ (depending on $b$ and the cylinders $\cc_m$) as follows. For any 
$u, u'\in \hU$, let $\ell(u,u') \geq 0$ be the length of the  smallest cylinder $Y(u,u')$ in $\hU$ containing $u$ and $u'$. Set
$\dd(u,u') = 0$ if $u = u'$. If $u \neq u'$ and  there exists  $p \geq 0$ with 
$$\sigma^p(Y(u,u')) \subset \cc'_m = \piU(\cc_m)$$
and $\ell(u,u') \geq p$ for some $m = 1, \ldots, m_0$, take the maximal $p$ with this property and the corresponding $m$ and set 
$$\dd(u,u') = \frac{\dte(u,u')}{\diamtef (\cc_m)} ,$$
where $\diamte (\cc_m)$ is the diameter of $\cc_m$ with respect to the
metric $D_\theta$. Finally, if no $p$ exists as in the previous sentence, set $\dd(u,u') = 1$. Then,
choosing appropriately a large constant $E > 0$, let $\kk_E$ be  the set of all functions $H \in \ff_\theta(\hU)$ 
such that  $H > 0$ on $\hU$ and 
$$\frac{|H(u) - H(u')|}{H(u')} \leq E\, \dd(u,u')$$ 
for all $u,u'\in \hU$ for which there exists an integer $p \geq 0$ with $\sigma^p(Y(u, u')) \subset \cc_m$ 
for some $m \leq m_0$ and $\ell(u,u') \geq p$. It turns out that $\nn_J(\kk_E) \subset \kk_{E}$ for any
of the contraction operators $\nn_J$ defined earlier, and this turns out to be very important for the estimates that follow.


Next, consider the following assumption for points $u,u' \in \hU$ contained in some cylinder 
$\cc'_m$ ($1\leq m \leq m_0$), an integer $p \geq 0$ and points $v,v'\in \hU$:
\be
u,u'\in \cc'_m \:  , \: \sigma^p(v) = \vli(u)\;, \; \sigma^p(v')  = \vli(u')\;,\; \ell(v,v') \geq N  ,
\ee
for some $i = 1,2$.
Assume $f \in \ff_{\theta_1}(\hU)$, and denote by $\kk_b$  the family  of all pairs $(h, H)$ such that $h \in \ff_{\theta}(\hU)$, 
$H \in \kk_E$, $|h| \leq H$ on $\hU$, and for any $u,u' \in \hU$ contained in a cylinder $\cc_m$ for some $m = 1, \ldots, m_0$, 
any integer  $p \geq 0$ and any points $v,v'\in \hU_1$
satisfying (1.7)  for some $i = 1,2$ and $\ell = 1, \ldots, \ell_0$ we have
$$|h(v) - h(v') | \leq E \, |b|\, \theta_2^{p+N} \, H(v') \, \diam(\tcc_m) .$$
Here $\theta_2 \in (\theta,1)$ is an appropriately chosen constant sufficiently close to $1$.

We then succeed to derive that some kind of cancellation occurs in the actions of the operators $\nn_J$.
More precisely, we prove that for small $|a|$ and large $|b|$, defining the cylinders $\cc_m$ and the metric $\dd$ as above,
for any $(h,H) \in \kk_b$  there exists a contraction operator $\nn_J$ such that $(\lab^{N} h, \nn_J H ) \in \kk_b$.

Perhaps the cancellation we have just mentioned looks a bit tiny and achieved under very special conditions, however it turns out
it is enough to prove (1.5) and then (1.6). This is done in Sect. 7.   Theorem 1.3 is derived from these in Sect. 8 using some
`standard procedures'. It seems the latter are difficult to find in the literature,
so in Sect. 8 we include detailed arguments on how exactly these procedures are done.

\bs

\noindent
{\bf Acknowledgements.} Thanks are due to Boris Hasselblatt, Dima Dolgopyat, Sebastian Gou\"ezel,
Carlangelo Liverani, Yakov Pesin, Vesselin Petkov, Mark Pollicott  and Amie Wilkinson who
provided me with various useful information and/or made valuable comments on this paper at various stages during the
work on it. Special thanks are due to Sebastian Gou\"ezel who pointed out some errors in the initial version of the paper. Finally,
thanks are due to the referee for his/her careful and thorough reading of the initial manuscript and numerous remarks and
suggestions that helped to improve the paper significantly. 


\newpage

\section{Preliminaries}
\setcounter{equation}{0}

Throughout this paper $M$ denotes a $C^2$ compact Riemannian manifold,  and  $\phi_t : M \longrightarrow M$ 
\index{$\phi_t$} ($t\in \R$) a  $C^2$  Anosov flow 
\index{Anosov flow} 
on $M$. That is, there exist   constants $C > 0$ and $0 < \lambda < 1$ such that there exists a $d\phi_t$-invariant decomposition  
$$T_xM = E^0(x) \oplus E^u(x) \oplus E^s(x)$$ 
of $T_xM$ ($x \in M$) 
\index{$E^u(x)$, $E^s(x)$} 
into a direct  sum of non-zero linear subspaces, where $E^0(x)$ is the one-dimensional subspace determined by the direction of the flow at $x$, 
$$\| d\phi_t(u)\| \leq C\, \lambda^t\, \|u\| \quad, \quad u\in E^s(x)\:, \:\: t\geq 0 ,$$
and
$$\| d\phi_t(u)\| \leq C\, \lambda^{-t}\, \|u\| \quad , \quad u\in E^u(x) \: , \:\: t\leq 0 .$$
Here $\|\cdot \|$ is the {\it norm determined by the Riemannian metric on $M$.}

For $x\in M$ and a sufficiently small $\epsilon > 0$ let 
$$\wloc^s(x) = \{ y\in M : d (\phi_t(x),\phi_t(y)) \leq \epsilon \: \mbox{\rm for all } \: t \geq 0 \; , \: d (\phi_t(x),\phi_t(y)) \to_{t\to \infty} 0\: \}\; ,$$
$$\wloc^u(x) = \{ y\in M : d (\phi_t(x),\phi_t(y)) \leq \epsilon \: \mbox{\rm for all } \: t \leq 0 \; , \: d (\phi_t(x),\phi_t(y)) \to_{t\to -\infty} 0\: \}$$
be the (strong) {\it stable} and {\it unstable manifolds} 
\index{strong stable/unstable manifolds} 
of size $\epsilon$.  Then $E^u(x) = T_x \wloc^u(x)$ and $E^s(x) = T_x \wloc^s(x)$.  Given $\delta > 0$, set 
$$E^u(x;\delta) = \{ u\in E^u(x) : \|u\| \leq \delta\} ;$$
$E^s(x;\delta)$ is defined similarly. For every $x \in M$ we will denote by
$$\exp^u_x : E^u(x;\ep) \longrightarrow \wloc^u(x) \quad \mbox{\rm and }\quad \exp^s_x : E^s(x;\ep) \longrightarrow \wloc^s(x)$$
the {\it exponential maps} acting on unstable and stable fibbers, respectively.

It follows from the hyperbolicity of the flow on $M$  that if  $\epsilon_0 > 0$ is sufficiently small,
there exists $\ep_1 > 0$ such that if $x,y\in M$ and $d (x,y) < \ep_1$, 
then $W^s_{\ep_0}(x)$ and $\phi_{[-\ep_0,\ep_0]}(W^u_{\ep_0}(y))$ intersect at exactly  one point $[x,y]$ 
\index{$[x,y]$}  
(cf. \cite{KH}). That is, there exists a unique  $t\in [-\ep_0, \ep_0]$ such that $\phi_t([x,y]) \in W^u_{\ep_0}(y)$. Setting $\Delta(x,y) = t$, 
\index{$\Delta(x,y)$} 
defines the so called {\it temporal distance
function\footnote{In fact in \cite{D1} and \cite{L1} a different definition for $\Delta$ is given, however
in the important case (the only one considered below) when $x\in W^u_\ep(z)$ and 
$y \in W^s_\ep(z)$ for some $z \in M$,  these definitions coincide with  the present one.}} 
\index{temporal distance function} 
(\cite{KB},\cite{D1}, \cite{Ch1}, \cite{L1}). For $x, y\in M$ with $d (x,y) < \ep_1$, define
$$\pi_y(x) = [x,y] = W^s_{\ep}(x) \cap \phi_{[-\ep_0,\ep_0]} (W^u_{\ep_0}(y)) .$$
Thus, for a fixed $y \in M$, $\pi_y : W \longrightarrow \phi_{[-\ep_0,\ep_0]} (W^u_{\ep_0}(y))$ is the
{\it projection} along local stable manifolds defined on a small open neighbourhood $W$ of $y$ in $M$.
Choosing $\ep_1 \in (0,\ep_0)$ sufficiently small,  the restriction 
$$\pi_y: \phi_{[-\ep_1,\ep_1]} (W^u_{\ep_1}(x)) \longrightarrow \phi_{[-\ep_0,\ep_0]} (W^u_{\ep_0}(y))$$
is called a {\it local stable holonomy map\footnote{In a similar way one can define
holonomy maps between any two sufficiently close local transversals to stable laminations; see e.g. \cite{PSW}.}.} 
\index{local holonomy map} 
Combining such a map with a shift along the flow we get another  local stable holonomy  map 
$$\hh_x^y : W^u_{\ep_1}(x)  \longrightarrow W^u_{\ep_0}(y) .$$
In a similar way one defines local holonomy maps along unstable laminations.
\index{$\pi_y$}
\index{$\hh_x^y$}
\index{$W^u(x)$, $W^s(x)$}

We will say that $A$ is an {\it admissible subset} 
\index{admissible set} 
of $W^u_{\ep}(z)$
if $A$ coincides with the closure of its interior in $W^u_\ep(z)$. Admissible subsets of  $W^s_\ep(z)$ are defined similarly.

Let $D$ be a submanifold of $M$ of codimension one such that $\diam(D) \leq \ep$ and
$D$ is transverse to the flow $\phi_t$. Assuming that $\ep > 0$ is sufficiently small, the projection 
$$\pr_D : \phi_{[-\ep,\ep]}(D) \longrightarrow D$$ 
 \index{$\pr_D$} 
 along the flow is well-defined and smooth. Given $x,y\in D$, set $\la x, y\ra_D = \pr_D([x,y])$. 
A subset $\tR$ of $D$ is called a {\it rectangle} 
\index{rectangle} 
if $\la x, y\ra_D \in \tR$ for all $x,y\in \tR$.  The rectangle $\tR$ is called {\it proper} if $\tR$ coincides with the closure of its
interior in $D$. For any $x\in \tR$ define the stable and unstable leaves through $x$ in $\tR$ by 
$$W^s_{\tR}(x) = \pr_D(W^s_\ep(x)\cap \phi_{[-\ep,\ep]}(D)) \cap \tR$$ 
and
$$W^u_{\tR}(x) = \pr_D(W^u_\ep(x)\cap \phi_{[-\ep,\ep]}(D)) \cap \tR .$$
For a subset $A$ of $D$ we will denote by $\Int_D(A)$
 \index{$\Int_D(A)$} 
 the {\it interior}  of $A$ in $D$. 

Let $\tRR = \{ \tR_i\}_{i=1}^{k_0}$ 
\index{$\tRR$, $\tR_i$} 
be a family of proper rectangles, where each $\tR_i$ is contained
in a submanifold $D_i$ of $M$ of codimension one. We may assume that each $\tR_i$ has the form 
$$\tR_i = \la U_i  , S_i \ra_{D_i} = \{ \la x,y\ra_{D_i} : x\in U_i, y\in S_i\}\;,$$
where $U_i \subset \wloc^u(z_i)$ and $S_i \subset \wloc^s(z_i)$, respectively, 
for some $z_i\in M$.  Moreover, we can take $D_i$ so that $U_i \cup S_i \subset D_i$. 
Set 
$$\tR =  \cup_{i=1}^{k_0} \tR_i .$$
We will denote by $\Int(\tR_i)$ the {\it interior} of the set $\tR_i$  in the topology of the disk $D_i$.
The family $\tRR$ is called {\it complete}
\index{complete family of rectangles} 
if  there exists $\chi > 0$ such that for every $x \in M$, $\phi_{t}(x) \in \tR$ for some  $t \in (0,\chi]$.  The {\it Poincar\'e map}
$$\tpp: \tR \longrightarrow \tR$$
related to a complete family $\tRR$ is defined by $\tpp(x) = \phi_{\ttau(x)}(x) \in \tR$, where
$\ttau(x) > 0$ is the smallest positive time with $\phi_{\ttau(x)}(x) \in \tR$. The function $\ttau$ 
\index{$\ttau$} 
is called the {\it first return time}  associated with $\tRR$. 
A complete family $\tRR = \{ \tR_i\}_{i=1}^{k_0}$ of rectangles in $M$ is called a  {\it Markov family} 
\index{Markov family} 
of size $\chi > 0$ for the  flow $\phi_t$ if: 
                                         
(a) $\diam(\tR_i) < \chi$ for all $i$; 

(b)  for any $i\neq j$ and any $x\in \Int_{D}(\tR_i) \cap \tpp^{-1}(\Int_{D}(\tR_j))$ we have   
$$W_{\tR_i}^s(x) \subset \overline{ \tpp^{-1}(W_{\tR_j}^s(\tpp(x)))} \quad , \quad                                                         
\overline{\tpp(W_{\tR_i}^u(x))} \supset W_{\tR_j}^u(\tpp(x)) ;$$

(c) for any $i\neq j$ at least one of the sets $\tR_i \cap \phi_{[0,\chi]}(\tR_j)$ and $\tR_j \cap \phi_{[0,\chi]}(\tR_i)$ is empty. 

The existence of a Markov family $\tRR$ of an arbitrarily small size $\chi > 0$ for $\phi_t$ follows from the construction of Bowen \cite{B}.
 \index{$\tpp$}
 \index{Poincar\'e map}

Following  \cite{R2} and \cite{D1}, we will now slightly change the Markov family $\tRR$ to a  {\it pseudo-Markov family} 
\index{pseudo-Markov family}  
$\rr = \{ R_i\}_{i=1}^{k_0}$ 
\index{$\rr$} 
of {\it pseudo-rectangles} 
\index{pseudo-rectangles}
$$R_i = [U_i  , S_i ] =  \{ [x,y] : x\in U_i, y\in S_i\}\,.$$
where $U_i$ and $S_i$ are as above. Set 
$$R =  \cup_{i=1}^{k_0} R_i .$$
Notice that  $\pr_{D_i} (R_i) = \tR_i$ for all $i$. Given $\xi = [x,y] \in R_i$, set 
$$W^u_{R}(\xi) = W^u_{R_i}(\xi) = [U,y] = \{ [x',y] : x'\in U_i\}$$ 
and
$$W^s_{R}(\xi) = W^s_{R_i}(\xi) = [x,S_i] = \{[x,y'] : y'\in S_i\} \subset W^s_{\ep_0}(x) .$$
The corresponding {\it Poincar\'e map} 
$$\pp: R \longrightarrow R$$ is defined by  
$\pp(x) = \phi_{\tau(x)}(x) \in R$, where $\tau(x) > 0$ 
\index{$\tau(x)$}  
is the smallest positive time with $\phi_{\tau(x)}(x) \in R$.  The function $\tau$  is the {\it first return time}  
\index{first return time functions} 
associated with $\rr$. 
The {\it interior} $\Int(R_i)$ of a rectangle $R_i$ is defined by  $\pr_D(\Int(R_i)) = \Int_D(\tR_i)$. In a similar way one can define
$\Intu(A)$ for a subset $A$ of some  $W^u_{R_i}(x)$ and $\Ints(A)$ for a subset $A$ of some  $W^s_{R_i}(x)$.
 \index{$\pp$}

We may and will assume that the family $\rr = \{ R_i\}_{i=1}^{k_0}$ has the same properties as $\tRR$, namely: 

$(a')$ $\diam(R_i) < \chi$ for all $i$; 

$(b')$ for any $i\neq j$ and any  $x\in \Int(R_i) \cap \pp^{-1}(\Int(R_j))$ we have  
$$\pp(\Int (W_{R_i}^s(x)) ) \subset \Ints (W_{R_j}^s(\pp(x))) \quad , \quad   \pp(\Int(W_{R_i}^u(x))) \supset \Int(W_{R_j}^u(\pp(x))) ;$$

$(c')$ for any $i\neq j$ at least one of the sets $R_i \cap \phi_{[0,\chi]}(R_j)$ and  $R_j \cap \phi_{[0,\chi]}(R_i)$ is empty. 

Define the matrix $\aa = (\aa_{ij})_{i,j=1}^k$ 
\index{$\aa$}  
by $\aa_{ij} = 1$ if $\pp(\Int (R_i)) \cap \Int(R_j) \neq  \e$ 
and $\aa_{ij} = 0$ otherwise.  According to \cite{BR} (see section 2 there), we may assume that $\rr$ 
is chosen in such a way that $\aa^{M_0} > 0$ (all entries of the $M_0$-fold product of $\aa$  by itself are 
positive) for some integer $M_0 > 0$. 
\index{$M_0$} In what follows we assume that the matrix $\aa$ has this property.

Notice that in general $\pp$ and $\tau$ are only (essentially) H\"older continuous. However there is an
obvious relationship between $\pp$ and the (essentially) Lipschitz map $\tpp$, and this will be used below.

From now on we will assume that $\tRR = \{ \tR_i\}_{i=1}^{k_0}$ is a fixed Markov family for  
$\phi_t$ of size $\chi < \ep_0/2 < 1$ and that  $\rr = \{ R_i\}_{i=1}^{k_0}$ is the related pseudo-Markov family. Set  
$$U = \cup_{i=1}^{k_0} U_i $$
and $\Intu (U) = \cup_{j=1}^{k_0} \Intu(U_j)$. 
\index{shift map} 
\index{$U$, $U_i$}
\index{$\piU$}
\index{projection along stable leaves}
\index{$W^s_R(x)$}

The {\it shift map} $\sigma : U   \longrightarrow U$ is given by $\sigma  = \piU \circ \pp$, where 
$$\piU : R \longrightarrow U$$
is the {\it projection} along stable leaves.  Notice that  $\tau$ is constant on each stable leaf $W_{R_i}^s(x) = W^s_{\ep_0}(x) \cap R_i$. 
For any integer $m \geq 1$ and any function $h : U \longrightarrow \C$ define $h_m : U \longrightarrow \C$ by
$$h_m(u) = h(u) + h(\sigma(u)) + \ldots + h(\sigma^{m-1}(u)) .$$
\index{$h_m$}

Denote by $\widehat{U}$ (or $\hR$) 
\index{$\hU$, $\hR$}
the {\it core} 
\index{core} 
of  $U$ (resp. $R$), i.e. the set  of those $x\in U$ (resp. $x \in R$)
such that  $\pp^m(x) \in \Int(R) = \cup_{i=1}^k \Int(R_i)$  for all $m \in \Z$.  It is well-known (see \cite{B}) that $\hU$ is a residual subset 
of $U$ (resp. $R$) and has full measure with respect to any Gibbs measure on $U$ (resp. $R$).
Clearly in general $\tau$ is not continuous on $U$, however $\tau$ is {\it essentially H\"older}  on $U$, i.e.
there exist constants $L>0$ and $\alpha >0$ such that $|\tau(x) - \tau(y)| \leq L\, (d(x,y))^\alpha$
whenever $x,y \in U_i$ and $\sigma(x), \sigma(y)\in U_j$ for some $i, j$.
The same applies to $\sigma : U \longrightarrow U$.  Throughout we will mainly 
work with the restrictions of $\tau$ and $\sigma$ to $\hU$. Set $\hU_i = U_i \cap \hU$.
For any $A \subset M$, let $\hA $  be {\it the set of all $x \in A$ whose trajectories 
do not pass through boundary points} of $R$.
\index{essentially H\"older continuous function}
 \index{$\|g\|_0$}
\index{$L_g$}
\index{Ruelle transfer operator}

Let $B(\hU)$ be the {\it space of  bounded functions} $g : \hU \longrightarrow \C$ with its standard norm  $\|g\|_0 = \sup_{x\in \hU} |g(x)|$.
Given a function $g \in B(\hU)$, the  {\it Ruelle transfer operator }  $L_g : B(\hU) \longrightarrow B(\hU)$  is defined by 
$$\di (L_gh)(u) = \sum_{\sigma(v) = u} e^{g(v)} h(v) .$$
Given $\alpha > 0$, let $\clip(\hU)$ 
\index{$\clip(\hU)$}
denote the {\it space of the essentially $\alpha$-H\"older continuous functions}
$h :\hU \longrightarrow \C$, i.e. such that there exists $L \geq 0$ with  $|h(x) -  h(y) | \leq L \, (d(x,y))^\alpha$
for all $i = 1, \ldots, k_0$ and all $x,y\in \hU_i$. The smallest $L > 0$ 
with this property is called the {\it $\alpha$-H\"older exponent} of $h$ and is denoted $|h|_\alpha$.   Set 
$$\|h\|_\alpha = \|h\|_0 + |h|_\alpha .$$

The hyperbolicity of the flow  implies the existence of constants $c_0 \in (0,1]$ and $\gamma_1 > \gamma > 1$ such that
\be
c_0 \gamma^m\; d (x,y) \leq  d (\tpp^m(x)), \tpp^m(y)) \leq \frac{\gamma_1^m}{c_0} d (x,y)
\ee
for all $x,y\in \tR$ such that $\tpp^j(x), \tpp^j(y)$ belong to the same $\tR_{i_j}$ for all $j = 0,1, \ldots, m$. 

Throughout this paper $\alpha_1 \in (0,1]$ 
\index{$\alpha_1$}
will denote the {\it largest constant } such that  $\tau\in C^{\alpha_1}(\hU)$ and the local stable/unstable
holonomy maps are uniformly $\alpha_1$-H\"older. We will also need to fix a constant $\talpha_1 \in (0,1)$
(take e.g. the largest again)  such that the projection $\tPsi : R \longrightarrow \tR$ along stable leaves is $\talpha_1$-H\"older. 
\index{$\talpha_1$}
\index{$\|h\|_\alpha$}
\index{$|h|_\alpha$}
\index{$c_0$} 
\index{$\gamma$} 
\index{$\gamma_1$}

 
\section{Lyapunov exponents and Lyapunov regularity functions}
\setcounter{equation}{0}

Let $M$ be a $C^2$ Riemannian manifold, and let $\phi_t$ be a $C^2$ Anosov  flow on $M$. Set $\varphi = \phi_1$. 
\index{$\varphi$}
Let $F_0$ 
\index{$F_0$}be a H\"older continuous real-valued function on $M$ and let $\m$ 
\index{$\m$}
be the Gibbs measure  generated by $F_0$ on $M$. It follows from the well-known {\it Oseledets Multiplicative Ergodic Theorem} (\cite{Os}) 
that in the situation considered here  there exists a $\phi_t$-invariant subset $\ll$ of $M$ with $\m(\ll) = 1$ such that for every  $x \in \ll$ 
\index{$\ll$} there exists numbers
$$0 < \chi_1  < \chi_2 < \ldots < \chi_{\tk}$$
and a $d\phi_t$-invariant decomposition  
$$E^u(x) = E^u_1(x) \oplus E^u_2(x) \oplus \ldots \oplus E^u_{\tk}(x)$$
of $E^u(x)$ into subspaces of constant dimensions such that
\begin{equation}
\lim_{t\to \infty} \frac{1}{t}\, \log \|d \phi_t(x)\cdot v\| = \chi_i \quad , \quad v \in E_i(x)\setminus \{0\} ,
\end{equation}
for all $i = 1, \ldots,\tk$ (see \cite{BP}, \cite{PS} or \cite{V}).  

As in Sect. 2,  here and in what follows we denote by $\|\cdot \|$ the {\it norm on the tangent spaces
$T_xM$ ($x\in M$) induced by the Riemannian metric} on $M$.

The points of $\ll$ are called {\it Lyapunov regular  points}, 
\index{set of Lyapunov regular points} while the numbers
$\chi_i > 0$ are called (the positive) {\it Lyapunov exponents} 
\index{Lyapunov exponents} of $\phi_t$.
As we stated above, the dimension $n_i$ of $E^u_i(x)$ is constant on $\ll$. Clearly 
$$n_1 + n_2 + \ldots + n_{\tk} = n^u = \dim(E^u(x))$$
for all $x\in \ll$.
We have a similar decomposition for $E^s(x)$, $x\in \ll$, and corresponding negative Lyapunov exponents. If the flow is contact, we have 
$n^s = \dim(E^s(x))  = n^u$ for all $x \in \ll$, and the negative Lyapunov exponents are $-\chi_i$, $i = 1, \ldots, \tk$.
\index{$\tk$}  
\index{$n_i$} 
\index{$n^u$, $n^s$} 
\index{$E^u_i(x)$}
\index{$\chi_i$}

Given $\ep > 0$, it follows from (3.1) that for every $x\in \ll$ and every $i = 1, \ldots, \tk$ we have
$$\lim_{n\to\infty} \frac{\|d\varphi^n(x)_{|E_i(x)}\|}{e^{(\chi_i+\ep) n}} = 0 ,$$
therefore
$$G(x) = \max_{1\leq i \leq \tk} \: \sup_{n \geq 0} \frac{\|d\varphi^n(x)_{|E_i(x)}\|}{e^{(\chi_i+\ep) n}} < \infty$$
The function $G(x)$ just defined is an example of a Lyapunov $\ep$-regularity function.

More generally, a Borel function $R: \ll \longrightarrow (1,\infty)$ such that
\be
\frac{1}{R (x)} \leq \frac{\|d\varphi^n(x)\cdot v\|}{e^{(\chi_i + \ep)n}\, \|v\|} 
\leq R (x) \quad , \quad x\in \ll \;, \; v\in E^u_i(x)\setminus \{0\} \;, \; n \geq 0 ,
\ee
for all $i = 1, \ldots, \tk$, and
\be
e^{-\ep} \leq \frac{R (\varphi(x))}{R (x)} \leq  e^{\ep} \quad , \quad x\in \ll ,
\ee
is called a {\it Lyapunov $\epsilon$-regularity function}.  
In particular, a positive Borel function $R$ on $\ll$ satisfying (3.3) is called a {\it slowly varying function}.
\index{slowly varying function}
\index{Lyapunov regularity function}

Notice that the sum of two Lyapunov $\ep$-regularity functions is also a Lyapunov $\ep$-regularity function.
Also, if $R_i$ is a  Lyapunov $\ep_i$-regularity function for all $i = 1, \ldots , k$, then $R_1 \ldots R_k$ is a 
Lyapunov $(\ep_1+ \ldots + \ep_k)$-regularity function.

In what follows we will frequently use the notation $R_\ep$  just
to indicate that this is a Lyapunov $\ep$-regularity function, however  we have to emphasise that in different situations
this notation will be used for (possibly) different functions. Following \cite{PS}, by an {\it $\ep$-slowly varying radius function}
 \index{$\ep$-slowly varying radius function}
we mean a function of the form $r_\ep(x) = 1/R_\ep(x)$, $x \in \ll$, where $R_\ep$ is a Lyapunov $\ep$-regularity function on $\ll$.
For such a function $r_\ep$ and every $x \in \ll$, on the open ball with centre $0$ and
radius $r_\ep(x)$ in $E^u(x)$, the linear maps $d\varphi^n (x)$ behave as in the case of an uniformly hyperbolic flow --
see the relations (3.12) - (3.15) below. Again, in what follows the notation $r_\ep$ will be `generic', i.e. in different situations 
it may denote different $\ep$-slowly varying radius functions. 

Notice that if $r_1 (x) = 1/R_1(x)$ and $r_2(x) = 1/R_2(x)$ are $\ep$-slowly varying radius functions, then 
$r(x) = 1/(R_1(x) + R_2(x))$ is also an $\ep$-slowly varying radius function and $r(x) \leq \min\{ r_1(x), r_2(x) \}$ for all $x \in \ll$.
As in the remark  above about Lyapunov regularity functions, if $r_i$ is an $\ep_i$-slowly varying radius function for all $i = 1, \ldots, k$,
then $r_1 \ldots r_k$ is an $(\ep_1 + \ldots+\ep_k)$-slowly varying radius function.

A subset $P$ of $\ll$ is called a {\it Pesin set} 
\index{Pesin set} if there exist constants $\ep > 0$ and $C > 0$ and a Lyapunov
$\ep$-regularity function $R_\ep(x)$ on $\ll$ such that $R_\ep(x) \leq C$ for all $x \in P$.
We will say that the {\it Pesin set $P$  has exponentially small tails}\footnote{A more straightforward definition would require 
the existence of constants $C > 0$ and $c > 0$ such that  $\m (M \setminus \cup_{j=0}^{n-1} \varphi^{-j}(P))\leq C e^{- c n}$ for all integers 
$n \geq 1$. However in this work  we need a stronger condition.}
 if for every $\delta > 0$ there exist constants $C > 0$ and $c > 0$ such that
 $$\m \left( \left\{ x\in \ll : \:\: \sharp \, \left\{ j : 0 \leq j \leq n-1 \: 
\mbox{\rm and } \:  \varphi^j(x) \notin P \right\} \geq \delta n \right\}\right) \leq C e^{- c n}  $$ 
for all integers $n \geq 1$. The latter says that for $x\in \ll$ outside a set with exponentially small measure,  `most' points $\pp^j(x)$ on 
the orbit $x, \pp(x), \ldots, \pp^{n-1}(x)$ (more precisely, at least $(1-\delta)n$ of them) belong to the Pesin set $P$. 
Consequently, for such points $\pp^j(x)$ we have a relatively good control over the norms $\|d\varphi^n(\pp^j(x))_{|E^u(\pp^j(x))}\|$.

\bs

Set  $\lambda_i = e^{\chi_i}$ 
\index{$\lambda_i$} 
for all $ i = 1, \ldots,\tk$. {\bf Fix an arbitrary constant} $\beta \in (0,1]$ 
\index{$\beta$} such that 
$$\lambda_j^\beta < \lambda_{j+1} \quad , \quad 1 \leq j < \tk .$$
Take $\hep > 0$ 
\index{$\hep$} so small that 
\begin{eqnarray}
e^{8\hep} < \lambda_1  \:\: , \:\: e^{8\hep} < \lambda_{j}/\lambda_{j-1} \:\: (j = 2, \ldots,\tk) .
\end{eqnarray}
Some further assumptions about $\hep$ will be made later.
Set
\be
1 < \nu_0 = \lambda_1 e^{-8\hep} <  \mu_{j} = \lambda_{j} e^{-\hep} < \lambda_{j} < \nu_{j} = \lambda_{j} e^{\hep}
\ee
for all $j = 1, \ldots,\tk$. 
\index{$\nu_i$}  
\index{$\mu_i$}

{\bf Fix $\hep > 0$ with the above properties}.

\index{$\hWuj(x)$}
\index{$\tWuj(x)$}   
\index{$\Phi_x^u$} 
\index{$\hE^u_j(x)$}  
\index{$\tE^u_j(x)$}

For $x \in \ll$ and $1 \leq j \leq \tk$ set 
$$\hE^u_j(x) = E^u_1(x) \oplus \ldots \oplus E^u_{j-1}(x) \quad , \quad \tE^u_{j}  = E^u_j(x) \oplus \ldots \oplus E^u_{\tk}(x) .$$
Also set $\hE^u_1(x) = \{0\}$ and $\hE^u_{\tk+1}(x) = E^u(x)$.
For any $x\in \ll$ and any $u\in E^u(x)$ we will write $u = (\uo,\ut , \ldots, \utk)$,  where $\ui \in E^u_i(x)$ for all $i$.

It follows from the general theory of non-uniform hyperbolicity (see \cite{P1}, \cite{BP}) 
that for any $j = 1, \ldots,\tk$ the invariant bundle  $\{\tE^u_{j}(x)\}_{x\in \ll}$ is uniquely integrable over $\ll$, i.e. 
there exists a measurable $\varphi$-invariant family $\{ \tWuj_{\tr(x)}(x)\}_{x\in \ll}$  
of $C^2$ submanifolds $\tWuj(x) = \tWuj_{\tr(x)}(x)$ of $M$ tangent to the bundle $\tE^u_{j}$ for some 
 {\it $\hep$-slowly varying radius function} $\tr = \tr_{\hep} : \ll \longrightarrow (0,1)$.
Moreover, with $\beta \in (0,1]$ as in the beginning of this section, for $j > 1$ it follows from Theorem 6.6 in 
\cite{PS} and (3.4) that there exists a $\varphi$-invariant family $\{ \hWuj_{\tr(x)}(x)\}_{x\in \ll}$ of
$C^{1+\beta}$ submanifolds $\hWuj(x) = \hWuj_{\tr(x)}(x)$ 
of $M$ tangent to the bundle $\hE^u_j$.
(However this family is not unique in general.) For each $x\in \ll$ and each $j = 2, \ldots, \tk$ fix
an $\varphi$-invariant family $\{ \hWuj_{\tr(x)}(x)\}_{x\in \ll}$ with the latter properties.
Then we can find an {\it $\hep$-slowly varying radius function}  $r = r_{\hep}: \ll \longrightarrow (0,1)$ 
and for any $x\in \ll$ a $C^{1+\beta}$ diffeomorphism
$$\Phi_x^u: E^u (x; r(x)) \longrightarrow \Phi_x (E^u (x; r (x)) \subset W^u_{\tr (x)}(x)$$
such that 
\be
\Phi^u_x(\hE^u_j(x; r(x))) \subset \hWuj_{\tr(x)}(x)\quad , \quad
\Phi^u_x(\tE^u_{j}(x; r(x))) \subset \tWuj_{\tr (x)}(x) 
\ee
for all $x \in \ll$ and $j = 2, \ldots , \tk$. Moreover, since for each $j > 1$ the submanifolds
$ \hWuj_{r(x)}(x)$ and $\exp^u_x(\hE^u_j(x;r(x)))$ of $W^u_{\tr(x)}(x)$ are tangent at $x$ of order 
$1+\beta$, we can choose  $\Phi^u_x$ so that the diffeomorphism
$$\Psi^u_x =  (\exp^u_x)^{-1} \circ \Phi^u_x   : E^u(x: r(x)) \longrightarrow \Psi^u_x (E^u(x: r(x))) \subset E^u(x; \tr(x))$$ 
is $C^{1+\beta}$-close to identity. Thus, we can choose a Lyapunov $\hep$-regularity function $R(x)$ such that
\be
\| \Psi^u_x(u) - u\| \leq R(x) \|u\|^{1+\beta} \quad , \quad \| (\Psi^u)^{-1}_x(u) - u\| \leq R(x) \|u\|^{1+\beta}
\ee
for all $ x\in \ll$ and  $u \in E^u(x; r(x))$, and also  that 
\be
\|d \Phi^u_x(u)\| \leq R(x) \quad ,\quad  \|(d\Phi^u_x(u))^{-1}\| \leq R(x) \quad,  
\quad x\in \ll \:, \: u \in E^u (x ; r(x)) .
\ee
Finally, replacing $R(x)$ with a larger Lyapunov $\hep$-regularity function if necessary (e.g take the 
sum of the initial one with another appropriately chosen Lyapunov $\hep$-regularity function\footnote{As remarked earlier, the sum of two
Lyapunov $\hep$-regularity functions is also a Lyapunov $\hep$-regularity function.}), we may assume that
\be
\quad\|\Phi^u_x(v) - \Phi^u_x(u) - d\Phi^u_x(u) \cdot (v-u)\| 
\leq R(x) \, \|v-u\|^{1+\beta} \:\:,\:  x\in \ll\:, \: u,v \in E^u (x ; r(x)) ,
\ee
and
\be
\|d\Phi^u_x(u) - \id\| \leq R (x)\, \|u\|^\beta \quad , \quad x\in \ll \:, \: u \in E^u (x; r(x)) .
\ee
In a similar way one defines the maps $\Phi^s_x$ and we will assume that $r(x)$ is chosen so that 
these maps satisfy the analogues of the above properties.
\index{$\Psi^u_x$}
\index{$\hf_x$} 
\index{$\hf_y^j$}
\index{$\Gamma = \Gamma_{\hep}$}  
\index{$r = r_{\hep}$} 
\index{$\|\cdot\|'_x$}

For any $x\in \ll$ consider the $C^{1+\beta}$ map (defined locally near $0$)
$$\hf_x = (\Phi^u_{\varphi(x)})^{-1} \circ \varphi \circ \Phi^u_x : E^u (x)  \longrightarrow E^u (\varphi(x))\;.$$
It is important to notice that
$$\hf^{-1}_x(\hE^u_j (\varphi(x) ; r(\varphi(x))) \subset \hE^u_j (x; r(x)) \quad , \quad 
\hf^{-1}_x(\tE^u_j (\varphi (x) ; r(\varphi (x))) \subset \tE^u_j (x; r(x)) \;$$
for all $x \in \ll$ and $j > 1$.

Given $y \in \ll$ and any integer $j \geq 1$ we will use the notation
$$\hf_y^j = \hf_{\varphi^{j-1}(y)} \circ \ldots \circ \hf_{\varphi (y)} \circ \hf_y\quad,
\quad \hf_y^{-j} = (\hf_{\varphi^{-j}(y)})^{-1} \circ \ldots \circ (\hf_{\varphi^{-2}(y)})^{-1}  \circ (\hf_{\varphi^{-1}(y)})^{-1} \;,$$
at any point where these sequences of maps are well-defined.

It is well known (see  e.g. the Appendix in \cite{LY1} or Sect. 3 in \cite{PS})  that there exists a 
Lyapunov $\hep$-regularity function $\Gamma = \Gamma_{\hep} : \ll \longrightarrow [1,\infty)$ and an $\hep$-slowly varying
radius function $r = r_{\hep} : \ll \longrightarrow (0,1)$ (we will assume that it is the same as the one chosen above)
and  for each $x\in \ll$ a norm $\| \cdot \|'_x$  on $T_xM$ such that
\be
\|v\| \leq \|v\|'_x \leq \Gamma (x) \|v\| \quad ,\quad x\in \ll\:,\: v \in T_xM ,
\ee
$$\angle(\hE^u_j (x), \tE^u_{j}(x)) \geq \frac{1}{\Gamma(x)} \quad , \quad x\in \ll\;, \; 2 \leq j \leq d ,$$
and for any $x \in \ll$ and any integer $m \geq 0$, assuming 
$\hf_x^j(u), \hf_x^j(v) \in E^u(\varphi^j(x), r (\varphi^j(x)))$ are well-defined for all  $j =1, \ldots,m$, the following hold:
\be
\mu^m_j\, \|u-v\|'_x  \leq \|\hf_x^m(u) - \hf_x^m(v)\|'_{\varphi^m(x)} \quad , \quad u,v \in \tE^u_j(x; r(x)) ,
\ee
\be
\mu^m_1\, \|u-v\|'_x  \leq \|\hf_x^m(u) - \hf_x^m(v)\|'_{\varphi^m(x)}  \quad , \quad u,v \in E^u(x; r(x))  ,
\ee
\be 
\mu_1^m \, \|v\|'_x \leq \|d\hf_x^m(u)\cdot v\|'_{\varphi^m(x)} \leq \nu_{\tk}^m \, \|v\|'_x
\quad , \quad u \in E^u(x;r(x))\;,\; v \in E^u(x) ,
\ee
\be 
\mu_j^m \, \|v\|'_x \leq \|d\hf_x^m(0)\cdot v\|'_{\varphi^m(x)} \leq \nu_j^m \, \|v\|'_x
\quad , \quad v \in E^u_j(x) .
\ee

We will also use the norm  
$$|u| = \max\{ \|\ui\| : 1 \leq i \leq \tk\} .$$ 
Clearly, 
$$\|u\| = \| u^{(1)} + \ldots + u^{(\tk)}\| \leq \sum_{i=1}^{\tk} \|\ui\| \leq \tk \, |u| .$$
Taking the regularity function $\Gamma(x)$ appropriately (as explained above),   we have $|u|\leq \Gamma(x) \|u\|$, so
\be
\frac{1}{\tk}\, \|u\| \leq |u| \leq \Gamma (x) \|u\| \quad , \quad x\in \ll\;,\; u\in E^u(x) .
\ee
\index{$|u|$} 

Next, Taylor's formula (see also Sect. 3 in \cite{PS}) implies that there exists 
a Lyapunov $\hep$-regularity function  $D = D_{\hep} : \ll \longrightarrow  [1,\infty)$ such that  for any $i = \pm 1$ 
and any $x \in \ll$ we have
\be
\quad\|\hf_x^i(v) - \hf_x^i(u) - d\hf^i_x(u) \cdot (v-u)\| \leq D(x) \, \|v-u\|^{1+\beta} \:\:, \: u,v \in E^u (x ; r(x)) ,
\ee
and
\be
\|d\hf^i_x(u) - d\hf^i_x(0)\| \leq D (x)\, \|u\|^\beta \quad , \quad  u \in E^u (x; r(x)) .
\ee

Finally, we state here a Lemma from \cite{St4} which will be used several times later.

\bs

\noindent
{\bf Lemma 3.1.} (Lemma 3.3 in \cite{St4}) {\it There exist a Lyapunov $\hep$-regularity function
$$L = L_{\hep}: \ll \longrightarrow [1,\infty)$$
and an $\hep$-slowly varying radius function 
$$r: \ll \longrightarrow (0,1)$$
such that for any $x\in \ll$, any integer $p \geq 1$ and any $v\in E^u(z, r (z))$ with  $\|\hf_z^p(v)\|\leq r(x)$,  where $z = \varphi^{-p}(x)$, we have
$$\|\wo_p - \vo_p\| \leq L(x) |v_p|^{1+\beta} ,$$
where $v_p = \hf^p_z(v) \in E^u (x)$ and $w_p = d\hf_z^p(0)\cdot v \in E^u (x)$. Moreover, if $|v_p| = \|\vo_p\| \neq 0$, then}  
$1/2 \leq \|\wo_p\|/\|\vo_p\| \leq 2$.

\index{$L$} 

\bs

\noindent
{\bf Remark.} Notice that if $v\in E^u_1(z, r (z))$ in the above lemma, then $v_p, w_p \in E^u_1(x)$, so $\|w_p - v_p \| \leq L(x)\, \|v_p\|^{1+\beta}$.

\ms

\index{$Q_p(\ep)$} 
\index{$\mu$} 
\index{$\tau_0$, $\htau_0$}

Let $F_0 : M \longrightarrow \R$ be a H\"older continuous functions as in Sect. 1.3 and let $\m$ be the  {\it Gibbs measure} 
\index{Gibbs measure} determined by $F_0$ (see Sect. 5.1 for the definition).
Let $\rr = \{ R_i\}_{i=1}^{k_0}$  be a pseudo-Markov family for $\phi_t$ as in Ch. 2, and let
$$\tau: R = \cup_{i=1}^{k_0} R_i \longrightarrow [0,1/2]$$ 
and $\pp : R\longrightarrow R$ be the 
corresponding first return time and the Poincar\'e map. As before fix constants  $0 < \tau_0 < \htau_0 \leq 1/2$ so that 
$$\tau_0 \leq \tau(x) \leq \htau_0 \quad , \quad x\in R .$$
The Gibbs measure $\m$ induces a {\it Gibbs measure} 
\index{Gibbs measure} $\mu$ 
on $R$ (with respect to the Poincar\'e map $\pp$) for the function
$$F(x) = \int_0^{\tau(x)} F_0 (\phi_s(x))\, ds \quad, \quad x\in R .$$
The function $F$ is H\"older and, using Sinai's Lemma, it is cohomologous to a {\it H\"older function $f: R \longrightarrow \R$  
\index{$f$} 
which is constant on stable leaves in rectangles} $R_i$ in $R$. Thus, $\mu$ 
coincides with the Gibbs measure determined by $f$.
For every continuous function $H$ on $M$ we then have (see e.g. \cite{PP})
\be
\int_M H\, d\m = \frac{\int_R \left( \int_0^{\tau(x)} H(\phi_s(x))\, ds\right) d\mu(x)}{\int_R \tau\, d\mu } .
\ee

Let $R_\ep$ be  a Lyapunov $\ep$-regularity function on $\ll$. Then any set of the form
$$Q_{p}(\ep) = \{ x\in \ll : R_{\ep}(x) \leq e^{p} \} $$
is a Pesin set.  Given  $p > 0$, $\ep > 0$, $\delta  > 0$ and an integer $n \geq 1$,  set
\be
\Xi_n(p,\ep,\delta) = \left\{ x\in \ll \cap R: \:\: \sharp \, 
\left\{ j : 0 \leq j \leq n-1 \: \mbox{\rm and } \:  \pp^j(x) \notin Q_{p}(\ep) \right\} 
\geq \delta\, n \: \: \right\}.
\ee
As we remarked earlier in the definition of a Pesin set with exponentially small tails,
when $x \in \Xi_n(p,\ep,\delta)$ `most' points $\pp^j(x)$ on the orbit $x, \pp(x), \ldots, \pp^{n-1}(x)$ (more precisely,
at least $(1-\delta)n$ of them) belong to the Pesin set $Q_p(\ep)$. Consequently, for such points $\pp^j(x)$ we have relatively good
control over the norms $\|d\varphi^n(\pp^j(x))_{|E^u(\pp^j(x))}\|$. This will be significantly used in the estimates of temporal distances in Sect. 5.3.

\index{$\Xi_n = \Xi_n(p, \ep, \delta)$}

\bs

\noindent
{\bf Definition.} (\cite{GSt})
Consider a log-integrable linear cocycle $M$ abovea transformation $(T,\mu)$, with Lyapunov exponents $\lambda_1 \geq \ldots \geq \lambda_d$. 
We say that $M$ has {\it exponential large deviations for all exponents}
\index{exponential large deviations for all exponents} if, for
any $i\leq d$ and any $\epsilon>0$, there exists $C>0$ such that, for all $n\geq 0$,
\begin{equation}
\label{eq:exp_dev_all_exp}
  \mu\{  x : |\log \|\Lambda^i M^n(x)\| - n (\lambda_1+\ldots+\lambda_i)| \geq n \epsilon\} \leq C e^{- n/C}.
\end{equation}

The following theorem, which is a special case of Theorem 1.7 in \cite{GSt},  shows that if $d \varphi$ has exponential large deviations 
for all exponents, then most points in $\ll$ return exponentially often to some Pesin set. 

\bs

\noindent
{\bf Theorem 3.2.} (\cite{GSt})
{\it  Assume that $d \pp$ has exponential large deviations for all exponents with respect to $\mu$.
Let $\hep > 0$ and $\hd_0 > 0$. Then there exist constants $p_0 > 0$, $C > 0$ and $c  > 0$ such that
\begin{equation}\label{eq:tailestimate}
\mu (\Xi_n(p_0,\hep,\hd_0))\leq C e^{- c n}
\end{equation}
for all $n \geq 1$. }

\bs

Clearly, if (3.22) holds for $p_0$, then it will hold with $p_0$ replaced by any $p \geq p_0$.
The relation (3.22) shows that, as a subset of $\ll\cap R$ and with respect to the Poincar\'e map $\pp$, the set $Q_p(\hep)$ is 
a Pesin set with exponentially small tails.

As established in \cite{GSt} (see Theorem 1.5 there), for a transitive subshift of finite type
$T$ on a space $\Sigma$, if $\mu$ is a Gibbs measure for a H\"older-continuous
potential and $M$ is a continuous linear cocycle on a vector bundle $E$
above $T$, each of the following conditions is sufficient for $M$ to have exponential large deviations 
for all exponents: 

(i) if all its Lyapunov exponents coincide;

(ii) if there is a continuous decomposition of $E$ as a direct sum of
    subbundles $E=E_1 \oplus \ldots \oplus E_k$ which is invariant under
    $M$, such that the restriction of $M$ to each $E_i$ has exponential
    large deviations for all exponents;
    
(iii) more generally, if there is an invariant continuous flag
    decomposition 
    $$\{0\}=F_0 \subseteq F_1 \subseteq \ldots \subseteq F_k = E$$
    such that the cocycle induced by $M$ on each $F_i/F_{i+1}$ has
    exponential large deviations for all exponents;
    
(iv) if the cocycle $M$ is locally constant in some trivialization of the
    bundle $E$ (this is equivalent to the existence of invariant continuous
    holonomies which are commuting);
    
(v) if the cocycle $M$ admits invariant continuous holonomies, and if it
    is pinching and twisting in the sense of Avila-Viana \cite{AV};
    
(vi) if the cocycle $M$ admits invariant continuous holonomies, and the
    bundle is $2$-dimensional.

It follows from the above and Theorem 9.18 in \cite{V} that generic linear cocycles have exponential large 
deviations for all exponents. Moreover,  amongst fiber bunched cocycles\footnote{Which are the most frequently 
met H\"older cocycles.}, generic cocycles in the H\"older topology also have exponential large deviations for all exponents.

\bs

In what follows sometimes it will be more convenient to use sets of the form
$$\hXi_k(p,\ep,\delta) = \left\{ x\in \ll : \:\: \sharp \, \left\{ s \in \Z : 0 \leq s < k \: \mbox{\rm and } \:  \varphi^s(x) \notin Q_{p}(\ep) \right\} 
\geq \delta\, k \: \: \right\}.$$
These sets are naturally related to sets of the form $\Xi_n(p,\ep,\delta)$. More
precisely, using the fact that $\tau_0 \leq \tau(x) \leq \htau_0$ for all $x \in R$, it is easy to check that
\begin{equation}
\Xi_m(p_0,\hep,\hd_0) \subset \hXi_{m\htau_0} (p_0 - \hep, \hep, \hd_0/(\tau_0 \htau_0) )
\end{equation}
and
\begin{equation}
\hXi_{m \htau_0} (p_0 + \hep ,\hep,\hd_0/\tau_0) \subset \Xi_{m} (p_0 , \hep, \hd_0)
\end{equation}
for all integers $m \geq 1$. Indeed, to check (3.23), let $z_0 \in \Xi_m(p_0,\hep,\hd_0)$.
Then 
$$\sharp \,  \left\{ j : 0 \leq j \leq m-1 \: , \:  R_{\hep}(\pp^j(z_0)) > e^{p_0} \right\}  \geq \hd_0\, m .$$
For every $j < m$ with $ R_{\hep}(\pp^j(x)) > e^{p_0}$, let $s = [\tau_j(z_0)]$. Then $s \leq j \htau_0 < m\htau_0$
and 
$$R_{\hep}(\varphi^s(z_0)) \geq e^{-\hep} R_{\hep}(\pp^j(z_0)) > e^{-\hep} e^{p_0} = e^{p_0- \hep} .$$
The number of $j$'s with $j < m$ and $s = [\tau_j(z_0)]$ is not more than $1/\tau_0$, so
\begin{eqnarray*}
&         &\sharp \,  \left\{ s : 0 \leq s < m \htau_0 \: , \:  R_{\hep}(\varphi^s(z_0)) > e^{p_0- \hep} \right\}\\
& \geq &  \frac{1}{\tau_0} \; \sharp \,  \left\{ j : 0 \leq j \leq m-1 \: , \:  R_{\hep}(\pp^j(z_0)) > e^{p_0} \right\}  
\geq \frac{\hd_0}{\tau_0}\, m  = \frac{\hd_0}{\tau_0 \htau_0}\, m \htau_0 ,
\end{eqnarray*}
which shows that $z_0$ belongs to the right-hand-side of (3.23). 

In a similar way one proves (3.24).

\index{$\hXi_k\hXi_k(p,\ep,\delta)$}

\newpage

\section{Non-integrability of contact Anosov flows}
\setcounter{equation}{0}

The main aim of this section is to present some kind of a non-integrability property of the flow
due to the preservation of a contact form. An important role here is played by Liverani's Lemma B.7 from \cite{L2} 
which we state as Lemma 4.1. We then choose appropriate Lyapunov $\hep$-regularity functions and an $\hep$- slowly varying radius
function, and define a specifically designed Pesin set $P_0$  with exponentially small tails.
We then state three important, although technical, lemmas -- Lemmas 4.2 - 4.4. All of them deal with unstable cylinders $\cc$ 
 in $R$ that have common points with a specially defined `large part' of $P_0$.
For such cylinders $\cc$ we estimate the diameters of their projections $\tcc$  to true unstable manifolds and certain type of
increments along local stable manifolds of related temporal distance functions (Lemma 4.2).
Lemma 4.3 is a non-integrability lemma which plays a fundamental role in the proof of the central result Theorem 1.3.
Both Lemmas 4.2 and 4.3 are proved in Sect. 9 using Liverani's Lemma 4.1 and the choice of the Pesin set $P_0$.

Lemma 4.4 states that unstable cylinders $\cc$ of length $m$ with $\cc \cap P_0 \neq \e$
and $\pp^m(\cc) \cap P_0 \neq \e$ have properties similar to  the so called {\it regular distortion along unstable manifolds}
considered in \cite{St2}. These properties are used many times in the following chapters. The proof of this lemma is given in Sect. 10. 

Throughout this section we assume that the Lyapunov $\hep$-regularity functions $R(x)$, $\Gamma(x)$, $D(x)$ and $L(x)$,
the $\hep$-slowly varying radius function $r(x)$ and the $\hep/2$-slowly varying radius function $\tr(x)$ are chosen as in Sect. 3, 
so that they satisfy (3.6) -- (3.18)  and the conclusions of  Lemma 3.1.

\subsection{Choice of constants, sets of Lyapunov regular points} 

In what follows we assume that $\tRR = \{ \tR_i\}_{i=1}^{k_0}$  is a fixed Markov partition for 
$\phi_t$ on $M$ of size $< 1/2$ and $\rr = \{ R_i\}_{i=1}^{k_0}$  is the related pseudo-Markov 
family as in Sect. 2. We will use the notation associated with these from Sect. 2, and
we will assume that for any $i = 1, \ldots,k_0$, $z_i$ is chosen so that $z_i \in \Intu(W^u_{R_i}(z_i))$.
For any $x\in R$, any $y \in \tR$ and $\delta > 0$ set 
$$B^u(x,\delta) = \{ y\in W^u_{R_i}(x) : d(x,y) < \delta\}  \quad , \quad \tB^u(y,\delta) = \{ z\in W^u_{\tR_i}(z) : d(z,y) < \delta\} .$$
In a similar way define $B^s(x,\delta)$. For brevity sometimes we will use the notation 
$$U_i(z) = W^u_{R_i}(z)$$ 
for $z\in  R_i .$ 

Fix constants $0 < \tau_0 < \htau_0 < 1$ so that 
$$\tau_0 \leq \tau(x) \leq \htau_0 \quad, \quad x\in R ,$$ 
and 
$$\tau_0 \leq \ttau(x) \leq \htau_0 \quad, \quad x\in \tR .$$

Let $\alpha_1 > 0$ be as in Sect. 2, and let $f$ be an essentially $\alpha_1$-H\"older continuous potential on $R$.  Set $g = f - P_f \tau$, 
where $P_f \in \R$ is chosen so that the topological pressure of $g$ with respect to the Poincar\'e map 
$\pp : R = \cup_{i=1}^{k_0} R_i \longrightarrow R$  is $0$.
Let $\mu = \mu_g$ be the Gibbs measure on $R$ determined by $g$; then $\mu(R) = 1$.
We will assume that $f$ (and therefore $g$) depends on forward coordinates only, i.e. it is constant on stable leaves of  $R_i$ for each $i$. 
\index{$B_u(x,\delta)$} 
\index{$\tB^u(y,\delta) $} 
\index{$g$}  
\index{$P_f$} 
\index{$C_[\ii]$} 
\index{$C_W[\ii]$} 

\bs

\noindent
{\bf Remark.} If the initial potential $F$ on $R$ is $\alpha^2$-H\"older, applying Sinai's Lemma (see e.g. Proposition 1.2 in \cite{PP}) produces 
an $\alpha$-H\"older potential $f$ depending on forward coordinates only. More precisely, there exists an $\alpha$-H\"older
function $h$ so that $F = f + h - h\circ \pp$ (i.e. $F$ is co-homological to $f$). Moreover, $F$ and $f$ have the same topological pressure
and the same equilibrium state (see e.g. Proposition 3.6 in \cite{PP}).

\bs

Given an unstable leaf $W = W^u_{\tR_i}(z)$ in some rectangle $\tR_i$ 
and an admissible sequence $\ii = i_0, \ldots,i_m$ of integers $i_j \in \{ 1,\ldots, k_0\}$, the set
$$C_W[\ii] = \{ x\in W : \tpp^j(x) \in \tR_{i_j} \:, \: j = 0,1,\ldots,m\}$$
will be called a {\it cylinder of length} $m$ in $W$ (or an {\it unstable cylinder} in $\tR$ in general). 
When $W = U_i$ we will simply write $C[\ii]$. 
In a similar way one defines cylinders $C_V[\ii]$, 
where $V = W^u_{R_i}(z)$ is an unstable leaf in some rectangle $R_i$. 

Let
$$\pr_D : \cup_{i=1}^{k_0} \phi_{[-\ep,\ep]}(D_i) \longrightarrow \cup_{i=1}^{k_0} D_i$$
be the {\it projection along the flow}, i.e. for all $i = 1,\ldots,k_0$ and all 
$x\in \phi_{[-\ep,\ep]}(D_i) $ we have $\pr_D(x) = \pr_{D_i}(x)$ (see Sect. 2). For any $z\in R$ denote by 
$\chU(z)$ the {\it part of the unstable manifold} $W^u_{\ep_0}(z)$ such that $\pr_D(\chU(z)) = W^u_{\tR_i}(z)$.
The shift along the flow determines bi-H\"older continuous bijections 
$$\T_z :  W^u_{R}(z) \longrightarrow \chU(z) \quad , \quad \tPsi : W^u_{R}(z) \longrightarrow W^u_{\tR}(z)$$
for all $i$.  These define bi-H\"older continuous bijections   
$$\Psi: R \longrightarrow \check{R} = \cup_{i=1}^{k_0} \check{R}_i  \quad \mbox{\rm and} \quad  \tPsi : R \longrightarrow \tR ,$$ 
where $\check{R}_i = \cup_{z\in S_i} \chU(z)$ and
$\Psi_{| W^u_R(z)} = (\T_z)_{|W^u_R(z)}$ for $z\in S_i$. 
Notice that there exists a global constant $C > 1$ such that 
$$\frac{1}{C} d(x,y) \leq d(\T_z(x), \T_z(y)) \leq C\, d(x,y)$$ 
for all $z\in \tR$ and all $x,y \in W^u_{\tR}(z)$.
\index{$\chU(z)$} 
\index{$\T_z$} 
\index{$\tPsi$} 
\index{$\Psi$} 
\index{$\check{R}_i$} 
\index{$\omega$}

\subsection{Non-integrability}

Throughout we assume that $\phi_t$ is a $C^2$ contact Anosov flow on $M$ with a $C^2$ invariant  contact form $\omega$. 
Then the two-form $d\omega$ is $C^1$, so there exists a constant $C_0 > 0$ such that
\be
|d\omega_x(u,v)| \leq C_0 \|u\|\, \|v\| \quad , \quad u, v \in T_xM\:, \: x\in M .
\ee
Moreover, there exists a constant $\theta_0 > 0$ such that for any $x\in M$ and any $u\in E^u(x)$ with $\|u\| = 1$
there exists $v\in E^s(x)$ with $\|v\| = 1$ such that $|d\omega_x(u,v)| \geq \theta_0$. 
\index{$\theta_0$}

The main ingredient in this section is the following lemma of Liverani (Lemma B.7 in \cite{L2})
which significantly strengthens a lemma of Katok and Burns (\cite{KB}).

\bs

\noindent
{\bf Lemma 4.1.} (\cite{L2}) {\it Let $\phi_t$ be a $C^2$ contact flow on $M$ with a 
$C^2$ contact form $\omega$.
Then there exist constants $C_0  > 0$, $\vartheta > 0$ and $\hep_0 > 0$ such that
for any $z\in M$, any $x\in  W^u_{\hep_0}(z)$ and any  $y\in W^s_{\hep_0}(z)$ we have
\be\label{eq:liverbound}
|\Delta(x,y)  - d\omega_z(u,v)| \leq C_0\,\left[ \|u\|^2\, \|v\|^\vartheta + \|u\|^\vartheta \|v\|^2 \right]\;,
\ee
where $u \in E^u(z)$ and $v \in E^s(z)$ are such that $\exp^u_z(u) = x$ and $\exp^s_z(v) = y$.}

\index{$ C_0$} 
\index{$\exp^u_z$, $\exp^s_z$}

 \bs

\noindent
{\bf Note.}
Actually Lemma B.7 in \cite{L2} is more precise with a particular choice of the constant $\vartheta$
determined by the (uniform) H\"older exponents of the stable/unstable foliations and the corresponding
local holonomy maps. However in this paper we do not need this extra information.

\bs

We will show below that for Lyapunov regular points $x\in \ll$ the estimate (4.2) can be 
improved what concerns the involvement of $u$ for certain choices of $u$ and $v$. More precisely, 
we will show that choosing $v$ in a special way, $\Delta(x,y)$ becomes a $C^1$ function of 
$x = \exp^u_z(u)$ with a non-zero uniformly bounded derivative in a certain direction. 

\bs

{\bf From now on we will assume that $C_0 > 0$, $\vartheta > 0$ and $\hep$ from Theorem 3.2 satisfy
 (4.1) and (4.2)}. As in Sect. 3, set
$$Q_p(\hep) = \{ x\in \ll  : R_{\hep}(x) \leq e^p \}$$
for all  $p > 0$. Then 
$$Q_0(\hep) \subset Q_1(\hep) \subset \ldots \subset Q_n(\hep) \subset \ldots $$
and $\cup_{p=0}^\infty Q_p(\hep) = \ll$.
{\bf Fix an integer} $p_0 \geq 1$ so large that $\mu(Q_{p_0}(\hep)) > 1- \hd_0$ for some small appropriately 
chosen $\hd_0 > 0$ (to be determined later). 
Set
\be \label{eq:pesinset}
P_0 = Q_{p_0}(\hep) \subset \hP_0 = Q_{p_0 + \hep}(\hep) .
\ee
Then the Lyapunov $\hep$-regularity functions $R(x)$, $\Gamma(x)$, $D(x)$ and $L(x)$ are bounded above by some constants on $\hP_0$ 
(and therefore on $P_0$ as well). Similarly, the $\hep$-slowly varying radius function $r(x)$ and the $\hep/2$-slowly varying radius 
function $\tr(x)$ are bounded  below by some constant on $\hP_0$. Thus, we may assume that
\be\label{eq:constbounds}
R (x) \leq R_0 \:\:, \:\: r(x), \tr(x) \geq r_0 \:\: , \:\: \Gamma(x) \leq \Gamma_0 \:\: , \:\:
L(x) \leq L_0 \: \: , \:\: D(x) \leq D_0 
\ee
for all $x\in \hP_0$ {\bf for some positive constants $R_0, \Gamma_0, L_0, D_0 \geq 1$ and $r_0 > 0$.  We fix $r_0 > 0$}
so that $r_0 \leq \frac{1}{R_0}$.
\index{$P_0$, $\hP_0$} 
\index{$\hep_0$}
\index{$\hd_0$}
\index{$R_0$} 
\index{$r_0$} 
\index{$\Gamma_0$}   
\index{$L_0$} 
\index{$D_0$} 

It follows easily from the properties of Markov  families\footnote{Easy proof by contradiction.} that there exists a constant 
$r_1 > 0$  such that for every $i$ and every $x\in \partial R_i$  
there exists $y \in R_i$ such that $\dist(y, \partial R_i) \geq r_1$ and $d(x,y) < r_0/2$.
{\bf Fix a constant} 
$$r_1 < \frac{r_0}{2R_0}$$ 
with this property.
\index{$r_1$} 
\index{$C_1$}  
\index{$\beta_1$} 

Using the notation in (3.20), set
$$\Xi_m = \Xi_m(p_0,\hep,\hd_0)$$
where the constants $p_0$, $\hep$ and $\hd_0$ are as above.

We will now state two Main Lemmas. Their proofs, both using Liverani's Lemma 4.1, are given in Sect. 9.

\bs

\noindent
{\bf Lemma 4.2.}  {\it  There exist constants $C_1 > 0$ and $\beta_1 \in (0,1)$ with the 
following properties:}

(a)  {\it For any unstable cylinder $\cc$ in $R$ of length $m$ with $\cc  \cap P_0\setminus \Xi_m \neq \e$ 
and any $z\in \cc$ we have
\be\label{eq:cylbound}
\frac{1}{C_1 \lambda^p_1 } \leq \diam(\tPsi(\cc)) \leq \frac{C_1 e^{2\hep p}}{\lambda_1^p} ,
\ee
where $p = [\tau_m(z)]$.}

\ms 

(b) {\it For any unstable cylinder $\cc$ of length $m$ in $R$ with $\cc \cap P_0\setminus \Xi_m \neq \e$, 
any $\hx_0,\hz_0 \in \cc$ and any $\hy_0, \hb_0 \in W^s_{R}(\hz_0)$ we have
$$|\Delta(\hx_0 ,\hy_0)  - \Delta(\hx_0, \hb_0)| \leq C_1 \, \diam (\tPsi(\cc)) \, (d(\hy_0,\hb_0))^{\beta_1} .$$
In particular,}
$$|\Delta(\hx_0 , \hy_0)| \leq C_1 \, \diam (\tPsi(\cc)) \, (d(\hy_0,\hz_0))^{\beta_1}  \leq C_1 \, \diam (\tPsi(\cc)) .$$

\bs

{\bf Fix a constant} $C_1 > 0$ with properties in Lemma 4.3. We take $C_1 \geq C_0$.
Set  
\be
\di \beta_0 = \frac{1}{\sqrt{1+ \theta_0^2/(64 C_1^2)}} .
\ee
Next, {\bf fix an integer} $\ell_0 = \ell_0(\beta_0) \geq 1$ so large that we can find
unit vectors $\eta_1, \eta_2 , \ldots,\eta_{\ell_0}$ in $\R^{n_1}$ such that for any unit vector $\xi \in \R^{n_1}$ 
there exists $j$ with $\la \xi , \eta_j\ra \geq \beta_0$. Then {\bf fix measurable families}
$$\eta_1(x), \eta_2(x), \ldots, \eta_{\ell_0}(x) \quad , \quad x\in \ll , $$
of unit vectors in $E^u_1(x)$ such that for any
$x\in \ll$ and any $\xi\in E^u_1(x)$ with $\|\xi\| = 1$ there exists $j$ with  $\la \xi, \eta_j(x) \ra \geq \beta_0$.
\index{$\ell_0$} 
\index{$\eta_j(x)$} 
\index{$\beta_0$} 

Recall the projections $\T_z : W^u_R(z) \longrightarrow \chU(z) \subset W^u_{\ep_0}(z)$ for $z\in R$ from Sect. 4.1.


The following lemma is derived from the non-integrability of the flow which stems from the fact that the
flow is contact. It will play a very essential role in proving that our contraction operators actually do have some
contraction properties due to certain  cancelations provided by property (4.9) below.

\bs

\noindent
{\bf Lemma 4.3.} {\it Let $\phi_t$ be a $C^2$ contact Anosov flow on $M$. Let  
$\eta_1(x), \eta_2(x), \ldots, \eta_{\ell_0}(x)$ ($x\in \ll_0$) be families of unit vectors in $E^u_1(x)$ as above, 
and let $\kappa \in (0,1)$ be a constant.  Then there exist constants  $\ep'' > 0$, 
$0 < \delta'' < \delta'$ (depending on $\kappa$ in general) and $\delta_0 \in (0,1)$
 with the following properties:}

\ms

(a) {\it  For any integer $m \geq 1$ and any $Z \in P_0 \setminus \Xi_m$ there exist families of points 
$$y_\ell(Z) \in  B^s (Z,\delta') \quad , \quad  \ell = 1, \ldots,\ell_0 ,$$
such that 
if $\cc$ is a cylinder of length $m$ in $W^u_R(Z)$ with $Z \in \cc$, then for any $x_0  \in \T_Z(\cc)$ and
$z_0  \in \T_Z(\cc \cap P_0)$ of the form $x_0 = \Phi^u_{Z}(u_0)$,  $z_0 = \Phi^u_{Z}(w_0)$ such that
\be
d(x_0,z_0) \geq \kappa \, \diam (\cc'') , 
\ee 
where $\cc'' = \T_z(\cc)$, and
\be
\left\langle \frac{w_0 - u_0}{\|w_0 - u_0\|} , \eta_\ell (Z) \right\rangle \geq \frac{\beta_0}{2R_0} 
\ee
for some $\ell = 1, \ldots, \ell_0$, then we have
\be \label{eq:deltabound}
\frac{\beta_0 \delta_0 \kappa}{R_0} \, \diam(\cc'') 
\leq |\Delta(x_0,\pi_{d_1}(z_0))  - \Delta(x_0, \pi_{d_2}(z_0))| 
\ee
for any $d_1 \in B^s(y_\ell (Z), \delta'')$ and  $d_2 \in B^s(Z, \delta'')$.}

\ms

(b)  {\it There exists an integer $N_0 \geq 1$ such that for any integer $N \geq N_0$, 
any integer $m \geq 1$ 
and any $Z \in P_0 \setminus \Xi_m$ there exist families of points 
$$\ylo(Z), \ylt(Z) \in \pp^N(B^u(Z;\ep'')) \cap B^s (Z,\delta')  \quad , \quad \ell = 1, \ldots,\ell_0 ,$$
such that if $\cc$ is a cylinder of length $m$ in $W^u_R(Z)$ with $Z \in \cc$, $x_0  \in \T_Z(\cc)$ and
$z_0  \in \T_Z(\cc \cap P_0)$ have the form $x_0 = \Phi^u_{Z}(u_0)$,  $z_0 = \Phi^u_{Z}(w_0)$ and
{\rm (4.7)} and {\rm (4.8)} hold for some $\ell = 1, \ldots, \ell_0$, then {\rm (4.9)} holds
for any $d_1 \in B^s(\ylo (Z), \delta'')$ and  $d_2 \in B^s(\ylt(Z), \delta'')$.}

\index{$\ylo(Z)$, $\ylt(Z)$} 
\index{$d_1$, $d_2$}

\subsection{Regular distortion of cylinders}

In \cite{St5} we established some nice properties concerning diameters of cylinders for 
Axiom A flows on basic sets satisfying a pinching condition which we called {\it regular distortion 
along unstable manifolds}. In \cite{St4} something similar was established for Anosov flows with 
Lipschitz local stable holonomy maps. It seems unlikely that every Anosov flow will 
have such properties, however it turns out that for general Anosov flows something similar holds 
for `sufficiently regular' cylinders in $R$.

If $\cc$ is an unstable cylinder of length $m$ and $\cc'$ is a sub-cylinder of $\cc$ of length $m+k$,
we will say that $\cc'$ has {\it co-length} $k$ in $\cc$.
\index{co-length}


\bs

\noindent
{\bf Lemma 4.4.} 
(a) {\it   There exists a constant $0 < \rho_1  < 1$ 
such that for any unstable leaf $W$ in $R$, any cylinder $C_W[\ii] = C_W[i_0, \ldots,i_m]$  in $W$  and any sub-cylinder 
$$C_W[\ii'] = C_W[i_0,i_{1}, \ldots, i_{m+1}]$$ 
of $C_W[\ii]$ of co-length $1$ 
such that $C_W[\ii'] \cap P_0 \neq \e$ and $\pp^{m+1}(C_W[\ii']) \cap P_0 \neq \e$  we have}
$$\rho_1 \; \diam ( \tPsi(C_W [\ii]) ) \leq  \diam ( \tPsi(C_W [\ii'] )) .$$

(b) {\it For any constant $\rho' \in (0,1)$ there exists an integer $q' \geq 1$ such that
for  any unstable leaf $W$ in $R$, any cylinder $C_W[\ii] = C_W[i_0, \ldots,i_m]$ of length $m$ in $W$  and any sub-cylinder 
$$C_W[\ii'] = C_W [i_0,i_{1}, \ldots, i_{m+1}, \ldots, i_{m+q'}]$$ 
of $C_W[\ii]$ of co-length $q'$  such that $C_W[\ii'] \cap P_0 \neq \e$ and $\pp^{m+q'}(C_W[\ii']) \cap P_0 \neq \e$  we have}
$$\diam(\tPsi(C_W[\ii']) ) \leq \rho' \, \diam (\tPsi(C_W[\ii])) .$$

(c)  {\it There exist an integer $q_0 \geq 1$ and a constant $\rho_1 \in (0,1)$ 
such that for any unstable leaf $W$ in $R$  and any cylinder $C_W[\ii] = C_W[i_0, \ldots,i_m]$ in $W$ 
such that 
$$C_W[\ii] \cap P_0 \neq \e \quad \mbox{\it  and} \quad \pp^m(C_W[\ii]) \cap P_0 \neq \e$$
there exist points $z, x\in C_W[\ii]$ such that for the sub-cylinder
$$C_W[\ii'] = C_W[i_0,i_{1}, \ldots, i_{m+1}, \ldots, i_{m+q_0}]$$
of  $C_W[\ii]$ of co-length $q_0$ containing $x$ we have 
$$d(z,y) \geq \rho_1 \, \diam (\tPsi(C_W[\ii]))$$
for all $y \in C_W[\ii']$.}

\ms

(d) {\it We can choose the constant $\rho_1 \in (0,1)$ from part {\rm (a)} 
such that for any unstable leaf $W$ in $R$, any cylinder $C_W[\ii] = C_W[i_0, \ldots,i_m]$  in $W$   and any sub-cylinder 
$C_W[\ii'] = C_W[i_0,i_{1}, \ldots, i_{m+1}]$ of $C_W[\ii]$ of co-length $1$ 
such that there exists $x \in C_W[\ii']$ with $\pp^m(x) \in P_0$  we have}
$$\rho_1 \; \diam ( (\Phi_x)^{-1}(C_W [\ii] )) \leq  \diam ( (\Phi_x)^{-1}(C_W [\ii'] )) .$$

\ms

Clearly, all statements in Lemma 4.4 remain true replacing $P_0$ by $\hP_0$ and slightly changing the
constants $\rho'$, $q'$, $q_0$ and $\rho_1$. 

Notice that in the last part (d) we do not require $C_W[\ii'] \cap P_0 \neq \e$, however this time we
measure diameters `upstairs' in the tangent bundle $E^u(x)$.

Lemma 4.4 will be used essentially in the proof of the main result in Sections 5-8 below. 
Its proof is given in Sect. 10.

\index{$q_0$} 
\index{$\rho_1$}

\section{Main estimates for temporal distances}
\setcounter{equation}{0}

This section contains the main application of Lemma 4.3, namely Lemma 5.4 which provides estimates from below of 
differences of temporal distances in a form convenient for the estimates of contraction operators in Sect. 6. 
As in \cite{D1} (and \cite{St2}), given a large integer $N$, the construction of the contraction operators involves the choice of some particular
inverses, or rather, families of inverses (as in \cite{St2}) of the map $\sigma^N$. Here however the construction of these inverses is much more
complicated and involves the Pesin set $P_0$ and the set $\Xi_m$ from Sect. 4.2. This is done in the proof of Lemma 5.4 which is the basis for
the construction of the contraction operators in Sect. 6.

\subsection{Normalized Ruelle operators and the metric $\dte$}

Let the constants $C_0 > 0$, $c_0 > 0$, $1 < \gamma < \gamma_1$  be as in Sections 2 and 4, so that
(2.1) and (4.1) hold. {\bf Fix a constant} $\theta$ such that
\be
\frac{1}{\gamma^{\alpha_1 }}  = \hat{\theta} \leq \theta < 1,
\ee
where $\alpha_1 > 0$ is the constant chosen at the end of Ch. 2.
\index{$\theta$} 
\index{$\hat{\theta}$}
 
Recall the metric $\dte$ on $\hU$  and the space $\ff_\theta(\hU)$ from Sect. 1.1.
In the same way we define the distance $\dte(x,y)$ for $x,y$ in $W \cap \hR$. Lemma 5.2 below shows that $\tau\in \ff_\theta(\hU)$.
For a non-empty subset  $A$ of $U$ (or some $W^u_R(x)$) let $\diamte(A)$ be the {\it diameter}  of $A$ with respect to $\dte$.

Let $f \in \ff_\theta (\hU)$ be a fixed real-valued function  and let  $g = f - P_f\, \tau$, where $P_f\in \R$ is such that $\Pr_{\sigma} (g) = 0$. 
Since $f$ is a H\"older continuous function on $\hU$, it can be extended to a H\"older continuous function on $R$ which is constant on stable leaves. 
Set 
$$\Fa = f - (P_f + a) \tau .$$
By Ruelle-Perron-Frobenius' Theorem (see e.g. chapter 2 in \cite{PP}) for any real number $a$  
with $|a|$ sufficiently small, as an operator on $\ff_\theta(\hU)$, $L_{\Fa}$ has a  {\it largest eigenvalue} 
$\lambda_{a}$ and there exists a (unique) regular probability measure $\hnu_a$ on $\hU$ with 
\index{$\hnu_a$} 
$L_{\Fa}^*\hnu_a = \lambda_a\, \hnu_a$, i.e.  
$$\int L_{\Fa} H \, d\hnu_a = \lambda_a\, \int H\ d\hnu_a$$
for every $H \in \ff_\theta(\hU)$.
Fix a corresponding (positive) eigenfunction $h_{a} \in \ff_\theta(\hU)$ such that $\int h_{a} \, d\hnu_a = 1$. 
\index{$h_a$} 
Then $d\nu = h_0\, d\hnu_0$ defines a {\it $\sigma$-invariant  probability measure} $\nu$ on $U$, called the {\it Gibbs measure}
determined by the function $F^{(0)}$.
Since $\Pr_\sigma (f- P_f\tau) = 0$, it follows from the main properties of pressure (cf. e.g. chapter 3 in \cite{PP}) 
that $|\Pr_\sigma(\Fa)| \leq \|\tau \|_0 \, |a|$.  Moreover, for small $|a|$ the maximal eigenvalue
$\lambda_{a}$ and the eigenfunction $h_{a}$ are Lipschitz in  $a$, so
there exist constants $a_0 > 0$ and $C > 0$ such that $|h_{a} - h_0| \leq C\, |a|$ on $\hU$ and
$|\lambda_a - 1| \leq C |a|$ for  $|a| \leq a_0$.
\index{$\nu$} 
\index{$\Fa$}

For $|a|\leq a_0$,  as in \cite{D1}, consider the function
$$\fa(u) = f (u) - (P_ f+ a) \tau(u) + \ln h_{a}(u) -  \ln h_{a}(\sigma(u)) - \ln \lambda_{a}$$
\index{$\fa$} 
and the operators 
$$ \lab = L_{\fa - \i\,b\tau} : \ff_\theta(\hU) \longrightarrow \ff_\theta (\hU)\:\:\: , \:\:\: \ma = L_{\fa} : \ff_\theta(\hU) \longrightarrow \ff_\theta(\hU) .$$
The Ruelle operator $\ma$ is the standard normalisation of $L_f$ which has largest eigenvalue $1$ and  satisfies $\ma \; 1 = 1$ (see e.g. Ch. 2 in \cite{PP}).

One checks that  
$$\di |(\lab^m h)(u)| \leq (\ma^m |h|)(u) \quad , \quad u\in \hU ,$$
for all $h\in \ff_\theta (\hU)$ and $m \geq 0$. Moreover, $L_{f^{(0)}}^*\nu = \nu$,  i.e.  
$$\int L_{f^{(0)}} H \, d\nu = \int H\, d\nu$$ 
for any $H \in \ff_\theta (\hU)$.
\index{$\lab$}
\index{$\ma$}  

Since $g$ has zero topological pressure with respect to the shift map $\sigma : U \longrightarrow U$, there exist constants $0 < c_1 \leq c_2$  
such that for any cylinder $\cc = \cc^u[i_0, \ldots, i_m]$  of length $m$ in $U$ we have
\be\label{eq:Gibbs}
c_1 \leq  \frac{\nu(\cc)}{e^{g_m(y) }} \leq c_2 \quad , \quad y \in \cc ,
\ee
(see e.g. \cite{PP} or \cite{P2}).
\index{$c_1$, $c_2$} 
We now state some basic properties of the metric $D_\theta$ that will be needed later.

\bs

\noindent
{\bf Lemma 5.1.}  
(a) {\it For any cylinder $\cc$ in $U$ the characteristic function $\chi_{\cc}$ of
$\cc$ on $U$ is Lipschitz with respect to $\dte$ and $\Lip_\theta(\chi_\cc) \leq 1/\diamte(\cc)$.}

\ms

(b) {\it There exists a constant $C_2 > 0$ such that if $x,y\in \hU_i$ for some $i$,  then 
$$|\tau(x) - \tau(y) |  \leq  C_2\, \dte(x,y) .$$
That is, $\tau\in \ff_\theta(\hU)$. Moreover, we can choose $C_2 > 0$ so that
$$|\tau_m(x) - \tau_m(y)| \leq C_2\,  \dte(\sigma^m(x),\sigma^m(y))$$
whenever $x,y\in \hU_i$ belong to the same cylinder  of length $m$.}

\ms

(c) {\it There exist  constants $C_2 > 0$ and $\alpha_2 > 0$  such that for any $z\in R$, any cylinder 
$\cc$ in $W^u_R(z)$ and any $x,y\in \cc$ we have 
$$d(\tPsi(x),\tPsi(y))  \leq  C_2\, \dte(x,y)$$ 
and
$$\diamte(\cc) \leq C_2 (\diam(\tPsi(\cc)))^{\alpha_2} .$$
Moreover, we can take $\alpha_2 > 0$ so that $1/(\gamma_1)^{\alpha_2} = \theta$.}
\index{$\alpha_2$} 

\bs

\noindent
{\it Proof.} (a) Let $\cc$ be a cylinder in $U$ and let $x,y\in \hU$.  If $x,y\in \cc$ or $x\notin \cc$ and $y\notin \cc$,
then $\chi_\cc(x) - \chi_\cc(y) = 0$. Assume that $x\in \cc$ and $y \notin \cc$.
Let $\dte(x,y) = \theta^{N+1}$ and let $\cc'$ be a cylinder of length $N$ containing both $x$
and $y$. Since $x\in \cc$, as well, and $x$ is an interior point of $\cc$, we must have
$\cc\subset \cc'$. Thus, $\diamte(\cc) \leq \dte(x,y)$. This gives
$$|\chi_\cc(x) - \chi_\cc(y)| = 1 = \frac{\diamte (\cc)}{\diamte(\cc)} \leq \frac{1}{\diamte(\cc)}\, \dte(x,y) ,$$
which proves the assertion.

\ms

(b), (c)  Assume $x\neq y$ and let $\cc$ be the cylinder of largest length $m$
containing both $x$ and $y$. Set $\tx = \tPsi(x), \ty = \tPsi(y) \in \tR$. Then $\dte(x,y) = \theta^{m+1}$. On the other hand, (2.1) and (5.1)  imply
\begin{eqnarray*}
|\tau(x) - \tau(y)| 
& \leq & |\tau|_{\alpha_1}\, (d(\tx,\ty))^{\alpha_1}   \leq \frac{\Con}{(\gamma^{\alpha_1})^{m}} \leq  \Con \theta^m \leq C_2\, \dte(x,y)
\end{eqnarray*}
for some global constant $C_2 > 0$. The above also shows that 
$$ d(\tx ,\ty ) \leq   \Con \theta^m \leq C_2\, \dte(x,y) ,$$
which proves half of  part (c). The second part of (c) follows by using a similar estimate and the other half of (2.1).

Next, assume that $x,y$ belong to the same cylinder $\cc$ of length $m$. Let 
$$\pp_j(x), \pp^j(y) \in R_{i_j} \quad , \quad j = 0,1, \ldots, m .$$
Assume that  $\dte(x',y') = \theta^{p+1}$, where
$x' = \sigma^m(x)$ and $y' = \sigma^m(y)$. Then $\dte (x,y) = \theta^{m+p+1}$
and moreover 
$$\dte(\sigma^j(x), \sigma^j(y)) = \theta^{m-j+p+1} \quad , \quad j = 0,1, \ldots,m-1 .$$
Now (2.1) and (5.1)  imply
\begin{eqnarray*}
&        & |\tau(\sigma^j(x)) - \tau( \sigma^j(y))| \\
& \leq & |\tau|_{\alpha_1} (d (\sigma^j(\tx), \sigma^j(\ty))^{\alpha_1} 
\leq  \Con\,  (d (\tPsi(\sigma^j(x)), \tPsi(\sigma^j(y)))^{\alpha_1} \\
& \leq & |\tau|_{\alpha_1} \left(\frac{1}{c_0 \gamma^{m-j+p}} 
d(\tpp^{m+p-j}(\tPsi(\sigma^j x)), \tpp^{m+p-j}(\tPsi(\sigma^j y)))  \right)^{\alpha_1 }\\
& \leq & \frac{\Con}{(\gamma^{\alpha_1 })^{(m-j+p)}} \leq \Con\, \theta^{m-j+p} 
\leq \Con \theta^{m-j-1} \dte(x',y') .
\end{eqnarray*}
Therefore
\begin{eqnarray*}
|\tau_m(x) - \tau_m(y)| 
& \leq & \sum_{j=0}^{m-1}|\tau(\sigma^j(x)) - \tau( \sigma^j(y))|  \\
& \leq & \Con \dte(x',y')\, \sum_{j=0}^{m-1} \theta^{m-j+1}  \leq  \Con \, \dte(x',y') ,
\end{eqnarray*}
which proves the statement.
\endofproof

\bs

It follows from Lemma 5.1 that $\tau \in \ff_\theta(\hU)$, so assuming $f \in \ff_\theta(\hU)$, we have
$h_a\in \ff_\theta(\hU)$ for all $|a|\leq a_0$. Then $\fa \in \ff_\theta(\hU)$ for all such $a$.
Moreover, using the analytical dependence of $h_a$ and $\lambda_a$ on $a$ and assuming that the 
constant $a_0 > 0$ is sufficiently small, there exists $T_0 = T_0(a_0)$ such that
\begin{equation}
T_0 \geq \max \{ \, \|\fa \|_0 \, , \, |\fa_{|\hU}|_\theta \, , |\tau_{|\hU}|_\theta \, \}\;
\end{equation}
for all $|a| \leq a_0$. Fix $a_0 > 0$ and $T_0 > 0$ and with these properties.
Taking the constant $T_0 > 0$ sufficiently large, we have  $\|\fa - f^{(0)}\|_0 \leq T_0\, |a|$ 
on $\hU$ for $|a| \leq a_0$.

\index{$T_0$} 
\index{$a_0$} 

The following Lasota-Yorke type inequality is similar to that in \cite{D1}, and the corresponding one
in \cite{St2} (although we now use a different metric) and its proof is also very similar. We include 
a proof in the Appendix (Sect. 11) for completeness.

\bs

\noindent
{\bf Lemma 5.2.}  {\it There exists a constant $A_0 > 0$  such that for all $a\in \R$ with $|a|\leq a_0$ the following hold:
 If the functions $h$ and  $H$ on $\hU$ and the constant $B > 0$ are such that $H > 0$ on $\hU$ and 
$$|h(v) - h(v')| \leq B\, H(v')\, \dte (v,v')$$
for any $i$ and any $v,v'\in \hU_i$, 
then for any $b\in \R$ with  $|b|\geq 1$ and any integer $m \geq 1$ we have 
$$ |\lab^m h(u) - \lab^m h(u')| \leq  A_0\,\left[ B\,\theta^m \, (\ma^m H)(u')  + |b|\, (\ma^m |h| )(u')\right]\, \dte (u,u')$$
whenever $u,u'\in \hU_i$ for some $i = 1, \ldots,k_0$.} 
\index{$A_0$} 

\bs

\noindent
{\bf Remark.} It follows from the proof of this lemma  that the constant
$A_0$ depends only on $\|f\|_\theta$ and some global constants, e.g. $c_0$ and $\gamma$ in {\rm (2.1)}.

\subsection{First step -- using Lemmas 4.2 and 4.3}

Let the constants $c_1$ and $c_2$ be as in (5.2). 
{\bf Fix constants} $\rho_1 \in (0,1)$ and $q_0 \geq 1$ such that Lemma 4.4(a), (b), (c) and (d)
hold with $\rho' = \rho_1/8$ and $q' = q_0$, {\bf both for the set $P_0$ and when $P_0$ is replaced by $\hP_0$.}
{\bf In what follows we will use the entire set-up and notation from Sect. 4}, 
e.g. the subsets $P_0$ and $\hP_0$ of $\ll$, the numbers $r_0 \geq r_1 > 0$, $R_0 > 1$, etc.,  satisfying  (4.3), (4.4), etc.
Let $\eta_1(x), \eta_2(x), \ldots, \eta_{\ell_0}(x)$ ($x\in M$) 
be families of unit vectors in $E^u_1(x)$ as in the text just  before Lemma 4.3, and let $\ep'' \in (0,\ep')$, 
$0 < \delta'' < \delta' $, $\delta_0 > 0$ (depending on the  choice of $\kappa$), $\beta_1 \in (0,1)$,  $C_1 > 0$
 be constants with the properties described in Lemmas 4.3 and 4.2.  
 \index{$\ep''$} 
 \index{$\delta''$} 
 \index{$\delta'$} 
 \index{$\delta_0$}

Let $E > 1$ be a constant -- we will see later how large it should be, and  let $\ep_1 > 0$ be a constant with
\be\label{eq:epone}
0 < \ep_1 \leq  \min\left\{ \;\frac{1}{32 C_0R_0 }\;, \;  \frac{1}{4EC_1 R_0^2} \, \right\} .
\ee
\index{$E$} 
\index{$\ep_1$} 

Having fixes $\theta \in [\htheta,1)$ in (5.1), now {\bf fix a constant} $\theta_1 \in (0,\theta)$ with
$$0 < \theta_1 = \frac{1}{\gamma_1} = \theta^{1/\alpha_2}  < \theta ,$$
\index{$\theta_1$} 
recalling the choice of $\alpha_2 > 0$ in Lemma 5.1(c),  and set
$$\theta_2 = \max\{ \theta , 1/ \gamma^{\alpha_1 \beta_1}\} .$$
\index{$\theta_2$} 

Recall from Sect. 4.1 that for any $z \in R_i$ ($i = 1, \ldots, k_0$), $\chU(z)$ denotes the part of the unstable manifold $W^u_{\ep_0}(z)$ such that $\pr_D(\chU(z)) = W^u_{\tR_i}(z)$
and the map $\T_z :  W^u_{R}(z) \longrightarrow \chU(z) $ is defined using the shift along the flow.

{\bf Fix a constant} $A > 0$ such that the maps $\tPsi^{-1}\circ \T_z$ and $\T_z^{-1} \circ \tPsi$ are $A$-Lipschitz for any
$z \in R$. We will also assume that $A > 0$ is so large that for any $i =1, \ldots, k_0$ and any $x,x'\in R_i$ we have 
$$\diam(\tPsi(W^u_{R_i}(x))) \leq A\, \diam(\tPsi(W^u_{R_i}(x'))) .$$
\index{$A$} 
\index{$N_0$}

Fix an integer  $N_0 \geq 1$ with the property described in Lemma 4.3(b), and then
take $N \geq N_0$ such that
\be\label{eq:Ncond}
\gamma^N > \frac{1}{\delta''}   \quad , \quad 
\theta^N < \frac{\rho^2_1 \beta_0 \delta_0\ep_1}{256  E} 
\quad , \quad \theta_2^N < \frac{\delta_0 \hrho \ep_1}{64 E A R_0^2} ,
\ee
where  $\beta_1 > 0$ is the constant from Lemma 4.2. 
A few additional conditions on $N$ will be imposed later. Set 
\be \label{eq:hddefn}
\di \hd = \frac{\beta_0 \delta_0 \rho_1}{16 R_0 C_0^3} ,
\ee
where $\beta_0 > 0$ is defined by (4.6).

\bs

\noindent
{\bf Lemma 5.3.} {\it Let $\cc$ be an unstable cylinder in $R$ of length 
$m \geq 1$  such that $\cc \cap P_0 \neq \e$ and $\pp^m(\cc) \cap P_0 \neq \e$. }

\ms

(a) {\it There exist sub-cylinders $\dd$ and $\dd_1$ of $\cc$ of co-length $q_0$ such that 
$$d(\tPsi(y),\tPsi(x)) \geq \frac{\rho_1}{2}  \, \diam(\tPsi(\cc))$$ 
for all $y \in \dd_1$ and $x\in \dd$.}

\ms

(b)   {\it  There exists an integer $q_1 \geq q_0$ such that for any  sub-cylinder 
$\cc_1$ of $\cc$ of co-length $q_1$ with $\pp^m(\cc_1)\cap P_0 \neq \e$ 
we have} 
$$\diam(\tPsi(\cc_1)) \leq \min\left\{ \frac{\rho_1}{8} , \frac{\hd }{8 C_1} \right\}\, \diam (\tPsi(\cc)) .$$

\ms

\noindent
{\it Proofs.} (a) Take $z,x  \in \cc$ as in Lemma 4.4(c), and let $\dd$ and $\dd_1$ be
the sub-cylinders of $\cc$ of co-length $q_0$ containing $z$ and $x$, respectively. By Lemma 4.4(b)
and the choice of $q_0$ and $\rho' = \rho_1/8$  it follows that 
$$\diam(\tPsi(\dd)) \leq \frac{\rho_1}{8} \diam(\tPsi(\cc)) .$$

Next, by the choice of $z,x$ in Lemma 4.4(c), for any $y \in \dd_1$ we have 
$$d(\tPsi(y),\tPsi(z)) \geq \rho_1\, \diam(\tPsi(\cc)) .$$
Then
\begin{eqnarray*}
d(\tPsi(x),\tPsi(y)) 
& \geq & d(\tPsi(y), \tPsi(z)) - d(\tPsi(x), \tPsi(z))  \\
& \geq & \rho_1 \diam(\tPsi(\cc)) - \frac{\rho_1}{8} \diam (\tPsi(\cc))
  >  \frac{\rho_1}{2}\,  \diam (\tPsi(\cc))  
\end{eqnarray*}
for any $y\in \dd_1$ and any $x\in \dd$.

\ms

(b) This follows from Lemma 4.4(b): take $q_1 = q_0^r$ for some sufficiently large integer $r \geq 1$.
\endofproof

\subsection{Main estimates for temporal distances}

We will use Lemma 4.3 with 
$$\kappa = \hrho/2 \quad \mbox{\rm  where} \quad  \hrho = \frac{\rho_1}{8C_0^2} ,$$ 
$\rho_1$ being the constant from Lemma 4.4. Define $\hd$ by (5.6). 
We will also use the integers $N_0  \geq 1$ and the constants $\ep'' > 0$ and $\delta' > \delta'' > 0$ from Lemma 4.3.
\index{$\hrho$} 
\index{$\hd$} 

Assume the {\it integer $n_0 \geq 1$ is chosen so large that} for any $z\in R$ and any
unstable cylinder $\cc$ of length $\geq n_0$ in $R$ we have $\diam(\tPsi(\cc)) \leq \ep''$
and $\diam(\T_z(\cc)) \leq \ep''$ for any $z\in \cc$.
\index{$n_0$} 

In what follows we will use the constants $p_0$, $\hep$ and $\hd_0$, the Pesin set $P_0$ and the
set $\Xi_m = \Xi_m(p_0,\hep, \hd_0)$ from Sect. 4.2.

Given $m \geq n_0$ and $Z\in P_0 \setminus \Xi_m$, let
$$\ylo(Z) \in B^s (Z,\delta') \cap \pp^N(B^u(Z,\ep'')) \quad , \quad \ylt(Z) \in B^s (Z,\delta') \cap \pp^N(B^u(Z,\ep'')) ,$$
($\ell = 1, \ldots,\ell_0$)  be families of points satisfying the requirements of Lemma  4.3(b). 
\index{$\ylo(Z)$, $\ylt(Z)$} 

\bs

\noindent
{\bf Lemma 5.4.} 
{\it For any  $m \geq n_0$, any point $Z \in P_0\setminus \Xi_m$, any integer $N \geq N_0$, any 
$\ell = 1, \ldots, \ell_0$ and any $i = 1,2$ there exists a (H\"older) continuous map 
$$B^u(Z',\ep'') \ni x \mapsto v^{(\ell)}_i (Z, x) \in U ,$$
where $Z' = \piU(Z) \in U$, such that $\sigma^N(\vl_i(Z, x)) = x$  for all $x\in B^u(Z,\ep'')$
and the following property holds:

For  any cylinder $\cc$ in $W^u_R(Z)$ of length $m$ with $Z\in \cc$ and $\pp^m(Z) \in P_0$
there exist sub-cylinders $\dd$ and $\dd_1$ of $\cc$ of co-length $q_1$ 
and  $\ell = 1, \ldots, \ell_0$ such that $Z \in \dd$ and for any points $x\in \dd$ 
and $z\in \dd_1$, setting $x' = \piU(x)$, $z' = \piU(z)$,  we have 
$$d(\T_Z(x), \T_Z(z)) \geq \frac{\hrho}{2} \, \diam(\T_Z(\cc))$$ 
and 
$$I_{N, \ell}(x',z') =  \left|  \varphi_{\ell}(Z,x') - \varphi_{\ell}(Z,z') \right|   \geq \hd \, \diam(\tPsi(\cc)) ,$$
where 
$$\varphi_{\ell}(Z, x) = \tau_{N}(\vl_1(Z,x)) - \tau_{N}(\vl_2(Z,x)) .$$
Moreover, 
$$I_{N,\ell}(x',z') \leq C_1\, \diam(\tPsi(\cc))$$ 
for any $x,z \in \cc$,  where $C_1 > 0$ is the constant from Lemma} 4.2.
\index{$I_{N,\ell}(x,z)$} 
\index{$\varphi_{\ell}(Z,x)$} 

\bs

\noindent
{\it Proof.}  Fix for a moment  $Z \in P_0 \setminus \Xi_m$, $N \geq N_0$ and $\ell = 1, \ldots, \ell_0$. 
Assume $Z\in R_{i_0}$. Using  Lemma 4.3, there exist points  
$$\ylo = \ylo(Z) \quad , \quad \ylt = \ylt(Z) \in W^s_{\delta'} (Z)$$ 
such that  the property (b) in Lemma 4.3 holds.

Given $i = 1,2$, there exists  a cylinder $L^{(\ell)}_i = L^{(\ell)}_i(Z)$ of length $N$ in $W^u_{R_{i_0}}(Z)$  so that 
$$\pp^{N} : L^{(\ell)}_i \longrightarrow W^u_{R_{i_0}}(\yli)$$
is a bijection; then it is a bi-H\"older homeomorphism. Consider its inverse and its H\"older continuous
extension  $\pp^{-N} : W^u_{R_{i_0}}(\yli)   \longrightarrow L^{(\ell)}_i$.

Set 
$$M^{(\ell)}_i = M^{(\ell)}_i(Z) = \piU(L^{(\ell)}_i(Z)) \subset U ;$$ 
this is then a cylinder  of length $N$ in $U_{i_0}$. Define the maps
$$\tvli (Z, \cdot) : U_{i_0} \longrightarrow L^{(\ell)}_i \subset B^u(Z,\ep'')  
\quad , \quad \vli (Z, \cdot)   : U_{i_0} \longrightarrow M^{(\ell)}_i  \subset U$$
by 
$$\tvli(Z, y) = \pp^{- N} (\pi_{\yli}(y)) \quad ,\quad  \vli(Z, y) = \piU(\tvli(Z, y)) .$$
Then  
\be
\pp^{N}(\tvli(Z, y)) = \pi_{\yli}(y) = W^s_{\ep_0}(y) \cap W^u_{R_{i_0}}(\yli) ,
\ee
and
\be
\pp^{N}(\vli(Z, y)) = W^s_{\ep_0}(y) \cap \pp^{N}(M^{(\ell)}_i) = \pi_{\dli}(y) ,
\ee
where $\dli = \dli(Z) \in W^s_{R}(Z)$ is such that $\pp^{N}(M^{(\ell)}_i) = W^u_{R}(\dli)$.
Next, there exist $x' \in M^{(\ell)}_i$ and $y'\in L^{\ell}_i$ with $\pp^{N} (x') = \dli$ and
$\pp^{N}(y') = \yli$. Since stable leaves shrink exponentially fast, using (2.1) and (\ref{eq:Ncond}) we get
\be
d(\dli , \yli) \leq \frac{1}{c_0 \gamma^{N}} d(x',y') \leq  \frac{1}{\gamma^{N}} < \delta'' .
\ee
{\bf Thus, $\dlo, \dlt$ satisfy the assumptions and therefore the conclusions of Lemma 4.3(b).}
\index{$\dli(Z)$}

Let $\cc$ be a cylinder in $W^u_{R_{i_0}}(Z)$ of length $m$ with 
$$Z \in \cc \cap P_0 \setminus \Xi_m \quad , \quad \pp^m(Z) \in P_0 .$$
Set $\cc'' = \T_Z(\cc)$ and $\tcc = \tPsi(\cc)$. By the choice of the constant $C_0$, we have 
$$\frac{1}{C_0} \diam(\tcc) \leq \diam(\cc'') \leq C_0\, \diam (\tcc) .$$
Let $\dd$ be the sub-cylinder of $\cc$ of co-length $q_1$ containing $Z$.

Next, by Lemma 5.3(a), there exists a sub-cylinder $\dd_1$ of $\cc$ of co-length  $q_1$ such that
$$d(\tPsi(y), \tPsi(x)) \geq \frac{\rho_1}{2} \diam(\tcc)$$ 
for all $y \in \dd_1$ and $x\in \dd$. Thus, 
$$d(\T_Z(y), \T_Z(x)) \geq \frac{\rho_1}{2C_0^2} \diam(\cc'') \geq \hrho \, \diam(\cc'')$$ 
for all $y \in \dd_1$ and $x\in \dd$.

Let $x \in \dd$, $x \neq Z$. Set $x_0 = \T_Z(x)$ and let $x_0 = \Phi^u_{Z}(u_0)$, where $u_0 \in E^u(Z)$.
Since $x \neq Z$ we have $u_0 \neq 0$.
By the choice of the constant $\beta_0$ and the family of unit vectors 
$\{\eta_\ell(Z)\}_{\ell=1}^{\ell_0}$, there exists some $\ell = 1, \ldots, \ell_0$ such that 
$$ \left\langle \frac{u_0}{\|u_0\|}, \eta_\ell(Z) \right\rangle \geq \beta_0 .$$
Moreover, $d(x_0, Z) \geq \hrho\, \diam(\T_Z(\cc))$.
It then follows from Lemma 4.3(b) with $\kappa = \hrho$ and (5.9) that 
\be
\frac{\beta_0 \delta_0 \hrho}{R_0 } \, \diam(\T_z(\cc)) \leq \left| \Delta(x_0, \dlo) - \Delta(x_0, \dlt)\right| .
\ee
(In the present situation, since $\dlo , \dlt \in W^s_R(Z)$, we have $\pi_{\dli}(Z) = \dli$ for $i = 1,2$.)

Consider the projections of $x, Z$ to $U$ along stable leaves: $x' = \piU(x) \in U_i$, $ Z' = \piU(Z) \in U_{i_0}$, 
where as before $R_{i_0}$ is the rectangle containing $Z$ (and therefore $\cc$).  We have
\begin{eqnarray*}
&     & I_{N, \ell}(x',Z')\\
& = & \left| \; [\tau_{N}(\vl_1(Z,x')) - \tau_{N}(\vl_2(Z,x'))] - [\tau_{N}(\vl_1(Z,Z')) - \tau_{N}(\vl_2(Z,Z'))] \; \right|\\
& = & \left| \; [\tau_{N}(\vl_1(Z, x')) - \tau_{N}(\vl_1(Z,Z'))]  - [\tau_{N}(\vl_2(Z,x')) - \tau_{N}(\vl_2(Z,Z'))] \; \right|\\
& = & \left| \, \Delta(\pp^{N} (\vl_1(Z, x')), \pp^{N}(\vl_1(Z, Z')))  -  \Delta(\pp^{N} (\vl_2(Z, x')), \pp^{N}(\vl_2(Z, Z'))) \,\right|\\
& = & \left| \, \Delta(\pi_{\dlo}(x'), \pi_{\dlo}(Z') ) -  \Delta(\pi_{\dlt}(x') , \pi_{\dlt} (Z') ) \,\right|\\
& = & \left| \, \Delta( x' , \dlo ) -  \Delta( x', \dlt ) \,\right|  = \left| \, \Delta( x , \dlo ) -  \Delta( x  , \dlt ) \,\right|.
\end{eqnarray*}
We claim that the latter is the same as the right-hand-side of (5.10). Indeed, let
$\Delta(x, \dlo) = s_1$ and $\Delta(x, \dlt) = s_2$. 
Then $\phi_{s_1}([x, \dlo]) \in W^u_{\ep_0}(\dlo)$ and $\phi_{s_2}([x, \dlt]) \in W^u_{\ep_0}(\dlt)$. 
Let $\phi_s(x_0) = x$. It is then straightforward to see that
$\Delta(x_0, \dlo) = s + s_1$ and $\Delta(x_0, \dlt) = s+ s_2$. Thus,
$$\left| \Delta(x_0, \dlo) - \Delta(x_0, \dlt)\right| = |(s+s_1) - (s+s_2)| = |s_1-s_2| = \left| \Delta(x, \dlo) - \Delta(x, \dlt)\right| .$$
Combining this with (5.10) gives
$$I_{N,\ell}(x' ,Z') \geq \frac{\beta_0 \delta_0 \hrho}{R_0} \, \diam(\T_z(\cc)) \geq 2\hd \, \diam(\tPsi(\cc)) .$$

For arbitrary $x,z\in \cc$, setting $x' = \piU(x)$, $z' = \piU(z)$,
the above calculation and Lemma 4.2 give
$$ I_{N,\ell}(x' ,z') = \left| \Delta(x, \pi_{\dlo}(z))  - \Delta(x, \pi_{\dlt}(z))\right| \leq C_1\, \diam(\tPsi(\cc)) .$$
The same argument shows that for any $z\in \dd$, using Lemma 5.3(b) and the fact that $Z \in \dd$,  we have
\begin{eqnarray*}
I_{N,\ell}(z',Z') =  \left| \Delta(z, \dlo) 
- \Delta(z, \dlt)\right| \leq C_1\, \diam(\tPsi(\dd))  \leq \frac{\hd}{8} \,\diam (\tPsi(\cc)) .
\end{eqnarray*}
Similarly, for any $z\in \dd_1$ and $z_0 = \T_Z(z)$ we have
$$I_{N,\ell}(z',Z') = \left| \Delta(z_0, \dlo) - \Delta(z_0, \dlt)\right| \leq C_1\, \diam(\tPsi(\dd'))  \leq \frac{\hd}{8} \,\diam (\tPsi(\cc)) .$$
Since 
$\Delta(x,\pi_{y}(z)) = \Delta(x, y) - \Delta(z, y) $
for any $y \in W^s_\ep(Z)$, it follows that
\begin{eqnarray*}
I_{N,\ell}(x' ,z') 
& =      &  I_{N,\ell}(x' , Z') -  I_{N,\ell}(z' , Z') \\
& \geq &  2\hd \, \diam(\tPsi(\cc)) - \hd \,\diam (\tPsi(\cc)) = \hd\, \diam(\tPsi(\cc)) .
\end{eqnarray*}
This completes the proof of the lemma.
\endofproof

\newpage

\section{Contraction operators}
\setcounter{equation}{0}

In this section we consider a fixed Pesin set $P_0$ with exponentially small tails, defined as in Sect. 3, and 
for  some sufficiently large value of the parameter $b$ we construct several objects that depend on $b$. First, we construct in a special
way a `sufficiently representative' compact subset $P'_0$ of $P_0$  and set $K_0 = \piU(P'_0)$, 
where $\piU : R \longrightarrow U$ is the projection along stable leaves in $R$. Then we choose a family of
unstable maximal cylinders $\cc_1, \ldots, \cc_{m_0}$ such that their projections $\cc'_m = \piU(\cc_m)$ cover $K_0$
and the {\it projections $\tcc_m = \tPsi(\cc_m)$ along the flow to $\tR$} satisfy $\diam(\tcc_m) \leq \frac{C_1}{|b|}$ 
for some constant $C_1 > 0$. The set $K_0$ plays a very important role in this and subsequent sections -- roughly speaking it is the behaviour 
`around' $K_0$ of the operators $\nn_J$ defined later which  ensures that they have some kind of contraction properties.

Next, given a large integer $N$, independent of $b$, we  construct a family of pairs of local inverses of the map $\sigma^N$ 
and use them to define `contraction' operators $\nn_J(a,b)$ for a family of symbols $J$, small $|a|$ and large $|b|$, very much as in \cite{St2}.
To prove some  contraction properties of these operators we  define a special metric $\dd$ on $\hU$, 
depending on $b$ and on the cylinders $\cc_m$ as well. It turns out that using bounded distortion and uniform expansion in this metric $\dd$,
some kind of a (partial) Lasota-Yorke type inequality can be derived for the log-H\"older constants of the operators $\nn_J$  (Lemma 6.5 below).  

We then succeed to derive that some kind of cancellation occurs in the actions of the operators $\nn_J$ (nearby the set $K_0$), and so they are indeed
contraction operators in some sense. Finally, we establish some estimates of Ruelle operators of the form $\lab^N$ by  means of the contraction operators $\nn_J$.

Although the contraction properties we derive in this section look rather weak, in Sect. 7 below 
we succeed to use  these to derive true $L^1$-contraction estimates.

Throughout this section we use the notation and the entire the set-up from Sect. 5. 

\subsection{Choice of cylinders, definition of the contraction operators}

Throughout the rest of Sect. 6, $b$ will be a {\bf fixed real number}  with $|b| \geq b_0$,
where the constant $b_0$ is defined by (6.8) below.

Recall the constants $R_0$ from (4.4) and $\ep_1$ with (5.4).
Below we will consider certain unstable cylinders $\cc$ in $R$ such that
\be
\hrho\, \frac{\ep_1}{R_0 |b|} \leq \diam(\tPsi(\cc)) \leq \frac{R_0\ep_1}{|b|}  ,
\ee
where, as before, 
$$\hrho = \frac{\rho_1}{8 C_0^2} \in (0,1) ,$$ 
and $\rho_1 \in (0,1)$ is the constant from Lemma 4.4.
Then by Lemma 5.1, if $\ell$ is the length of such $\cc$, then
$$\frac{-\log C_2 - \alpha_2 \log (R_0 \ep_1)}{|\log \htheta|} + \frac{\alpha_2}{|\log\htheta|}\, \log |b| \leq
\ell \leq \frac{\log C_2 - \log (\hrho \ep_1/R_0)}{|\log \htheta|} + \frac{1}{|\log\htheta|}\, \log |b| ,$$
where $\htheta \in (0,1)$ is given by (5.1).
Thus, there exists a global constant $B > 1$ (independent of $b$) such that if $|b| \geq 10$, then
\be\label{eq:Bcond}
\frac{1}{B} \log |b| \leq \ell \leq B\, \log |b|  .
\ee
{\bf Fix a constant } $B > 1$ with this property. Later we may have to impose some further requirements on $B$. 
Then take $A_0 > 0$ as in Lemma 5.2 such that 
$$A_0 \geq \frac{2B}{|\log \htheta|} .$$
\index{$\hrho$} 
\index{$B$} 
\index{$\Omega^{(n)}_B$}

Next, {\bf fix constants} $p_0 > 0$, $\hep, \hd_0 \in (0,1)$ as in Ch. 3 so that 
$$ \lambda_k^{\frac{4\htau_0}{\tau_0^2} \hd_0} \leq e^{\hep} .$$
(This particular choice will be used later in Sect. 9.3.2.) 
Recall the sets $\Xi_n = \Xi_n(p_0,\hep, \hd_0)$ defined by (3.20), and set
\be\label{eq:Xin}
\Omega^{(n)}_B= \ll\setminus  \bigcup_{n/B \leq \ell \leq B n} \Xi_\ell   .
\ee 
It follows from  Theorem 3.2 that, choosing the  constants $C_3, c_3 > 0$ appropriately, we have
\be\label{eq:XinB}
\mu(R\setminus \Omega^{(n)}_B) \leq C_3 \, e^{-c_3 n/B}
\ee
 for all  $n \geq 1$.

{\bf Fix integers} $d \geq 1$ and $t_0$ such that
\begin{eqnarray}\label{eq:t0condition}
c_0 \gamma^d & > & \frac{1}{\hrho} \quad, \quad 
t_0 \geq \frac{1}{\beta_1 \log \gamma} \left|\log  \frac{4 C_1 A^2 R_0^2}{\beta_0 \delta_0 \hrho^{q_1}}\right| 
+ \left|\frac{\log c_0}{\log \gamma} \right|\\  
s_0 & = & \frac{1}{|\log \rho_1|}\, \left|\log \frac{\beta_0\delta_0 \hrho^{q_1}}{2 C_1 A^5 R_0^4} \right| ,\nonumber
\end{eqnarray}
where $q_1 \geq 1$ is the constant from Lemma 5.3(b), while $c_0 > 0$ and 
$\gamma > 1$ are the constants from the end of Sect. 2. 
\index{$d$} 
\index{$t_0$} 
\index{$s_0$} 

Let $N \geq N_0$ be as in Sect. 5.  {\bf Assume also that $N > t_0 + s_0$.} 
Choose two other constants $0 < \delta_1 < \delta_2$  so small that
\be\label{eq:dcond}
\delta_1 < \frac{\hd_0}{100 d B} \quad ,\quad 
\delta_2 = (2d+3) B \delta_1 < \frac{1}{2}.
\ee
Set
\be\label{eq:mucond}
\mu_0 = \mu_0(N,\theta) = \min \left\{\;
\frac{\theta^{2N + 2d} }{6 e^{T_0/(1-\theta)}}\; , \; 
\frac{1}{10\,e^{2 T_0 N}} \, \sin^2\left(\frac{\hd \hrho \ep_1}{8 R_0}\right) \; , \; \frac{\theta^{s_0+t_0}}{100 e^{T_0/(1-\theta)}}\; \right\} ,
\ee
and
\be\label{eq:b0cond}
b_0 = b_0(N, \theta) = \max\left\{ \theta^{-N}\,  , \, \left(\frac{2 C_0 \gamma_1^{d}}{c_0 \hd_0} \right)^{1/\alpha_1} 
\, , \, R_0  \left(\frac{3C_2 T_0 e^{T_0/(1-\theta)}}{(1-\theta) } \right)^{1/\alpha_2}\right\} ,
\ee
where $\alpha_2 > 0$ is as in Lemma 5.1(c).
\index{$\delta_1$} 
\index{$\delta_2$} 
\index{$\mu_0$} 
\index{$b_0$} 

As we stated in the beginning of this section, we assume that the parameter $b \in \R$ satisfies $|b| \geq b_0$.
 Set 
\be\label{eq:hbdef}
\hb = \lceil \log |b| \rceil  .
\ee
\index{$\hb$} 

For every $z\in \ll \cap R$ denote by $\cc(z)$ {\bf the maximal cylinder} in $W^u_{R}(z)$ with
\be
\diam ((\Phi_z)^{-1}(\cc(z))) \leq \ep_1/|b| . 
\ee
\index{$\cc(z)$} 

If $\pp^{m_z}(\cc(z)) \cap P_0 \neq \e$, where $m_z$ is the length of $\cc(z)$,
then the maximality of $\cc(z)$ and Lemma 4.4(d) imply 
$$\diam ((\Phi_z)^{-1}(\cc(z)))\geq \hrho \,\ep_1/|b| .$$
Define the subset $P_1 = P_1 (b)$ of  $\ll$ by
$$P_1 = \{ z\in \ll\cap R : \pp^{m_z}(z) \in P_0 \} ,$$
and set
$$K_0 = \piU(P_1 \cap P_0 \cap \Omb_B) .$$ 
\index{$K_0$} 
\index{$P_1$} 
\index{$Z(u)$} 

Next, we define an important family of cylinders in  $R$ and $U$ and some sub-cylinders of
theirs that will play an important role throughout Sect. 6 and 7.

\bs

\noindent
{\bf Definition 6.1(Choice of cylinders):} 
For any $u \in K_0$ amongst the cylinders 
$\cc(z)$ with $z\in P_1\cap P_0 \cap \Omb_B$ and $\piU(z) = u$, there 
is {\bf one of maximal length} (and so of smallest $\diam_\theta$). Choose one of these -- it has the form $C(Z(u))$ for some 
$Z(u) \in P_1 \cap P_0\cap \Omb_B$ with $\piU(Z(u)) = u$. Set $\cc'(u) = \piU(C(Z(u)))$. 
It follows from this choice that for any $z\in R$ with $\piU(z) = u$ we have $\cc'(u) \subset \piU(\cc(z))$.

Since the lengths of the cylinders $C'(u)$ are bounded above and 
$K_0 \subset \cup_{u\in K_0} \cc'(u)$, there exist
finitely many different cylinders $\cc'_m = \cc'(u_m)$ for some $m =1, \ldots, m_0$ such that
$$K_0 \subset \cup_{m=1}^{m_0} \cc'_m .$$
Different cylinders have no common interior points, so
$\cc'_m \cap \cc'_{m'} \cap \hU = \e$ for $m \neq m'$. For each $m$, set $\cc_m = \cc(Z(u_m))$; then
$$Z_m = Z(u_m) \in \cc_m \cap P_1 \cap P_0 \cap \Omb_B$$ 
is so that $\piU(Z_m) = u_m$. 
According to the definitions of the cylinders $\cc(z)$, $\cc_m = \cc(Z_m)$ is a 
{\it maximal  closed cylinder in $W^u_R(Z_m)$ with} 
$$\diam ((\Phi_z)^{-1}(\cc_m)) \leq \ep_1/|b| .$$
Let 
$\dd_1, \ldots, \dd_{\tilde{j}}$
be the list of {\it all closed unstable cylinders  in $R$ which are sub-cylinders of co-length} $q_1$ of some $\cc_m$. 
Here $q_1 \geq 1$ is the constant from Lemma 5.3(b). Set $\dd_j' = \piU(\dd_j) \subset U$.
Re-numbering the cylinders  $\dd_j$ if necessary, we may assume there exists $j_0 \leq \tilde{j}$ such that
$$\dd_1, \ldots, \dd_{j_0}$$ 
is {\it the list of all sub-cylinders $\dd_j$ } with 
$$\dd_j \cap P_1 \cap P_0 \cap \Omb_B \neq \e .$$
\index{$\dd_j$} 
\index{$\cc_m$}
\index{$V_b$} 

\ms

From the choice of the cylinders $\cc_m$  and Lemmas 4.4 and  5.3, and using (3.7) as well and the
fact that $R(x) \leq R_0$ on $P_0$,  we get:
\be
\hrho\, \frac{\ep_1}{R_0 |b|} \leq \diam(\tPsi(\cc_m)) \leq \frac{R_0 \ep_1}{|b|}  \quad , \quad 1 \leq m \leq m_0 .
\ee
If $\ell_m$ is {\it the length of the cylinder} $\cc_m$, it follows from  (\ref{eq:Bcond}) that
\be
\frac{1}{B} \log |b| \leq \ell_m \leq B\, \log |b| \quad , \quad m = 1, \ldots, m_0 .
\ee

\index{$\vli = \vli(Z_m, \cdot)$} 
\index{$\xijl$}

Set 
\be
V_b = \cup_{j=1}^{j_0} \dd'_j  \subset U .
\ee
It follows from the construction that
$$K_0 \subset V_b  .$$

We are now ready to define an important family of {\it contraction operators}.
For any $\ell = 1, \ldots, \ell_0$, $i = 1,2$ and $j = 1, \ldots, j_0$, consider the unique 
$m = 1, \ldots, m_0$ with $\dd_j \subset \cc_m$, and set 
$$\vli = \vli(Z_m, \cdot) \quad , \quad \xijl = \vli(\dd'_j) \subset U ,$$
where $\vli(Z_m, \cdot)$ is the map from Lemma 5.4 for the integer $N$. We will consider this map only on $\cc'_m$.
By Lemma 5.1(a), the {\it characteristic function} 
$$\omijl = \chi_{\xijl} : \hU \longrightarrow [0,1]$$ 
of $\xijl$ belongs to $\ff_\theta(\hU)$ and  
$$\Lip_\theta(\omijl) \leq 1/\diamte(\xijl) .$$

\index{$\omijl$} 


A subset $J$ of the set
$$\Pi(b)  = \{\; (i, j , \ell) \; :  \; 1\leq i \leq 2\; ,\;  1\leq j\leq j_0 \;,  \; 1 \leq \ell \leq \ell_0 \;\} $$
will be called {\it representative} 
\index{representative set} if for every $j = 1, \ldots, j_0$ there exists 
at most one pair $(i,\ell)$ such that
$(i,j,\ell) \in J$, and for any $m = 1, \ldots, m_0$ there exists $(i,j, \ell) \in J$ 
such that $\dd_j \subset \cc_m$.
Let $\jj(b)$ be the {\it family of all representative subsets} $J$  of $\Pi(b)$.
\index{$\jj(b)$} 
\index{$\Pi(b)$}

Given $J  \in \jj(b)$, define the function    $\omega_{J}  : \hU \longrightarrow [0,1]$ by
$$\di \omega_J  = 1- \mu_0 \,\sum_{(i, j,\ell) \in J} \omijl .$$
Clearly $\omega_J \in \ff_\theta(\hU)$ and 
$$\frac{1}{2} \leq 1-\mu_0 \leq \omega_J (u) \leq 1 \quad , \quad u \in \hU .$$
Define the {\it contraction operator} 
\index{contraction operator}
$$\nn = \nn_J(a,b) : \ff_\theta (\hU) \longrightarrow \ff_\theta (\hU) \quad  \mbox{\rm by } \quad \nn h = \ma^{N} (\omega_J \cdot h) ,$$
where $\ma = L_{\fa}$ is defined in Sect. 5.1.
\index{$\nn_J = \nn_J(a,b)$} 
\index{$\omega = \omega_J$}

\subsection{Main properties of the contraction operators}

First, we  derive an important consequence of the construction of the cylinders $\cc_m$ and $\dd_j$.

\bs

\noindent
{\bf Lemma 6.2.} {\it If $\sigma^p(\dd'_j) \subset \cc'_{k}$  for some $p\geq 0$, $j \leq j_0$ and  
$k \leq m_0$, then $p \leq t_0$, where $t_0$ is as in {\rm (6.5)}. Moreover, the co-length of $\sigma^p(\dd'_j)$ in $\cc'_k$ 
does not exceed the constant $s_0$ from {\rm (6.5)}.}

\bs

\noindent
{\it Proof.} Assume  $\dd = \sigma^p(\dd'_j) \subset \cc'_{k}$  for some $p > 0$, $j \leq j_0$ and  $k \leq m_0$. 
From the assumptions we get $\pi_{Z_k}(\dd) \subset \cc_{k} \subset W^u_R(Z_k)$.

According to the choice of the sub-cylinders $\dd_j$, there exists 
$$Z \in \dd_j \cap P_1 \cap P_0\cap \Omb_B .$$
Then using Lemma 4.3(a) with $\kappa = 1/2$ and Lemma 4.4(c) and an appropriately chosen 
$X \in \T_Z(\dd_j) \subset W^u_{\ep}(Z)$ with 
\be
d(X, Z) \geq \frac{1}{2} \diam(\T_z(\dd_j)) \geq \frac{\hrho^{q_1}\ep_1}{2A R_0 |b|} ,
\ee
we can find points $d_1,d_2 \in W^s_R(Z)$ such that
$$|\Delta(X, \pi_{d_1}(Z)) - \Delta (X, \pi_{d_2}(Z))| \geq \frac{\beta_0 \delta_0}{2R_0}\, \diam (\T_Z(\dd_j)) ,$$
that is
\be 
|\Delta(X, d_1) - \Delta (X, d_2)|
\geq \frac{\beta_0 \delta_0}{4 A R_0^2}\, \frac{\hrho^{q_1}\ep_1}{|b|} .
\ee

Let $T = \tau_p(Z)$; then $z = \pp^p(Z) = \phi_T(Z)$. 
Next, consider the points 
$$d'_i = \phi_T(d_i) \in W^s_R(z) \quad , \quad i = 1,2 ,$$
and
$$x = \phi_T(X) \in W^u_{\ep}(z) \subset \phi_T(W^u_{\ep}(X)) .$$
It follows from (6.15) and the properties of temporal distance that
\be
|\Delta(x, d'_1) - \Delta (x, d'_2)|
\geq \frac{\beta_0 \delta_0}{4 A R_0^2}\, \frac{\hrho^{q_1}\ep_1}{|b|} ,
\ee
while (2.1) yields
\be
d(d'_1, d'_2) = d(\pp^p(d_1), \pp^p(d_2)) \leq \frac{1}{c_0 \gamma^p}\, d(d_1,d_2)
\leq \frac{1}{c_0 \gamma^p} .
\ee
Let $y \in W^s_R(z)$ be such that $\cc_k \subset W^u_R(y)$. Since $\pp^p(\dd'_j) \subset \cc'_k$,
we have $y \in \cc_k$. For the point $x' = \pi_y(x) \in W^{uc}_{\ep}(y)$ we have $\phi_t(x') \in \cc_k$ for some $t\in \R$, so 
$$x'' = \T_y(\phi_t(x')) \in \T_y(\cc_k) .$$
Moreover it is easy to see, using just the definition of the temporal distance function
and the fact that $x'' = \phi_s(x')$ for some $s\in \R$, that
$$|\Delta(x, d'_1) - \Delta(x, d'_2)| = |\Delta(x', d'_1) - \Delta(x', d'_2)| = |\Delta(x'', d'_1) - \Delta(x'', d'_2)| .$$
This and (6.16) give
$$|\Delta(x'', d'_1) - \Delta(x'', d'_2)|  \geq \frac{\beta_0 \delta_0 \hrho^{q_1}}{4A R_0^2}\, \frac{\ep_1}{|b|} .$$
Combining the latter with 
$$\diam(\T_y(\cc_k)) \leq A R_0\ep_1/|b| ,$$ 
$y, x'' \in \T_y(\cc_k)$, and $d'_1, d'_2 \in W^s_R(y)$, and using the fact that $\cc_k \cap P_1 \cap \Omb_B \neq \e$, Lemma 4.3(b) implies
\begin{eqnarray*}
\frac{\beta_0 \delta_0 \hrho^{q_1}}{4 A R_0^2}\, \frac{\ep_1}{|b|}
&  \leq & 
|\Delta(x'', d'_1) - \Delta(x'', d'_2)| \leq C_1 \diam(\T_y(\cc_k)) \, (d(d'_1, d'_2))^{\beta_1}\\
& \leq & C_1 A R_0 \frac{\ep_1}{|b|} \, (d(d'_1, d'_2))^{\beta_1} .
\end{eqnarray*}
This and (6.17) give
$$\left(\frac{\beta_0 \delta_0 \hrho^{q_1}}{4 C_1 A^2 R_0^3}\right)^{1/\beta_1} \leq 
\frac{1}{c_0 \gamma^p} ,$$
so $p \leq t_0$, where $t_0 > 0$ is the integer from (6.5).

Next, let $s$ be the co-length of $\sigma^p(\dd'_j)$ in $\cc'_k$. Denote by $Q$ the cylinder
in $W^u_R(z)$ such that $Q \parallel \cc_k$, i.e. $\piU(Q) = \piU(\cc_k)$. Then $\pp^p(\dd_j)$ is a sub-cylinder of $Q$
of co-length $s$, so $\dd_j$ is a sub-cylinder of co-length $s$ of $Q' = \pp^{-p}(Q)$. 
Since 
$$Z \in \dd_j \cap P_1 \cap P_0 \cap \Omb_B ,$$ 
it follows from Lemma 4.3(a) with $\kappa = 1/2$ that there exist $x_0 \in \T_Z(Q')$ and $y_1,y_2 \in W^s_R(Z)$ such that
$$\frac{\beta_0\delta_0}{2 A R_0^2} \, \diam(\T_Z(Q')) \leq |\Delta(x_0, y_1) - \Delta(x_0,y_2)| .$$
Setting $x'_0 = \phi_T(x_0) \in \T_z(Q)$ and $y'_i = \phi_T(y_i) \in W^s_R(z)$, $i = 1,2$, we have 
$$|\Delta(x_0, y_1) - \Delta(x_0,y_2)| = |\Delta(x'_0, y'_1) - \Delta(x'_0,y'_2)| , $$ 
so
$$\frac{\beta_0\delta_0}{2 A R_0^2} \, \diam(\T_Z(Q')) \leq |\Delta(x'_0, y'_1) - \Delta(x'_0,y'_2)| .$$
As above, denoting by $x''_0 \in \T_{Z_k}(\cc_k)$ the unique point such that $x''_0 \in W^{sc}(x'_0)$, and using Lemma 4.3(b), we get
\begin{eqnarray*}
\frac{\beta_0\delta_0}{2 A R_0^2} \, \diam(\T_Z(Q'))
& \leq & |\Delta(x'_0, y'_1) - \Delta(x'_0,y'_2)| = |\Delta(x''_0, y'_1) - \Delta(x''_0,y'_2)|\\
& \leq & C_1 \diam(\T_{Z_k}(\cc_k)) ,
\end{eqnarray*}
so
$$\frac{\beta_0\delta_0}{2 A R_0^2} \, \diam(\T_Z(Q')) \leq  \frac{C_1 A R_0 \ep_1}{|b|} 
\leq  \frac{C_1A^2 R_0}{\hrho^{q_1}} \frac{\hrho^{q_1} \ep_1}{A |b|} \leq \frac{C_1 A^2 R^2_0}{\hrho^{q_1}} \, \diam(\T_Z(\dd_j)) .$$
On the other hand, it follows from Lemma 4.4(b) with $\rho' = \rho_1$  that
$$\diam(\T_Z(\dd_j)) \leq A^2  \rho_1^s\, \diam(\T_Z(Q')) .$$
Thus,
$$\frac{\beta_0\delta_0}{2 A R_0^2} \, \diam(\T_Z(Q')) 
\leq \frac{C_1 A^4 R_0^2}{\hrho^{q_1}} \,\rho_1^s\,  \diam(\T_Z(Q')) ,$$
so $\frac{\beta_0\delta_0 \hrho^{q_1}}{2 C_1 A^5 R_0^4} \leq \rho_1^s$, which implies
$s \leq s_0 $, the constant from (6.5).
\endofproof

\index{$\ell(u,u')$} 
\index{$Y(u,u')$}
\index{$\dd$}  

\bs

Given $u, u'\in \hU$, we will denote by $\ell(u,u') \geq 0$  {\it the length of the smallest cylinder $Y(u,u')$   in $\hU$ containing} $u$ and $u'$.

\bs

\noindent
{\bf Definition 6.3.} 
Define the {\bf distance} $\dd(u,u')$  for $u,u' \in \hU$ by\footnote{Clearly $\dd$ depends on  the cylinders $\cc_m$ and therefore
on the parameter $b$ as well.}: 

(i) $\dd(u,u') = 0$ if $u = u'$; 

(ii)  Let $u \neq u'$, and assume there exists 
$p \geq 0$ with $\sigma^p(Y(u,u')) \subset \cc'_m$, $\ell(u,u') \geq p$, for some $m = 1, \ldots, m_0$. 
Take the maximal $p$ with this property and the corresponding $m$ and set 
$$\dd(u,u') = \frac{\dte(u,u')}{\diamtef (\cc_m)} ;$$

(iii) Assume $u \neq u'$,
however there is no $p \geq 0$ with the property described in (ii). Then set $\dd(u,u') = 1$.

Notice that $\dd(u,u') \leq 1$ always. 
Some other properties of $\dd$ are contained in the following, part (b)  of which needs Lemma 6.2.

\bs


\noindent
{\bf Lemma 6.4.} {\it Assume that $u,u' \in \hU$, $u \neq u'$, and $\sigma^N(v) = u$, $\sigma^N(v') = u'$
for some $v,v'\in \hU$ with $\ell(v,v') \geq N$. Assume also that there exists 
$p \geq 0$ with $\sigma^p(Y(u,u')) \subset \cc'_m$, $\ell(u,u') \geq p$, for some $m = 1, \ldots, m_0$. }

\ms

(a) {\it We have $\dd(v,v') = \theta^N\, \dd (u,u')$.}

\ms

(b) {\it  Assume in addition that $\omega_J(v) < 1$ and $\omega_J(v') = 1$ for some $J \in \jj(b)$. 
Then $p \leq t_0$ and}
$$|\omega_J (v) - \omega_J (v')| \leq \frac{\mu_0}{\theta^{t_0 + s_0 }}\, \dd(u,u') .$$

\ms

\noindent
{\it Proof.} 
(a) Let $p$ be the maximal integer with the given property and let $m \leq m_0$ correspond to $p$. 
Then $\sigma^{p+N}(Y(v,v')) \subset \cc'_m$, $\ell(v,v') \geq p+N$, and $p+N$ is the maximal 
integer with this property. Thus,
$$\dd(v,v') = \frac{\dte(v,v')}{\diamtef (\cc_m)} = \theta^N \, \frac{\dte(u,u')}{\diamtef (\cc_m)} 
= \theta^N \, \dd(u,u') .$$

(b)  $\omega_J(v) < 1$ means that $v\in \Xijl$ for some $(i,j,\ell) \in J$, and so 
$v = \vl_i(u)$ for some $u \in \dd'_j$. Then  $u = \sigma^N(v)$.  
If $u' \in \dd'_j$, then $v'' = \vl_i(u') \in \Xijl$ and $\sigma^N (v'') = u'$, so we must have
$v'' = v'$, which implies $\omega_J(v') = \omega_J(\vl_i(u')) < 1$, a contradiction.
This shows that $u'\notin \dd'_j$, and so 
$$\dte(u,u') \geq \diamte(\dd'_j) .$$

Since $u\in \dd'_j$, $u' \notin \dd'_j $ and $\ell(u,u') \geq p$, it follows that 
$\sigma^p(u) \in \sigma^p(\dd'_j)$ and $\sigma^p(u') \notin \sigma^p(\dd'_j) $. 
On the hand, by assumption, $\sigma^p(u), \sigma^p(u') \in \cc'_m$.  Thus, the cylinder
$\sigma^p(\dd'_j)$ must be contained in $\cc'_m$. Now
Lemma 6.2 gives $p \leq t_0$ and the co-length $s$ of $\sigma^p(\dd'_j)$ in $\cc'_m$ is
$s \leq s_0$. If $\ell_m = \length(\cc_m)$, and $\ell = \length(\dd_j)$ we have 
$$\ell-p -s = \length (\sigma^p(\dd'_j)) -s = \length (\cc'_m) = \ell_m .$$ 
Hence
$$\ell = \ell_m + p + s \leq \ell_m + t_0 + s_0 ,$$
and using
$$\dd(u,u') = \frac{\dte(u,u')}{\diamtef (\cc_m)} ,$$ 
we get
\begin{eqnarray*}
|\omega (v) - \omega (v')| 
& = &  \mu_0 = \mu_0 \, \frac{\dte(u,u')}{\dte(u,u')}  \leq \mu_0 \frac{\dte(u,u')}{\diamtef(\dd'_j) } 
=  \mu_0 \frac{\dte(u,u')}{  \theta^{\ell} }\\
& \leq & \mu_0 \frac{\dte(u,u')}{\theta^{\ell_m + t_0 + s_0}}   = \mu_0 \frac{\dte(u,u')}{\theta^{t_0 + s_0}\diamtef (\cc_m) } 
= \frac{\mu_0}{ \theta^{t_0 + s_0}} \, \dd(u,u').
\end{eqnarray*}
This proves the lemma. 
\endofproof

\index{$\kk_E$}

\bs

Given $E > 0$ as in Sect. 5.1, let $\kk_E$ 
be {\it the set of all functions $H \in \ff_\theta(\hU)$  such that } $H > 0$ on $\hU$ and 
$$\frac{|H(u) - H(u')|}{H(u')} \leq E\, \dd(u,u')$$ 
for all $u,u'\in \hU$ for which there exists an integer $p \geq 0$ with $\sigma^p(Y(u, u')) \subset \cc_m$ 
for some $m \leq m_0$ and $\ell(u,u') \geq p$.

Using Lemma 6.4 we will now prove the main lemma in this section, which makes it possible to use an
inductive procedure involving the contraction operators $\nn_J$. As the reader would observe, in the proof below
we essentially derive a Lasota-Yorke type inequality, using bounded distortion and uniform expansion in the metric $\dd$.
If one disregards the assumption that we only consider pairs $u,u'\in \hU$  with $\sigma^p(Y(u, u')) \subset \cc_m$ for some 
$p$ and $m$, this would look like an inequality for the log-H\"older constant of $\nn H$.

\bs

\noindent
{\bf Lemma 6.5.} {\it For any $J \in \jj(b)$ we have  $\nn_J(\kk_E) \subset \kk_{E}$.} 

\bs

\noindent
{\it Proof.} 
Let  $u, u' \in \hU$  be such that  there exists
an integer $p \geq 0$ with $\sigma^p(Y(u, u')) \subset \cc_m$ for some $m = 1, \ldots,  m_0$ and $\ell(u,u') \geq p$.

Given $v \in \hU$ with $\sigma^{N}(v) = u$, let $C[\ii] = C[i_0, \ldots,i_{N}]$ 
be the  cylinder of length $N$ containing $v$.  
Set $\hC[\ii] = C[\ii]\cap \hU$. Then $\sigma^{N}(\hC[\ii]) = \hU_i$. Moreover, 
$\sigma^N : \hC[\ii] \longrightarrow \hU_i $ is a homeomorphism, so
there  exists a unique $v' = v'(v)\in \hC[\ii]$ such that $\sigma^N(v') = u'$.  Then 
$$\dte (\sigma^j(v),\sigma^j(v'(v))) = \theta^{N-j}\, \dte (u,u')$$ 
for all $j = 0,1, \ldots, N-1$. 
Also $\dte(v,v'(v))= \theta^N \dte (u,u')$ and $\dd(v,v'(v)) = \theta^N \dd(u,u')$. Using (5.3), we get
\begin{eqnarray}
&        & |\fa_N(v) - \fa_N(v')| 
 \leq  \sum_{j=0}^{N-1} |\fa(\sigma^j(v)) - \fa (\sigma^j(v'))|\\
& \leq & \sum_{j=0}^{N-1} |\fa|_\theta \,\theta^{N-j}\, \dte (u,u')
 \leq  \frac{T_0}{1-\theta}\, \dte (u,u'). \nonumber
\end{eqnarray}

 Let $J \in \jj(b)$ and let $H \in \kk_E$. Set $\nn = \nn_J$. We will show that $\nn H \in \kk_{E }$.

Using the above  and the definition of $\nn = \nn_J$, and setting  $v' = v'(v)$ and $\omega = \omega_J$ for brevity, we get
\begin{eqnarray*}
&        & \frac{|(\nn H)(u) - (\nn H)(u')|}{\nn H (u')} \\
& =     & \frac{\di \left| \sum_{\sigma^N v = u} e^{\fa_N(v)}\, \omega(v) H(v) -  
\sum_{\sigma^N v = u} e^{\fa_N(v'(v))}\, \omega(v'(v)) H(v'(v)) \right|}{\nn H (u')} \\
& \leq & \frac{\di \left| \sum_{\sigma^N v = u}  e^{\fa_N(v)}\, [\omega(v) H(v) -  \omega(v') H(v')]\right|}{\nn H (u')}\\ 
&        &  + \frac{\di \sum_{\sigma^N v = u}   \left|e^{\fa_N(v)}-  e^{\fa_N(v')}\right| \, \omega(v') H(v')}{\nn H (u')} \\
& \leq & \frac{\di \sum_{\sigma^N v = u} e^{\fa_N(v) - \fa_N(v')}  e^{\fa_N(v')}\, |\omega(v) - \omega(v')| H(v')}{\nn H (u')}\\
&        &    + \frac{\di \sum_{\sigma^N v = u} e^{\fa_N(v)}\, \omega(v) | H(v) -   H(v')|}{\nn H (u')} \\
&       & + \frac{\di \sum_{\sigma^N v = u}  \left| e^{\fa_N(v) - \fa_N(v')} - 1\right| 
\,  e^{\fa_N(v')} \omega(v') H(v')}{\nn H (u')} .
\end{eqnarray*}

By the definition of $\omega$, either $\omega(v) = \omega(v')$ or at least one of these numbers is  $ < 1$.
Using Lemma 6.4 we then get 
$|\omega(v) - \omega(v')| \leq \frac{\mu_0}{\theta^{t_0 + s_0 }} \dd(u,u')$. Apart from that
$H \in \kk_E$ implies 
$$|H(v) - H(v')| \leq E H(v') \dd(v,v') = E H(v') \theta^N \dd(u,u') ,$$
while 
$$\left| e^{\fa_N(v) - \fa_N(v')} - 1\right| \leq e^{T_0/(1-\theta)} \frac{T_0}{1-\theta}\, \dte (u,u') .$$
Thus,
\begin{eqnarray*}
         \frac{|(\nn H)(u) - (\nn H)(u')|}{\nn H (u')} 
& \leq & e^{T_0/(1-\theta)} \frac{\mu_0}{\theta^{t_0 + s_0 }} \,  \frac{\di \sum_{\sigma^N v = u} e^{\fa_N(v')}\, \dd(u,u')\, H(v')}{\nn H (u')} \\
&        &  + \frac{\di \sum_{\sigma^N v = u}  e^{\fa_N(v)- \fa_N(v')}\,e^{\fa_N(v')}\,  2\omega(v')\, E H(v') \theta^N \dd(u,u')}{\nn H (u')} \\
&        &  + e^{T_0/(1-\theta)} \frac{T_0}{1-\theta}\, \dte (u,u')\\
& \leq & 2 e^{T_0/(1-\theta)} \frac{\mu_0}{\theta^{t_0 +s_0 }} \, \dd(u,u')  + 2 e^{T_0/(1-\theta)} E \theta^N \, \dd(u,u')\\
&       &  +  e^{T_0/(1-\theta)} \diamte(\cc_m) \frac{T_0}{(1-\theta) }\, \dd (u,u')
 \leq E \, \dd(u,u') , 
\end{eqnarray*}
where we used (\ref{eq:mucond}) and Lemma 5.1, and assumed
$$2 e^{T_0/(1-\theta)} \theta^N \leq 1/3 \quad , \quad 
e^{T_0/(1-\theta)} C_2 (R_0 \ep_1/|b|)^{\alpha_2} \frac{T_0}{(1-\theta) } \leq \frac{1 }{3}\leq \frac{E}{3} .$$
The latter follows from $|b| \geq b_0$ and (\ref{eq:b0cond}). Hence $\nn H \in \kk_E$.
\endofproof


\def\tkk{\widetilde{\kk}}

\subsection{Return times estimates}

Recall the numbers $\theta_1, \theta_2 \in (0,1)$ defined in the beginning of Sect. 5.2.
Then using the proof of Lemma 5.1(c) and taking $C_2 > 0$ sufficiently large we have\footnote{Notice that
for (6.19) choosing $\theta_1$ with $\theta_1 ^{\alpha_2} \leq \theta$ would be enough. However in the beginning
of Sect. 6.1 we imposed a stronger condition on $\theta_1$ which will be used later on (see the end of the proof
of Theorem 1.3 in Ch. 8).}
\be
\diam_{\theta_1}(\cc) \leq C_2\, \diam(\tPsi(\cc))
\ee
for any cylinder $\cc$ in $U$.

Throughout the rest of this section {\bf we assume that} $f \in \ff_{\theta_1}(\hU)$.

\medskip

Given a cylinder $\cc_m$, recall the point $Z_m \in \cc_m$ chosen in Sect. 6.1.
For points $u,u' \in U$ we will denote $\tu = \T_{Z_m}(\pi_{Z_m}(u))$ and $\tu' = \T_{Z_m}(\pi_{Z_m}(u'))$; 
these are then points on the true unstable manifold $W^u_{\ep_0}(Z_m)$. In this section we will 
frequently work under the following assumption for points $u,u' \in \hU$ contained in some cylinder 
$\cc'_m = \piU(\cc_m)$ ($1\leq m \leq m_0$), an integer $p \geq 0$ and points $v,v'\in \hU$:
\be\label{eq:u-cond}
u,u'\in \cc'_m \:  , \: \sigma^p(v) = \vli(u)\;, \; \sigma^p(v')  = \vli(u')\;,\; \ell(v,v') \geq N  ,
\ee
for some $i = 1,2$.
From (\ref{eq:u-cond}) we get $\ell(v,v') \geq N + p$ and 
$$\sigma^{N+p}(v) = u \quad , \quad \sigma^{N+p}(v') = u' .$$
We will use the notation 
$$\tcc_m = \tPsi(\cc_m) \subset \tR .$$
\index{$\tc_m$} 

The following estimate plays a central role in this section.

\bs

\noindent
{\bf Lemma 6.6.} {\it There exists a global constant $C_4 > 0$ independent of $b$ and $N$ 
such that if the points  $u,u' \in \hU$, the  cylinder $\cc_m$, the integer 
$p \geq 0$ and the points $v,v'\in \hU$ satisfy {\rm (\ref{eq:u-cond})}  for some 
$i = 1,2$ and $\ell = 1, \ldots, \ell_0$, and $w,w'\in \hU$ are such that $\sigma^N w = v$, 
$\sigma^N w' = v'$ and $\ell(w,w') \geq N$, then}
$$|\tau_N(w) - \tau_N(w')| \leq C_4 \, \theta_2^{p+N} \, \diam(\tcc_m) .$$
\index{$C_3$} 

\ms

\noindent
{\it Proof.} Assume that the points $u,u',v,v',w,w'$ and the cylinder $\cc$
satisfy the assumptions in the lemma. Clearly, $\ell(w,w') \geq p+2N$ and
\be
\tau_N(w) - \tau_N(w') = [\tau_{p+2N}(w) - \tau_{p+2N}(w')] - [\tau_{p+N}(v) - \tau_{p+N}(v')] .
\ee

Recall the construction of the map $\vl_i$ from the proof of Lemma 5.4. In particular by (5.8),
$\pp^N(\vl(u)) = \pi_{\dli} (u)$, where we set $\dli = \dli(Z_m) \in W^s_{R}(Z_m)$ for brevity. 
Since $\sigma^{p}(v) = \vl_i(u)$ and $\sigma^p(v') = \vl_i(u')$, we have $\sigma^{p+N}(v) = u$ and 
$\sigma^{p+N}(v') = u'$, so 
$$\pp^{p+N}(v) \: , \: \pp^{p+N}(v') \in W^u_{R}(d')$$ 
for some $d' \in W^s_{R}(Z_m)$.
Moreover, $\pp^p(v) \in W^s_{R_1}(\vl_i(u))$ and the choice of $N$ imply (as in the proof of Lemma 5.4) that 
$d(\dli, d') < \delta''$, the constant from Lemma 4.3.
Similarly, 
$$\pp^{p+2N}(w) \: , \: \pp^{p+2N}(w') \in W^u_{R}(d'')$$ 
for some $d'' \in W^s_{R}(Z_m)$ with
$d(\dli, d'') < \delta''$.  Moreover, since the local stable/unstable holonomy maps are uniformly 
$\alpha_1$-H\"older (by the choice of $\alpha_1$), there exists a global constant $C'_4 > 0$ such that 
$$d(d', d'') \leq C'_4 (d(\pp^{p+N}(v), \pp^{p+2N}(w)))^{\alpha_1} .$$
Using this and (2.1) for points on local stable manifolds, i.e. going backwards along the flow, we get
$$d(d', d'') \leq C'_4 (d(\pp^{p+N}(v), \pp^{p+2N}(w)))^{\alpha_1} 
\leq C'_4 \left( \frac{d(v, \pp^N(w))}{c_0 \gamma^{p+N}} \right)^{\alpha_1} 
\leq \frac{C'_4}{c_0^{\alpha_1} \gamma^{\alpha_1 (p+N)}}  .$$
Hence
$$(d(d', d''))^{\beta_1} \leq  (C'_4/c_0^{\alpha_1})^{\beta_1} (1/\gamma^{p+N})^{\alpha_1 \beta_1} 
\leq  C''_4 \theta_2^{p+N} .$$

We are preparing to use Lemma 4.2. Set $\hu = \pi_{Z_m}(u)$ and $\hu' = \pi_{Z_m}(u')$,
and let $\tu \in R$ and $\tu' \in R$ be the shifts along the flow of the points $\hu$ and $\hu'$.
Then we have $\tu = \phi_{t(u)} (\hu)$ and  $\tu' = \phi_{t(u')} (\hu')$ 
for some small $t(u), t(u') \in \R$. So
\begin{eqnarray*}
\tau_{p+N}(v) - \tau_{p+N}(v') 
& = & \Delta(\pp^{p+N}(v), \pp^{p+N}(v')) = \Delta(u, \pi_{d'}(u')) = \Delta(\hu, \pi_{d'}(\hu'))\\
& = & \Delta(\tu, \pi_{d'}(\tu')) + t(u) - t(u') ,
\end{eqnarray*}
and similarly
\begin{eqnarray*}
\tau_{p+2N}(w) - \tau_{p+2N}(w') 
& = & \Delta(\pp^{p+2N}(w), \pp^{p+2N}(w')) = \Delta(\hu, \pi_{d''}(\hu'))\\
& = & \Delta(\tu, \pi_{d''}(\tu')) + t(u) - t(u') .
\end{eqnarray*}
This, (6.21), Lemma 4.2  and the above estimate yield
\begin{eqnarray*}
|\tau_N(w) - \tau_N(w')|
& =     & |[ \Delta(\tu, \pi_{d'}(\tu')) + t(u) - t(u')] - [ \Delta(\tu, \pi_{d''}(\tu')) + t(u) - t(u')]|\\
& =     & |  \Delta(\tu, \pi_{d'}(\tu')) -  \Delta(\tu, \pi_{d''}(\tu'))| 
\leq C_1 \diam (\tc_m) \, (d(d',d''))^{\beta_1}\\
& \leq & C_1 C''_4 \theta_2^{p+N} \, \diam(\tc_m) .
\end{eqnarray*}
This proves the lemma.
\endofproof

\bs

Set $M_1 = M_0 + a_0$ (see Sect. 2 for the choice of $M_0$). Let  $T_0 $ be as in (5.3), and let 
$$E_1 = 2C_5 e^{C_5} \quad \mbox{\rm where } \quad C_5 = \frac{T_0 C_2}{1-\theta} + M_1 C_4 ,$$ 
and $C_4 > 0$ is the constant from Lemma 6.6. Assume $N$ is so large that
$$\theta_2^N e^{C_5} \leq \frac{1}{2} .$$
\index{$M_1$} 
\index{$E_1$} 
\index{$C_4$} 

Denote by $\kk_0$
\index{$\kk_0$} 
{\it the set of all $h \in \ff_{\theta}(\hU)$ such that $h \geq 0$ on $\hU$ and for any 
$u,u' \in \hU$ contained in some cylinder $\cc'_m = \piU(\cc_m)$ ($1\leq m \leq m_0$), any integer 
$p \geq 0$ and any points $v,v'\in \hU$ satisfying {\rm (\ref{eq:u-cond})} for some $i = 1,2$ and 
$\ell = 1, \ldots, \ell_0$ we have}
\be
|h(v) - h(v') | \leq E_1 \, \theta_2^{p+N} \, h(v')\, \diam(\tcc_m) .
\ee

We are going to show that the eigenfunctions $h_a \in \kk_0$ for $|a| \leq a_0$ (see Sect. 5.1). 
This will be derived from the following.

\bs

\noindent
{\bf Lemma 6.7.} {\it For any real constant $s$ with $|s| \leq M_1$ we have 
$L^{qN}_{f - s \tau}(\kk_0) \subset \kk_0$ for all integers $q \geq 1$.} 

\bs

\noindent
{\it Proof.} We will use Lemma 6.6 and a standard argument.

Assume that $u,u' \in \hU$, the  cylinder $\cc'_m$ in $U$, the integer $p \geq 0$ 
and the points $v,v'\in \hU$ satisfy (\ref{eq:u-cond}) for some $i = 1,2$ and $\ell = 1, \ldots, \ell_0$, and 
$w,w'\in \hU$ are such that $\sigma^N w = v$, $\sigma^Nw' = v'$ and $\ell(w,w') \geq N$; then $w' = w'(w)$ is
uniquely determined by $w$. 

Using $f \in \ff_{\theta_1}(\hU)$, the choice of $\theta_1$ and (6.19), we get 
\begin{eqnarray*}
|f_N (w) - f_N(w')| 
& \leq & \frac{T_0}{1-\theta_1} D_{\theta_1}(v,v') = \frac{T_0}{1-\theta_1} \theta_1^{p+N}\, D_{\theta_1}(u,u')\\
& \leq &  \frac{T_0}{1-\theta_1} \theta_1^{p+N}\, \diam_{\theta_1}(\cc_m) \leq C'_5 \theta_2^{p+N} \, \diam(\tcc_m) ,
\end{eqnarray*}
where $C'_5 = C_2 T_0/ (1-\theta_1)$. This and Lemma 6.6 imply
$$|(f-s\tau)_N(w) - (f- s\tau)_N(w')| \leq C_5 \theta_2^{p+N} \, \diam(\tcc_m) $$
for all $s\in \R$ with $|s| \leq M_1$, where $C_5 > 0$ is as above.

Thus, given $s$ with $|s| \leq M_1$ and $h \in \kk_0$ we have:
\begin{eqnarray*}
&        & |(L^N_{f-s\tau}h)(v) - (L^N_{f- s\tau} h)(v')|\\
& =     &  \left| \sum_{\sigma^N w = v} e^{(f-s \tau)_N(w)}\, h(w)   -  \sum_{\sigma^N w = v} e^{(f - s \tau)_N(w'(w))}\, h(w'(w)) \right| \\
& \leq & \left| \sum_{\sigma^N w = v} e^{(f- s\tau)_N(w)}\, [h(w) -  h(w')]\right| \\ 
&       & + \sum_{\sigma^N w = v}  \left|e^{(f- s\tau)_N(w)} -  e^{(f- s\tau)_N(w')}\right| \, h(w') \\
& \leq & \sum_{\sigma^N w = v} e^{(f- s\tau)_N(w) - (f- s\tau)_N(w')}  e^{(f- s \tau)_N(w')}\,E_1 \theta_2^{p+2N} \diam (\tcc_m)\,  h(w') \\  
&       & + \sum_{\sigma^N w = v}  \left| e^{(f- s \tau)_N(w) - (f- s \tau)_N(w')} - 1\right|  \,  e^{(f- s \tau)_N(w')} h(w') \\
& \leq & E_1 \theta_2^{p+2N} \diam (\tcc_m)\,e^{C_5} \, (L^N_{f-s \tau} h)(v') +  e^{C_5} \,C_5 \theta_2^{p+N} \, \diam(\tcc_m)\, (L^N_{f-s \tau} h)(v')\\
& \leq & E_1 \,\theta_2^{p+N}\,   \diam(\tcc_m)\, (L^N_{f-s \tau} h)(v') , 
\end{eqnarray*}
since $e^{C_5} \,C_5 \leq E_1/2$ and $\theta_2^N \,e^{C_5} \leq 1/2$ by the choice of $N$.
Hence $L^N_{f-s \tau} h \in \kk_0$.

A simple induction implies now that $L^{qN}_{f - s \tau} h \in \kk_0$ for all $h \in \kk_0$ and all integers $q \geq 1$.
\endofproof

\bs

Lemma 6.7 is required to prove the following.

\bs

\noindent
{\bf Corollary 6.8.} {\it For any real constant $a$ with $|a| \leq a_0$ we have $h_a \in \kk_0$.} 

\bs

\noindent
{\it Proof.} Let $|a| \leq a_0$. Since the constant  function $h = 1 \in \kk_0$, it follows from 
Lemma 6.7 that $L^{qN}_{f-(P+a)\tau} 1 \in \kk_0$ for all $q \geq 1$. Now the Ruelle-Perron-Frobenius 
Theorem (see e.g. \cite{PP}) and the fact that $\kk_0$ is closed in $\ff_\theta(\hU)$ imply  $h_a\in \kk_0$. 
\endofproof

\subsection{Estimates of $\lab^N$ by contraction operators}

We will now define a class of pairs of functions similar to $\kk_0$ however involving the parameter $b$. 
We continue to {\bf assume that} $f \in \ff_{\theta_1}(\hU)$.

Denote by $\kk_b$  {\it the set of all pairs $(h, H)$ such that $h \in \ff_{\theta}(\hU)$, 
$H \in \kk_E$, and the following two conditions are  satisfied:}
\index{$\kk_b$} 

\ms

$(T1)$\quad  $|h| \leq H$ on $\hU$,

\ms

$(T2)$ \quad for any $u,u' \in \hU$ contained in a cylinder $\cc'_m = \piU(\cc_m)$ for some $m = 1, \ldots, m_0$, 
any integer  $p \geq 0$ and any points $v,v'\in \hU_1$
satisfying (\ref{eq:u-cond})  for some $i = 1,2$ and $\ell = 1, \ldots, \ell_0$ we have
\be
|h(v) - h(v') | \leq E \, |b|\, \theta_2^{p+N} \, H(v') \, \diam(\tcc_m) .
\ee
Recall that here $\tcc_m = \tPsi(\cc_m)$.

\ms

Our aim in this section is to prove the following.

\ms

\noindent
{\bf Lemma 6.9.} {\it Choosing $E > 1$ and $\mu_0$ as in Sect. {\rm 5.2} and assuming $N$ is sufficiently large, 
for any $|a| \leq a_0$, any $|b| \geq b_0$ and any $(h,H) \in \kk_b$ 
there exists $J \in \jj(b)$ such that $(\lab^{N} h, \nn_J H ) \in \kk_b$.}

\bs

To prove this  we need the following lemma, whose proof is essentially the same as that of 
Lemma 14 in \cite{D1}.  For completeness we prove it in the Appendix.

\bs

\noindent
{\bf Lemma 6.10.}  {\it Let $(h,H) \in \kk_b$. Then for any $m \leq m_0$,  any $j = 1, \ldots, j_0$
with $\dd_j \subset \cc_m$, any $i = 1,2$ and  $\ell = 1,\ldots,  \ell_0$ we have:}

(a) {\it $\di\frac{1}{2} \leq \frac{H(\vli (u'))}{H(\vli (u''))} \leq 2$ for all} $u', u'' \in \dd'_j$;

(b) {\it Either for all $u\in \dd'_j$ we have
$|h(\vli (u))|\leq \frac{3}{4}H(\vli (u))$, or $|h(\vli (u))|\geq \frac{1}{4}H(\vli (u))$ 
for all $u\in \dd'_j$.}

\bs

\noindent
{\it Proof of Lemma} 6.9. 
The constant $E_1 > 1$ from Sect. 6.3 depends only on $C_5$, and we take $N$ so large that
$E_1 \theta_2^N \leq 1/4$; then  $C_5 \theta_2^N \leq 1/2$ holds, too.

Let $|a| \leq a_0$, $|b| \geq b_0$ and $(h, H) \in \kk_b$. We will construct a representative set 
$J \in \jj(b)$ such that  $(\lab^N h , \nn_J H) \in \kk_b$. 

Consider for a moment an arbitrary (at this stage) representative set $J$. 
We will first show that  $(\lab^N h , \nn_J H) $ has property $(T2)$. 

Assume that the points $u,u'$, the cylinder $\cc_m$ in $U$, the integer $p \geq 0$ and the points 
$v,v'\in \hU$ satisfy (\ref{eq:u-cond})  for some $i = 1,2$ and $\ell = 1, \ldots, \ell_0$. 

From the definition of $\fa$, for any $w, w'$ with $\sigma^N w = v$, $\sigma^N(w') = v'$ 
and $\ell(w,w') \geq N$ we have 
\begin{eqnarray*}
\fa_N(w) 
& =  & f_N(w) - (P+a)\tau_N(w) + (\ln h_a - \ln h_a\circ \sigma)_N(w) - N \lambda_a\\
& =  & f_N(w) - (P+a)\tau_N(w) + \ln h_a (w)  - \ln h_a (v) - N \lambda_a .
\end{eqnarray*}
Since $h_a \in \kk_0$ by Corollary 6.8,
$$|\ln h_a(w) - \ln h_a(w') | \leq \frac{| h_a(w) - h_a(w')|}{\min \{ |h_a(w)|, |h_a(w')|\}} 
\leq E_1 \,\theta_2^{p+2N} \, \diam(\tc_m) , $$
and similarly, $|\ln h_a(v) - \ln h_a(v') | \leq E_1 \,\theta_2^{p+N} \, \diam(\tc_m)$.
Using this and Lemma 6.6, as in the proof of Lemma 6.7 we get
\begin{eqnarray}
&        & |\fa_N(w) - \fa_N(w')| 
 \leq  C_5 \theta_2^{p+2N} \, \diam(\tc_m) + 2E_1 \, \theta_2^{p+N}\, \diam(\tc_m)\\
& \leq & (C_5 + 2 E_1)  \, \theta_2^{p+N}\, \diam(\tc_m)  \leq 1 , \nonumber
\end{eqnarray}
by the choice of $N$.

Hence for any $a$ and $b$ with $|a| \leq a_0$ and $|b| \geq b_0$, using (6.24) and Lemma 6.6, we get  
\begin{eqnarray*}
&        &|(\lab^N h)(v) - (\lab^N  h)(v')|\\
& =    &   \left| \sum_{\sigma^N w = v} e^{(\fa_N - \i b \tau_N) (w)}\, h(w)   -  \sum_{\sigma^N w = v} e^{(\fa_N - \i b \tau_N)(w'(w))}\, h(w'(w)) \right| \\
& \leq & \left| \sum_{\sigma^N w = v} e^{(\fa_N - \i b \tau_N)(w)}\, [h(w) -  h(w')]\right| \\
&       &  + \sum_{\sigma^N w = v}  \left|e^{(\fa_N - \i b \tau_N)(w)} -  e^{(\fa_N - \i b \tau_N) (w')}\right| \, |h(w')| \\
& \leq & \sum_{\sigma^N w = v} e^{(\fa_N (w) - \fa_N(w')} e^{\fa_N (w')}\,E |b|  \theta_2^{p+2N} \diam (\tcc_m)\,  H(w') \\  
&       & + \sum_{\sigma^N w = v}  \left| e^{(\fa_N - \i b \tau_N)(w) -  (\fa_N - \i b  \tau_N)(w')} - 1\right| \,  e^{\fa_N(w')} H(w') \\
& \leq & e\, E |b|\, \theta_2^{p+2N} \diam (\tc_m) \, (\ma^N H)(v') \\
&       & +  e \,(C_5 + 2E_1 +  C_4 |b| ) \theta_2^{p+N} \, \diam(\tc_m)\, (\ma^N H)(v')\\
& \leq & [2 e\, E \theta_2^N  + 2 e (C_5 + 2E_1 + C_4) ]  \,|b|\, \theta_2^{p+N}\,    \diam(\tc_m)\, (\nn_J H)(v')\\
& \leq &  E |b|\, \theta_2^{p+N}\,   \diam(\tc_m)\, (\nn_J H)(v') , 
\end{eqnarray*}
assuming $2 e \theta^N \leq 1/2$ and $2 e (C_5 + 2E_1 + C_4)  \leq E/2$. 
Thus, $(\lab^N h , \nn_J H)$ has property $(T2)$.

So far the choice of $J$ was not important. We will now construct a representative set $J$ so that 
$(\lab^N h , \nn_J H)$ has property $(T1)$, namely 
\be\label{eq:Lab}
|\lab^N h|(u) \leq (\nn_J H) (u) \quad , \quad u \in \hU . 
\ee

Define the functions $\psi_\ell, \gao_\ell, \gat_\ell : \hU  \longrightarrow \C$  by
$$\di \psi_\ell(u) = e^{(\fa_{N}+\i b\tau_{N})(\vl_1(u))} h(\vl_1(u)) 
+ e^{(\fa_{N}+\i b\tau_{N})(\vl_2(u))} h(\vl_2(u)) ,$$
$$\di \gao_\ell(u) = (1-\mu_0)\, e^{\fa_{N} (\vl_1(u))} H(\vl_1(u)) + e^{\fa_{N}(\vl_2(u))} H(\vl_2(u)) ,$$
while $\gat_\ell(u)$ is defined similarly with a coefficient $(1-\mu_0)$ 
in front of the second term. 

Notice that (\ref{eq:Lab}) is trivially satisfied for $u \notin V_b$ for any choice of $J \in \jj(b)$.

Consider an arbitrary  $m = 1, \ldots, m_0$. We will construct  $j \leq j_0$ with $\dd_{j} \subset \cc_m$,  
and a pair  $(i, \ell)$ for which $(i,j,\ell)$ will be included in $J$.

Recall the functions $\varphi_{\ell,m}(u) = \varphi_\ell(Z_m, u)$, $u \in U$,  from Lemma 5.4.

\ms

\noindent
{\bf Case 1.}
There exist $j\leq j_0$ with $\dd_j \subset \cc_m$, $i = 1,2$ and $\ell \leq \ell_0$   such that the 
first alternative in Lemma 6.10(b) holds for $\dd_j$,  $i$ and $\ell$. 
For such $j$, choose $i= i_j$ and $\ell = \ell_j$ with this property and include $(i,j,\ell)$ in $J$. 
Then $\mu_0 \leq 1/4$ implies $|\psi_\ell(u)| \leq \gamma^{(i)}_\ell(u)$ for all 
$u\in \dd'_j$, and regardless how the rest of $J$ is defined, (\ref{eq:Lab}) holds for all $u\in \dd'_j$, since
\begin{eqnarray}
&        & \left| (\lab^N h)(u)\right| 
 \leq  \left| \sum_{\sigma^N v = u, \;v\neq \vl_1(u),\vl_2(u)} e^{(\fa_N+\i b\tau_N)(v)} h(v) \right|  + |\psi_\ell(u)| \\
& \leq & \sum_{\sigma^N v = u, \;v\neq \vl_1(u),\vl_2(u)} e^{\fa_N(v)} |h(v)|  + \gamma^{(i)}_\ell(u)\nonumber\\
& \leq & \sum_{\sigma^N v = u, \;v\neq \vl_1(u),\vl_2(u)} e^{\fa_N(v)} \omega_J(v) H(v)\nonumber\\
&      & + \left[e^{\fa_N(\vl_1(u))} \omega_J(\vl_1(u)) H(\vl_1(u)) 
+  e^{\fa_N(\vl_2(u))} \omega_J(\vl_2(u)) H(\vl_2(u))\right]  \nonumber\\
& \leq &  (\nn_J H) (u) . \nonumber
\end{eqnarray}

\noindent
{\bf Case 2.}
For all $j\leq j_0$ with $\dd_j \subset \cc_m$, $i = 1,2$ and $\ell \leq \ell_0$ the 
second alternative in Lemma 6.10(b) holds for $\hdd_j$,  $i$ and $\ell$, i.e.
\be
|h(\vl_i(u))|\geq \frac{1}{4}\, H(\vl_i(u)) > 0
\ee
for any  $u \in \cc'_m$.

Let $u,u' \in \cc'_m$, and let $i = 1,2.$ 
Using (6.23) and the assumption that $(h,H) \in \kk_b$, and in particular property 
$(T2)$ with $p = 0$, $v = \vl_i(u)$ and $v' = \vl_i(u')$, and assuming e.g.
$$\min\{ |h(\vl_i(u))| , |h(\vl_i(u'))|\}  = |h(\vl_i(u'))| ,$$
we get
\begin{eqnarray*}
\frac{|h(\vl_i(u)) - h(\vl_i(u'))|}{\min\{ |h(\vl_i(u))| , |h(\vl_i(u'))| \}}
& \leq &  \frac{E|b|\,\theta_2^N H(\vl_i(u'))}{|h(\vl_i(u'))| } \diam (\tPsi(\cc_m)) \\
& \leq & 4 E|b|\,   \theta_2^{N} \frac{\ep_1 A R_0}{|b|} = 4 E A \theta_2^N \ep_1 .
\end{eqnarray*}
So, the difference between the arguments of  the complex numbers 
$h(\vl_i(u))$ and $h(\vl_i(u'))$ (regarded as vectors in $\R^2$)  is  
$< 8 E \theta_2^N \ep_1 A R_0  < \pi/6$, assuming $N$ is sufficiently large.  
In particular,  for each $i = 1,2$ we can choose a real continuous
function $\thetam_i(u)$, $u \in  \cc'_m$, with values in $[0,\pi/6]$  and a constant $\lambdam_i$ such that
$$\di h(\vl_i(u)) = e^{\i(\lambdam_i + \thetam_i(u))}|h(\vl_i(u))| \quad , \quad u\in \cc'_m .$$
Fix an arbitrary $u_0\in \cc'_m$ and set $\lambdam = |b| \varphi_{\ell,m}(u_0)$. 
Replacing e.g $\lambdam_2$ by $\lambdam_2 +  2r \pi$ for some integer $r$, 
we may assume that 
$$|\lambdam_2 - \lambdam_1 + \lambdam | \leq \pi .$$

Using the above, $\theta \leq 2 \sin \theta$ for $\theta \in [0,\pi/3]$,  
and some elementary geometry  yields
$$|\thetam_i(u) - \thetam_i(u')|\leq 2 \sin |\thetam_i(u) - \thetam_i(u')| < 16 E A R_0 \theta_2^N \ep_1$$ 
for all $u,u' \in \cc'_m$.

The difference between the arguments of the complex numbers
$$e^{\i \,b\,\tau_N(\vl_1(u))} h(\vl_1(u)) \quad \mbox{\rm and} \quad e^{\i \,b\, \tau_N(\vl_2(u))} h(\vl_2(u))$$
is given by the function
\begin{eqnarray*}
\Gl(u) 
& = & [b\,\tau_N(\vl_2(u)) + \thetam_2(u) + \lambdam_2] -  
[b\, \tau_N(\vl_1(u)) + \thetam_1(u) + \lambdam_1]\\
& = & (\lambdam_2-\lambdam_1) + |b| \varphi_{\ell,m}(u) + (\thetam_2(u) - \thetam_1(u)) .
\end{eqnarray*}

It follows from Lemma 5.4 that there exist $j, j' \leq j_0$ with $j \neq j'$, 
$\dd_j, \dd_{j'} \subset \cc_m$, and  $\ell = 1, \ldots, \ell_0$ such that
for all $u \in \dd'_j $ and $u'\in \dd'_{j'}$  we have 
\begin{eqnarray}
\frac{\hd \hrho \ep_1}{R_0 |b|} 
& \leq & \hd \, \diam(\tPsi(\cc_m)) \leq |\varphi_{\ell,m}(u) - \varphi_{\ell,m}(u')|\\ 
& \leq & C_1 \, \diam(\tPsi(\cc_m)) \leq  C_1 \frac{R_0 \ep_1}{|b|} .\nonumber
\end{eqnarray}
Fix $ \ell$ with this property.
Then for $u \in \dd'_j$ and $u'\in \dd'_{j'}$ we have
\begin{eqnarray*}
|\Gl(u)- \Gl(u')|
& \geq & |b|\, |\varphi_{\ell,m}(u) - \varphi_{\ell,m}(u')| - |\thetam_1(u) - \thetam_1(u')|\\
&        & - |\thetam_2(u) - \thetam_2(u')|\\
& \geq & \frac{\hd \hrho \epsilon_{1}}{R_0} - 32 E \theta_2^N A R_0 \ep_1 >  2\ep_3 ,
\end{eqnarray*}
since $32 E A \theta_2^N  < \frac{\hd \hrho}{2R_0}$ by (5.5), (5.6), where 
$$\di \ep_3 = \frac{ \hd \hrho \epsilon_{1}}{4 R_0} .$$

Thus,  $|\Gl(u)- \Gl(u')|\geq 2\epsilon_{3}$ for all  $u\in \dd'_j$ and all $u'\in \dd'_{j'}$. Hence either 
$|\Gl(u)| \geq \ep_3$ for all $u\in \dd'_j$ or $|\Gl(u')| \geq \ep_3$ for all $u'\in \dd'_{j'}$. Indeed, if
$|\Gl(u')| < \ep_3$ for some $u'\in \dd'_{j'}$, then for every $u\in \dd'_{j}$ we get
$$|\Gl(u)| = |(\Gl(u)- \Gl(u')) + \Gl(u')| \geq |\Gl(u)- \Gl(u')| - |\Gl(u')| > \ep_3 . $$
Similarly, if $|\Gl(u)| < \ep_3$ for some $u\in \dd'_{j}$, then $|\Gl(u')| \geq \ep_3$ for every $u'\in \dd'_{j'}$.

Assume for example that $|\Gl(u)| \geq \ep_3$ for all $u\in \dd'_j$. 
On the other hand, (6.28) and the choice of $\ep_1$ imply that
for any $u \in \cc'_m$ we have
\begin{eqnarray*}
|\Gamma_\ell(u)| 
& \leq & |\lambdam_2-\lambdam_1 +\lambdam | +  |b|\, |\varphi_\ell(u) -\varphi_\ell(u_0)| + |\thetam_2(u) - \thetam_1(u)| \\
& \leq & \pi + C_1 \ep_1R_0 + 16 E A \theta_2^N \ep_1 R_0 < \frac{3\pi}{2} .
\end{eqnarray*}
Thus, $\ep_3 \leq |\Gl(u)| <  \frac{3\pi}{2}$ for all $u \in \dd'_j$.

Hence, we see that for $u\in \dd'_j$ the difference $\Gl(u) $ between the  arguments of the complex numbers
$$e^{\i \,b\,\tau_N(\vl_1(u))} h(\vl_1(u)) \quad \mbox{\rm and} \quad  e^{\i \,b\, \tau_N(\vl_2(u))} h(\vl_2(u)) ,$$
defined as a number in the interval $[0, 2\pi)$, satisfies
$$\Gl(u) \geq \ep_3 \quad , \quad u\in \dd'_j .$$

\def\halpha{\hat{\alpha}}
\def\hgamma{\hat{\gamma}}

As in \cite{D1} it follows from Lemma 6.10 that either $H(\vl_1(u)) \geq H(\vl_2(u))/4$ for all
$u \in \dd'_j$ or $H(\vl_2(u)) \geq H(\vl_1(u))/4$ for all $u \in \dd'_j$. Indeed, fix
an arbitrary $u'\in \dd'_j$ and assume e.g. $H(\vl_1(u')) \geq H(\vl_2(u'))$. Then for any $u \in \dd'_j$
using Lemma 6.10(a) twice we get 
$$H(\vl_1(u)) \geq H(\vl_1(u'))/2 \geq H(\vl_2(u'))/2 \geq H(\vl_2(u))/4 .$$
Similarly, if $H(\vl_2(u')) \geq H(\vl_1(u'))$, then $H(\vl_2(u)) \geq H(\vl_1(u))/4$ for all $u \in \dd'_j$.

Now assume e.g. that $H(\vl_1(u)) \leq H(\vl_2(u))/4$ for all $u \in \dd'_j$.
As in \cite{D1} (see also \cite{St2}) we will show that 
$|\psi_\ell(u)| \leq \gamma^{(1)}_\ell (u)$ for all $u \in \dd'_j$.
Given such $u$, consider the points 
$$z_1 = e^{(\fa_{N}+\i b\tau_{N})(\vl_1(u))} h(\vl_1(u)) \:\:\: ,  \:\:\: z_2 =  e^{(\fa_{N}+\i b\tau_{N})(\vl_2(u))} h(\vl_2(u))$$
in the complex plane $\C$, and let $\varphi$ be the smaller angle between the arguments of $z_1$ and $z_2$.
It then follows from the above estimates for $\Gl(u)$ that $\epsilon_3 \leq \varphi \leq 3\pi/2$. Moreover,  (6.24) and (6.27)  imply
\begin{eqnarray*}
\frac{|z_1|}{|z_2|} 
=     e^{\fa_N(\vl_1(u)) - \fa_N(\vl_2(u))}   \frac{|h(\vl_1(u))|}{|h(\vl_2(u))|}
  \leq    \frac{H(\vl_1(u))}{H(\vl_2(u))/4} \leq 16  .
 \end{eqnarray*}
This yields
\begin{equation}
|z_1 + z_2| \leq (1- t)  |z_1| + |z_2| ,
\end{equation}
where we can take e.g. 
$$t = \frac{1 - \cos (\epsilon_3)}{20} .   $$
Indeed, we have
$$|z_1 + z_2|^2 = |z_1|^2 + |z_2|^2 + 2 \langle z_1, z_2\rangle \leq  |z_1|^2 + |z_2|^2 + 2 |z_1|\,  |z_2| (1- \alpha) ,$$
where $\alpha = 1- \cos \ep_3$. Thus, (6.29) will hold if
$$|z_1|^2 + |z_2|^2 + 2 |z_1|\,  |z_2| (1- \alpha)  \leq  (1-t)^2 |z_1|^2 + |z_2|^2 + 2 (1-t) |z_1|\,  |z_2| ,$$
that is if 
$$(1 - (1-t)^2) |z_1| + 2   |z_2| (1- \alpha)  \leq  2 (1-t)   |z_2| ,$$
which equivalent to
$$|z_1| \leq 2 \frac{\alpha - t}{t(2-t)} \, |z_2| .$$
If $t = \alpha/20$, then $16 <  2 \frac{\alpha - t}{t(2-t)} = \frac{38}{2- \alpha/20}$, so the above inequality holds.
This proves (6.29) with the given choice of $t$.

Since $\mu_0 \leq t$ by (\ref{eq:mucond}), we have $|\psi_\ell(u)| \leq \gamma^{(1)}_\ell (u)$  
for all $u \in \dd'_j$.  Now set $j_m = j$, $\ell_m = \ell$ and  $i_m = 1$,
and include  $(i_m,j_m,\ell_m)$ in the set $J$.  
Then $\dd_{j_m} \subset \cc_m$ and as in the proof of (6.26) we deduce that  (\ref{eq:Lab}) holds on $\dd'_{j_m}$.

This completes the construction of the set 
$$J = \{ (i_m, j_m, \ell_m) : m = 1, \ldots, m_0\} \in \jj(b)$$
and also the proof  of (\ref{eq:Lab}) for all $u \in V_b$. As we mentioned in the beginning of the proof, (\ref{eq:Lab}) 
always holds for $u \in \hU\setminus V_b$.
\endofproof

\section{$L^1$ contraction estimates}
\setcounter{equation}{0}

Here we obtain $L^1$-contraction estimates for large powers of the contraction operators
using the properties of these operators on $K_0$, the strong mixing
properties of the shift map $\pp : R \longrightarrow R$ and the Pesin set $P_0$ with exponentially
small tails. We continue to use the notation from Sections 5 and 6. We will also use the constants $p_0$, $\hep$ and $\hd_0$, 
the Pesin set $P_0$ and the set $\Xi_m = \Xi_m(p_0,\hep, \hd_0)$ from Sect. 4.2.

For any $J \in \jj(b)$ set
$$W_J = \cup\{ \dd'_j : (i,j,\ell) \in J \; \mbox{\rm for some }\; i, \ell \} \subset V_b .$$
Using Lemma 6.4 and the class of functions $\kk_E$ we will now prove the following 
important estimates\footnote{This should be regarded as the analogue of 
Lemma 12 in \cite{D1} (and Lemma 5.8 in \cite{St2}).}.
\index{$W_J$}

\bs

\noindent
{\bf Lemma 7.1.}  {\it Let $f \in \ff_{\theta_1}(\hU)$. }

\medskip

(a) {\it There exists a global constant $C_6 > 0$, independent of $b$ and $N$, such that for any $H\in \kk_E$  
and any $J\in \jj(b)$ we have}
\be
\int_{V_b} H^2 \, d\nu \leq C_6\, \int_{W_J} H^2 \, d\nu .
\ee

\ms

(b) {\it Assuming that $a_0 > 0$ is sufficiently small, for any $H\in \kk_E$   and any $J\in \jj(b)$ we have
\be
\int_{V_b} (\nn_J H)^2 \, d\nu \leq \rho_3\, \int_{V_b} L^N_{\f0} (H^2) \, d\nu ,
\ee
where} 
$$\di \rho_3 = \rho_3(N) = \frac{e^{a_0NT_0}}{1+ \frac{\mu_0 e^{-NT_0}}{C_6}} < 1 .$$

\ms

\noindent
{\it Proofs.} (a) Let $H \in \kk_E$ and let $J \in \jj(b)$. Consider an arbitrary 
$m = 1, \ldots,m_0$. There exists $(i_m,j_m,\ell_m) \in J$ such that $\dd_{j_m} \subset \cc_m$.
It follows from (\ref{eq:Gibbs}) that there exists a global constant $\omega_0 \in (0,1)$ such that
 $$\frac{\nu(\dd'_{j_m})}{\nu(\cc'_m)} \geq 1- \omega_0 .$$
Since $H \in \kk_E$, for any $u,u'\in \cc'_m$ we have 
$$\frac{|H(u) - H(u')|}{H(u')} \leq E\dd(u,u') \leq E , $$ so
$$H(u)/H(u') \leq 1+ E \leq 2E .$$
Thus, if 
$$L_1 = \max_{\cc'_m} H \quad , \quad L_2 = \min_{\cc'_m}H ,$$
we have $L_1/L_2 \leq 2E$. This gives
$$\int_{\cc'_m} H^2 \, d\nu \leq L_1^2 \nu(\cc'_m) 
\leq  \frac{4E^2}{1-\omega_0}\, \int_{\dd'_{j_m}} H^2\, d\nu . $$
Hence
$$\int_{V_b} H^2 \, d\nu \leq \sum_{m=1}^{m_0} \int_{\cc'_m} H^2 \, d\nu  
\leq  \frac{4E^2}{1-\omega_0}\, \sum_{m=1}^{m_0}  \int_{\dd'_{j_m}} H^2\, d\nu
\leq C_6\, \int_{W_J} H^2 \, d\nu ,$$
with $C_6 =  \frac{4E^2}{1-\omega_0}$, since  $\cup_{m=1}^{m_0} \dd'_{j_m} = W_J$. 
This proves (7.1).

\bs

(b) Let again $H \in \kk_E$ and $J \in \jj(b)$. By Lemma 6.4, $\nn_J H \in \kk_E$, while
the Cauchy-Schwartz  inequality implies
$$(\nn_J H)^2 = (\ma^N \omega H)^2 \leq (\ma^N \omega_J^2)\, (\ma^N H^2) \leq 
(\ma^N \omega_J)\, (\ma^N H^2) \leq \ma^N H^2 .$$
Notice that if $u \notin W_J$, then $\omega_{J}(u) = 1$. Let $u \in W_J$; then $u \in \dd'_j$ 
for some (unique) $j\leq j_0$, and there exists a unique $(i(j), j, \ell(j)) \in J$. 
Set $ i = i(j)$, $\ell = \ell(j)$ for brevity.
Then $\vl_i(u) \in \Xijl$, so $\omijl(\vl_i(u)) = 1$, and therefore  
$$\omega(v^{(\ell)}_i(u)) \leq 1-\mu_0\, \omijl(\vl_i(u)) = 1-\mu_0 .$$
In fact,  if $\sigma^N(v) = u$ and $\omega(v) < 1$, then
$\omega^{(\ell')}_{i',j'}(v) = 1$ for some $(i',j',\ell') \in J$, so $v \in X^{(\ell')}_{i',j'}$. This yields
$$u = \sigma^N(v)  \in \sigma^N(X^{(\ell')}_{i',j'}) = \dd'_{j'} .$$ 
Thus, we must have $j' = j$, and since for a given $j$, there is 
only one element $(i(j) , j , \ell(j) )$ in $J$, we must have also 
$i' = i(j)$ and $\ell' = \ell(j)$. Assuming e.g. that $i = 1$,  this implies $v = \vl_1(u)$. Thus, 
\begin{eqnarray*}
(\ma^N \omega_{J})(u)
& = &  \sum_{\sigma^N v = u, \;v\neq \vl_1(u)} e^{\fa_N(v)}   + e^{\fa_N(\vl_1(u))} \omega_{J} (\vl_1(u))  \\
& =  & \sum_{\sigma^N v = u, \;v\neq \vl_1(u)} e^{\fa_N(v)}   + (1-\mu_0)  e^{\fa_N(\vl_1(u))} \\ 
&  =   & \sum_{\sigma^{N}v = u} e^{\fa_{N}(v)} - \mu_0 \,e^{\fa_{N}(\vl_1(u))}  \\
& \leq & (\ma^{N}\; 1)(u) - \mu_0 \, e^{-N T_0} = 1 - \mu_0 \, e^{-NT_0} .
\end{eqnarray*}
This holds for all $u \in W_J$, so  
$$(\nn_J H)^2 \leq (1- \mu_0 e^{-NT_0})\, (\ma^N H^2)$$ 
on $W_J$. Using this and part (a) we get:
\begin{eqnarray*}
\int_{V_b} (\nn_J H)^2\, d\nu
& =     & \int_{V_b \setminus W_J} (\nn_J H)^2 \, d\nu +  \int_{W_J} (\nn_J H)^2 \, d\nu\\
& \leq &  \int_{V_b \setminus W_J} (\ma^N H)^2 \, d\nu  + (1-\mu_0 e^{-NT_0})\,  \int_{W_J} (\ma^N H)^2 \, d\nu\\
& =     & \int_{V_b} (\ma^N H)^2 \, d\nu  -\mu_0 e^{-NT_0}\,  \int_{W_J} (\ma^N H)^2 \, d\nu\\
& \leq & \int_{V_b} (\ma^N H)^2 \, d\nu  - \mu_0 e^{-NT_0}\,  \int_{W_J} (\nn_J H)^2 \, d\nu \\
& \leq & \int_{V_b} (\ma^N H)^2 \, d\nu  - \frac{\mu_0 e^{-NT_0}}{C_6} \,  \int_{V_b} (\nn_J H)^2 \, d\nu .
\end{eqnarray*}
From this and 
$$(\ma^N H)^2 \leq (\ma^N 1)^2 (\ma^NH^2) \leq \ma^N H^2 
= L^N_{\f0} (e^{\fa_N - \f0_N} H^2)\leq e^{a_0 NT_0} (L^N_{\f0} H^2) ,$$ 
we get
$$(1+\mu_0 e^{-NT_0}/C_6)\, \int_{V_b} (\nn_J H)^2\, d\nu \leq  \int_{V_b} (\ma^N H)^2 \, d\nu
\leq e^{a_0NT_0}\, \int_{V_b} L^N_{\f0} H^2 \, d\nu .$$
Thus, (7.2) holds with the given $\rho_3 $ in the statement of the lemma.
Taking $a_0 = a_0(N) > 0$ sufficiently small, we have $\rho_3 < 1$.
\endofproof

\bs

We can now prove that iterating sufficiently many contraction operators provides an $L^1$-contraction on $U$.

Define $\hb$ by (6.9), and let $0 < \delta_1 < \delta_2$ be as in (6.6). Set
\be\label{eq:Lambda}
\Lambda_N (b) = \left\{ x\in \ll \cap R : \:\: \sharp \, \left\{ j : 0 \leq j <  N\, \hb \; , \;  
\pp^j(x) \notin P_1 \cap P_0 \cap \Omb_{B}\right\}  \geq \frac{\delta_2}{N}\,  N \, \hb \: \: \right\}.
\ee
\index{$\Lambda_N(b)$}

In what follows we will use the estimate (6.4). Apart from that it follows from Theorem 3.2 that there
exist constants $C'_7 = C'_7(N) > 0$ and $c'_7 = c'_7(N) > 0$ such that
\be
\mu (\Xi_{n}(p_0,\hep, \delta_1/N)) \leq C'_7\, e^{-c'_7 n} \quad, \quad n  \geq 1.
\ee

\ms

The following lemma will play a significant role in showing that some cancellation occurs in the actions
of the contraction operators $\nn_J$.

\bs

\noindent
{\bf Lemma 7.2.} {\it Let $s > 0$ be a constant. There exist constants $C_7 > 0$ and $N_0 \geq 1$, independent of $b$, such that
$$\mu (\Lambda_N (b)) \leq \frac{C_7}{|b|^{s}}$$
for all $b$ with $|b| \geq b_0$ and all integers $N \geq N_0$.}

\bs

\noindent
{\it Proof.} 
Set $m = \hb$. We claim that 
\be
\Lambda_N (b) \subset \Xi_{Nm}(p_0,\hep, \delta_1/N) \cup Y ,
\ee
where
\be
Y = \cup_{j=0}^{Nm} \pp^{-j}(R\setminus \Omb_{B}) .
\ee

Set $\Xi = \Xi_{Nm}(p_0,\hep, \delta_1/N)$ for brevity.
 Assume that there exists
$$x \in \Lambda_N(b) \setminus (\Xi \cup Y) .$$
Then $x \notin Y$, so $\pp^j(x) \notin R\setminus \Omb_{B}$ for all $j \in [0,Nm]$, i.e.
\be
\pp^j(x) \in \Omb_{B} \quad , \quad 0 \leq j \leq Nm .
\ee
Next, $x \notin \Xi$ gives 
\be
\# \{ j :  0 \leq j <  N m \:, \:  \pp^j(x) \notin P_0\} < \frac{\delta_1}{N} N m .
\ee
Now $x \in \Lambda_N (b)$ and (7.7) imply $\pp^j(x) \notin P_1 \cap P_0$ for at least $\frac{\delta_2}{N}\,  N \, m$ values of 
$j = 0,1, \ldots, N m -1$. Notice that by (\ref{eq:dcond}), 
$$\frac{\delta_2}{N}\,  N \, m = \frac{(2d+3) B \delta_1}{N} N m .$$
This and (7.8) now yield 
\be
\# \{  j :  0 \leq j <  N m \:, \:  \pp^j(x) \notin P_1\} \geq  \frac{(2d+2) B \delta_1}{N} N m .
\ee

We need to extend (7.8) a bit. Using (7.7) for $j = Nm$ and (6.3), we get $\pp^{Nm}(x) \notin \Xi_{Bm}(p_0, \hep_0, \delta_1)$, 
which means that
$$\# \{ i :  0 \leq i <  B m \:, \:  \pp^{Nm+i}(x) \notin P_0\} < \delta_1 B m = \frac{\delta_1 B}{N} Nm .$$
Combining the latter with (7.8) gives
\be
\# \{ j :  0 \leq j <  (N+B) m \:, \:  \pp^j(x) \notin P_0\} < \frac{2B \delta_1}{N} N m .
\ee

Given $y = \pp^j(x)$ for some $j = 0,1,\ldots, N m-1$, let $m_j = m_{y}$ be the length of the maximal cylinder $\cc(y)$ in $W^u_{R}(y)$ with
$$\diam ((\Phi_y)^{-1}\cc(y))\leq \ep_1/|b| .$$ 
Then by (\ref{eq:Bcond}), $m_j \leq B \hb = Bm$. Moreover, if $y \notin P_1$, 
then $\pp^{m_j}(y) \notin P_0$, i.e. $\pp^{m_j + j} (x) \notin P_0$. By (7.7), $y = \pp^j(x) \in \Omb_B$, and then by (6.3), 
$y \notin \Xi_\ell$ for all $\ell \in [m/B, Bm]$. Using the latter with $\ell = m/B$, the definition of $\Xi_\ell$ yields $\pp^i(y) \in P_0$ for
some $i \in [(1-\delta_0)m/B, m/B]$. Since $m_j \geq m/B$ by (6.2), we have $(1-\hd_0)m/B \leq i \leq m_j$. 
Now we can apply Lemma 4.4(d) using the maximality of the cylinder $\cc(y)$, and the points $y$ and $\pp^i(y) \in P_0$
and derive that 
$$\diam (\Phi_y^{-1}(\cc(y)))\geq \rho_1 \ep_1/|b| .$$
More precisely, let $\cc_1$ be the unstable cylinder 
in $W^u_{R}(y)$ of length $i -1$; this is then $ < m_j$, the length of $\cc(y)$. Thus $\cc(y) \subset \cc_1$ and the maximality 
of $\cc(y)$ implies $\diam (\Phi_y^{-1}(\cc_1)) > \ep_1/|b|$. Since $\pp^i(y) \in P_0$, it follows from Lemma 4.4(d) 
and $\hrho < \rho_1$  that
$$\diam (\Phi_y^{-1}(\cc(y))) \geq \rho_1\, \diam (\Phi_y^{-1}(\cc_1)) > \rho_1 \frac{\ep_1}{|b|} .$$
On the other hand, by the choice of the constant $d$ we have 
$c_0\gamma^d \rho_1 > 1$, therefore 
$$\diam (\Phi_{\pp^d(y)}^{-1}(\pp^d(\cc(y))) > \ep_1/|b| .$$ 
Thus,
$\cc(\pp^d(y)) = \cc(\pp^{j+d}(x))$ is a proper sub-cylinder of $\pp^d(\cc(y))$ and therefore its length $m_{j+d}$ is
strictly larger that the length $m_j-d$ of $\pp^d(\cc(y))$, i.e. $m_{j+d} > m_j-d$. Thus we have 
$$m_j+j < m_{j+d} + (j+d) .$$
Moreover, for $j \in [0,Nm)$ we have 
$$m_j + j \in [0, (N+B)m) .$$

By (7.9) there are at least $\frac{(2d+2) B \delta_1}{N}  N m$ values of $j = 0,1, \ldots, Nm -1$ with $\pp^j(x) \notin P_1$;
then $\pp^{m_j+j}(x) \notin P_0$. From those $j's$ we get a strictly increasing sub-sequence with at least
 $$\frac{(2d+2) B \delta_1}{dN}  N m > \frac{2B \delta_1}{N}  N m$$ 
 members which is a contradiction with (7.10).
 
 This proves (7.5). Now using (6.4), (7.4) and $m = \hb = \lceil \log |b| \rceil $ we get
 \begin{eqnarray*}
 \mu(\Lambda_N(b))
 & \leq & C'_7\, e^{-c'_7 N\hb} + N \hb\, \mu(R \setminus \Omb_B)
 \leq \frac{C'_7}{|b|^{N c'_7}} + N\hb\, C_3 e^{-c_3 N \hb/B}  ,
 \end{eqnarray*}
 where $C_3$ and $c_3$ are the constants from (6.4).
 Since $\hb > 1$, assuming $N_0 \geq 1$ is sufficiently large and $N \geq N_0$, we have
 $$N\hb\, e^{-c_3 N \hb/B}  \leq e^{-c_3 N \hb/(2B)} .$$
 Thus, for such $N$,
 $$ \mu(\Lambda_N(b)) \leq  \frac{C'_7}{|b|^{N c'_7}} + C_3 e^{-c_3 N \hb/(2B)} 
 = \frac{C'_7}{|b|^{N c'_7}} + \frac{C_3}{|b|^{N c_3/(2B)}} \leq \frac{C_7'+ C_3}{|b|^s} ,$$
 assuming $N_0 c'_7 \geq s$ and $\frac{N_0 c_3}{2B} \geq s$.
  This proves the lemma.
\endofproof

\index{$\rho_3$}
\index{$S_0$} 
\index{$a_0$}

\bs

Consider the constants
$$ \rho_3 = \frac{e^{a_0NT_0}}{1+ \frac{\mu_0 e^{-NT_0}}{C_6}} < 1 \quad ,  \quad S_0 = e^{a_0 N T_0} > 1 ,$$
and notice that $\rho_3$ is as in Lemma 7.1. Then define the function $\hath$ by
$$\hath = \rho_3\, \chi_{V_b} + S_0\, \chi_{U\setminus V_b} .$$

We will assume $a_0 = a_0(N) > 0$ is chosen so small that
\be
8 a_0NT_0 < \log\left(1+ \frac{\mu_0 e^{-NT_0}}{2C_6} \right) .
\ee
Recall the constants $\theta_1$ and $\theta_2$ from the beginning of Sect. 5.2.

After the comprehensive study of contraction operators in Sections 6 and 7 so far, we can now prove that the
contraction operators $\nn_J$ do have some contraction properties.

\bs

\noindent
{\bf Lemma 7.3.}  {\it Let $f \in \ff_{\theta_1}(\hU)$ and let $s \geq 1$ be a constant. There exist global constants
$N_0 \geq 0$,  $C_8 > 0$ and $C_9 > 0$ such that for any $N \geq N_0$ there exist  constants $k = k(N) \geq 1$,
$a_0 = a_0(N) > 0$  and $b_0 = b_0(N) \geq 1$ such that for any  $|a| \leq a_0$ and $|b| \geq b_0$ we have the following:}

\ms

(a) {\it  For any sequence $J_1, J_2, \ldots, J_r \ldots $ of elements  of $\jj(b)$, setting
$H^{(0)} = 1$ and $H^{(r+1)} = \nn_{J_r} (H^{(r)})$ ($r \geq 0$) we have}
\be
\di\int_{U} (H^{(k \hb)})^2 \, d\nu \leq \frac{C_8}{|b|^{10 s}} .
\ee

(b) {\it For all $h \in \ff_{\theta_1}(U)$ we have}
\be
\|\lab^{2k N \hb} h \|_0 \leq  \frac{C_9}{|b|^{s}} \, \|h\|_{\theta_1,b} .
\ee

\ms

\noindent
{\bf Remark.} Notice that in general the operator $\lab$ does not have to preserve the space $\ff_{\theta_1}(\hU)$.
Indeed, the function $\fa$ involves $\tau$ which is in $\ff_\theta(\hU)$, however not necessarily in $\ff_{\theta_1}(\hU)$.
So, in the left-hand-side of (7.13) we just have the sup-norm of a function in $\ff_\theta(\hU)$.

\bs

\noindent
{\it Proof of Lemma} 7.3. (a) Set  $\omega_r = \omega_{J_r}$, $W_r = W_{J_r}$ and $\nn_r = \nn_{J_r}$.
Since $H^{(0)} = 1 \in \kk_E$, it follows from Lemma 6.5 that $H^{(r)} \in \kk_E$ for all $r \geq 1$. 

Using 
$$L^N_{\fa}((\hath\circ \sigma^N) \, H) = \hath\, (L_{\f0}^N H)$$ 
and Lemma 7.1(b) we get
\begin{eqnarray*}
\int_{U} (H^{(m)})^2 \, d\nu
& =      & \int_{V_b} (H^{(m)})^2 \, d\nu + \int_{U\setminus V_b} (H^{(m)})^2 \, d\nu\\
& \leq  & \rho_3\, \int_{V_b} L^{N}_{\f0} (H^{(m-1)})^2\, d\nu 
+ e^{a_0 NT}\, \int_{U\setminus V_b} L^N_{\f0} (H^{(m-1)})^2\, d\nu\\
& =     & \int_{U} \hath\,  (L_{\f0}^N (H^{(m-1)})^2)\, d\nu 
= \int_{U}   L_{\f0}^N ( (\hath \circ \sigma^N)\, (H^{(m-1)})^2)\, d\nu\\
& =     &  \int_{U} (\hath\circ \sigma^N)\,  (H^{(m-1)})^2\, d\nu .
\end{eqnarray*}
Similarly,
\begin{eqnarray*}
\int_{U} (\hath\circ \sigma^N)\,  (H_{m-1})^2 \, d\nu
& \leq  &  \int_{U}  (\hath\circ \sigma^{2N})\, (\hath\circ \sigma^N)\,  (H^{(m-2)})^2\, d\nu .
\end{eqnarray*}
Continuing by induction and using $H^{(0)} = 1$, we get
\be
\int_{U} (H^{(m)})^2 \, d\nu \leq   \int_{U}  (\hath\circ \sigma^{mN})\, (\hath\circ \sigma^{(m-1)N}) 
 \ldots  (\hath\circ \sigma^{2N})\, (\hath\circ \sigma^N)\, d\nu .
\ee

Let $s > 0$ be a constant. Using Lemma 7.2, choose the constants $C_7 > 0$ and $N_0 \geq 1$ so that
for $N \geq N_0$ we have 
$\mu(\Lambda_N(b)) \leq C_7/|b|^{12s} .$

Given an integer $k \geq 1$, set $m = k \hb$ and let  $0 < \delta_1 < \delta_2$ be as in (6.6). Let $N \geq N_0$.
Set 
$$W = \{ x \in U : x\in \sigma^{-jN}(U\setminus V_b) \: \mbox{\rm for at least $\delta_2 m$ values of }\: j = 0,1, \ldots,m-1 \} .$$
Since $K_0 \subset V_b$, $x\in \sigma^{-jN}(U\setminus V_b)$ implies  $x\in \sigma^{-jN}(U\setminus K_0)$.

Notice that 
\be
(\piU)^{-1}(W) \subset \Lambda_{N} (b) ,
\ee
the set defined by (\ref{eq:Lambda}). Indeed, if $x\in W$ and $y \in W^s_R(x)$, then for any $j = 0,1, \ldots, m-1$ with
$$\sigma^{jN}(x) \notin K_0 = \piU(P_1\cap P_0 \cap \Omb_B) ,$$
it follows from $\piU(\pp^{jN}(y)) = \sigma^{jN}(x)$ that  
$\pp^{jN}(y) \notin P_1 \cap P_0 \cap \Omb_B .$
Thus, the latter holds for at least $\delta_2 m$ values of $j = 0, 1, \ldots, m-1$, so 
$$\pp^i(y) \notin P_1 \cap P_0 \cap \Omb_B$$ 
for at least $\delta_2 m = \frac{\delta_2}{N} N  \hb$ values of
$i = 0,1, \ldots, N \hb-1$. It follows from (\ref{eq:Lambda}) that $y \in \Lambda_{N}(b)$. This proves (7.15), and now Lemma 7.2 implies
\be
\nu(W) \leq \frac{C_7}{|b|^{12 s}} .
\ee

Notice that if $x\in U\setminus W$, then $x\in \sigma^{-j N}(V_b)$ for at least $(1-\delta_2)m$ values of
$j = 0,1, \ldots, m-1$, so $(\hath\circ \sigma^{jN})(x) = \rho_3$ for that many $j's$.  Thus, 
using $\log(1+x) > 1 + x/2$ for $0 < x < 1$, (7.11) and (7.16) yield
\begin{eqnarray*}
\int_U (H^{(m)})^2\, d\nu 
& \leq & \int_{U\setminus W} \prod_{j=1}^{m} (\hath\circ \sigma^{jN}) \, d\nu  + \int_{W} \prod_{j=1}^{m} (\hath\circ \sigma^{jN}) \, d\nu \\
& \leq & \rho_3^{(1-\delta_2)m}S_0^{\delta_2 m} + S_0^{m} \nu(W) 
\leq (\rho_3^{1-\delta_2} S_0^{\delta_2})^m + \frac{C_7 \, S_0^m}{|b|^{12s}} \\
& \leq & \left(e^{(1-\delta_2)\log \rho_3 + \delta_2 a_0NT_0}\right)^{k \log |b| } + \frac{C_7 \, e^{a_0 NT_0 k \log |b|}}{|b|^{12 s}}\\
& \leq & \left(e^{(1-\delta_2) a_0NT_0 - (1-\delta_2)\log (1+\mu_0e^{-NT_0}/C_6) + \delta_2 a_0NT_0}\right)^{k \log |b| } \\
&       & + \frac{C_7}{|b|^{12s - a_0 NT_0k}}\\
& \leq & \left(e^{ a_0NT_0 - \frac{1}{2}\log (1+\mu_0e^{-NT_0}/C_6)}\right)^{k \log |b| } + \frac{C_7}{|b|^{12 s - a_0NT_0k}}\\
& \leq & \left(e^{a_0 N T_0 - \frac{1}{4}\mu_0e^{-NT_0}/C_6}\right)^{k \log |b| } + \frac{C_7}{|b|^{12 s- a_0 N T_0 k}}.
\end{eqnarray*}
Now choose 
$$k = \left\lceil  \frac{100 s C_6 e^{NT_0}}{\mu_0} \right\rceil ,$$
and assume that
$$a_0 \leq \frac{\mu_0}{100 C_6 NT_0 e^{NT_0} } ,$$
Then the above yields
\begin{eqnarray*}
\int_U (H^{(m)})^2\, d\nu 
& \leq &\left(e^{- \frac{\mu_0}{8 C_6 e^{NT_0}} }\right)^{k \log |b| } + \frac{C_7}{ |b|^{12s- \frac{\mu_0 k}{100 C_6 e^{NT_0}} }}\\
& \leq & \frac{1}{ |b|^{\frac{k \mu_0}{8 C_6 e^{NT_0}} }   } 
+ \frac{C_7}{ |b|^{12s- \frac{\mu_0}{100 C_6 e^{NT_0}} ( \frac{100 s C_6 e^{NT_0}}{\mu_0}+1)}}\\
& \leq & \frac{1}{|b|^{10 s}} + \frac{C_7}{ |b|^{11s- \frac{\mu_0}{100 C_6 e^{NT_0}} }} \leq \frac{1 + C_7}{|b|^{10 s}} .
\end{eqnarray*}
Thus, we can take $C_8 = 1 + C_7$.

\ms

(b) Let $h\in \ff_{\theta_1}(\hU)$ be such that 
$\| h\|_{\theta_1,b} \leq 1$.  Then  $|h(u)| \leq 1$ for all $u\in \hU$ and  $|h|_{\theta_1} \leq |b|$.  

Assume that the points $u,u'$, the cylinder $\cc = \cc_m$ in $U$, the integer $p \geq 0$ and the points $v,v'\in \hU_1$
satisfy (6.20)  for some $i = 1,2$. Then, using (6.19)  and $|h|_{\theta_1} \leq |b|$ we get
\begin{eqnarray*}
|h(v) - h(v')| 
& \leq & |b|\, D_{\theta_1} (v,v') = |b|\, \theta_1^{p+N}\, D_{\theta_1} (u,u')
\leq  |b|\, \theta_1^{p+N}\, \diam_{\theta_1}(\cc) \\
& \leq & |b|\, \theta_1^{p+N}\, C_2 \, \diam(\tPsi(\cc)) \leq E  |b|\, \theta_2^{p+N} \, \diam(\tPsi(\cc)) ,
\end{eqnarray*}
since $C_2 \leq E$. Thus, $(h, 1) \in \kk_b$.
Set $h^{(m)} = \lab^{m N} h$ for $m \geq 0$. Define the sequence of functions $\{ H^{(m)}\}$ recursively by
$H^{(0)} = 1$ and $H^{(m+1)} = \nn_{J_m} H^{(m)}$, where $J_m \in \jj(b)$ is chosen by induction as follows.
Since $(h^{(0)}, H^{(0)}) \in \kk_b$,  using Lemma 6.9  we find $J_0 \in \jj(b)$ such that
for $h^{(1)} = \lab^{N} h^{(0)}$ and $H^{(1)} = \nn_{J_0} H^{(0)}$ we have
$(h^{(1)} , H^{(1)}) \in \kk_b$. Continuing in this way we construct by induction
an infinite sequence of functions $\{ H^{(m)}\}$ with
$H^{(0)} = 1$, $H^{(m+1)} = \nn_{J_m} H^{(m)}$ for all $m \geq 0$, such that $(h^{(m)} , H^{(m)}) \in \kk_b$.

Next, a choose $a_0$ and $k \geq 1$ as in part (a) and set $m = k \hb$.  Then part (a) implies
$$\int_U (H^{(m)})^2 d\nu \leq \frac{C_8}{|b|^{10 s}} .$$
Hence
$$\int_{U} |\lab^{m N} h|^2 \; d\nu  = \int_{U} |h^{(m)}|^2\; d\nu \leq \int_{U}  (H^{(m)})^2\; d\nu  
\leq \frac{C_8}{|b|^{10s}}.$$
From this it follows that for any $h \in \ff_{\theta_1}(\hU)$ we have 
$$\int_{U} |\lab^{m N} h|^2 \; d\nu \leq  \frac{C_8}{|b|^{10s}} \, \|h\|^2_{\theta_1, b} ,$$
and so
\be
\int_{U} |\lab^{m N} h| \; d\nu \leq  \frac{\sqrt{C_8}}{|b|^{5 s}} \,\|h\|_{\theta_1,b} .
\ee

We will now use a standard procedure (see \cite{D1}) to derive an estimates of the form (7.13) from (7.17).

\bs

\noindent
{\bf\large 7.4 Standard procedure:} First, recall from the Perron-Ruelle-Frobenius Theorem (see e.g. \cite{PP}) 
that there exist global constants $C'_9 \geq 1$ and $\rho_4 \in (0,1)$, independent of $b$ and $N$, such that 
\be
\| L_{\f0}^n w - h_0 \, \int_{U} w \, d\nu \| \leq C'_9\, \rho_4^n \, \|w\|_{\theta}
\ee
for all $w \in \ff_{\theta}(\hU)$ and all integers $n \geq 0$, where $h_0 > 0$ is the normalised eigenfunction of 
$L_{f-P_f \tau}$ in $\ff_{\theta}(\hU)$ (see the beginning of Sect. 5.1).

Given $h \in \ff_{\theta_1}(\hU)$ with  $\|h\|_{\theta_1,b} \leq 1$, we have
$|h|_\theta \leq |h|_{\theta_1} \leq |b|$, so using Lemma 5.2 with $H = 1$ yields
\be
|\lab^{r} h|_{\theta} \leq A_0 [ |b| \theta^{r} + |b| ]\leq 2A_0 |b|
\ee
for any integer $r \geq 0$.

Choose again $a_0$ and $k \geq 1$ as in the proof of part (a) and set $m = k \hb$.  
Then (7.17) holds.
Write $\rho_4 = e^{-\beta_3}$ for some global constant $\beta_3 > 0$. Given $h \in \ff_{\theta_1}(\hU)$
with $\|h\|_{\theta_1,b} \leq 1$, we have
\begin{eqnarray*}
|\lab^{2mN} h|
& =    & |\lab^{mN} \left(|\lab^{mN} h |\right)| \leq \ma^{mN} |\lab^{mN} h| 
= L_{\f0}^{mN} \left(e^{\fa_{mN} - \f0_{mN}}\, |\lab^{mN} h|\right)\\
& \leq & \left(L_{\f0}^{mN}\left(e^{\fa_{mN} -\f0_{mN}}\right)^{2}\right)^{1/2}\,
\left( L_{\f0}^{mN} |\lab^{mN} h|^{2} \right)^{1/2} .
\end{eqnarray*}
For the first term in this product (5.3) implies 
$$\left(L_{\f0}^{mN} \left(e^{\fa_{mN} - \f0_{mN}})\right)^{2}\right)^{1/2} \leq e^{a_0 NT_0 m} 
\leq  e^{a_0 N T_0 k \log |b|} = |b|^{a_0 NT_0 k}  .$$
By the choice of $k$ and $a_0$,
$$a_0 NT_0 k \leq a_0 N T_0 \left(\frac{100 s C_6 e^{NT_0}}{\mu_0} +1 \right)
\leq  \frac{\mu_0 N T_0}{100 C_6 NT_0 e^{NT_0} }  \, \frac{100 s C_6 e^{NT_0}}{\mu_0}   + a_0 N  T_0 < s + \frac{1}{2} ,$$
assuming $a_0 NT_0 < 1/2$. Thus,
$$\left(L_{\f0}^{mN} \left(e^{\fa_{mN} - \f0_{mN}})\right)^{2}\right)^{1/2} \leq |b|^{s+1/2} .$$

For the second term, using (7.18) with $w = |\lab^{mN} h|$, we get
$$L_{\f0}^{mN} |\lab^{mN} h|^2 \leq L_{\f0}^{mN} |\lab^{mN} h| \leq \|h_0\|\, \int_U |\lab^{mN} h|\, d\nu 
+ C'_9\, \rho_4^{mN} \, \|\lab^{mN} h\|_{\theta} $$
for some global constant $C'_9 > 0$.
By (7.19), 
$$\|\lab^{mN} h\|_{\theta} \leq 2 A_0 |b| ,$$
and  (7.17) implies
$$L_{\f0}^{mN} |\lab^{mN} h|^2 \leq \frac{C''_9}{|b|^{5s}} + 2 A_0 C''_9 |b| \rho_4^{mN} $$
for another global constant $C''_9 > 0$.
Since
$$\rho_4^{mN} \leq e^{-\beta_3 N \log |b|} = \frac{1}{|b|^{N \beta_3}} < \frac{1}{|b|^{5s+1}} ,$$
assuming $N \beta_3 > 5s+1$,  we get
$$L_{\f0}^{mN} |\lab^{mN} h |^2 \leq \frac{C'''_9}{|b|^{5s}} .$$
Combining the estimates of the two terms and using $s \geq 1$, it follows that
$$|\lab^{2mN} h | \leq  |b|^{s+1/2} (C'''_9/|b|^{5s})^{1/2}
\leq \frac{C_9}{|b|^{s}} ,$$
for some global constant $C_9 > 0$,
assuming $a_0 = a_0(N) > 0$ is chosen sufficiently small, e.g. $a_0NT_0 < s/2$.

The estimate (7.13) now follows immediately.
\endofproof

\newpage

\section{Proofs of the main results}
\setcounter{equation}{0}

Here we prove Theorems 1.3 and 1.1 and Corollary 1.4.

\bs

\noindent
{\it Proof of Theorem } 1.3. Let $\theta_1 = \theta_1(\theta) \in (0,\theta]$ and  
$\theta_2 = \theta_2(\theta)\in [\theta,1)$  be as in the beginning of Sect. 5.2. 

We will again assume that $f \in \ff_{\theta_1}(\hU)$; the general case $f \in \ff_\theta(\hU)$
will be done later using an approximation procedure.

Let $\hat{\theta} \leq \theta < 1$, where $\hat{\theta}$ is as in (5.1). Set 
$$ s =  \frac{2}{\alpha_2} ,$$
where $\alpha_2 > 0$ is the constant from Lemma 5.1(c), and recall that $\theta_1^{\alpha_2} = \theta$.
Next, choose $N_0 \geq 1$ as in Lemma 7.2, replacing $s$ by $12s$, so that
$\mu(\Lambda_N(b)) \leq C_7/|b|^{12s}$, as in the proof of Lemma 7.3.
Let $N \geq N_0$. Choose $k = k(N)$, $a_0 = a_0(N)$, $b_0 = b_0(N)$, $\rho_4 = \rho_4(N) \in (0,1)$, 
$C_8, C_9 > 0$ as in Lemma  7.3. Then (7.12) and (7.13) hold.

Let $|a|\leq a_0$ and  $|b|\geq b_0$, and let $h\in \ff_{\theta}(\hU)$ be such that 
$\| h\|_{\theta,b} \leq 1$.  Then  $|h(u)| \leq 1$ for all $u\in \hU$ and  $|h|_{\theta} \leq |b|$.

Take the smallest integer $p$ so that $\theta^p \leq 1/|b|^2$. It is known  (see e.g. the end of Ch. 1 in \cite{PP}) that there exists 
$h' \in \ff_{\theta_1}(\hU)$ which is constant on cylinders of  length $p$ so that $\|h - h'\|_0 \leq |h|_\theta \, \theta^p$. Then 
$\|h-h'\|_0  \leq 1/|b|$ and so $\|h'\|_0 \leq 2$, and it follows easily from this that 
$$|h'|_{\theta_1} \leq \frac{4}{\theta_1^{p-1}} \leq \frac{4}{\theta^{(p-1)/\alpha_2}}  \leq C'_{10} |b|^{2/\alpha_2} .$$
Thus, $\|h'\|_{\theta_1,b} \leq 2 C'_{10} |b|^{2/\alpha_2 - 1}$ and (7.13) gives 
$$\|\lab^{2k N \hb} h'\|_{\theta_1,b} \leq \frac{C_9}{|b|^{s}} \; 2 C'_{10} |b|^{2/\alpha_2 - 1} \leq
\frac{2 C'_{10} C_9}{|b|^{2/\alpha_2 - 2/\alpha_2+1}}  = \frac{2C'_{10} C_9}{|b|} .$$
Hence
$$\|\lab^{2k N  \hb} h'\|_{0} \leq \frac{C_{10}}{|b|} .$$
for  some global constant $C_{10} > 0$.  Therefore
$$\|\lab^{2 k N \hb} h\|_0 \leq \|\lab^{2 k N \hb} h'\|_0 + \|\lab^{2 k N \hb} (h-h')\|_0
\leq \frac{C_{10}}{|b|} + \frac{1}{|b|} \leq \frac{2C_{10}}{|b| } .$$

Next, using Lemma 5.2, and writing $\theta = e^{-\beta_5}$ for some constant $\beta_5 > 0$, we get
\begin{eqnarray*}
|\lab^{4 k N \hb} h|_\theta 
 &   =  &  |\lab^{2k N \hb} ( \lab^{2 k N \hb} h) |_\theta 
 \leq   A_0 \left[  |b| \, \theta^{2 k N \hb}  + |b| \, \| \lab^{2 k N \hb} h\|_0 \right] \\
 & \leq & A_0\left[ |b| \; \frac{1}{|b|^{2k N \beta_5} } + |b| \; \frac{2C_{10}}{|b|}\right]\leq C'_{11}  ,
\end{eqnarray*}
assuming $2 N\beta_5 \geq 1$. 
This yields
$$\|\lab^{4 k N \hb} h\|_{\theta,b} \leq \frac{C_{11}}{|b|}\, \|h\|_{\theta,b}$$
for all $h \in \ff_\theta(\hU)$, where $C_{11} > 0$ is a global constant.

Let $n \geq 4k N \hb$ be an arbitrary integer. Writing $n = r( 4k N \hb)+ \ell$ for some
$\ell  = 0,1, \ldots, 4k N \hb-1$, and using the above $r$ times we get
$$\|\lab^{r 4k N \hb} h\|_{\theta,b} \leq \frac{1}{|b|^{r}}\, \|h\|_{\theta,b} .$$
As before, using Lemma 5.2 with $H = 1$ and $B = |\lab^{ r4k N \hb} h|_\theta$, implies
\begin{eqnarray*}
|\lab^{n} h|_\theta = |\lab^{\ell} ( \lab^{  r 4 k N \hb} h) |_\theta 
 \leq   A_0 \left[ |\lab^{r 4 k  N \hb} h|_\theta \, \theta^{\ell} + |b| \, \| \lab^{r 4 k N \hb} h\|_0 \right] ,
\end{eqnarray*}
so 
$$\frac{1}{|b|} |\lab^{n} h|_\theta \leq 2A_0 \|\lab^{r 4k  N \hb} h\|_{\theta,b} \leq \frac{2A_0 }{|b|^{r}}\, \|h\|_{\theta,b} .$$
This and 
$$\|\lab^{n} h\|_0 \leq \|\lab^{r 4k  N\hb} h|\|_0 \leq \frac{1}{|b|^{r}}\, \|h\|_{\theta,b}$$ 
give
\begin{eqnarray*}
\|\lab^{n} h\|_{\theta,b} 
& \leq & \frac{3A_0}{|b|^{r}}\, \|h\|_{\theta,b} = 3A_0  e^{- r \log |b|} \|h\|_{\theta,b} .
\end{eqnarray*}
We have $r \geq (r+1)/2$ for all $r \geq 1$, so the above implies
\begin{eqnarray} 
\|\lab^{n} h\|_{\theta,b} 
&  \leq &  3A_0  e^{-\frac{(r+1)  \log |b|}{2}} \|h\|_{\theta,b}
 \leq  3A_0  e^{-\frac{(r+1)4 k N  \hb }{8 k N}} \|h\|_{\theta,b} \\
& \leq & 3A_0  \rho_6^n \|h\|_{\theta,b} ,\nonumber
\end{eqnarray}
where $\rho_6 = e^{-1/(8 k  N)} \in (0,1)$. 

Thus, (8.1) holds for all $h \in \ff_\theta(\hU)$ and all integers $n \geq 4 k N \hb = 4 k  N \lceil\log |b|\rceil$. 
Finally, recall the eigenfunction $h_a \in \ff_\theta(\hU)$ for the operator $L_{f- (P_f+a)\tau}$
from Sect. 5.1. It is known that $\|h_a\|_\theta \leq \Con$ for bounded $a$, e.g. for $|a|\leq a_0$.
It now follows from
$$\lab^n(h/h_a) = \frac{1}{\lambda_a^n h_a}\, L^n_{f- (P+a+\i b) \tau} h$$
and the above estimate that there exist constants $0 < \rho < 1$, $a_0 > 0$, 
$b_0 \geq 1$ and  $C > 0$ such that if $a,b\in \R$  satisfy $|a| \leq a_0$ and $|b| \geq b_0$, 
then 
\be
\|L_{f -(P_f+a+ \i b)\tau}^n h \|_{\theta,b} \leq C \;\rho^n \; \| h\|_{\theta,b}
\ee
for any integer $n \geq 4 k N \log |b|$ and any  $h\in \ff_\theta (\hU)$. So, we can just set $T = 4k N$.

\bs

This completes the proof of Theorem 1.3 under the assumption that $f \in \ff_{\theta_1}(\hU)$. The case $f \in \ff_{\theta}(\hU)$ 
follows by using an approximation procedure. To our knowledge this has not been done anywhere in details, and the argument 
involved is  not trivial, so we will sketch it for completeness.

\bs

\noindent
{\bf\large Sketch of the proof of Theorem 1.3 for arbitrary $f \in \ff_\theta(\hU)$:}\\
We will use again the constants from the beginning of Sect. 6.1, including $\theta_1$, $\theta_2$, etc.
Fix $B$, $N$ as in Sect. 6.1 and define $\hb$ by (\ref{eq:hbdef}). 
Let $|a| \leq a_0$ and $|b| \geq b_0$, where $b_0$ is given by (\ref{eq:b0cond}). 

Let $f \in \ff_\theta(\hU)$ be an arbitrary real-valued function. Take the minimal integer $t = t(b)  > 0$ so that
\be
\theta^{t+1} \leq  \frac{2A_0}{\log|b|}  ,
\ee
where $A_0$ is the constant from Lemma 5.2. There exists a real-valued $\ftt$ depending only on $t$ coordinates
such that 
$$\|f- \ftt\|_0 \leq |f|_\theta \, \theta^t \leq T_0\, \theta^t $$
(see the end of Ch. 1 in \cite{PP}), where $T_0$ is as in (5.3). Then $\ftt \in \ff_{\theta_1}(\hU)$,
$\|\ftt\|_0 \leq 2T_0$ and 
$$|\ftt|_\theta \leq \frac{4}{\theta^t} \leq \frac{2 \log |b|}{A_0} \quad, \quad |\ftt|_{\theta_1} \leq \frac{4}{\theta_1^t}.$$

Let $\lambda_{at}$ be the {\it largest eigenvalue} of
$$\Fat =  \ftt - (P_t + a)\tau,$$
where $P_t = P_{\ftt}$, and let $h_{at} \in \ff_\theta(\hU)$ be a corresponding (positive) eigenfunction such that 
$$\int_{\hU} h_{at} \, d\hnu_{at} = 1 ,$$
 where $\hnu_{at}$ is the unique regular probability measure on $\hU$ with  $(\Fat)^* \hnu_{at} = \hnu_{at}$.

\def\nnt{\nn^{(t)}}

For $|a|\leq a_0$,  as in \cite{D1}, consider the function
$$\fat(u) = \ftt (u) - (P_t + a) \tau(u) + \ln h_{at}(u) -  \ln h_{at}(\sigma(u)) - \ln \lambda_{at}$$
and the operators 
$$ \labt = L_{\fat - \i\,b\tau} : \ff_\theta(\hU) \longrightarrow \ff_\theta (\hU)\:\:\: , 
\:\:\: \mat = L_{\fat} : \ff_\theta(\hU) \longrightarrow \ff_\theta(\hU) .$$
Then $\mat \; 1 = 1$ and $\di |(\labt^m h)(u)| \leq (\mat^m |h|)(u)$ for all $u\in \hU$.

Using part of the proof of Lemma 4.1 in \cite{PeS5}, one shows that $|h_{at}|_\theta \leq \Con\, |\ftt|_\theta$ for some global constant
$\Con > 0$. Thus, $|\fat|_\theta \leq \Con\, |\ftt|_\theta$, and it is also clear that $\|\fat\|_0 \leq \Con$.

Next, define the set $K_0$, cylinders $\cc_m$ and their sub-cylinders $\dd_j$ and the function
$\omega_J$ as in Sect. 6.1 and consider the operator $\nnt = \nnt_J$ on $\ff_\theta(\hU)$
defined by
$$\nnt(h) = \mat^N (\omega_J\, h) = L^N_{\fat} (\omega_J\, h) .$$
It is important to notice that 
\be
e^{|\ftt|_\theta} \; \diam_\theta(\cc_m) \leq \frac{1}{|b|^{1/A_0} },
\ee
provided we took the constants $A_0$ in Lemma 5.2 and  $B$ in Sect. 6.1 so that $A_0 \geq \frac{4B}{|\log \theta|}$.
Indeed, for the length $\ell_m$ of $\cc_m$ we have (\ref{eq:Bcond}), so
\begin{eqnarray*}
e^{|\ftt|_\theta}\; \diam_\theta(\cc_m) 
& \leq &  e^{\frac{2\log |b|}{A_0}} \theta^{\ell_m} =
|b|^{2/A_0} e^{-\ell_m |\log\theta| } \leq  |b|^{2/A_0} e^{-(|\log\theta|/B) \log|b| }\\
& =    & |b|^{2/A_0 - \frac{1}{B} |\log \theta|} \leq |b|^{-1/A_0} ,
\end{eqnarray*}
which proves (8.4).

Then we define the metric $\dd(u,u')$ on $\hU$ and the class of positive functions $\kk_E$ as in Sect. 6.2.
Now with the above one easily shows that Lemma 6.5 is valid in the form $\nnt (\kk_E) \subset \kk_E$.
Indeed, the main observation to make to prove this is that, given $u, u' \in \hU$ such that  there exists
an integer $p \geq 0$ with $\sigma^p(Y(u, u')) \subset \cc'_m$ for some $m \leq m_0$ and $\ell(u,u') \geq p$,
then for any integer $k \geq 1$, if $v, v'(v) \in U$ satisfy $\sigma^{k}(v) = u$, $\sigma^k(v') = u'$
and belong to the same cylinder of length $k$, then
\begin{eqnarray}
&        & |\fat_k(v) - \fat_k(v')|  \leq  \sum_{j=0}^{m-1} |\fat|_\theta \,\theta^{m-j}\, \dte (u,u')\\
& \leq  & \Con\; |\ftt|_\theta\,D_\theta(u,u') 
 \leq  \Con\; |\ftt|_\theta\, \diam_\theta(\cc_m) \leq \Con .\nonumber
\end{eqnarray}
With this observation, a simple modification of the proof of Lemma 6.5 gives $\nnt H \in \kk_E$ for
every $H \in \kk_E$.

Next, we define the class of functions $\kk_0$ as in Sect. 6.3 and prove the analogue of Lemma 6.7:
$L_{\ftt - s \tau}(\kk_0) \subset \kk_0$ for all $s$ with $|s| \leq M_1$ and all integers $q\geq N$.
To prove this, the choice of $\theta_1$ is important; it implies (using Lemma 5.1 and $\theta_1^{\alpha_2} = \theta$)
$$\diam_{\theta_1}(\cc_m) = \theta_1^{\ell_m} \leq \theta_1^{\ell_m/2} \; (\theta^{\ell_m})^{1/\alpha_2}
\leq \theta_1^{\ell_m/2} \; \diam(\tcc_m) .$$
Then, assuming $u, u', v, v', w, w'$ are as in the proof of Lemma 6.7, we derive
\begin{eqnarray}
|\fat_N(w) - \fat_N(w')| 
& \leq & \Con\, |\ftt|_{\theta_1} \theta_1^{p+N} \diam_{\theta_1}(\cc_m)\\
& \leq & \Con \; \frac{4}{\theta_1^t} \theta_1^{p+N} \, \theta_1^{\ell_m/2} \; \diam(\tcc_m) \nonumber\\
& \leq & \Con \, \theta_2^{p+N} \theta_1^{\ell_m/2 - t} \; \diam (\tcc_m) \nonumber\\
& \leq &  \Con \, \theta_2^{p+N} \; \diam (\tcc_m)  \leq 1,\nonumber
\end{eqnarray}
since $t < < \ell_m/2$. Now the rest of the proof of Lemma 6.7 is the same, and as a consequence one
gets (as in Corollary 6.8) that the eigenfunctions $h_{at}$ belong to $ \kk_0$.

Finally, the arguments in Sect. 6.4 can be repeated with very little change -- the main one is that
in the first estimate of $|(\labt^Nh)(v) - (\labt^N)(v')|$ one has to use (8.5) again. This proves the
analogue of Lemma 6.9, where the operator $\lab^N$ is replaced by $\labt$.

We will now prove Lemma 6.9 in its original form under the present assumption that $f \in \ff_\theta(\hU)$.

\bs

\noindent
{\bf Lemma 8.1.}  {\it Assume $f\in \ff_\theta(\hU)$.
Choosing $E > 1$ and $\mu_0$ as in Sect. {\rm 5.2} and assuming $N$ is sufficiently large, 
for any $|a| \leq a_0$, any $|b| \geq b_0$ and any $(h,H) \in \kk_b$ 
there exists $J \in \jj(b)$ such that $(\lab^{N} h, \nn_J H ) \in \kk_b$.}

\bs

\noindent
{\it Proof.} Consider the function
$$\zeta = \fa_N - \fat_N \in \ff_\theta(\hU) .$$
Notice that for any $u \in \hU$ and any function $h$ on $\hU$ we have
\begin{eqnarray*}
(\ma^N h)(u) 
& =  & \sum_{\sigma^Nv = u} e^{\fa_N(v)} h(v)
=   \sum_{\sigma^Nv = u} e^{\fat_N(v)} e^{(\fat_N - \fa_N)(v)} h(v)\\
& = & (\mat^N (e^{\zeta} h))  (u) .
\end{eqnarray*}
Thus, $\ma^N h = \mat^N (e^\zeta h)$, and similarly one observes that
 $\lab^N h = \labt^N (e^\zeta h)$.
 
 We will now repeat the argument from the proof of Lemma 6.9.

Let $|a| \leq a_0$, $|b| \geq b_0$ and $(h, H) \in \kk_b$. We will construct a representative set 
$J \in \jj(b)$ such that  $(\lab^N h , \nn_J H) \in \kk_b$.  Given  an arbitrary representative set $J$,
we will first show that  
$$(\lab^N h , \nn_J H) = (\labt^n (e^\zeta h), \nnt_J (e^\zeta H)) $$ 
has property $(T2)$ from Sect. 6.4. 

Assume that the points $u,u'$, the cylinder $\cc_m$ in $U$, the integer $p \geq 0$ and the points 
$v,v'\in \hU$ satisfy (\ref{eq:u-cond})  for some $i = 1,2$ and $\ell = 1, \ldots, \ell_0$. 
Since $h_a, h_{at} \in \kk_0$, we have
$$|\ln h_a(w) - \ln h_a(w') | \leq \frac{| h_a(w) - h_a(w')|}{\min \{ |h_a(w)|, |h_a(w')|\}} 
\leq E_1 \,\theta_2^{p+2N} \, \diam(\tcc_m) , $$
and similarly, $|\ln h_{at}(v) - \ln h_{at}(v') | \leq E_1 \,\theta_2^{p+N} \, \diam(\tcc_m)$.
By (8.6),
\be
|\fat_N(w) - \fat_N(w')| \leq \Con\, \theta_2^{p+N}\, \diam(\tcc_m)  \leq 1,
\ee
assuming $N$ is chosen appropriately. Thus, using (8.6) and a 
similar but simpler estimate for $| \fa_N(w) - \fa_N (w')|$, we get
\begin{eqnarray}
|\zeta(w) - \zeta (w')|
& =    & | (\fa - \fat)_N(w) - (\fa - \fat)_N (w')|\\
& \leq &  C'_{12} \, \theta_2^{p+N} \; \diam (\tcc_m) < C'_{12}\nonumber
\end{eqnarray}
for some global constant $C'_{12} > 0$.
This implies
\begin{eqnarray}
\qquad |e^{\zeta (w) - \zeta (w')} - 1|
 \leq  \Con\, \left| \zeta (w) - \zeta (w') \right|  \leq  C_{12} \, \theta_2^{p+N} \; \diam (\tcc_m) < C_{12}
\end{eqnarray}
for some global constant $C_{12} > 0$.

Hence for any $a$ and $b$ with $|a| \leq a_0$ and $|b| \geq b_0$   we have:
\begin{eqnarray*}
&        &|(\lab^N h)(v) - (\lab^N  h)(v')| = |(\labt^N (e^\zeta h))(v) - (\labt^N  (e^\zeta h))(v') |\\
&  =    &      \left| \sum_{\sigma^N w = v} e^{(\fat_N - \i b \tau_N) (w)}\, e^{\zeta (w)} h(w) 
 -  \sum_{\sigma^N w = v} e^{(\fat_N - \i b \tau_N)(w'(w))}\, e^{\zeta (w')} \, h(w'(w)) \right| \\
& \leq & \left| \sum_{\sigma^N w = v} e^{(\fat_N - \i b \tau_N)(w)}\, e^{\zeta (w)}\, [h(w) -  h(w')]\right|\\
&        &  +
\left| \sum_{\sigma^N w = v} e^{(\fat_N - \i b \tau_N)(w)}\, \left[e^{\zeta (w)} - e^{\zeta (w')}  \right]\,   h(w')\right|  \\
&       & + \sum_{\sigma^N w = v}  \left|e^{(\fat_N - \i b \tau_N)(w)} -  e^{(\fat_N - \i b \tau_N) (w')}\right| \,e^{\zeta (w')}\, |h(w')| \\
& \leq & \sum_{\sigma^N w = v} e^{(\fat_N (w) - \fat_N(w'))} e^{\fat_N (w')}\,E |b| 
\theta_2^{p+2N} \diam (\tcc_m)\,  e^{C'_{12}}\, e^{\zeta (w')} \, H(w') \\  
&      & + \sum_{\sigma^N w = v} e^{(\fat_N (w) - \fat_N(w'))} e^{\fat_N (w')}\, C_{12} \, \theta_2^{p+N} \; \diam (\tcc_m) \, e^{\zeta (w')} H(w') \\
&       & + \sum_{\sigma^N w = v}  \left| e^{(\fat_N - \i b \tau_N)(w) - 
(\fat_N - \i b  \tau_N)(w')} - 1\right| \,  e^{\fat_N(w')} e^{\zeta (w')}\, H(w') \\
& \leq & e^{1+C'_{12}}\, E |b|\, \theta_2^{p+2N} \diam (\tcc_m) \, \mat^N (e^\zeta  H)(v') 
+  C_{13}   \theta_2^{p+N} \, \diam(\tcc_m)\, \mat^N (e^\zeta H)(v')\\
&      & +  C_{13} |b|  \theta_2^{p+N} \, \diam(\tcc_m)\, \mat^N (e^\zeta H)(v')\\
& \leq & [2 e^{1+ C'_{12}}\, E \theta_2^N  + 2 C_{13} ]  \,|b|\, \theta_2^{p+N}\,    \diam(\tcc_m)\, \nnt_J (e^\zeta H)(v') 
\end{eqnarray*}
for some global constant $C_{13} > 0$.
Assuming $2 e^{1+C'_{12}} \theta^N \leq 1/2$ and $2 C_{13}  \leq E/2$, we get
\begin{eqnarray*}
|(\lab^N h)(v) - (\lab^N  h)(v')| 
& \leq &  E \,|b|\, \theta_2^{p+N}\,    \diam(\tcc_m)\, \nnt_J (e^\zeta H)(v') \\
& =     & E \,|b|\, \theta_2^{p+N}\,    \diam(\tcc_m)\, (\nn_J H)(v') , 
\end{eqnarray*}
so, $(\lab^N h , \nn_J H)$ has property $(T2)$ from Sect. 6.4.

Now  we will construct $J$ so that  $|\lab^N h|(u) \leq (\nn_J H) (u)$ for all $u \in \hU$, which is equivalent to
\be
|\labt^N (e^\zeta h)|(u) \leq (\nnt_J (e^\zeta H)) (u)
\ee
for all $u \in \hU$. 

Define the functions 
$$\psi_\ell, \gao_\ell, \gat_\ell : \hU  \longrightarrow \C$$ 
 as in the proof of Lemma 6.9. Notice that
$$\di \psi_\ell(u) = e^{(\fat_{N}+\i b\tau_{N})(\vl_1(u))} (e^{\zeta} h) (\vl_1(u)) 
+ e^{(\fat_{N}+\i b\tau_{N})(\vl_2(u))} (e^\zeta h) (\vl_2(u)) ,$$
$$\di \gao_\ell(u) = (1-\mu_0)\, e^{\fat_{N} (\vl_1(u))} (e^\zeta H) (\vl_1(u)) + e^{\fat_{N}(\vl_2(u))} (e^\zeta H) (\vl_2(u)) ,$$
and
$$\di \gat_\ell(u) = e^{\fat_{N} (\vl_1(u))} (e^\zeta H) (\vl_1(u)) + (1-\mu_0)\,  e^{\fat_{N}(\vl_2(u))} (e^\zeta H) (\vl_2(u)) .$$
We will use again the functions $\varphi_{\ell,m}(u) = \varphi_\ell(Z_m, u)$, $u \in U$,  from Sect. 5.3. 

As before (8.10) is trivially satisfied for $u \notin V_b$ for any choice of $J$.

Consider an arbitrary  $m = 1, \ldots, m_0$. We will construct  $j \leq j_0$ with $\dd_{j} \subset \cc_m$,  
and a pair  $(i, \ell)$ for which $(i,j,\ell)$ will be included in $J$.

\ms

\noindent
{\bf Case 1.}
There exist $j\leq j_0$ with $\dd_j \subset \cc_m$, $i = 1,2$ and $\ell \leq \ell_0$   such that the 
first alternative in Lemma 6.10(b) holds for $\hdd_j$,  $i$ and $\ell$. This case is dealt with exactly
as in the proof of Lemma 6.9.
 
\medskip

\noindent
{\bf Case 2.}
For all $j\leq j_0$ with $\dd_j \subset \cc_m$, $i = 1,2$ and $\ell \leq \ell_0$ the 
second alternative in Lemma 6.10(b) holds for $\hdd_j$,  $i$ and $\ell$, i.e.
\be
|h(\vl_i(u))|\geq \frac{1}{4}\, H(\vl_i(u)) > 0
\ee
for any  $u \in \hcc'_m$.

Let $u,u' \in \hcc'_m$, and let $i = 1,2.$ 
Using (6.20) and the assumption that $(h,H) \in \kk_b$, and in particular property 
$(T2)$ with $p = 0$, $v = \vl_i(u)$ and $v' = \vl_i(u')$, and also (8.8) and (8.9) with $p =0$, 
and assuming e.g.  
$$\min\{ |e^{\zeta(\vl_i(u))} h(\vl_i(u))| , |e^{\zeta(\vl_i(u'))} h(\vl_i(u'))|\}  = |e^{\zeta(\vl_i(u'))} h(\vl_i(u'))| ,$$ 
we get
\begin{eqnarray*}
&        & \frac{|e^{\zeta(\vl_i(u))} h(\vl_i(u)) - e^{\zeta(\vl_i(u'))} h(\vl_i(u'))|}{\min\{ |e^{\zeta(\vl_i(u))} h(\vl_i(u))| , |e^{\zeta(\vl_i(u'))} h(\vl_i(u'))| \}}\\
& \leq &   \frac{|e^{\zeta(\vl_i(u))}  - e^{\zeta(\vl_i(u'))}|\; | h(\vl_i(u'))|}{|e^{\zeta(\vl_i(u'))} h(\vl_i(u'))|}
+  \frac{e^{\zeta(\vl_i(u))} | h(\vl_i(u)) -  h(\vl_i(u'))|}{ |e^{\zeta(\vl_i(u'))} h(\vl_i(u'))| }\\
& \leq &   |e^{\zeta(\vl_i(u))  - \zeta(\vl_i(u'))} - 1|
+  C_{14} \frac{ | h(\vl_i(u)) -  h(\vl_i(u'))|}{ | h(\vl_i(u'))| }\\
& \leq & C_{14} \theta_2^N \diam(\tcc_m) + \frac{E|b|\,\theta_2^N H(\vl_i(u'))}{|h(\vl_i(u'))| } \diam (\tcc_m) \\
&  \leq &  (C_{14} + 4 E|b|) \,   \theta_2^{N} \frac{A R_0 \ep_1}{|b|} < 5  E \theta_2^N \ep_1 ,
\end{eqnarray*}
assuming $E\geq C_{14} A R_0$. So, the difference  between the arguments of the vectors
$$e^{\zeta(\vl_i(u))} h(\vl_i(u))\quad \mbox{\rm  and} \quad e^{\zeta(\vl_i(u'))} h(\vl_i(u'))$$ 
in $\R^2$  is   $< 10 E A R_0\theta_2^N \ep_1 < \pi/6$ by (\ref{eq:Ncond}).  
Since $e^{\zeta(\vl_i(u)} $ and $e^{\zeta(\vl_i(u))} $ are real numbers, the arguments of the complex numbers
$$e^{\zeta(\vl_i(u))} h(\vl_i(u)) \quad \mbox{\rm and} \quad e^{\zeta(\vl_i(u'))} h(\vl_i(u'))$$
are the same as those of  $h(\vl_i(u))$ and $h(\vl_i(u'))$. 

As before,  for each $i = 1,2$ we can choose a real continuous
function $\thetam_i(u)$, $u \in  \cc'_m$, with values in $[0,\pi/6]$  and a constant $\lambdam_i$ such that
$$\di h(\vl_i(u)) = e^{\i(\lambdam_i + \thetam_i(u))}|h(\vl_i(u))| \quad , \quad u\in \cc'_m .$$
Fix an arbitrary $u_0\in \cc'_m$, set $\lambdam = |b| \varphi_{\ell,m}(u_0)$, and assume again that 
 $$|\lambdam_2 - \lambdam_1 + \lambdam | \leq \pi .$$
 Then
$$|\thetam_i(u) - \thetam_i(u')|\leq 2 \sin |\thetam_i(u) - \thetam_i(u')| < 16 E A R_0 \theta_2^N \ep_1$$
for all $u,u' \in \cc'_m$.

As in the proof of Lemma 6.9, 
the difference between the arguments of the complex numbers
$$e^{\i \,b\,\tau_N(\vl_1(u))} (e^\zeta h) (\vl_1(u)) \quad \mbox{\rm and} \quad e^{\i \,b\, \tau_N(\vl_2(u))} (e^\zeta h) (\vl_2(u))$$
is given by the function
\begin{eqnarray*}
\Gl(u) 
& = & [b\,\tau_N(\vl_2(u)) + \thetam_2(u) + \lambdam_2] -  
[b\, \tau_N(\vl_1(u)) + \thetam_1(u) + \lambdam_1]\\
& = & (\lambdam_2-\lambdam_1) + |b| \varphi_{\ell,m}(u) + (\thetam_2(u) - \thetam_1(u)) ,
\end{eqnarray*}
and as before we prove that there exist $j \leq j_0$ and $\ell \leq \ell_0$ such that
$\ep_3 \leq |\Gl(u)| <  \frac{3\pi}{2}$ for all $u \in \dd'_j$.

As in the proof of Lemma 6.9 either $H(\vl_1(u)) \geq H(\vl_2(u))/4$ for all
$u \in \dd'_j$ or $H(\vl_2(u)) \geq H(\vl_1(u))/4$ for all $u \in \dd'_j$.
Assume e.g. that $H(\vl_1(u)) \leq H(\vl_2(u))/4$ for all $u \in \dd'_j$.
We will show that  $|\psi_\ell(u)| \leq \gamma^{(1)}_\ell (u)$ for all $u \in \dd'_j$.
Given such $u$, let $\varphi$ be the smaller angle between the vectors
$$z_1 = e^{(\fat_{N}+\i b\tau_{N})(\vl_1(u))} (e^\zeta h)(\vl_1(u)) \:\:\: , 
\:\:\: z_2 =  e^{(\fat_{N}+\i b\tau_{N})(\vl_2(u))} (e^\zeta h) (\vl_2(u))$$
in the complex plane $\C$; then $\epsilon_3 \leq \varphi \leq 3\pi/2$. Moreover, 
(8.8) and  (8.9)   imply
$$\frac{|z_1|}{|z_2|} = e^{\fat_N(\vl_1(u)) - \fat_N(\vl_2(u))} \; e^{\zeta(\vl_1(u)) - \zeta(\vl_2(u))} \;
 \frac{|h(\vl_1(u))|}{|h(\vl_2(u))|} \leq(2C_{12}) 16 = 32 C_{12} .$$
As in the proof of Lemma 6.9 this  yields
\begin{equation}
|z_1 + z_2| \leq (1- t) |z_1| + |z_2| ,
\end{equation}
for some constant $t \in (0,1)$ depending on $\ep_3$ and $C_{11}$. Assuming that the constant 
$\mu_0 \leq t$, we have $|\psi_\ell(u)| \leq \gamma^{(1)}_\ell (u)$  
for all $u \in \dd'_j$.  Now set $j_m = j$, $\ell_m = \ell$ and $i_m = 1$,
and include  $(i_m,j_m,\ell_m)$ in the set $J$.  
Then $\dd_{j_m} \subset \cc_m$ and we deduce that  (8.10) holds on $\dd'_{j_m}$.
\endofproof

\bs

Next, we proceed with what is done in Sect. 7. First, we prove parts (a) and (b)  of Lemma 7.1 assuming $f \in \ff_\theta(\hU)$.
Part (a) goes without a change. In part (b) one proves that 
\be
\int_{V_b} (\nn_J H)^2 \, d\nu \leq \rho_3\, \int_{V_b} L^N_{\f0} (H^2) \, d\nu ,
\ee
for any $H\in \kk_E$  and any $J\in \jj(b)$,  where $\di \rho_3 = \rho_3(N) < 1$ is possibly a slightly larger constant, and again
$a_0 = a_0(N) > 0$ is chosen sufficiently small. The proof of this uses the same lines as the ones in the proof of Lemma 7.1(b)
combined with the fact that $\|\zeta_N\|_0 = \|\fa_N - \fat_N\| \leq N \, \Con$ for some global constant $\Con > 0$. 

Then, using the analogue of Lemma 7.1 with $f \in \ff_\theta(\hU)$ and  Lemma 7.2 one proves Lemma 7.3 in the same form -- the 
difference is that now $f \in \ff_\theta(\hU)$ compared to the previous stronger assumption $f \in \ff_{\theta_1}(\hU)$. This gives
the estimate (7.13) in exactly the same form under this more general assumption.  And then one just needs to repeat the argument
from the proof of Theorem 1.3 (the same as under the assumption $f \in \ff_{\theta_1}(\hU)$). 
\endofproof

\bigskip

\noindent
{\it Proof of Theorem} 1.1. This follows from the procedure described in \cite{D1} (see Sect. 4 and Appendix 1  there).
\endofproof

\bs

\noindent
{\it Proof of Corollary 1.4.} Let again $\hat{\theta}$ be as in (5.1). 
Choose the constants $C > 0$, $\rho \in (0,1)$, $a_0 > 0$ and $b_0 \geq 1$, $T > 0$  as in the proof of Theorem 1.3.
Let $\hat{\theta} \leq \theta < 1$. As in the proof of Lemma 5.1,
$(d(x,y))^{\alpha} \leq \Con \dte (x,y)$ will always hold assuming $1/\gamma^{\alpha} \leq \theta$, 
i.e. $\alpha \geq \frac{|\log\theta|}{\log \gamma}$. Here 
$1 < \gamma < \gamma_1$ are the constants from (2.1).
Then for such $\alpha$ we have $|h|_\theta \leq \Con |h|_{\alpha}$. 

Recall that $\alpha_1 = \frac{|\log \hat{\theta}|}{\log \gamma} > 0$ by (5.1), and that
$\alpha_1 \in (0,1]$ is chosen such  that the local stable holonomy maps on $\tR$ are uniformly 
$\alpha_1$-H\"older, i.e. there exists a constant $C_{15} > 0$ such that for any $z,z' \in \tR_i$ 
for some $i = 1, \ldots, k_0$ and any 
$x,y \in W^u_{\tR}(z)$ for the projections $x',y'\in W^u_{\tR}(z')$ of $x, y$ along stable leaves we have
$d(x',y') \leq C_{15} \, (d(x,y))^{\alpha_1}$.

Let $\alpha \in (0,\alpha_1]$; then  $\alpha = \frac{|\log \theta|}{\log \gamma} $
for some $\theta \in [\hat{\theta},1)$. As above this gives $|h|_\theta \leq C'_{16} |h|_\alpha$ for any
$h \in C^\alpha(\hU)$. 

Assume that for a given $h \in C^\alpha(\hU)$ we have $\|h\|_{\alpha,b} \leq 1$; then $\|h\|_0 \leq 1$
and $|h|_\alpha \leq |b|$, so $|h|_\theta \leq C'_{16}|b|$ and therefore $\|h\|_{\theta,b} \leq C'_{16}+1$.
As in (8.1),
$$\|\lab^{n} h\|_{\theta,b} \leq 3A_0  \rho_{6}^n \quad, \quad n \geq 4 k N \hb = 4 k N \lceil \log |b|\rceil  ,$$
so in particular 
\be
\|\lab^{n} h \|_0 \leq 3A_0 \rho_{6}^n
\ee
for all $n \geq  4 k N \hb$.

Next, one needs to repeat part of the arguments from the proof of Theorem 1.3 above. 

First, one needs a version of Lemma 5.2(b) for functions $w\in C^\alpha(\hU)$. Given an integer
$m \geq 0$ and $u,u' \in U_i$ for some $i = 1, \ldots,k_0$, notice that if $\sigma^m(v) = u$,
$\sigma^m(v') = u'$ and $v' = v'(v)$ belongs to the cylinder of length $m$ containing $v$, then
\begin{eqnarray}
|w(\sigma^jv) - w(\sigma^j(v')| 
& \leq & |w|_\alpha (d(\sigma^jv), \sigma^j(v'))^\alpha \\
& \leq & \frac{|w|_\alpha}{c_0^\alpha \gamma^{\alpha(m-j)}}  (d(\tpp^{m-j}(\sigma^jv),\tpp^{m-j}(\sigma^j v'))^\alpha \nonumber\\
& \leq & C''_{16} \frac{|w|_\alpha}{\gamma^{\alpha(m-j)}} (d(u,u'))^{\alpha \alpha_1} .\nonumber
\end{eqnarray}
This implies 
\be
|w_m(v) - w_m(v')| \leq C_{17} |w|_\alpha (d(u,u'))^{\alpha \alpha_1} .
\ee
In particular this holds for $w = f$, $w = \tau$. Now repeating the argument in the proof of Lemma 5.2(b), 
for $|a| \leq a_0$ and $w \in C^\alpha(\hU)$ we get
\begin{eqnarray*}
&        &|(L^m_{f-(P+a)\tau} w)(u) - (L^m_{f- (P+a)\tau} w)(u')|\\
&  =    &   | \sum_{\sigma^m v = u} e^{(f- (P+a) \tau)_m(v)}\, w(v)
 -\sum_{\sigma^m v = u} e^{(f- (P+a)\tau)_m(v'(v))}\, w(v'(v)) | \\
& \leq & \left| \sum_{\sigma^m v = u} e^{(f- (P+a) \tau)_m(v)}\, [w(v) -  w(v')]\right| \\
&       &  +
\sum_{\sigma^m v = u}  \left|e^{(f- (P+a)\tau)_m(v)} -  e^{(f- (P+a)\tau)_m(v')}\right| \, |w(v')| \\
& \leq &  C_{17} \left[ \frac{|w|_\alpha}{\gamma^{\alpha m}} 
   + \|w\|_0 \right]\, (L^m_{f- (P+a)\tau} 1)(u)\, (d(u,u'))^{\alpha \alpha_1} .
\end{eqnarray*}
In particular this shows that $L^m_{f-(P+a)\tau} w \in C^{\alpha \alpha_1}(\hU)$ for all $w\in C^\alpha(\hU)$
and all integers $m \geq 0$. Since $w = 1 \in C^\alpha(\hU)$, it now follows from the Perron-Ruelle-Frobenius
Theorem that the eigenfunction $h_a \in C^{\alpha \alpha_1}(\hU)$ and so $\fa\in C^{\alpha\alpha_1}(\hU)$
for all $|a| \leq a_0$. Moreover, taking $a_0$ sufficiently small, we may assume that 
$\|h_a\|_{\alpha\alpha_1} \leq C'_{18} = \Con$ for all $|a| \leq a_0$. 
Using (8.16) with $w = f_a$ and $\alpha$ replaced by $\alpha\alpha_1$, we get
$$|\fa_m(v) - \fa_m(v')| \leq C_{18}''  (d(u,u'))^{\alpha \alpha_1^2}, $$
and also 
$$|\fa(v) - \fa(v')| \leq C_{18}'' \, \rho_{14}^m (d(u,u'))^{\alpha \alpha_1^2} .$$

Now, using standard arguments, for $h \in C^\alpha(\hU)$ we get
\begin{eqnarray*}
&         & |\lab^m h(u) - \lab^m h(u')|\\
& =      &  \di \left| \sum_{\sigma^m v = u} \left( e^{\fa_m(v) - \i b \tau_m(v)}\, h(v) -
  e^{\fa_m(v') - \i b \tau_m(v')}\, h(v')\right)\right|\\
& \leq & \di \left| \sum_{\sigma^m v = u} e^{\fa_m(v) - \i b \tau_m(v)}\, [h(v) - h(v')]\right|\\
&        & \di  +  \left| \sum_{\sigma^m v = u} \left( e^{\fa_m(v) - \i b \tau_m(v)} -
  e^{\fa_m(v') - \i b \tau_m(v')}\right)\, h(v')\right|\\
& \leq & \di \sum_{\sigma^m v = u} e^{\fa_m(v) }\, | h(v) - h(v')|\\
&        & +   \sum_{\sigma^m v = u} \left| e^{\fa_m(v) - \i b \tau_m(v) -
  \fa_m(v') + \i b \tau_m(v')} - 1\right|\, e^{\fa_m(v') }\, |h(v')| .
\end{eqnarray*}
Using (8.15) with $w = h$ and $j = 0$ gives 
$$|h(v) - h(v')| \leq C_{17} \frac{|h|_\alpha}{\gamma^{\alpha m}} (d(u,u'))^{\alpha \alpha_1} .$$
Moreover, 
$$\left| e^{\fa_m(v) - \i b \tau_m(v) -\fa_m(v') + \i b \tau_m(v')} - 1\right|
\leq C'_{18}  |b|  \,(d(u,u'))^{\alpha \alpha_1^2} .$$
Thus,
\be
|\lab^m h(u) - \lab^m h(u')| \leq C_{18} \left[\rho_{7}^m\, |h|_\alpha + |b|\, \|h\|_0\right]\, (d(u,u'))^{\alpha \alpha_1^2} .
\ee
Since $|h|_\alpha \leq |b|$ and $\|h\|_0 \leq 1$, this gives 
$|\lab^mh|_{\alpha \alpha_1^2} \leq \Con\, |b|$ for all $m \geq 0$. Using (8.14) and (8.17) 
with $h$ replaced by $\lab^m h$ and $\alpha$ replaced by $\alpha \alpha_1^2 \leq \alpha_0$, we get
\begin{eqnarray*}
|(\lab^{2m} h) (u) - (\lab^{2m} h)(u')| 
& \leq & \Con \left[\rho_{7}^m\, |\lab^m h|_{\alpha \alpha_1^2}  + |b|\, \|\lab^m h\|_0\right]\, (d(u,u'))^{\alpha \alpha_1^4} \\
& \leq & C'_{19} \left[\rho_{7}^m\, |b| + |b|\, \rho_{6}^m \right]\, (d(u,u'))^{\alpha \alpha_1^4} 
\end{eqnarray*}
for $m \geq 4 k N \hb$. 
Thus, 
$$\|\lab^{2m} h\|_{\alpha \alpha_1^4,b} \leq C_{19}  \rho_{7}^m$$ 
for all $m \geq  4 k N \hb$ and all
$h \in C^\alpha(\hU)$ with $\|h\|_\alpha \leq 1$.  
Since 
$$L^m_{f- (P_f+a + \i b)\tau} h = \frac{1}{h_a}\, \lab^m (h_a h) ,$$ 
it is now easy to get
$$\|L^{m}_{f-(P+a - \i b)} h\|_{\alpha \alpha_1^4,b} \leq C_{19} \rho_{7}^m\, \|h\|_{\alpha,b}$$
for all $m \geq 4 k N \hb$ and all $h \in C^\alpha(\hU)$.  
Setting $\alpha_0 = \alpha_1$ and $\hbeta = \alpha_1^4$ proves the assertion.
\endofproof

\newpage

\section{Temporal distance estimates on cylinders}
\setcounter{equation}{0}

Here we prove Lemmas 4.2 and 4.3.

Apart from the constants $\hep > 0$, $\beta > 0$, etc. from Sect. 3 (see p. 20) we will also need  the constant $\hep_1 = \frac{\beta \hep}{6}$. 
Moreover, in order to be able to use Lemma 10.7 below, throughout this chapter it will be convenient to assume that 
$R(x)$, $\Gamma(x)$, $D(x)$ and $L(x)$ are Lyapunov $\hep_1$-regularity functions, while $r(x)$ is an
$\hep_1$-slowly varying radius function and $\tr(x)$ is an $\hep_1/2$-slowly varying radius function with the properties in Sect. 3, 
so that they satisfy (3.6) -- (3.18)  and the conclusions of  Lemma 3.1 with $\hep$ replaced by $\hep_1$. Since $\hep_1 < \hep$, 
$R(x)$, $\Gamma(x)$, $D(x)$ and $L(x)$ are still Lyapunov $\hep$-regularity functions, $r(x)$ is still an
$\hep$-slowly varying radius function and $\tr(x)$ is an $\hep/2$-slowly varying radius function.
We will also assume that the set  $P_0 = Q_{p_0}(\hep)$ is as in (4.3), the constants $R_0$, $\Gamma_0$, $L_0$, $D_0$ and  $r_0$ 
satisfy (4.4), and that $r_1$, $p_0$, $\hd_0$ and $\Xi_m = \Xi_m(p_0, \hep, \hd_0)$ are as in Sect. 4.2.

\subsection{A  technical lemma}

Notice that in Lemma 4.1 the exponential maps are used to parametrize $W^u_\ep(z)$ and $W^s_\ep(z)$.
The particular choice of the exponential maps is not  important, however it is important that these
maps are $C^2$. So, we cannot use the maps $\Phi^u_z$ and $\Phi^s_z$ defined in Sect. 3. 
In order to use Lemma 4.1 we will need in certain places to replace the  local lifts $\hf^p_z$ of the 
iterations $\varphi^p$ of the map $\varphi$ by slightly different maps.

For any $x\in \ll$ consider the $C^{2}$ map 
$$\tvarphi_x = (\exp^u_{\varphi(x)})^{-1} \circ \varphi \circ \exp^u_x : E^u (x; r(x))  \longrightarrow E^u (\varphi(x), \tr (\varphi(x)))\;.$$
It is well-defined assuming that the $\hep$-slowly varying radius function $r(x)$ and the $\hep/2$-slowly varying radius function $\tr(x)$
are chosen appropriately as in Sect. 3.
As with the maps $\hf$, for $y \in \ll$ and an integer $j \geq 1$ we will use the notation
$$\tvarphi_y^j = \tvarphi_{\varphi^{j-1}(y)} \circ \ldots \circ \tvarphi_{\varphi(y)} \circ \tvarphi_y\quad,
\quad \tvarphi_y^{-j} = (\tvarphi_{\varphi^{-j}(y)})^{-1} \circ \ldots \circ (\tvarphi_{\varphi^{-2}(y)})^{-1}  \circ (\tvarphi_{\varphi^{-1}(y)})^{-1} \;,$$
at any point where these sequences of maps are well-defined.
In a similar way one defines the maps $\tvarphi_x$ and their iterations on $E^s(x;r(x))$.
\index{$\tvarphi_x$}
\index{$\tvarphi_y^j$}

Following the notation in Sect. 3 and using the fact that the flow $\phi_t$ is contact,
the negative Lyapunov exponents over $\ll$ are 
$$-\log \lambda_1 = - \chi_1 > -\log \lambda_{2} = - \chi_2 > \ldots > - \log \lambda_{\tk} = \chi_{\tk} .$$
Fix $\hep > 0$ as in Sect. 3, assuming in addition that 
$$\hep \leq \frac{\log \lambda_1}{100}\, \min\{ \beta, \vartheta\}  \quad , \quad \hep \leq \frac{\log \lambda_1 \, (\log \lambda_2 - 
\log \lambda_1)}{4 \log \lambda_1  + 2\log \lambda_2} .$$
Set 
$$\hep_1 = \frac{\beta\, \hep}{6} ,$$
where $\beta > 0$ is as in Sect. 3 ( see p. 20 above).

For $x\in \ll$ we have an $\varphi$-invariant decomposition
$$E^s(x) = E^s_1(x) \oplus E^s_2(x) \oplus \ldots \oplus E^s_{\tk}(x)$$
into subspaces of dimensions $n_1, \ldots, n_{\tk}$, where $E^s_i(x)$ $(x\in \ll$) is the
$d\varphi$-invariant subbundle corresponding to the Lyapunov exponent $- \log \lambda_i = - \chi_i$.
As we stated above, we assume that the {\it Lyapunov $\hep_1$-regularity function}  $R = R_{\hep_1} : \ll \longrightarrow (1,\infty)$ is
chosen as in  Sect. 3. This is also a Lyapunov $\hep$-regularity function, so it follows from (3.4) that
\be
\frac{1}{R (x)\, e^{m\hep}} \leq \frac{\|d\varphi^m(x)\cdot v\|}{\lambda_i^{-m}\|v\|} 
\leq R (x)\, e^{m\hep} \quad , \quad x\in \ll  \;, \; v\in E^s_i(x)\setminus \{0\} \;, \; m \geq 0.
\ee

For the contact form $\omega$ it is known (see e.g. \cite{KH} or Appendix B in
\cite{L1}) that $\omega$ vanishes on every stable/unstable manifold of a point on $M$,
while $d\omega$ vanishes on every weak stable/unstable manifold.
For Lyapunov regular points we get a bit of extra information.

\bs

\noindent
{\bf Lemma 9.1.} {\it For every $x\in \ll$ and every
$u = (\uo, \ldots, u^{(\tk)}) \in E^u(x;r(x))$ and $v = (\vo, \ldots,\vk)\in E^s(x;r(x))$
we have}
\be
d\omega_x(u,v) = \sum_{i=1}^{\tk} d\omega_x(\ui,\vi) .
\ee

\noindent
{\it Proof.} It is enough to show that $d\omega_x(\ui, \vj) = 0$ if $i \neq j$. Let e.g. $i < j$.  
Using (3.4), (4.1), (9.1) and  the fact that $d\omega$ is $d\varphi$-invariant, for $m \geq 0$ and $x_m = \varphi^m(x)$ we get
\begin{eqnarray*}
 |d\omega_x(\ui,\uj)|  
& =    & |d\omega_{x_m}(d \varphi^m(x)\cdot \ui, d \varphi^m(x)\cdot \vj)|\\
& \leq & C \|d \varphi^m(x)\cdot \ui\| \,\|d \varphi^m(x)\cdot \vj)\|\\
& \leq & C R^2(x) \|\ui\|\, \|\vj\|\, \frac{(\lambda_i e^{2\hep})^m}{\lambda_j^m} .
\end{eqnarray*}
Since $\lambda_i e^{2\hep} < \lambda_j$, the latter converges to $0$ as $m \to \infty$,  so $d\omega_x(\ui,\vj) = 0$.

The case $i > j$ is considered similarly by taking $m \to -\infty$.
\endofproof

\def\tyj{\tilde{y}^{(j)}}
\def\cvj{\check{v}_j}
\def\heta{\hat{\eta}}

\def\etao{\eta^{(1)}}
\def\etai{\eta^{(i)}}
\def\zetao{\zeta^{(1)}}
\def\hbeta{\hat{\beta}}
\def\tdelta{\tilde{\delta}}

\def\tf{\tilde{\varphi}}

\subsection{Proof of Lemma 4.2(a)} 

We will consider cylinders $\cc$ of length $m \geq 1$  in $\check{R}$ with 
$\cc \cap P_0 \setminus \Xi_m \neq \e$ (instead of considering cylinders $\cc$ in $R$) 
with corresponding obvious changes in the estimates we need to prove.

Let $\cc$ be a cylinder of length $m$ in  $\check{R}$. 
Fix an arbitrary $z_0 \in \cc \cap P_0 \setminus \Xi_m$. Given $x_0\in \cc$, write 
$$x_0 = \Phi^u_{z_0}(\xi_0) = \exp^u_{z_0}(\txi_0)$$
for some $\xi_0, \txi_0 \in E^u(z_0)$ with $\txi_0 = \Psi^u_{z_0}(\xi_0)$.  
Then $ \|\xi_0\| \leq R_0 \, \diam(\cc)$. 
Set $\cc' = \tPsi \circ \Psi^{-1}(\cc) \subset \tR$, $T = \ttau_m(z_0)$ and $p = [T]$, so that $p \leq T < p+1$.
Moreover, $m\tau_0 \leq \ttau_m(z_0) \leq m \htau_0$ gives $\frac{p}{\htau_0} \leq m \leq \frac{p+1}{\tau_0} \leq \frac{2p}{\tau_0}$.

Since $m$ is the length of  $\cc'$, $\tpp^m(\cc')$ contains a whole unstable leaf
of a proper rectangle $\tR_j$.  Moreover, $z_0\in \cc \setminus \Xi_m$ shows that there exists an integer
$m'$ with $m (1-\hd_0)  \leq m' \leq m $ such that  $z = \tpp^{m'}(z_0) \in P_0$. 
Let $z \in \tR_i$. By the choice of the constant $r_1 > 0$ (see Sect. 4.2), there exists $y \in R_i$ such that
$B^u(y,r_1) \subset W^u_{\tR_i}(z)$ and $d(z,y) < r_0/2$. In particular, for every point $b'\in B^u(y,r_1)$
there exists $b \in \cc$ with $\tpp^{m'}(b) = b'$. Set $p' = [\ttau_{m'}(z_0)]$; then $\varphi^{p'}(z_0) \in \hP_0$.

To estimate $p-p'$, notice that, as above,
$\frac{p'}{\htau_0} \leq m' \leq \frac{2p'}{\tau_0}$, so the relationship between $m'$ and $m$ implies
$$p - p' \leq \ttau_m(z_0) - \ttau_{m'}(z_0) + 1 \leq (m-m') \htau_0 +1 \leq \hd_0 m \htau_0 +1 \leq \frac{3\hd_0 \htau_0}{\tau_0}\, p ,$$
assuming that $m$ (and so $p$) is sufficiently large. Thus, 
$$p (1-3\hd_0 \htau_0/\tau_0)  \leq p' \leq p ,$$ 
where $3\hd_0 \htau_0/\tau_0 < 1,$
 and $z_{p'} = \varphi^{p'}(z_0) \in  \hP_0$, so $r(z_{p'}) \geq r_0$ by (4.3). Clearly, $p' \geq \ttau_{m'}(z_0)$.
Then for every $b\in W^u_{r_1}(z_{p'})$ there exists $b \in \cc$ with $\varphi^{p'}(b) = b'$. Consider an arbitrary
$\zeta_{p'} \in E^u(z_{p'}; r_1/R_0)$ such that $\|\zetao_{p'}\| \geq r_1/R_0$, 
and set $\zeta = \hf^{-p'}_{z_{p'}}(\zeta_{p'})$.  Then $x = \Phi^u_z(\zeta) \in \cc$, so 
$$\diam(\cc) \geq d(z_0,x) \geq \frac{\|\zeta\|}{R_0} \geq \frac{\|\zetao\|}{\Gamma_0 R_0} .$$
On the other hand, Lemma 3.5 in \cite{St4} (see Lemma 10.1 below)  and $\mu_1^{p'} \leq \mu_1^p < \lambda_1^p$ give
$$\|\zetao\|
 \geq  \frac{1}{\Gamma_0} \|\zetao\|'_{z_0} \geq \frac{\|\zetao_{p'}\|}{\Gamma_0 \mu_1^{p'}}
 \geq  \frac{r_1/R_0}{\Gamma_0 \mu_1^{p}}  > \frac{r_1}{R_0\Gamma_0 \lambda_1^{p} } .$$
Hence $\diam(\cc) \geq  \frac{c_4}{\lambda_1^{p} }$, where $c_4 =  \frac{r_1}{R^2_0\Gamma^2_0} \geq 1$. 

This proves the left-hand-side inequality in (\ref{eq:cylbound}) with $C_1 = 1/c_4$. 
The other inequality in (\ref{eq:cylbound}) follows by a similar (in fact, easier) argument, using the above estimate of
$p-p'$. We omit the details.

\def\hb{\hat{b}}

\subsection{Proof of Lemma 4.2(b)}

Let $\cc$ be a cylinder  of length $m$ in $R$ such that there exists $\hz \in \cc \cap P_0 \setminus \Xi_m$.
Set $\tcc = \tPsi(\cc)$ as before.

Let  $\hx_0,\hz_0 \in \cc$, $\hy_0, \hb_0 \in W^s_{R_1}(\hz_0)$. 
We can assume that $\cc$ is the {\bf smallest cylinder
containing} $\hx_0$ and $\hz_0$; otherwise we will replace $\cc$ by a smaller cylinder.

It is enough to consider the case when $\hz_0 = \hz$. Indeed, assuming 
the statement is true with $\hz_0$ replaced by $\hz$, consider arbitrary points $\hx_0, \hz_0\in \cc$.
Set 
$$\{ y \} = W^u_R(\hy_0) \cap W^s_R(\hz)\quad \mbox{\rm and} \quad\{ b \} = W^u_R(\hb_0) \cap W^s_R(\hz) .$$
Since the local unstable holonomy maps are
uniformly H\"older, there exist (global) constants $C' > 0$ and $\beta' > 0$ 
such that $d(y,b) \leq C' (d(\hy_0,\hb_0))^{\beta'}$. Thus, using the assumption, 
$$|\Delta(\hx_0,y) - \Delta(\hx_0,b)| \leq C_1 \diam(\cc) (d(y,b))^{\beta_1} 
\leq C_1 (C')^{\beta_1} \diam(\cc) (d(\hy_0,\hb_0))^{\beta'\beta_1} .$$
A similar estimate holds for $|\Delta(\hz_0,y) - \Delta(\hz_0,b)|$, so
\begin{eqnarray*}
|\Delta(\hx_0,\hy_0) - \Delta(\hx_0,\hb_0)| 
& =    &  |(\Delta(\hx_0,y) - \Delta(\hz_0,y) - (\Delta(\hx_0,b) - \Delta(\hz_0,b))|\\
& \leq & |\Delta(\hx_0,y)  - \Delta(\hx_0,b)| + |\Delta(\hz_0,y) - \Delta(\hz_0,b)|\\
& \leq & 2 C_1 (C')^{\beta_1} \diam(\cc) (d(\hy_0,\hb_0))^{\beta'\beta_1} .
\end{eqnarray*}

So, from now on we will assume that $\hz_0 = \hz \in \cc \cap P_0 \setminus \Xi_m$.
Then $R(\hz_0) \leq R_0$, $r(\hz_0) \geq r_0$, etc.
Set 
$$x_0 = \Psi(\hx_0) \quad , \quad z_0 = \Psi(\hz_0) \quad , \quad y_0 = \Psi(\hy_0) \in \tR \quad , \quad b_0 = \Psi(\hb_0) ,$$ 
and then write 
$$x_0 = \Phi^u_{z_0}(\xi_0) = \exp^u_{z_0}(\txi_0)$$
for some $\xi_0, \txi_0 \in E^u(z_0)$ with $\txi_0 = \Psi^u_{z_0} (\xi_0)$.  
Then 
$$ \|\xi_0\|, \|\txi_0\| \leq R_0 \; \diam(\tcc) .$$
Similarly, write 
$$y_0 = \exp^s_{\tz_0}(\tv_0) = \Phi^s_{z_0}(v_0) \quad \mbox{\rm and} \quad b_0 = \exp^s_{z_0}(\teta_0) = \Phi^s_{z_0}(\eta_0)$$
for some 
$v_0, \tv_0, \eta_0, \teta_0 \in E^s(z_0)$ with $\tv_0 =  \Psi^s_{z_0} (v_0)$ and
$\teta_0 =  \Psi^s_{z_0} (\eta_0)$. It follows from  (3.6) that
\be
\|\tv_0 - v_0\| \leq R_0 \|v_0\|^{1+\beta} \:\: , \:\: \|\txi_0 - \xi_0\|\leq R_0\|\xi_0\|^{1+\beta} \:\:, \:\:
 \|\teta_0 - \eta_0\| \leq R_0\| \eta_0\|^{1+\beta} .
\ee

\subsubsection{Pushing forward}

Set $p = [\ttau_m(z_0)]$; then (\ref{eq:cylbound}) holds.
Set $q = [p/2]$. We will in fact assume that $q = p/2$; the difference with the case
when $p$ is odd is insignificant. For any integer $j \geq 0$ set
$$z_j = \varphi^j(z_0) \quad , \quad x_j = \varphi^j(x_0) \quad , \quad y_j = \varphi^j(y_0) ,$$ 
$$\hxi_j =  d\hf^j_{z_0}(0)\cdot \xi_0 \quad , \quad  \xi_j = \hf^j_{z_0}(\xi_0) \quad , \quad \txi_j = \tf^j_{z_0}(\txi_0) ,$$
$$\hv_j =  d\hf^j_{z_0}(0)\cdot v_0 \quad  , \quad v_j = \hf^j_{z_0}(v_0) \quad , \quad \tv_j = \tf^j_{z_0}(\tv_0) ,$$
$$b_j = \varphi^j (b_0) \quad , \quad  \heta_j =  d\hf^j_{z_0}(0)\cdot \eta_0 \quad  , \quad 
\eta_j = \hf^j_{z_0}(\eta_0) \quad  , \quad \teta_j = \tf^j_{z_0}(\teta_0) .$$

Notice that $\txi_0 = \Psi^u_{z_0}(\xi_0)$,  $\tv_0 = \Psi^s_{z_0}(v_0)$, and also
$$\txi_j = \Psi^u_{z_j} (\xi_j) \:\: , \:\: \Phi^u_{z_j}(\xi_j) = x_j \quad, \quad \tv_j = \Psi^s_{z_j} (v_j) 
\quad, \quad \teta_j = \Psi^s_{z_j} (\eta_j) ,$$
so it follows from (3.6) that
\be
\|\xi_j - \txi_j \| \leq R(z_j) \|\xi_j\|^{1+\beta} \: , \: \|v_j - \tv_j\| \leq R(z_j) \|v_j\|^{1+\beta} 
\: , \: \|\eta_j - \teta_j\| \leq R(z_j) \|\eta_j\|^{1+\beta} .
\ee
Moreover, 
$$\exp^u_{z_j}(\txi_j) = \varphi^j(\exp^u_{z_0}(\xi_0)) = \varphi^j(x_0) = x_j ,$$ 
$\exp^s_{z_j}(\tv_j) = y_j$ and $\exp^s_{z_j}(\teta_j) = b_j$, so Lemma 4.1 implies
\be
|\Delta(x_j,y_j)  - d\omega_{z_j}(\txi_j,\tv_j)| 
\leq  C_0 \,\left[ \|\txi_j\|^2\, \|\tv_j\|^\vartheta + \|\txi_j\|^\vartheta \|\tv_j\|^2 \right] .
\ee
Similarly
$$|\Delta(x_j,b_j)  - d\omega_{z_j}(\txi_j,\teta_j)| 
\leq  C_0 \,\left[ \|\txi_j\|^2\, \|\teta_j\|^\vartheta + \|\txi_j\|^\vartheta \|\teta_j\|^2 \right] $$
for every integer $j \geq 0$. From (9.4) one gets
$$|d\omega_{z_j} (\txi_j, \tv_j) - d\omega_{z_j}(\xi_j,v_j)| 
\leq 2 C_0 R(z_j)\, \|\xi_j\|\, \|v_j\| (\|\xi_j\|^\beta + \|v_j\|^\beta) ,$$
$$|d\omega_{z_j} (\txi_j, \teta_j) - d\omega_{z_j}(\xi_j,\eta_j)| 
\leq 2 C_0 R(z_j)\, \|\xi_j\|\, \|\eta_j\| (\|\xi_j\|^\beta + \|\eta_j\|^\beta) ,$$
and also $\|\txi_j\|\leq 2\|\xi_j\|$, $\|\tv_j\| \leq 2 \| v_j\|$ and $\|\teta_j\| \leq 2 \| \eta_j\|$.
Indeed, from (9.4), 
$$\|\txi_j\| \leq \|\xi_j\| (1+ R(z_j) \|\xi_j\|^\beta) \leq \|\xi_j\| (1 + R_0 e^{j\hep} r_0/\mu_1^{j} )
\leq \|\xi_j\| (1+  R_0 r_0 ) \leq 2 \|\xi_j\| ,$$
since $r_0 \leq 1/R_0$. Similarly, $\|\tv_j\| \leq 2 \| v_j\|$ and $\|\teta_j\| \leq 2 \| \eta_j\|$.

Using these, it follows from (9.5) that
\begin{eqnarray}
|\Delta(x_j,y_j)  - d\omega_{z_j}(\xi_j,v_j)| 
& \leq &  2 C_0 R(z_j)\, \|\xi_j\|\, \|v_j\| (\|\xi_j\|^\beta + \|v_j\|^\beta) \\ 
&       & + 8 C_0\,\left[ \|\xi_j\|^2\, \|v_j\|^\vartheta + \|\xi_j\|^\vartheta \|v_j\|^2 \right] .\nonumber 
\end{eqnarray}
and similarly
\begin{eqnarray}
|\Delta(x_j,b_j)  - d\omega_{z_j}(\xi_j,\eta_j)| 
& \leq &  2 C_0 R(z_j)\, \|\xi_j\|\, \|\eta_j\| (\|\xi_j\|^\beta + \|\eta_j\|^\beta) \\ 
&       & + 8 C_0\,\left[ \|\xi_j\|^2\, \|\eta_j\|^\vartheta + \|\xi_j\|^\vartheta \|\eta_j\|^2 \right] . \nonumber
\end{eqnarray}
for every integer $j \geq 0$.

Next, we will be estimating 
$|\Delta(x_0,y_0) - d\omega_{z_0}(\xi_0,v_0)|  .$
Since $\Delta$ is $\varphi$-invariant and $d\omega$ is $d\varphi$-invariant  we have
$$\Delta(x_0,y_0) = \Delta (x_j,y_j)  \quad , \quad  d\omega_{z_0}(\xi_0,v_0) = d\omega_{z_j}(\hxi_j , \hv_j) ,$$
and also $\Delta(x_0,b_0) = \Delta_{z_j}(x_j,b_j)$ and   
$d\omega_{z_0}(\xi_0,\eta_0) = d\omega_{z_j}(\hxi_j , \heta_j)$
for all $j$. (Notice that $d\hf_{x}(0) = d\varphi(x)$ for all $x\in M$.)

It follows from $z_0\notin \Xi_m$ and (3.23) that 
$z_0 \notin \hXi_{m\htau_0} (p_0+\hep, \hep, \hd_0/\tau_0) ,$
so there exists an integer $\ell$ with 
$$q - (m\htau_0) \frac{\hd_0}{\tau_0} \leq \ell \leq q \quad \mbox{\rm with} \quad \varphi^\ell(z_0) \in \hP_0 .$$ 
{\bf Fix an $\ell$ with this property}.
As in Sect. 9.2, we have $\frac{p}{\htau_0} \leq m \leq \frac{2p}{\tau_0}$.
Combined with $q = p/2$, this gives 
$$ (m\htau_0) \frac{\hd_0}{\tau_0} \leq  (4q \htau_0) \frac{\hd_0}{\tau^2_0} .$$
Thus, we have
\be
(1 - \tdelta_0) q \leq \ell \leq q \quad , \quad z_\ell = \varphi^\ell(z_0) \in P_0  ,
\ee
where
$$\tdelta_0 = \frac{4\htau_0}{\tau_0^2}\, \hd_0 < 1 .$$
It then follows from Lemma 3.1, the choice of $L_0$ and $\|\xi_\ell\| \leq r(z_\ell)$ (since $\ell \leq q= p/2$; see also Sect. 9.3.2 below) that
\be
\|\hxio_\ell - \xio_\ell\| \leq L_0 \|\xi_\ell\|^{1+\beta} .
\ee
Apart from that, using Lemma 10.7(b) below, backwards for stable manifolds, with 
$$a = d\hf^{-\ell}_{z_\ell}(0)\cdot v_\ell \in E^s(z_0) \quad \mbox{\rm and} \quad b = d\hf^{-\ell}_{z_\ell}(0)\cdot \eta_\ell \in E^s(z_0) ,$$ 
since
$v_0 = \hf^{-\ell}_{z_\ell}(v_\ell)$ and $\eta_0 = \hf^{-\ell}_{z_\ell}(\eta_\ell)$, it follows that
$$\|(\ao - \bo) - (\vo_0 - \etao_0)\| \leq L_0 \left[ \|v_0 - \eta_0\|^{1+\beta}  + \|\eta_0\|^\beta \|v_0-\eta_0\|\right]
\leq 2L_0 \| v_0 - \eta_0\| .$$
Thus,
\be
\| d\hf^{-\ell}_{z_\ell}(0) \cdot (\vo_\ell - \etao_\ell) - (\vo_0 - \etao_0)\| \leq 2L_0 \|v_0 - \eta_0\| .
\ee

In what follows we denote by $\Con$ a global constant (depending on constants like 
$C_0$, $L_0$, $R_0$ but independent of
the choice of the cylinder $\cc$, the points $x_0,z_0,y_0,b_0$, etc.) which may change from line to line.

Using  (9.9), (9.10) and the above remarks,  we obtain
\begin{eqnarray}
&        & |d\omega_{z_\ell} (\xi_\ell, v_\ell - \eta_\ell)| \\
& \leq & |d\omega_{z_\ell} (\xio_\ell, \vo_\ell - \etao_\ell)| 
        + C_0 \sum_{i=2}^{\tk} \|\xi_\ell^{(i)}\| \, (\|v_\ell^{(i)}\| + \| \eta_\ell^{(i)}\|)\nonumber\\
& \leq &  |d\omega_{z_\ell} (\hxio_\ell, \vo_\ell - \etao_\ell)| + \Con \|\xi_\ell\|^{1+\beta}  \|\vo_\ell - \etao_\ell\| \nonumber\\
&        & +  C_0 \sum_{i=2}^{\tk} \|\xi_\ell^{(i)}\| \, (\|v_\ell^{(i)}\|  + \| \eta_\ell^{(i)}\|)\nonumber\\
& \leq & |d\omega_{z_\ell} (d \hf^{\ell}_{z_0}(0) \cdot \xio_0, \vo_\ell - \etao_\ell)| 
+  \Con \|\xi_\ell\|^{1+\beta} \|\vo_\ell - \etao_\ell\|\nonumber\\
&       &  +  C_0 \sum_{i=2}^{\tk} \|\xi_\ell^{(i)}\| \,  (\|v_\ell^{(i)}\| + \| \eta_\ell^{(i)}\|)\nonumber\\
& =     & |d\omega_{z_0} (\xio_0, d\hf^{-\ell}_{z_\ell}(0)\cdot (\vo_\ell - \etao_\ell))|
+ \Con \|\xi_\ell\|^{1+\beta} \|\vo_\ell - \etao_\ell\| \nonumber\\
&        & +  C_0 \sum_{i=2}^{\tk} \|\xi_\ell^{(i)}\| \, (\|v_\ell^{(i)}\| + \| \eta_\ell^{(i)}\|)\nonumber\\
& \leq &  |d\omega_{z_0} (\xio_0, \vo_0 - \etao_0)| + 2C_0 L_0 \|\xi_0\|\, \|v_0 - \eta_0\| 
+ \Con \|\xi_\ell\|^{1+\beta} \|\vo_\ell - \etao_\ell\| \nonumber\\
&      & +  C_0 \sum_{i=2}^{\tk} \|\xi_\ell^{(i)}\| \, (\|v_\ell^{(i)}\| + \| \eta_\ell^{(i)}\|)\nonumber\\
& \leq & \Con \, \diam(\cc) \, \|v_0 - \eta_0\|  + \Con \|\xi_\ell\|^{1+\beta} \|\vo_\ell - \etao_\ell\| \nonumber\\
&        & +  C_0 \sum_{i=2}^{\tk} \|\xi_\ell^{(i)}\| \, (\|v_\ell^{(i)}\| + \| \eta_\ell^{(i)}\|) .\nonumber
\end{eqnarray}

\subsubsection{Estimates for $\|\xi_\ell\|$, $\|v_\ell\|$ and $\| \eta_\ell\|$}

We will now use the choice of $\ell$ to estimate $\|\xi_\ell\|$, $\|v_\ell\|$ and $\| \eta_\ell\|$ by means of $\|\xi_0\|$,
$\|v_0\|$ and $\|\eta_0\|$. We will first estimate $\|\xi_q\|$, $\|v_q\|$ and $\|\eta_q\|$.

Using the definition of $\xi_j$, $p = 2q$, $z_0\in P_0$  and (3.12) we get
$$\|\xi_{q}\| \leq \|\xi_q\|'_{z_q} \leq \frac{\|\xi_p\|'_{z_p}}{\mu_1^{p-q}}
\leq  \frac{\Gamma(z_p) e^{q \hep} \|\xi_p\|}{\lambda_1^{q}}
\leq  \frac{\Gamma_0 e^{3q \hep} \|\xi_p\|}{\lambda_1^{q}}\;.$$
Since $\Phi^u_{z_p}(\xi_p) = x_p$ and  $d(x_p,z_p) \leq \diam(\tR_i)$, we get
$$\|\xi_p\| \leq R(z_p) d(x_p,z_p) \leq R_0 e^{p \hep} r_1 < R_0 e^{p \hep} .$$
Thus,
\be
\|\xi_{q}\| \leq  \frac{R_0 \Gamma_0 e^{5q \hep}}{\lambda_1^{q}}\;.
\ee
Using (3.12) again (on stable manifolds) and $\|v_0\| \leq 2\delta'/R_0 < 1$, we get
\be
\|v_{q}\| = \|v_q\|'_{z_q} \leq  \frac{\|v_0\|'_z}{\mu_1^q}  \leq \frac{\Gamma_0 e^{q\hep} \|v_0\|}{\lambda_1^q} 
\leq \frac{\Gamma_0 e^{q\hep}}{\lambda_1^q} .
\ee
Similarly, $\|\eta_q\| \leq  \frac{\Gamma_0 e^{q\hep}}{\lambda_1^q}$.

Next, it follows from (\ref{eq:cylbound}) that $(\lambda_1)^{2q} \geq \frac{1}{C_1\, \diam(\cc)}$, so
\be
q \geq \frac{1}{2\log \lambda_1 }\, \log \frac{1}{C_1\, \diam(\cc)} .
\ee
This and (9.12) give
\begin{eqnarray}
\|\xi_{q}\|
& \leq & R_0 \Gamma_0 (\lambda_1 e^{-5\hep})^{-q} = R_0 \Gamma_0 e^{-q \log (\lambda_1 e^{-5\hep})} \\
& \leq & R_0 \Gamma_0 \, e^{-\frac{\log (\lambda_1 e^{-5\hep})}{2\log \lambda_1}
\log \left(\frac{1}{C_1\, \diamf(\cc)}\right) } \nonumber \\
& =    &  R_0 \Gamma_0 \,\left(\frac{ 1}{C_1\, \diam(\cc)}\right)^{- \frac{\log\lambda_1 - 5 \hep}{2\log \lambda_1 }}
\leq R_0 \Gamma_0 C_1\, \left(\diam(\cc)\right)^{\frac{\log\lambda_1 - 5 \hep}{2\log \lambda_1 }} ,\nonumber
\end{eqnarray}
since $\frac{\log\lambda_1 - 5 \hep}{\log \lambda_1 } < 1$.
Similarly, (9.13) yields
\begin{eqnarray*}
\|v_{q}\|
& \leq & \Gamma_0 (\lambda_1 e^{-\hep})^{-q}
 \leq  \Gamma_0 \, e^{-\frac{\log (\lambda_1 e^{-\hep})}{2\log \lambda_1 }  \log \left(\frac{1}{C_1\, \diamf( \cc)}\right) } 
\leq \Gamma_0 C_1 \, \left(\diam(\cc)\right)^{\frac{\log\lambda_1 -  \hep}{2\log \lambda_1}}\;.
\end{eqnarray*}
The same estimate holds for $\|\eta_q\|$.

We need similar estimates, however with $q$ replaced by $\ell$. 
Since $q-\ell \leq \tdelta_0 q$ by (9.8), as in (9.15) one obtains
\begin{eqnarray*}
\|\xi_\ell\| 
& \leq &  \|\xi_\ell\|'_{z_\ell} \leq \|\xi_q\|'_{z_q} \leq \Gamma_0 e^{q \hep} \|\xi_q\|
 \leq  R_0 \Gamma^2_0 C_1 \, \left(\diam(\cc)\right)^{\frac{\log\lambda_1  - 6 \hep}{2\log \lambda_1}} .
\end{eqnarray*}
Since $\lambda_{\tk}^{\tdelta_0} < e^{\hep}$ by the choice of $\hd_0$ in Sect. 6.1, we have
$\lambda_{\tk}^{q-\ell} \leq \lambda_{\tk}^{\tdelta_0 q} < e^{\hep q}$, and therefore
\begin{eqnarray*}
\|v_\ell\| 
 \leq  \Gamma(z_\ell) e^{(q-\ell)\hep} \lambda_{\tk}^{q-\ell}\|v_{q}\|
\leq \Gamma_0 e^{3q \hep} \|v_{q}\| \leq \Gamma^2_0 (\lambda_1 e^{-4\hep})^{-q}
\leq  \Gamma^2_0 C_1 \, \left(\diam(\cc)\right)^{\frac{\log\lambda_1 -  4\hep}{2\log \lambda_1 }} ,
\end{eqnarray*}
and again the same estimate holds for $\|\eta_\ell\|$.
Thus, taking the constant $C'' > 0$ so large that  $C'' \geq R_0 \Gamma_0^2 C_1$, we get
$$\|v_{\ell}\|, \|\eta_\ell\|  , \| \xi_\ell\|
\leq C'' \left(\diam(\cc)\right)^{\frac{\log\lambda_1 -  6\hep}{2\log \lambda_1}} .$$

Using these we get the following estimates for the terms in (9.11):
\begin{eqnarray*}
\|\xi_\ell\|\, \|v_\ell\| (\|\xi_\ell\|^\beta + \|v_\ell\|^\beta)
 \leq  2  (C'')^3 \left(\diam(\cc)\right)^{(2+\beta)\frac{\log\lambda_1 
 -  6\hep}{2\log \lambda_1}}  \leq  2 (C'')^3 (\diam(\cc))^{1+\hbeta} ,
\end{eqnarray*}
where we choose 
\be
0 < \hbeta  = \min\left\{  \frac{1}{4} \, \min \{\beta , \vartheta\} \; ,\:
\frac{\log\lambda_2 - \log \lambda_1}{2 \log \lambda_1}  \right\} ,
\ee
and we use the assumption $\hep  \leq \frac{\log \lambda_1}{100} \min\{ \beta , \vartheta \}$.
Then 
$$(2+\beta)\frac{\log\lambda_1 -  6\hep}{2\log \lambda_1} \geq 1+ \hbeta$$
and also
$$(2+\vartheta)\frac{\log\lambda_1 -  6\hep}{2\log \lambda_1 } \geq 1+ \hbeta , $$
which is used in the next estimate. Similarly,
$$\|\xi_\ell\|^{1+\beta} \|v_\ell\| \leq (C'')^3 (\diam(\cc))^{1+\hbeta} ,$$
and 
\begin{eqnarray}
\|\xi_\ell\|^2\, \|v_\ell\|^\vartheta +  \|\xi_\ell\|^\vartheta \|v_\ell\|^2  
\leq 2 (C'')^3 (\diam(\cc))^{1+\hbeta} .
\end{eqnarray}

Next, for any $\xi = \xio + \xit + \ldots + \xi^{(\tk)} \in E^u(z)$ or $E^s(z)$ for some $z\in M$ set
$\cxit = \xit + \ldots + \xi^{(\tk)}$, so that $\xi = \xio + \cxit$.
Using Lemma 3.5 in \cite{St4} (see Lemma 10.1 below), $p-\ell = 2q - \ell \geq q$
and the fact that 
$$\|\xi_{\ell}\| \leq \|\xi_p\| \leq R_0 r_1 \leq R_0 ,$$ 
we get
$$\| \cxit_{\ell}\|'_{z_\ell}  \leq \frac{\Gamma_0 \|\cxit_{\ell}\|}{\mu_2^{q}} 
\leq  \frac{\Gamma_0 \|\xi_{\ell}\|}{\mu_2^{q}} \leq  \frac{\Gamma_0 R_0 }{\mu_2^{q}} . $$
Similarly, using Lemma 3.5 in \cite{St4} (backwards for the map $f^{-1}$ on stable manifolds),
$z_0 \in P_0$, $v_0 = v_{j,1}(z_0) \in E^s(z_0,r'_0)$ and the fact that 
$\|v_0\| \leq \delta' < 1$, we get
$$\| \cvt_{\ell}\|'_{z_\ell}  \leq \frac{\Gamma_0 \|v_0\|}{\mu_2^{q (1-\tdelta_0)}} 
\leq  \frac{\Gamma_0 }{\mu_2^{q(1-\tdelta_0)}} . $$
Hence for $i \geq 2$ we have
$$\|\xii_\ell\| \leq |\cxit_\ell| \leq \|\cxit_\ell\| \leq \frac{\Gamma_0 R_0 }{\mu_2^{q}} , $$
and similarly 
$$\di \|\vi_\ell\| \leq \frac{\Gamma_0}{\mu_2^{q(1-\tdelta_0)}} .$$
From these estimates, (9.14), $\mu_2 = \lambda_2 e^{-\hep}$, and the assumptions about $\hep$, we get
\begin{eqnarray*}
\|\xii_\ell\|\, \| \vi_\ell\|
& \leq & \Gamma_0^2 R_0\,  (\lambda_2 e^{-2\hep})^{-  2q} 
= \Gamma_0^2 R_0\, e^{-2q \log(\lambda_2e^{-2\hep})}\\
& \leq &  \Gamma_0^2 R_0\, e^{\frac{- \log(\lambda_2e^{-2\hep})}{\log \lambda_1 }\, \log\left(\frac{1}{C_1\, \diamf(\cc)} \right)}\\
& \leq & \Gamma_0^2 R_0\,\left( C_1\, \diam(\cc)\right)^{\frac{\log \lambda_2 - 2\hep}{ \log\lambda_1} } 
\leq C'' (\diam(\cc))^{1+\hbeta} .
\end{eqnarray*}
Here we used 
$$\hbeta \leq \frac{\log\lambda_2 - \log \lambda_1}{2 \log \lambda_1}$$ 
which follows from  (9.16), and assumed
$$C'' \geq \Gamma_0^2 R_0 (C_1)^{\log\lambda_2/\log\lambda_1} .$$
Then
\begin{eqnarray*}
\frac{\log \lambda_2 - 2\hep}{ \log\lambda_1}  - 1
 =       \frac{\log \lambda_2  -  2\hep  - \log \lambda_1 }{ \log\lambda_1}   \geq \hbeta . 
\end{eqnarray*}

\subsubsection{Final estimate}

Using  (9.11) and the above estimates for $\|\xi_\ell\|$, $\|v_\ell\|$, $\|\xii_\ell\|\, \| \vi_\ell\|$,  we obtain
\begin{eqnarray*}
|d\omega_{z_\ell} (\xi_\ell, v_\ell - \eta_\ell)|
& \leq & \Con \, \diam(\cc) \, \|v_0 - \eta_0\| + \Con \, (\diam(\cc))^{1+ \hbeta} .
\end{eqnarray*}

It now follows from  (9.6) and (9.7) with $j = \ell$  and the previous estimate that
\begin{eqnarray}
&         & |\Delta(x_0, y_0) - \Delta (x_0, b_0)|
 =      |\Delta(x_\ell, y_\ell) - \Delta (x_\ell, b_\ell)|\\
& \leq & |d\omega_{z_\ell} (\xi_\ell, v_\ell - \eta_\ell)| + \Con \, (\diam(\cc))^{1+ \hbeta}\nonumber\\
& \leq & \Con \, \diam(\cc) \, \|v_0 - \eta_0\| + \Con \, (\diam(\cc))^{1+ \hbeta} .\nonumber
\end{eqnarray}

Next, we consider two cases.

\ms

\noindent
{\bf Case 1.} $\diam(\cc) \leq \| v_0 -\eta_0\|^{\vartheta/2}$. Then (9.18) immediately implies
$$|\Delta(x_0, y_0) - \Delta (x_0, b_0)| \leq \Con \, \diam(\cc) \, \|v_0 - \eta_0\|^{\hbeta \vartheta/2} .$$

\ms

\noindent
{\bf Case 2.} $\diam(\cc) \geq \| v_0 -\eta_0\|^{\vartheta/2}$. 
Set 
$$\{X'\} = W^u_R(y_{0}) \cap W^s_{R}(x_{0}) \quad \mbox{\rm  and} \quad X = \phi_{\Delta(x_{0},y_{0})}(X') .$$
Then $X \in W^u_{\ep_0}(y_{0})$ and it is easy to see that
$$|\Delta(x_0, y_0) - \Delta (x_0, b_0)|  = |\Delta(X,b_{0})| .$$
We have $X = \exp^u_{y_{0}}(\tt)$ and $b_{0} = \exp^s_{y_{0}}(\ts)$ for some $\tt \in E^u(y_{0})$ and 
$\ts \in E^s(y_{0})$. Clearly $\|\tt\| \leq \Con$.
Using Liverani's Lemma (Lemma 4.1) we get
$$|\Delta(X,b_{0})| \leq C_0 [ |d\omega_{y_{0}}(\tt, \ts)|  + \|\tt\|^2 \|\ts\|^{\vartheta} + \|\tt\|^{\vartheta} \|\ts\|^2]
\leq \Con \, \|\ts\|^{\vartheta} .$$
However, $\|\ts\| \leq \Con d(y_{0},b_{0}) \leq \Con \|v_{0}-\eta_{0}\|$, so
$$|\Delta(X,b_{0})| \leq \Con \|v_0 - \eta_0\|^{\vartheta}  \leq \Con \diam(\cc) \|v_0 - \eta_0\|^{\vartheta/2} .$$
This proves the lemma.
\endofproof

\subsection{Proof of Lemma 4.3} 

\subsubsection{Set-up -- choice of some constants and initial points}

Choosing a constant $\ep' \in (0,r_0/2)$ sufficiently small, for any $z \in M$ and any 
$z'\in B^u(z,\ep')$ the local unstable holonomy map 
$$\hh_{z}^{z'} : W^s_{\ep'}(z) \longrightarrow W^s_{\ep_0}(z')$$
is well defined and uniformly H\"older continuous. Replacing $\ep'$ by a smaller constant if necessary, by (3.7) 
for $z\in P_0$ and $z'\in P_0\cap B^u(z,\ep')$ the {\it pseudo-holonomy map} 
\index{pseudo-holonomy map}
$$\hhh_{z}^{z'} = (\Phi^s_{z'})^{-1} \circ \hh_{z}^{z'} \circ \Phi^s_{z} : E^s (z; \ep') \longrightarrow E^s (z'; r_0) $$
is uniformly H\"older continuous, as well. Thus, there exist constants
$C' > 0$ and $\beta'' > 0$ (depending on the set $P_0$) so that for $z$, $z'$ as above we have
\be
\|\hhh_z^{z'}(u) - \hhh_z^{z'}(v)\| \leq C' \|u-v\|^{\beta''} \quad , \quad u,v\in E^s(z;\ep') .
\ee
We will assume $\beta'' \leq \beta$, where $\beta \in (0,1]$ is the constant from Sect. 3.

Fix arbitrary constants $\delta' > 0$ and $s_1 > 0$ with
\be
(\delta')^{\beta''} < \frac{\beta_0 \kappa \theta_0}{16 L_0 C_1 R_0\Gamma_0^2}  \quad , \quad 0 < s_1 < \delta'/(2R^2_0),
\ee
 and $\delta''$ with 
\be
0 < \delta''  < \min\left\{\frac{\delta'}{3R_0} , \frac{\beta_0 \hd_0 \kappa}{100 R^3_0 L_0 C^2_1}  ,
 \frac{s_1 \theta_0 c_0}{4C_1 \gamma_1 R^2_0}  \right\}   ,
\ee  
Then set 
$$\delta_0 = \frac{s_1 \theta_0}{16R_0} > 0 .$$
\index{$\hhh_z^{z'}$}

Next, assuming $\beta'' > 0$ is taken sufficiently small and $C' > 0$ sufficiently large, for any 
$j = 1, \ldots, \ell_0$ there exists a Lipschitz\footnote{Uniform continuity is enough.} family of 
unit vectors\footnote{E.g. define
$\eta_j(Z,z) = \frac{\left((\Phi^u_z)^{-1}\circ \Phi^u_Z(r_0\eta_j(Z)/2)\right)^{(1)}}
{\|\left((\Phi^u_z)^{-1}\circ \Phi^u_Z(r_0\eta_j(Z)/2)\right)^{(1)}\|}$.}
$$\eta_j(Z,z) \in E^u_1(z) \quad, \quad  Z\in \tP_0 \:, \: z \in B^u(Z,r_0/2)\cap \tP_0 ,$$
such that $\eta_j(Z,Z) = \eta_j(Z)$ and for any $v \in E^s(Z)$ we have
$$| \omega_z(\eta_j(Z,z), \hhh_Z^z(v)) - \omega_Z(\eta_j(Z), v)| \leq C' d(Z,z) \|v\|^{\beta''} .$$
Fix a constant $\ep'' \in (0,\ep'/2)$ so small that $C' (\ep'')^{\beta''} < \delta'$. Then
\be
| \omega_z(\eta_j(Z,z), \hhh_Z^z(v)) - \omega_Z(\eta_j(Z), v)| \leq  \delta' \|v\|^{\beta''} 
\:\: , \:\: Z \in \tP_0\: , \: z \in B^u(Z,r_0)\cap \tP_0 .
\ee

Using the symbolic coding provided by the Markov family $\{R_i\}$ it is easy to see that
there exists an integer $N_0\geq 1$ such that for any integer $N \geq N_0$  we have
$$\pp^N(B^u_{\ep'}(z)) \cap B^s(z', \delta'') \neq \e$$ 
for any $z,z'\in R$ (see the notation in the beginning of Ch. 4).

Fix for a moment $Z\in P_0 $. 
Given $j = 1, \ldots, \ell_0$, since $\eta_j(Z) \in E^u_1(Z)$, by 
Lemma 9.1 and the choice of $\theta_0 > 0$ (see Sect. 4.2), there exists  $\cvj(Z)\in E^s_1(Z)$ with 
$$d\omega_{Z}(\eta_j(Z) ,\cvj(Z)) \geq \theta_0 \quad, \quad  \|\cvj(Z)\| = 1 .$$
{\bf Fix a vector $\cvj(Z)$} with the above property for every $j$.

Set
\be
v_j (Z) = \frac{s_1}{R_0}  \cvj(Z) \in E^s_1(Z)\quad , \quad  y_j(Z) = \Phi^s_{Z}(v_j (Z)) \in W^s_{s_1}(Z)  .
\ee
Then $s_1/R_0^2 \leq d(Z,y_j(Z)) \leq s_1$. Since   $d\omega_{Z}(\eta_j(Z) , v_j(Z)) \geq s_1\theta_0/R_0 ,$  by (4.1),
\be
|d\omega_{Z}(\eta_j(Z) , v ) | \geq \frac{s_1 \theta_0}{2R_0}  
\ee
for all $v\in E^s(Z)$ with $\|v - v_j(Z)\| \leq \frac{s_1 \theta_0}{2C_0 R_0}$.

{\bf Fix an arbitrary  $N \geq N_0$}. It follows from the above that  
for each  $Z\in P_0$, each $i = 1,2$ and each $j = 1,\ldots, \ell_0$ there exist
\be
\yjo(Z) \in \pp^N(B^u(Z,\ep')) \cap B^s(y_j (Z), \delta'')  \:\: , \:\:
\yjt(Z) \in \pp^N(B^u(Z,\ep')) \cap B^s(Z, \delta'') .
\ee
Fix points $\yji(Z)$ with these properties; then $\yji(Z) \in W^s_{\ep_0}(Z)$.  We have
$$\yji(Z) = \Phi^s_{Z}(\wji(Z))\quad \mbox{\rm  for some } \quad   \wji(Z) \in E^s(Z)$$
such that 
$$\wji(Z) \in (\Phi^s_{Z})^{-1}(B^s(\yji (Z), \delta'')) .$$   
For $z\in B^u(Z,\ep')$ set 
\be
\wji(Z,z) = \hhh_{Z}^{z} (\wji(Z)) \in E^s(z) .
\ee
Notice that  
\be
\Phi^s_{z}(\wji(Z,z)) = \pi_{\yji(Z)}(z) . 
\ee

Given $Z\in P_0$ and $z\in B^u(Z,\ep') \cap \tP_0$,  it follows from $d(y_j(Z), \yjo(Z)) \leq \delta''$ that
$$\|\wjo(Z) - v_j (Z)\| \leq \delta'' R_0  ,$$ 
In particular, 
$$\frac{s_1}{2R_0} \leq \|\wjo(Z)\| \leq \frac{2s_1}{R_0} .$$
Apart from that, $\|\wjt(Z)\| \leq \frac{\delta''}{R_0}$. Now (9.19) gives
\begin{eqnarray}
\|\wjt(Z,z)\| 
& =      &  \|\hhh_Z^z(\wjt(Z)) - \hhh_Z^z(0)\| \leq C' \|\wjt(Z)\|^{\beta''} \\
& \leq & C' \left(\frac{\delta''}{R_0}\right)^{\beta''} < \frac{s_1}{4 R_0^3} .\nonumber
\end{eqnarray}
A similar estimate holds for $\wjo(Z,z)$, so we get
\be
\|\wjt(Z,z)\| \leq  \frac{s_1}{2R_0^3} \leq \|\wjo(Z,z)\| \leq 2 s_1  R_0
\:\: , \:\:  Z\in P_0 \;,\;  z\in B^u(Z,\ep') \cap \tP_0 .
\ee

Next, (9.24) implies
$$|d\omega_{Z}(\eta_j(Z) , \wjo(Z) ) | \geq \frac{s_1 \theta_0}{2R_0} ,$$
while (9.22) yields 
$$|d\omega_{Z}(\eta_j(Z,z) , \wjo(Z,z) ) | \geq 
\frac{s_1 \theta_0}{2R_0} - \delta' s_1 \left(\frac{4s_1}{R_0} \right)^{\beta''}\geq \frac{s_1 \theta_0}{4R_0} ,$$
and therefore
\be
|d\omega_{Z}(\eta_j(Z,z) , \wjo(Z,z) ) | \geq 4 \delta_0 \quad, \quad Z\in P_0 \:,\: z\in B^u(Z, \ep'')\cap \tP_0 .
\ee

To finish with this preparatory section, let $\cc$ be a cylinder of length $m$ in $R$ such that
$\cc \cap P_0 \setminus \Xi_m \neq \e$.
Let 
$$Z \in \cc \cap P_0 \setminus \Xi_m \quad , \quad Z_0 = \Psi(Z) ,$$
and let $x_0  \in \Psi(\cc)$,  $z_0  \in \Psi(\cc)$   
have the form $x_0 = \Phi^u_{Z_0}(u_0)$,  $z_0 = \Phi^u_{Z_0}(w_0)$, where
\be
d(x_0,z_0) \geq \kappa \, \diam (\Psi(\cc)) 
\ee
for some $\kappa \in (0,1]$, and
\be
\left\langle \frac{w_0 - u_0}{\|w_0 - u_0\|} , \eta_j(Z_0) \right\rangle \geq \beta_0 
\ee
for some $j = 1, \ldots, \ell_0$. Fix $\kappa$ and $j$ with these properties. Set $\tcc = \tPsi(\cc)$.
Then $Z_0, z_0 \in \tcc \cap \hP_0$. By the assumption on $m$, $\diam(\tcc) < \ep''$, so
$z_0\in B^u(Z_0, \ep'')$. Let $z_0 = \phi_{t_0}(z)$ for some $z\in \cc$ and $t_0 \in (-\chi, \chi)$. Set
$$x_0 = \Phi^u_{z_0}(\xi_0) \quad , \quad v_0 = d\phi_{t_0}(z) \cdot \wjo(Z, z_0) \in E^u(z_0; r_0/R_0) ,$$ 
for some $\xi_0 \in E^u(z_0; r_0/R_0)$; then $\|\xi_0\| \leq R_0\, \diam(\tcc)$.

\def\qq{\mathcal Q}

\subsubsection{Estimates for $|d\omega_{z_0} (\xio_0 , \vo_0)| $}

Recall that $z_0 = \Phi^u_{Z_0} (w_0)$.
Since $Z_0, z_0  \in \hP_0$ and $\|w_0\| \leq R_0 \ep'' < < r_0/R_0$, the map
$$\qq = (\Phi^u_{z_0})^{-1}\circ \Phi^u_{Z_0} : E^u(Z_0;r_0/ R_0^2) \longrightarrow E^u(z_0)$$
is well-defined and $C^{1+\beta}$. Using $d(\Phi^u_{z_0})^{-1}(z_0) = \id$,  $\qq (w_0) = 0$ and $\qq (u_0) = \xi_0$, we get
$$d\qq (w_0) =   d(\Phi^u_{z_0})^{-1}(z_0) \circ d\Phi^u_{Z_0}(w_0) = d\Phi^u_{Z_0}(w_0) .$$
Now (3.8) implies
\be
\|\xi_0 - d\Phi^u_{Z_0}(w_0) \cdot (u_0 - w_0) \| \leq 10 R_0^3 \|u_0 - w_0\|^{1+\beta}  .
\ee

\noindent
[{\it Proof of } (9.33): Using $C^2$ coordinates in $W^u_{r_0}(Z_0)$, we can identify $W^u_{r_0}(Z_0)$ with an
open subset  $V$ of $\R^{n^u}$ and regard $\Phi^u_{Z_0}$ and $\Phi^u_{z_0}$ as $C^{1+\beta}$ maps on $V$
whose derivatives and their inverses are bounded by $R_0$. By Taylor's formula (3.8),
$$\Phi^u_{Z_0}(u_0) - \Phi^u_{Z_0}(w_0) = d\Phi^u_{Z_0}(w_0)\cdot (u_0-w_0) + \eta ,$$
for some $\eta \in \R^{n^u}$ with $\|\eta\| \leq R_0 \|u_0 - w_0\|^{1+\beta}$. Hence
$$d (\Phi^u_{z_0})^{-1}(z_0) \cdot (\Phi^u_{Z_0}(u_0) - \Phi^u_{Z_0}(w_0)) 
= d\Phi^u_{Z_0}(w_0)\cdot (u_0-w_0) + \eta .$$
Since $Z_0 \in P_0$, by (3.9),
$$\|d\Phi^u_{Z_0}(w_0) - \id\| = \|d\Phi^u_{Z_0}(w_0) -  d\Phi^u_{Z_0}(0)\|  \leq R_0 \|w_0\|^{1+\beta} ,$$
so $\| d\Phi^u_{Z_0}(w_0)\| \leq 2R_0$.
Using Taylor's formula again,
\begin{eqnarray*}
\qq (u_0) - \qq (w_0) 
& = & (\Phi^u_{z_0})^{-1} (\Phi^u_{Z_0}(u_0)) - (\Phi^u_{z_0})^{-1} (\Phi^u_{Z_0}(w_0))\\
& = & d (\Phi^u_{z_0})^{-1}(z_0) \cdot (\Phi^u_{Z_0}(u_0) - \Phi^u_{Z_0}(w_0)) + \zeta 
\end{eqnarray*}
for some $\zeta$ with 
\begin{eqnarray*}
\|\zeta\| 
& \leq & R_0 \|\Phi^u_{Z_0}(u_0) - \Phi^u_{Z_0}(w_0)\|^{1+\beta}
\leq R_0\left( 2R_0 \|w_0 - u_0\| + R_0 \|w_0 - u_0\|^{1+\beta} \right)^{1+\beta} \\
& \leq & 9 R_0^{3} \|u_0-w_0\|^{1+\beta}.
\end{eqnarray*}
Thus, 
$$\xi_0 = \qq (u_0) - \qq (w_0) = d\Phi^u_{Z_0}(w_0)\cdot (u_0-w_0) + \eta + \zeta ,$$
where $\|\eta + \zeta\| \leq (R_0+ 9R_0^3) \|u_0-w_0\|^{1+\beta} \leq 10 R_0^3 \|u_0-w_0\|^{1+\beta}$.]

\bs

Next, by (4.8) the direction of $w_0 - u_0$ is close to $\eta_j(Z_0)$. More precisely, let
$$w_0 - u_0 = t \eta_j(Z_0) + u$$
for some $t \in \R$ and $u \perp \eta_j(Z_0)$. 
Then for $s = t/\|w_0 - u_0\|$ we have
$$\frac{w_0 - u_0}{\|w_0 - u_0\|} = s \eta_j(Z_0) + \frac{u}{\|w_0 - u_0\|} \, ,$$
so 
$$\di s= \left\la \frac{w_0 - u_0}{\|w_0 - u_0\|} , \eta_j(Z_0) \right \ra  \geq \beta_0 ,$$ 
and therefore
$t = s \|w_0 - u_0\| \geq \beta_0 \|w_0 - u_0\|$. Moreover,
\begin{eqnarray*}
\|u\|^2 
& =    & \|w_0 - u_0 - t \eta_j(Z_0)\|^2 = \|w_0 - u_0\|^2 - 2t \la w_0 - u_0 , \eta_j(Z_0) \ra +t^2\\
& =    & \|w_0 - u_0\|^2 \left(1- 2s \left\la \frac{w_0 - u_0}{\|w_0 - u_0\|}, \eta_j(Z_0) \right\ra 
+ s^2\right)\\
&  = & \|w_0 - u_0\|^2 (1-2s^2 +s^2) = \|w_0 - u_0\|^2 (1-s^2) \leq (1 - \beta_0^2) \|w_0 - u_0\|^2\;,
\end{eqnarray*}
and therefore $\|u\| \leq \sqrt{1-\beta_0^2}\,  \|w_0 - u_0\|$. 

Since $v_0 = d\phi_{t_0}\cdot \wjo(Z, z_0) = \hhh_{Z}^{z_0} (\wjo(Z))$, it follows from (9.24) 
with $z = z_0$ and $w = v_0$ that
$$|d\omega_{z_0}(\eta_j(Z,z_0), v_0)| \geq 8 \delta_0 ,$$ 
while (9.29) gives 
$$s_1/(2R_0^3) \leq \|v_0\| \leq 2 s_1 R_0  \leq 2\delta'/R_0 .$$
Using $d\Phi^u_{Z_0}(0) = \id$ and (3.9), we have
$$\|d\Phi^u_{Z_0}(w_0) - \id\| \leq R_0 \|w_0\|^\beta \leq R_0 (R_0\ep'')^\beta \leq R_0^2 (\ep'')^\beta .$$
Moreover, $\beta_0^2  (1+ \theta_0^2/(64C_0)^2) = 1$, so
$ \beta_0^2 \theta_0^2 = (64C_0)^2 (1-\beta_0^2)$, and therefore 
$$4 C_0 \sqrt{1-\beta_0^2} = \beta_0 \theta_0/16 .$$
The above, (9.29), (9.23), (9.21),
(9.22), $\|\vo_0\| \leq |v_0| \leq \|v_0\|$, Lemma 9.1 and the fact  that $\eta_j(z_0) \in E^u_1(z_0)$ imply 
\begin{eqnarray*}
&        & |d\omega_{z_0} (\xio_0 , \vo_0)| \nonumber\\
& =      &  |d\omega_{z_0} (\xi_0 , \vo_0)|  \geq |d\omega_{z_0} (d \Phi^u_{Z}(w_0)\cdot (u_0 - w_0) , \vo_0 )|\\
&         &           - |d\omega_{z_0} (\xi_0 - d \Phi^u_{Z}(w_0)\cdot (u_0 - w_0) , \vo_0 )|\nonumber\\
& \geq & t  |d\omega_{z_0} (d \Phi^u_{Z}(w_0)\cdot \eta_j(Z_0) , \vo_0 )|  -  |d\omega_{z_0} (d \Phi^u_{Z}(w_0)\cdot u , \vo_0 )|\\
&         &   - 10 C_0 R_0^3 \|u_0 - w_0\|^{1+\beta}\|\vo_0\|\nonumber\\
& \geq &  \beta_0 \|u_0 - w_0\| \, [\;  |d\omega_{z_0} (\eta_j(z_0) , \vo_0 )| 
- |d\omega_{z_0} (d \Phi^u_{Z}(w_0)\cdot \eta_j(z_0) - \eta_j(z_0), \vo_0 )| \; ]\nonumber\\
&         & - C_0 (1+  R_0^2 (\ep'')^\beta)\sqrt{1-\beta_0^2}\,  \|u_0 - w_0\| \|\vo_0\|  
- 10 C_0 R_0^3 \|u_0 - w_0\|^{1+\beta}\|\vo_0\|\nonumber\\
& \geq & \|u_0 - w_0\| \, [ \; \beta_0  |d\omega_{z_0} (\eta_j(z_0) , v_0 )|  -  \beta_0 C_0 R_0^2 (\ep'')^\beta \|v_0\| 
           - 2 C_0 \sqrt{1-\beta_0^2}\,  \|v_0\| \\
&         &   - 10C_0 R_0^3 (2\ep'')^{\beta}\|v_0\| \;]\nonumber\\
& \geq & \|u_0-w_0\| \, [ \; 8 \beta_0 \delta_0 -  2 \beta_0 C_0 R_0  \delta'
          - 4 C_0 \sqrt{1-\beta_0^2}\,  \delta'/R_0 - 40 C_0 R_0^2 \delta' \;]\nonumber\\
& \geq & \|u_0-w_0\| \, [ \; 8 \beta_0 \delta_0 -  \beta_0 \delta_0
          - \beta_0 \delta_0 - \beta_0 \delta_0  \;] \geq 4 \|u_0-w_0\| \, \beta_0\delta_0 . 
\end{eqnarray*}
Combining this with  (4.7) and (3.7) gives
\be
|d\omega_{z_0} (\xio_0 , \vo_0)|  \geq \frac{4 \beta_0 \delta_0 \kappa }{R_0} \, \diam(\tcc)\, .
\ee

Next, set $\txi_0 = \Psi^u_{z_0} (\xi_0) \in E^u(z_0)$.  Then 
\be
\exp^u_{z_0}(\txi_0) = \Phi^u_{z_0}(\xi_0) = x_0 ,
\ee
and
\be
\frac{\kappa}{R_0}\, \diam(\tcc) \leq \|\xi_0\| \leq R_0  \diam(\tcc) . 
\ee

Setting $\tv_0 = \Psi_{z_0}^s (v_0)\in E^s(z_0)$ and $y_0 = \exp^s_{z_0}(\tv_0)$,
using $v_0 = \wjo (Z,z_0)$, (9.25) and (9.27), we get
\be 
y_0 = \exp^s_{z_0}(\tv_0) = \Phi^s_{Z}(\wjo(Z,z_0)) = \pi_{\yjo(Z)}(z_0)
\in  B^s(z_0, \delta'')  .
\ee

We will now prove that
\be
|\Delta(x_0,y_0)| \geq \frac{2 \beta_0 \delta_0\kappa}{R_0} \, \diam(\tcc) .
\ee
From this and Lemma 4.2(b), (4.9) follows easily for
$d_1 \in B^s(\yj_1(Z), \delta'')$ and $d_2 \in B^s (Z, \delta'')$,
using the choice of $\delta''$.

It follows from (3.6), $\|v_0\|\leq r_0/R_0$ and $\|\xi_0\|\leq r_0/R_0$ that
$$\|\tv_0 - v_0\| \leq R_0 \|v_0\|^{1+\beta} \quad , \quad \|\txi_0 - \xi_0\|\leq R_0\|\xi_0\|^{1+\beta} ,$$
and in particular $\|\tv_0\| \leq 2\|v_0\|$ and $\|\txi_0\| \leq 2\|\xi_0\| \leq 2 R_0 \diam(\tcc)$.

As in Sect. 9.3.1, set $p = [\ttau_m(z_0)]$, $q = [p/2]$, and for $j \geq 0$ define
$z_j = \varphi^j(z_0)$, $x_j = \varphi^j(x_0)$, $y_j = \varphi^j(y_0)$,
$\hxi_j =  d\hf^j_{z_0}(0)\cdot \xi_0$, etc. in the same way. By the choice of $\ep'' > 0$  all estimates in Sect. 
9.3.1 hold without change. Choosing an arbitrary $z\in \cc \cap P_0 \setminus \Xi_m$, as before we find
$j \geq 0$ with $\pp^j(z) \in P_0$ such that (9.8) holds for $\ell = [\ttau_j(\Psi(z))]$ and
$r(z_\ell) \geq r_0$. Fix $\ell$ with these properties; then (9.9) and (9.10) hold again.

We need an estimate from below for $|d\omega_{z_\ell}(\xi_\ell, v_\ell)|$ similar to (9.11).
Instead of using Lemma 10.7 this time it is enough to use Lemma 3.1.
Since $v_\ell = \hf^{\ell}_{z_\ell}(v_0)  \in E^s(z_\ell)$ and $z_0\in P$ implies  $L(z_0) \leq L_0$, for
$w = d\hf^{-\ell}_{z_\ell}(0)\cdot v_\ell$, using Lemma 3.1, we get
\be
\|\vo_0 - \wo \| \leq L_0(z) |v_0|^{1+\beta} \leq L_0  \|v_0\|^{1+\beta} .
\ee
As in the proof of (9.11) we will now use the estimates in Sect. 9.3.2.
It  follows from Lemma 9.1, (9.9) and (9.35) that
\begin{eqnarray*}
|d\omega_{z_\ell}(\xi_\ell, v_\ell)| 
& \geq & |d\omega_{z_\ell}(\xio_\ell, \vo_\ell)| - \sum_{i=2}^{\tk} |d\omega_{z_\ell}(\xii_\ell, \vi_\ell)| \\
& \geq & |d\omega_{z_\ell}(\hxio_\ell, \vo_\ell)| - C_0 L_0 \|\xi_\ell\|^{1+\beta} \|v_\ell\| 
- C_0 \sum_{i=2}^{\tk} \|\xii_\ell\|\,  \|\vi_\ell\| \\
& =      &  |d\omega_{z_0}(d\hf^{-\ell}_{z_\ell}(0)\cdot \hxio_\ell, d\hf^{-\ell}_{z_\ell}(0)\cdot\vo_\ell)| \\
&         &     - C_0 L_0  \|\xi_\ell\|^{1+\beta} \|v_\ell\| - C_0 \sum_{i=2}^{\tk} \|\xii_\ell\|\,  \|\vi_\ell\| \\
& =      &  |d\omega_{z_0}(\xio_0, \wo)| - C_0 L_0  \|\xi_\ell\|^{1+\beta} \|v_\ell\|  
- C_0 \sum_{i=2}^{\tk} \|\xii_\ell\|\,  \|\vi_\ell\| \\
& \geq &  |d\omega_{z_0}(\xio_0, \vo_0)| - C_0 L_0 R_0  \, \diam(\tcc) \|v_0\|^{1+\beta} - \Con \, (\diam(\tcc))^{1+\hbeta}  .
\end{eqnarray*}
Combining this with (9.6) and (9.34) gives
\begin{eqnarray*}
|\Delta(x_0,y_0)|
& =      &  |\Delta(x_\ell, y_\ell)| \geq |d\omega_{z_\ell}(\xi_\ell,v_\ell)| 
-  8 C_0 R_0\, \|\xi_\ell\|\, \|v_\ell\| (\|\xi_\ell\|^\beta + \|v_\ell\|^\beta)\nonumber\\ 
&         & - 8 C_0\,\left[ \|\xi_\ell\|^2\, \|v_\ell\|^\vartheta 
+ \|\xi_\ell\|^\vartheta \|v_\ell\|^2 \right]\nonumber\\
& \geq & |d\omega_{z_0}(\xio_0, \vo_0)| - C_0 L_0 R_0 \,  \diam(\tcc) \|v_0\|^{1+\beta}
- \Con\, (\diam(\tcc))^{1+\hbeta} \nonumber\\
& \geq & \frac{4\beta_0 \delta_0\kappa}{R_0} \, \diam(\tcc) - C_0 L_0 R_0  \, \diam(\tcc) \|v_0\|^{1+\beta}
- C''' (\diam(\tcc))^{1+\hbeta}          
\end{eqnarray*}
for some constant $C''' > 0$. Now assume
$$(2\ep'')^{\hbeta} \leq \frac{\beta_0 \delta_0 \kappa }{R_0 C''' } ,$$
and recall that $\|v_0\| \leq \delta'$ and $\diam(\tcc) \leq 2\ep''$. By (9.29), $\|v_0\|\leq 2s_1$, while
(9.20) implies $\|v_0\|^\beta \leq (\delta'')^\beta  < (\delta')^\beta < \frac{\beta_0\kappa \theta_0}{16 L_0 C_0 R_0}$. 
Thus, using (8.21),
$$ C_0 L_0 R_0 \,  \diam(\tcc) \|v_0\|^{1+\beta} \leq 
C_0 L_0 R_0 \diam (\tcc) \, 2s_1 \,\frac{\beta_0\kappa \theta_0}{16 L_0 C_0 R^2_0} 
\leq \diam (\tcc) \, \frac{\beta_0 \delta_0 \kappa}{R_0} ,$$
and therefore $\Delta(x_0, y_0 ) \geq  \frac{2\beta_0 \delta_0 \kappa}{R_0} \, \diam(\tcc)$.
This proves (9.38).
\endofproof

\newpage

\section{Regular distortion for Anosov flows}
\setcounter{equation}{0}

In this section we prove Lemma 4.4. Here we do not need to assume that the flow $\phi_t$ is contact.

As in Sect. 9,  here we will use  the constant $\hep_1 = \frac{\beta \hep}{6}$, where $\beta > 0$ is as in  Sect. 3 (see p. 19), and
throughout we will again assume that 
$R(x)$, $\Gamma(x)$, $D(x)$ and $L(x)$ are Lyapunov $\hep_1$-regularity functions, while $r(x)$ is an
$\hep_1$-slowly varying radius function and $\tr(x)$ is an $\hep_1/2$-slowly varying radius function with the properties in Sect. 3, 
so that they satisfy (3.6) -- (3.18)  and the conclusions of  Lemma 3.1 with $\hep$ replaced by $\hep_1$. 
Let us remark again that, since $\hep_1 < \hep$,  $R(x)$, $\Gamma(x)$, $D(x)$ and $L(x)$ are still Lyapunov $\hep$-regularity functions, $r(x)$ is still an
$\hep$-slowly varying radius function and $\tr(x)$ is an $\hep/2$-slowly varying radius function.

\subsection{Expansion along $E^u_1$}

Let again $M$ be a $C^2$ complete Riemannian manifold and let $\phi_t$  be a  $C^2$ Anosov flow on $M$.
Set
$$\hmu_2 = \lambda_1 + \frac{2}{3} (\lambda_2-\lambda_1)  \quad , \quad  \hnu_1 = \lambda_1 + \frac{1}{3} (\lambda_2 - \lambda_1) .$$
Then  $\hmu_2 < \mu_2 e^{-\hep}$ and 
$$\lambda_1 < \nu_1 < \hnu_1 < \hmu_2 < \mu_2 < \lambda_2 .$$
For $\hep > 0$, apart from (3.1), we assume in addition that 
$$e^{\hep} \leq \frac{2\lambda_2}{\lambda_2 + \hmu_2} .$$

For a non-empty set $X\subset E^u (x)$  set 
$$\ell(X) = \sup \{ \|u\|  : u \in X\} .$$
Given $z\in \ll$ and $p \geq 1$, setting $x = \varphi^p(z)$, define
$$\hB^u_p (z, \delta) = \{ u \in E^u (z) : \|\hf^p_z(u) \| \leq\delta\} .$$

Fix for a moment $x\in \ll$ and an integer $p \geq 1$, set $z = \varphi^{-p}(x)$ and given $v \in E^u (z; r(z))$, set 
$$z_j = \varphi^j(z)\quad, \quad v_j = \hf^j_z(v) \in E^u (z_j) \quad , \quad w_j = d\hf_z^j(0)\cdot v \in E^u (z_j)$$
for any $j = 0,1,\ldots, p$ (assuming that these points are well-defined).
\index{$\hmu_2$}
\index{$\hnu_1$}
\index{$\tvt$}
\index{$\ell(X)$}
\index{$\hB^u_p(z,\delta)$}

For any $v = \vo + \vt + \ldots+ v^{(\tk)} \in E^u(x)$ with $v^{(j)} \in E^u_j$, set
$$\tvt =  \vt + \ldots+ v^{(\tk)} \in \tE_2^u(x) .$$

\ms

\noindent
{\bf Lemma 10.1.}  {\it Assume that the $\frac{3\hep_1}{\beta}$-slowly varying radius function  $\hr(x) \leq r(x)$ on $\ll$ satisfies 
\be
\hr(x) \leq  \min \left\{  \left( \frac{1/\hmu_2 - 1/\lambda_2}{6  \Gamma^2 (x) D (x)} \right)^{1/\beta}  \;, \;
\left( \frac{1/\lambda_1 - 1/\hnu_1}{6 e^{3\hep} \Gamma^2(x) D(x)} \right)^{1/\beta} \right\}\;
\ee
for all $x\in \ll$. Then for any $x\in \ll$ and any 
$$V = \Vo + \tVt \in E^u(x; \hr(x)) , $$ 
setting $y = \varphi^{-1}(x)$ and $U = \hf_x^{-1}(V)$, we have
\be
\|\tUt\|'_{y} \leq \frac{\|\tVt\|'_x}{\hmu_2} ,
\ee
and
\be
\|\Uo\|'_{y} \geq \frac{\|\Vo\|'_x}{\hnu_1} \;.
\ee
Moreover, if $V, W \in E^u(x; \hr(x))$ and $\Wo = \Vo$, then for $S = \hf_x^{-1}(W)$ we have
\be
\|\tUt - \tSt\|'_{y} \leq \frac{\|\Vt - \tWt\|'_x}{\hmu_2} ,
\ee
and, if $\tWt =   \Vt \in E^u(x; \hr(x))$ and  $S = \hf_x^{-1}(W)$ again, then}
\be
\|\Uo - \So\|'_{y} \geq \frac{\|\Vo - \Wo\|'_x}{\hnu_1}  .
\ee

\ms

\noindent
{\it Proof.} The estimates (10.2) and (10.3) follow from Lemma 3.5 in \cite{St4}, while the proofs of (10.4)
and (10.5) are similar. We will prove (10.4) for completeness.

Assuming $V,W \in E^u(x; r(x))$ and $\Vo = \Wo$, by (3.16),
$$\|\hf_x^{-1}(V) - \hf_x^{-1}(W) - d\hf^{-1}(W)\cdot (V-W)\| \leq D(x) \|V-W\|^{1+\beta} ,$$
so
$$\hf_x^{-1}(V) - \hf_x^{-1}(W) = d\hf_x^{-1}(W)\cdot (\Vt - \tWt) + \xi$$
for some $\xi\in E^u(y)$ with  
$$\|\xi\| \leq D(x) \|\tVt - \tWt\|^{1+\beta} .$$
By (3.17),
$$\| d\hf_x^{-1}(W)\cdot (\tVt - \tWt)  - d\hf^{-1}_x(0)\cdot (\tVt- \tWt)\| \leq  D(x) \|W\|^{\beta}\, \|\tVt-\tWt\| ,$$
so
$$d\hf_x^{-1}(W)\cdot (\tVt - \tWt)  = d\hf^{-1}_x(0)\cdot (\tVt- \tWt) + \eta ,$$
with 
$$\|\eta\| \leq D(x) \|W\|^{\beta}\, \|\tVt-\tWt\| .$$
Now we get
$$U - S = \hf_x^{-1}(V) - \hf_x^{-1}(W) = d\hf_x^{-1}(0)\cdot (\tVt - \tWt) + \xi + \eta ,$$
which yields
$$\tUt - \tSt = d\hf_x^{-1}(0)\cdot (\tVt - \tWt) + \txi^{(2)} + \teta^{(2)} .$$
This, (3.10) and (3.14) imply
\begin{eqnarray*}
&       & \|\tUt - \tSt\|'_{y} 
 \leq  \|d\hf_x^{-1}(0)\cdot (\tVt - \tWt)\|'_y + \Gamma(y) (\|\xi\| + \|\eta\|)\\
& \leq & \|\tVt - \tWt\|'_x \; \left(\frac{1}{\mu_2} 
+ \Gamma(y) D(x)  (\|\tVt- \tWt\|^\beta + \|W\|^\beta)\right)\\
& \leq &  \|\tVt - \tWt\|'_x \; \left(\frac{1}{\mu_2}  + \Gamma(x) e^{\hep_1} D(x)  3\hr^\beta(x))\right) .
\end{eqnarray*}
By (10.1) and (3.2),
\begin{eqnarray*}
\frac{1}{\mu_2} + \Gamma(x) e^{\hep_1} D(x)  3 \hr^\beta (x)
 \leq  \frac{e^{\hep}}{\lambda_2} + e^{\hep}\, \frac{1/\hmu_2 - 1/\lambda_2}{2}
= e^{\hep} \frac{\lambda_2 + \hmu_2}{2 \lambda_2 \hmu_2} < \frac{1}{\hmu_2} ,
\end{eqnarray*}
since $e^{\hep} < \frac{2\lambda_2}{\lambda_2 + \hmu_2}$. 
This proves (10.4).
\endofproof

\bs

Next,  for any $y \in \ll$, $\ep\in (0,r(y)]$ and $p \geq 1$ set 
$$\hBo_p(y,\ep) =   \hB^u_p(u;\ep) \cap E^u_1(y) .$$
\index{$\hBo_p(y, \ep)$}

Replacing the $\frac{3\hep_1}{\beta}$-slowly varying radius function $\hr(x)$ by a smaller one\footnote{In this case we just need to multiply
the original $\hr(x)$ by a sufficiently small constantt.}, we may assume that
\be
(\Gamma(x))^{1+\beta} L(x) (\hr(x))^\beta \leq \frac{1}{100 n_1} \quad , \quad x\in \ll ,
\ee
where $n_1 = \dim (E^u_1(x))$.

In what follows we assume that $\hr(x)$ is a {\bf fixed $\frac{3\hep_1}{\beta}$-slowly varying radius function} on $\ll$ satisfying (10.1) and (10.6).

The proof of the following lemma is similar to the proof of Proposition 3.2 in \cite{St4}.

\bs

\noindent
{\bf Lemma 10.2.}  {\it 
Let  $z\in \ll$ and $x = \varphi^p(z)$ for some integer $p \geq 1$, and let $\ep \in (0, \hr(x)]$. Then
\be
\ell (\hB^{u}_p(z,\ep)) \leq  2 \tk \Gamma^4 (x) \,  \ell (\hBo_p(z,\ep) ) .
\ee
Moreover for any $\ep' \in (0, \ep]$ there exists $u \in \hBo_p(z,\ep')$ with}
\be
\|u\| \geq \frac{\ep'}{2 \tk \ep \Gamma^3 (x)}\ell (\hB^{u}_p(z,\ep))\quad \mbox{\rm and} \quad \|\hf^p_z(u)\| \geq \ep'/2 .
\ee

\ms

\noindent
{\it Proof.}   
Let $z\in \ll$ and  $x = \varphi^p(z)$ for some integer $p\geq 1$.
Let $v = (\vo,\tvt) \in \hB^u_p(z,\ep)$ be such that $\|v\|$ is the maximal possible,  i.e. $\|v\| = \ell(\hB^u_p(z,\ep))$. Set
$$V = (\Vo, \tVt) = \hf^p_z(v) \quad , \quad W = (\Wo, \tWt) = d\hf^p_z(0) \cdot v .$$
Then $|V| \leq \|V \| \leq \Gamma(x) \ep$, and by Lemma 3.1 and (10.6),
$$\|\Wo - \Vo\| \leq L(x) |V|^{1+\beta} \leq \frac{1}{100} |V| < \frac{\ep}{100} ,$$
so
$$\|\Wo\| \leq  \|\Vo\| + L(x) |V|^{1+\beta} \leq  \|\Vo\| + \frac{\ep}{100} .$$

\ms

\noindent
{\bf Case 1.} $\|\tvt\| \geq \|\vo\|$. Take $U = (\Uo,0)\in E^u_1(x,\ep)$ such that $\|\Uo\| \geq \ep/2$. 
Then $u = \hf_x^{-p}(U)\in \hB^u_p(z,\ep) \cap E^u_1(z)$. Applying (10.3) $p$ times gives
$$\|\uo\|'_z \geq \frac{\|\Uo\|'_x}{\nu_1^p} \geq \frac{\ep}{2\nu_1^p} .$$
Similarly, applying (10.2) $p$ times for $V = \hf^p_z(v)$ we get
\begin{eqnarray*}
\|\tvt \| 
& \leq & \|\tvt\|'_z \leq \frac{\|\tVt\|'_x}{\hmu_2^p} \leq \frac{\Gamma(x) \|\tVt\|}{\hmu_2^p} \leq \Gamma^2(x) \frac{\ep}{\hmu_2^p}\\
& \leq & 2 \Gamma^2 (x) \left(\frac{\hnu_1}{\hmu_2}\right)^p\, \frac{\ep}{2\hnu_1^p} 
\leq  2 \Gamma^2 (x) \left(\frac{\hnu_1}{\hmu_2}\right)^p\, \|\uo\|'_z .
\end{eqnarray*}
Thus,
\begin{eqnarray*}
\|v\| 
& \leq & \Gamma (z)\, |v| = \Gamma (z)\, \|\tvt\| \leq  2 \Gamma^2 (x) \Gamma(z) \left(\frac{\hnu_1}{\hmu_2}\right)^p\, \|\uo\|'_z\\
& \leq & 2 \Gamma^2 (x) \Gamma^2(z)   \left(\frac{\hnu_1}{\hmu_2}\right)^p\,\|\uo\|
\leq  2 \Gamma^4(x) \left(\frac{\hnu_1 \, e^{2\hep}}{\hmu_2}\right)^p\,  \|\uo\| ,
\end{eqnarray*}
which proves (10.7), since $\frac{\hnu_1 \, e^{2\hep}}{\hmu_2} < 1$.

To prove (10.8) in this case, we use the above argument however this time we choose 
$U = (\Uo,0)\in E^u_1(x,\ep')$ with $\|\Uo\| \geq \ep'/2$. 
Then as above we derive 
$$\|\uo\|'_z \geq  \frac{\ep'}{2\hnu_1^p}$$ 
and
$$\|\tvt \|  \leq 2 \Gamma^2 (x) \left(\frac{\hnu_1}{\hmu_2}\right)^p\, \frac{\ep}{2\hnu_1^p} 
\leq  2 \Gamma^2 (x) \frac{\ep}{\ep'} \left(\frac{\hnu_1}{\hmu_2}\right)^p\, \|\uo\|'_z 
\leq  2 \Gamma^3 (x)  \frac{\ep}{\ep'} \, \|\uo\| ,$$
which proves (10.8).

\bs

\noindent
{\bf Case 2.} $\|\tvt\| < \|\vo\|$. Set 
$\di u = \frac{\ep'}{2\ep \Gamma(x)} \vo \in E^u_1(z) .$
We will now check that $u \in \hB^u_p(z,\ep')$. Indeed, by Lemma 3.1 and (10.6),
$$\|\hf^p_z(u) - d\hf^p_z(0) \cdot u\| \leq L(x) |\hf^p_z(u)|^{1+\beta} \leq \frac{\|\hf^p_z(u)\|}{100} ,$$
therefore
\begin{eqnarray*}
\|\hf^p_z(u)\| 
& \leq & \|d\hf^p_z(0)\cdot u\| + L(x) |\hf^p_z(u)|^{1+\beta}\\
& \leq &  \|d\hf^p_z(0)\cdot u\|  + L(x) \, (\hr(x))^\beta \, (\Gamma(x))^{1+\beta} \|\hf^p_z (u)\| \\
& \leq & \|d\hf^p_z(0)\cdot u\|  + \frac{1}{100} \, \|\hf^p_z (u)\| .
\end{eqnarray*}
This gives
\begin{eqnarray*}
\frac{99}{100} \, \|\hf^p_z(u)\| 
& \leq & \|d\hf^p_z(0)\cdot u\| = \frac{\ep'}{2 \ep \Gamma(x)} \, \| d\hf ^p_z(0) \cdot \vo\|\\
& =    & \frac{\ep'}{2 \ep \Gamma(x)} \, \| \Wo\| \leq
\frac{\ep'}{2 \ep \Gamma(x)} \, \left( \| \Vo\| + \frac{\ep}{100}\right) \\
& \leq & \frac{\ep'}{2 \ep \Gamma(x)} \, \left( \Gamma (x) \ep + \frac{\ep}{100} \right)  < \ep' .
\end{eqnarray*}
Thus, $u \in \hB^u_p(z,\ep')$. Since 
$$\|u\| = \frac{\ep'}{2\ep \Gamma(x) } \|\vo\| = \frac{\ep'}{2\ep \Gamma(x)} |v| \geq \frac{\ep' \|v\|}{2 \ep \tk \, \Gamma (x)} ,$$
taking $\ep' = \ep$, proves (10.7).

To prove (10.8), let now $u \in \hBo_p(z,\ep')$ be such that $\|u\|$ is the maximal possible, and
set $U = \hf^p_z(u) \in E^u_1(x, \ep')$. It follows from the previous argument that 
$$\|u\| \geq \frac{\ep'}{2 \tk \ep \Gamma(x)}\, \ell(\hB^u_p(z,\ep)) .$$
It remains to show that $\|U\| \geq \ep'/2$. If $\|U\| < \ep'/2$, then by Lemma 3.1 and (10.6), 
$$\|U - d\hf^p_z(0) \cdot u\| \leq \frac{\|U\|}{100} < \frac{\ep'}{200} ,$$
so $\|d \hf^p_z(0) \cdot u \| < 2\ep'/3$. Setting $\hu = t u$ for some $t > 1$, $t$ close to $1$, we get
$$\|\hf^p_z(\hu) - d\hf^p_z(0)\cdot (\hu)\| \leq \frac{\|\hf^p_z(\hu)\|}{100} ,$$
so 
$$\|\hf^p_z(\hu)\| \leq \frac{100}{99} \|d \hf^p_z(0)\cdot 
\hu\| \leq \frac{100}{99}\cdot \frac{2t\ep'}{3} < \ep' ,$$
if $t$ is sufficiently close to $1$. Thus, $\hu \in \hBo_p(z,\ep')$ for $t > 1$, $t$ close to $1$.
However $\|\hu\| = t\|u\| > \|u\|$, contradiction with the choice of $u$. Hence we must have $\|U\| \geq \ep'/2$,
which proves (10.8).
\endofproof

\bs

To prove the main result in this section, it remains to compare diameters of sets of the form $\hBo_p(y,\ep)$. 

\bs  

\noindent
{\bf Lemma 10.3.} {\it There exists an $\frac{5\hep_1}{\beta}$-slowly varying radius function $\hr_1(x) \leq \hr(x)$ ($x\in \ll$)  such that:}

\ms

(a) {\it For any  $x\in \ll$ and any $0 < \delta \leq  \ep \leq \hr_1(x)$ we have
\be
\ell \left( \hBo_{p}(\varphi^{-p}( x),\ep)  \right) 
\leq  16 n_1 \frac{\ep}{\delta} \, \ell \left( \hBo_p (\varphi^{-p} ( x) , \delta) \right) \;
\ee
for any integer $p \geq  1$.}

\ms

(b) {\it For any  $x\in \ll$, any $0 < \ep \leq \hr_1(x)$  and any $\rho \in (0,1)$, for any $\delta$ with
$$0 < \delta \leq  \frac{\rho \, \ep}{16 n_1} $$
we have
\be
\ell \left( \hBo_p(\varphi^{-p}(x),\delta)  \right) \leq  \rho \,  \ell \left( \hBo_p(\varphi^{-p} (x) , \ep) \right) \;
\ee
for any  integer $p \geq 1$.}

\bs  

\noindent
{\bf Theorem 10.4.} {\it There exists an $\frac{6\hep_1}{\beta}$-slowly varying radius function $\hr_1(x) \leq \hr(x)$ ($x\in \ll$) such that:}

\ms

(a) {\it For any  $x\in \ll$ and any $0 < \delta \leq  \ep \leq \hr_1(x)$ we have
$$\ell \left( \hB^{u}_p(z,\ep)  \right) \leq  \frac{32 \tk n_1 \Gamma^4(x)  \ep}{\delta}  
\, \ell \left(\hB^{u}_p(z,\delta) \right) \;$$
for any integer $p \geq  1$, where $z = \varphi^{-p}(x)$.}

\ms

(b) {\it For any  $x\in \ll$, any $0 < \ep \leq \hr_1(x)$, any  $\rho \in (0,1)$ and any $\delta$ with
$$0 < \delta  \leq \frac{\rho \ep}{32 \tk n_1 \Gamma^3(x) }$$
we have
$$\ell \left( \hB^{u}_p(z ,\delta) \right) \leq  \rho \, \ell \left( \hB^{u}_p( z ,\ep)  \right) $$
for all integers $p \geq 1$, where $z = \varphi^{-p}(x)$.}

\ms

(c) {\it For any $x\in \ll$, any  $0 < \ep' < \ep \leq \hr_1(x)/2$, any 
$0 < \delta < \frac{\ep'}{100n_1}$ and any integer $p \geq 1$,
setting $z = \varphi^{-p}(x)$, there exists $u \in \hBo_p(z,\ep')$ such that for every $v \in E^u(z)$ with
$$\|\hf^p_z(u) - \hf^p_z(v)\| \leq \delta$$
we have} 
$$\|v\| \geq \frac{\ep'}{4 \ep \tk \Gamma^4(x)} \, \ell(\hB^u_p(z,\ep)) .$$

\bs

Using Lemma 10.3, we will now prove Theorem 10.4. The proof of Lemma 10.3 is given
in the next section. In fact, part (c) above is a consequence of Lemmas 3.1 and 
10.2 and does not require Lemma 10.3.

\bs

\noindent
{\it Proof of Theorem} 10.4. Choose the function $\hr_1(x)$ as in Lemma 10.3.

\ms

(a) Let $0 < \delta < \ep \leq \hr_1(x)$. 
Given an integer $p \geq 1$, set $z = \varphi^{-p}(x)$. Then Lemmas 10.2 and 10.3 and (10.7) imply
\begin{eqnarray*}
\ell (\hB^{u}_p(z,\ep)) 
& \leq & 2 \tk \Gamma^4(x) \, \ell (\hBo_p(z,\ep))
\leq  2 \tk \Gamma^4(x)  \, 16 n_1  \frac{\ep}{\delta}\, \ell (\hBo_p(z,\delta))\\
& \leq &  32 \tk n_1 \Gamma^4(x)  \frac{\ep}{\delta} \, \ell (\hB^{u}_p(z,\delta)) .
\end{eqnarray*}

(b) Let $x\in \ll$ and $0 < \ep \leq \hr_1(x)$. Given $\rho \in (0,1)$, set
$\rho' = \frac{\rho}{2\tk \Gamma^3(x)} < \rho$.  By Lemma 10.3(b), if
$0 < \delta \leq \frac{\rho' \ep}{16 n_1}$ then (10.10) holds with $\rho$ replaced by $\rho'$
for any  integer $p \geq 1$ with $z = \varphi^{-p}(x)$. Using this and Lemma 10.2 we get
\begin{eqnarray*}
\ell (\hB^{u}_p(z,\delta) ) \leq   2 \tk \Gamma^3(x) \, \ell (\hBo_p(z,\delta))
 \leq    2 \tk \Gamma^3(x) \,\rho'\,  \ell (\hB^{u}_p(z,\ep)) =  \rho\,  \ell (\hB^{u}_p(z,\ep)) ,
\end{eqnarray*}
which completes the proof.

\ms

(c) We will prove the following:

\ms

(d) {\it For any $x\in \ll$, any  $0 < \ep' < \ep \leq \hr_1(x)/2$, any 
$0 < \delta < \frac{\ep'}{100n_1}$ and any integer $p \geq 1$,
setting $z = \varphi^{-p}(x)$, there exists $u \in \hBo_p(z,\ep')$ such that for every $v \in E^u(z)$ with
$$\|\hf^p_z(u) - \hf^p_z(v)\| \leq \delta/\Gamma (x)$$ 
we have} 
$$\|v\| \geq \frac{\ep'}{4 \ep \tk \Gamma^3(x)} \, \ell(\hB^u_p(z,\ep)) .$$

As one can see, the only difference between (d) and (c) is in the condition involving $\delta$;
in (d) $\delta$ is replaced by $\delta/\Gamma(x)$. Once (d) is proved, replacing the $\frac{3\hep}{\beta}$-slowly varying radius  function
$\hr_1(x)$ by $\hr_1(x)/\Gamma(x)$, thus obtaining a $\frac{4\hep}{\beta}$-slowly varying radius  function,
and multiplying all constants appearing in (d) by $\Gamma(x)$ will prove (c).

\ms

So, we will prove (d) now.

Given $x\in \ll$, $z = \varphi^{-p}(x)$, let $\ep'$, $\ep$ and  $\delta$ be as in the assumptions in (d). 
Let $u \in \hBo_p(z,\ep')$ be such that $\|u\|$ is the maximal possible.
By Lemma 10.2, for $U = \hf^p_z(u) \in E^u_1(x)$ we have $\ep'/2 \leq \|U\| \leq \ep'$. 
Setting $W = d\hf^p_z(0) \cdot u \in E^u_1(x)$, Lemma 3.1 and (10.6) give
$$\| W - U\| \leq L(x) |U|^{1+\beta} \leq \frac{\|U\|}{100 n_1} \leq \frac{\ep'}{100 n_1} ,$$
so $\|W\| \leq \frac{101 \ep'}{100}$. 

Let $v = (\vo, \tvt) \in E^u(z)$ be such that for $V = \hf^p_z(v)$ we have $\|V - U\| \leq \delta/ \Gamma(x)$.
Then $|V-U| \leq \delta$, so $\|\Vo - \Uo\| \leq \delta$ and $\| \tVt \| \leq \delta$.  Moreover,
$$\|V\| \leq \|U\| + \delta/\Gamma(x) \leq \ep' + \delta/\Gamma (x) .$$

Set $S = d\hf^p_z(0) \cdot v$; then $\So = d\hf^p_z(0) \cdot \vo$. By Lemma 3.1 and (10.7),
$$\| \So - \Vo\| \leq L(x) |V|^{1+\beta} \leq L(x) (\Gamma(x))^{1+\beta} \|V\|^{1+\beta} \leq 
\frac{\|V\|}{100 n_1} \leq \frac{\ep' + \delta}{100 n_1} ,$$
so
\begin{eqnarray*}
\| \So - \Wo\| 
& \leq & \|\So - \Vo\| + \|\Vo - \Uo\| + \| \Uo - \Wo\|\\
& \leq & \frac{\ep' + \delta}{100 n_1} + \delta + \frac{\ep'}{100 n_1} < \frac{\ep'}{30 n_1} .
\end{eqnarray*}

Choose an orthonormal basis $e_1, \ldots, e_{n_1}$ in $E^u_1(x)$ such that $W = \Wo = c_1 e_1$ 
for some $c_1 \in [\ep'/3,\ep']$.
Let $\So = \sum_{i=1}^{n_1} d_i e_i$. Then the above implies 
$$|d_1-c_1| \leq \frac{\ep'}{30 n_1} \quad \mbox{\rm  and} \quad |d_i| \leq \frac{\ep'}{30 n_1} $$
for all $ i = 2,\ldots,n_1$.

Notice that for any $i = 1, \ldots, n_1$, 
$$u' = d\hf^{-p}_x(0) \cdot (\ep' e_i/2) \in \hBo_p(z,\ep') .$$
Indeed, by Lemma 3.1 and (10.6),
$$\|\hf^p_z(u') - d \hf^p_z(0) \cdot u' \| \leq \frac{\|\ep' e_i/2\|}{100n_1} = \frac{\ep'}{200n_1} ,$$ 
so
$$\| \hf^p_z(u')\| \leq \|d \hf^p_z(0) \cdot u' \| + \frac{\ep'}{200n_1}  = \frac{\ep'}{2} + \frac{\ep'}{200 n_1} < \ep' .$$
By the choice of
$u$, this implies $\|u'\| \leq \|u\|$, so $\| d\hf^{-p}_x(0) \cdot e_i\| \leq \frac{2 \|u\|}{\ep'}$ 
for all $i = 1, \ldots, n_1$.

The above yields
\begin{eqnarray*}
\|d_1 \, d\hf^{-p}_z(0) \cdot e_1 \| 
& \geq & \|c_1 d\hf^{-p}_z(0) \cdot e_1 \|  - \|(d_1- c_1) d\hf^{-p}_z(0) \cdot e_1 \|\\
& \geq & \|u\| - \frac{\ep'}{30 n_1} \cdot \frac{2\|u\|}{\ep'}  = \|u\| \left(1 - \frac{1}{15n_1}\right) .
\end{eqnarray*}
Moreover, for $i \geq 2$ we have
$$\|d_i \, d\hf^{-p}_z(0) \cdot e_i \| \leq \frac{\ep'}{30 n_1} 
\cdot \frac{2\|u\|}{\ep'}  = \frac{\|u\|}{15n_1} .$$
Hence
\begin{eqnarray*}
\|\vo\|
& =      & \| d \hf^{-p}_x(0) \cdot \So\| = \left\| \sum_{i=1}^{n_1} d_i \, d\hf^{-p}_z(0) \cdot e_i \right\|\\
& \geq & \|d_1 \, d\hf^{-p}_z(0) \cdot e_1 \|  - \sum_{i=2}^{n_1} \|d_i \, d\hf^{-p}_z(0) \cdot e_i \| \\
& \geq &   \|u\| \left(1 - \frac{1}{15n_1}\right) - n_1 \, \frac{\|u\|}{15n_1} > \frac{\|u\|}{2} .
\end{eqnarray*}
Combining this with Lemma 10.2 gives, 
$$\|v\| \geq |v| \geq \|\vo\| > \frac{\|u\|}{2} 
\geq \frac{\ep'}{4 \tk \ep \Gamma^3(x) }\, \ell( \hB^u_p(z,\ep)) ,$$
which proves (d).
\endofproof

\bs

What we actually need later  is the following immediate consequence of Theorem 10.4 which 
concerns sets of the form
$$B^u_T(z,\ep) = \{ y \in W^u_\ep(z) : d(\phi_T(y), \phi_T(z)) \leq \ep\}\;,$$
where $z\in \ll$, $\ep > 0$ and $T > 0$. 

\bs

\noindent
{\bf Corollary 10.5.} {\it 
There exist a $\frac{6\hep_1}{\beta}$-slowly varying radius function $\hr_1(x) \leq \hr(x)$ ($x\in \ll$) and a global constant $L_1 \geq 1$ such that:}

\ms

(a) {\it For any  $x\in \ll$ and any $0 < \delta \leq  \ep \leq \hr_1(x)$ we have
$$\diam \left( B^{u}_T(z,\ep)\right)  \leq  L_1\,\Gamma^3(x) R(z) 
\frac{\ep}{\delta}  \, \diam \left(B^{u}_T (z,\delta) \right) $$
for any $T > 0$, where $z = \phi_{-T}(x)$.}

\ms

(b) {\it For any  $x\in \ll$, any $0 < \ep \leq \hr_1(x)$, any $\rho \in (0,1)$ and any $\delta$ with
$$0 < \delta \leq \frac{\rho \ep}{L_1 \Gamma^3 (x) R(z)}$$ 
we have
$$\diam \left( B^{u}_T (z ,\delta)\right) \leq  \rho \, \diam \left( B^{u}_T ( z ,\ep)\right)$$
for all $T > 0$, where $z = \phi_{-T}(x)$.}

\ms

(c) {\it For any  $x\in \ll$, any $0 < \ep' < \ep \leq \hr_1(x)$, any $0 < \delta \leq \frac{\ep'}{100 n_1}$ and any
$T > 0$, for $z = \phi_{-T}(x)$ there exists $z' \in B^u_T(z,\ep')$ such that for every $y \in B^u_T(z', \delta)$
we have} 
$$d(z,y) \geq \frac{\ep'}{L_1 \ep \Gamma^3(x) R(z)} \, \diam(B^u_T(z,\ep)) .$$

\def\hq{\hat{q}}

\subsection{Linearization along $E^u_1$}

Here we prove Lemma 10.3 using arguments similar to these in the proofs of Theorem 3.1 and Lemma 3.2 in \cite{St5}.
 
We use the notation from Sect. 10.1. As in Sect. 10.1, $\hr(x)$, $x\in \ll$, will be an $\hep$-slowly varying radius
function satisfying (10.1) and (10.6).

\bs 

\noindent
{\bf Proposition 10.6.} {\it There exist a $\frac{5\hep_1}{\beta}$-slowly varying radius function $\hr_1(x) \leq \hr(x)$ 
and a Lyapunov $4 \hep_1$-regularity function $L_1(x)$, $x\in \ll$, such that:}

\ms

(a) {\it For every $x\in \ll$ and every $u\in E^u_1(x; \hr_1(x))$ there exists
$$F_x(u) = \lim_{p\to\infty} d\hf^p_{f^{-p}(x)}(0)\cdot \hf_x^{-p}(u) \in E^u_1(x; \hr(x)) .$$
Moreover, $\| F_x(u) - u\| \leq L_1(x)\, \|u\|^{1+\beta}$  for any $u\in E^u_1(x,\hr_1(x))$ and any integer $p \geq 0$.}

\ms 

(b) {\it The maps 
$$F_x : E^u_1( x; \hr_1(x)) \longrightarrow F_x (E^u_1( x; \hr_1(x)))  \subset E^u_1(x; \hr(x)) \quad , \quad x\in \ll ,$$
are uniformly Lipschitz. More precisely, 
$$\|F_x(u) - F_x(v) - (u-v)\| \leq C_1 \, \left( \|u-v\|^{1+\beta} + \| v\|^\beta\cdot \|u-v\| \right)$$
for all $x\in \ll$, $u,v\in E^u_1(x; \hr_1(x))$. Moreover,}
$$\frac{1}{2} \|u-v\| \leq \|F_x(u) -F_x(v)\| \leq 2\|u-v\| \quad, 
\quad x\in \ll\:\:, \: u,v\in E^u_1(x; \hr_1(x)) .$$

(c) {\it For any $x \in M$ and any integer $q\geq 1$, setting $x_q = \varphi^{-q}(x)$, we have
$$d\hf_{x_q}^q(0) \circ F_{x_q} (v) = F_x \circ \hf_{x_q}^q (v)$$ 
for any $v \in E^u_1 (x_q; \hr_1(x_q))$ with $\|\hf^q_{x_q}(v)\| \leq \hr_1(x)$.}

\bs

\def\ao{a^{(1)}}
\def\bo{b^{(1)}}

As in \cite{St5} this is derived from the following lemma. Part (b) below is a bit stronger than 
what is required here, however
we need it in this form for the proof of Lemma 4.2 in Sect. 8.

\bs

\noindent
{\bf Lemma 10.7.}  {\it There exist a $\frac{5\hep_1}{\beta}$-slowly varying radius functions $\hr_1(x) \leq \hr(x)$ 
and a Lyapunov $4 \hep_1$-regularity function $L_1(x)$, $x\in \ll$, with the following properties:}

\ms

(a) {\it If $x\in \ll$, $z = \varphi^{p}(x)$ and $\|\hf^p_z(v)\| \leq \hr_1(x)$ for some $v\in E^u_1 (z; \hr_1(z))$ and some
integer $p \geq 1$, then $\|d\hf^p_z(0)\cdot v\| \leq 2 \|\hf^p_z(v)\|$ and
$$\|d\hf_z^p(0)\cdot v - \hf^p_z(v)\| \leq L_1(x)\, \|\hf^p_x(v)\|^{1+\beta} .$$
Similarly, if $\|d\hf^p_z(0)\cdot v\|\leq \hr_1(x)$ for some $v\in E^u_1 (z)$ and some integer $p \geq 1$, then
$$\|\hf^p_x(v) \| \leq 2 \|d\hf^p_x(0)\cdot v\|$$ 
and}
$$\|\hf_x^p(v) - d\hf^p_x(0)\cdot v\| \leq L_1(x)\, \|d\hf^p_x(0)\cdot v\|^{1+\beta} .$$

\ms

(b)  {\it For any $x \in \ll$ and any integer $p \geq 1$, setting $z = \varphi^{-p}(x)$,  the map
$$F^p_x = d \hf^p_z(0) \circ (\hf^p_x)^{-1}:  E^u (x; \hr_1(x)) \longrightarrow E^u (x; \hr(x)) ,$$
is such that
\be
\left\| \left[ (F^p_x (a))^{(1)} - (F^p_x (b))^{(1)} \right] - [\ao - \bo ]\right\|  
\leq L_1(x) \, \left( \|a - b \|^{1+\beta} + \| b \|^\beta \cdot \|a - b\| \right)  
\ee
for all $a,b \in E^u (x ; \hr_1(x))$. Moreover,
\be
\frac{1}{2} \| a - b\| \leq \left\| d\hf^p_z(0)\cdot \left[ (\hf^p_x)^{-1}(a) - (\hf^p_x)^{-1}(b)\right] \right\|\leq 2  \| a - b\| 
\ee
for all $a,b \in E^u_1 (x ;  \hr_1(x))$. }

\ms

\noindent
{\it Proof of Lemma} 10.7. Set 
$$\hr_1(x) = \frac{\hr(x)}{(\Gamma(x)^{2+\beta} D(x))^{1/\beta}} $$
for all $x\in \ll$.

Part (a) follows from Lemma 3.1 (see also the Remark after the lemma).

\ms

(b)  Let $x\in \ll$ and $z = \varphi^{-p}(x)$ for some integer $p \geq 1$.
Given any  $a, b \in E^u (x; \hr_1(x))$, set $v = \hf_z^{-p}(a)$ and $\eta =  \hf_z^{-p}(b)$.
Then  $\|\hf^p_z(v)\| = \|a\| \leq \hr(x)$  and $\|\hf^p_z(\eta)\| = \|b\| \leq \hr(x)$.
Set  
$$z_j = \varphi^j(z) \quad, \quad v_j = \hf^j_z(v) \in E^u(z_j) \quad , \quad w_j = d\hf_z^j(0)\cdot v \in E^u(z_j) ,$$
$$\eta_j = \hf^j_z(\eta) \in E^u(z_j) \quad, \quad \zeta_j = d\hf_z^j(0)\cdot v \in E^u(z_j) .$$
Clearly $v_p = a$ and $\eta_p = b$. We need to prove that
$$\|(\wo_p - \zetao_p) - (\vo_p - \etao_p)\| \leq L_1(x) \, 
\left( \|v_p - \eta_p \|^{1+\beta} + \| \eta_p \|^\beta \cdot \|v_p - \eta_p\| \right) .$$  

Next, we use the argument from the proof of Lemma 3.1 (which is Lemma 3.3 
in \cite{St3}) with a small modification.

By (3.12),
$$\| v_k\| \leq \|\hf^k_z(v)\|'_{z_k} \leq \frac{1}{\mu_1^{p-k}}\, \|\hf^p_z(v)\|'_x \quad , \quad
\|\eta_k\| \leq \|\hf^k_z(\eta)\|'_{z_k}  \leq \frac{1}{\mu_1^{p-k}}\, \|\hf^p_z(\eta)\|'_x , $$
and also
\be
\|v_k - \eta_k\| \leq \|v_k - \eta_k\|'_{z_k}  \leq \frac{1}{\mu_1^{p-k}}\, \|v_p - \eta_p\|'_x
\ee
for all  $k = 0,1, \ldots,p$.

It follows from (3.16) that
$$\|\hf_z(v) - \hf_z(\eta) - d\hf_z(0)\cdot (v-\eta) \| \leq D(z) \, \left( \|v-\eta\|^{1+\beta} 
+ \|\eta\|^\beta\, \|v-\eta\| \right) ,$$
so 
$v_1 - \eta_1 = d\hf_z(0) \cdot (v-\eta)  + u_1$
for some $u_1 \in E^u(z_1)$ with 
$$\|u_1\|\leq D(z)\, [\|v-\eta\|^{1+\beta} + \|\eta\|^\beta \, \|v-\eta\| ]  .$$
Hence
$$d\hf_{z_1}(0)\cdot (v_1-\eta_1) = d\hf^2_{z}(0)\cdot (v-\eta) + d\hf_{z_1}(0)\cdot u_1 .$$
Using (3.16) again, 
$$\|\hf_{z_1}(v_1) - \hf_{z_1}(\eta_1) - d\hf_{z_1}(0)\cdot (v_1-\eta_1)\| \leq 
D(z_1) \, \left(\|v_1 - \eta_1\|^{1+\beta} + \|\eta_1\|^\beta\, \|v_1 - \eta_1\| \right) ,$$
so $v_2 - \eta_2 =  d\hf_{z_1}(0)\cdot (v_1-\eta_1) + u_2$
for some $u_2\in E^u (z_2)$ with 
$$\|u_2\| \leq D(z_1) \, [ \|v_1 - \eta_1\|^{1+\beta} + \|\eta_1\|^\beta\, \|v_1 - \eta_1 \| ] ,$$
which gives $v_2  - \eta_2 =( w_2 -\zeta_2) + u_2 + d\hf_{z_1}(0)\cdot u_1 .$

Assume that for some $k = 1, \ldots, p-1$ we have
\be
v_k - \eta_k = (w_k - \zeta_k) + u_k + d\hf_{z_{k-1}}(0)\cdot u_{k-1}  + \ldots + d\hf_z^{k-1}(0)\cdot u_1 ,
\ee
where  $u_j \in E^u (z_j)$ and 
$$\|u_j\| \leq D(z_{j-1}) \, \left( \|v_{j-1} - \eta_{j-1}\|^{1+\beta}  + \|\eta_{j-1}\|^\beta\, \|v_{j-1} - \eta_{j-1}\| \right)$$ 
for all $j = 1,\ldots,k-1$. Then 
$$d\hf_{z_k}(0)\cdot (v_k - \eta_k) = (w_{k+1} - \zeta_{k+1}) + d\hf_{z_k}(0)\cdot u_k + 
d\hf^2_{z_{k-1}}(0)\cdot u_{k-1}  + \ldots + d\hf_z^{k}(0)\cdot u_1 .$$
By (3.16), 
$$\|\hf_{z_k}(v_k) - \hf_{z_k}(\eta_k) - d\hf_{z_k}(0)\cdot (v_k- \eta_k)\| 
\leq D(z_k) \, \left[ \|v_k - \eta_k\|^{1+\beta} + \|\eta_k\|^\beta\, \|v_k - \eta_k\| \right] ,$$
therefore 
$$v_{k+1} - \eta_{k+1} = d\hf_{z_k}(0)\cdot (v_k- \eta_k) + u_{k+1}$$
for some $u_{k+1} \in E^u (z_{k+1})$ with 
$$\|u_{k+1}\| \leq D(z_k) \, \left( \|v_k - \eta_k\|^{1+\beta} + \|\eta_k\|^\beta\, \|v_k - \eta_k\| \right) .$$ 
Thus, (10.14)  holds for all $k = 1 , \ldots,p$.

It follows from (10.14) with $k=p$ that
$$\vo_p - \etao_p = (\wo_p - \zetao_p) + \uo_p + d\hf_{z_{p-1}}(0)\cdot \uo_{p-1}  
+ \ldots + d\hf_z^{p-1}(0)\cdot \uo_1 .$$
Next, (10.13) implies 
\begin{eqnarray*}
\|u_j\|'_{z_j} 
& \leq & \Gamma(z_j) \|u_j\|\leq \Gamma (z_j) D(z_{j-1})\, \left[ \|v_{j-1} - \eta_{j-1} \|^{1+\beta}  + \|\eta_{j-1}\|^\beta\, \|v_{j-1} - \eta_{j-1}\| \right]\\
& \leq & \frac{\Gamma(z_j) D(z_{j-1})}{\mu_1^{(1+\beta)(p-j+1)}}\, \left[ (\|v_p-\eta_p\|'_x)^{1+\beta}  + (\|\eta_p\|'_x)^\beta\, \|v_p - \eta_p\|'_x \right] \;,
\end{eqnarray*}
so by (3.14),
\begin{eqnarray*}
&        & \|d\hf_{z_j}^{p-j}(0)\cdot \uo_j\|\\
& \leq &  \|d\hf_{z_j}^{p-j}(0)\cdot \uo_j\|'_x \leq (\nu_1)^{p-j}\, \|\uo_j\|'_{z_j}\\
& \leq & \frac{\Gamma(z_j) D(z_{j-1})}{\mu_1^{1+\beta}}\, \left(\nu_1/ \mu_1^{1+\beta}\right)^{p-j} \, 
\left[(\|v_p-\eta_p\|'_x)^{1+\beta} + (\|\eta_p\|'_x)^\beta\, \|v_p - \eta_p\|'_x \right]\\
& =     & \Gamma(z_j) D(z_{j-1}) \Gamma(x)^{1+\beta} \,  \left(\nu_1/ \mu_1^{1+\beta}\right)^{p-j} \,
 \left[\|v_p-\eta_p\|^{1+\beta} + \|\eta_p\|^\beta\, \|v_p - \eta_p\| \right]\\
& \leq & \Gamma(x)^{2+\beta} D(x) \, 
\left(e^{2\hep}\nu_1/ \mu_1^{1+\beta}\right)^{p-j} \, \left[\|v_p-\eta_p\|^{1+\beta}  + \|\eta_p\|^\beta\, \|v_p - \eta_p\| \right]\\
\end{eqnarray*}
for all $j = 1, \ldots, p$.
Using the above, setting 
$$L_1(x) = \Gamma(x)^{2+\beta} D(x)$$ 
and assuming $\hep$ is chosen sufficiently small so that
$$e^{2\hep}\nu_1/\mu_1^{1+\beta}  < 1 ,$$
we get
$$\|(\vo_p - \etao_p) - (\wo_p - \zetao_p) \|    \leq L_1(x) \, \left[\|v_p-\eta_p\|^{1+\beta} + \|\eta_p\|^\beta\, \|v_p - \eta_p\| \right] ,$$
which proves (10.11).

Notice that when $a,b \in E^u_1(x)$ we have $v_p, \eta_p, w_p, \zeta_p \in E^u_1(x)$ as well, 
so (10.11) and the choice of $\hr(x)$ and $L_1(x)$ imply
$$\|(v_p - \eta_p) - (w_p - \zeta_p) \|   \leq \|v_p - \eta_p \|\, 3 L_1(x) (\hr(x))^\beta  \leq \frac{1}{2}\, \|v_p - \eta_p\| ,$$
which proves (10.13).
\endofproof

\bs

\noindent
{\it Proof of Proposition 10.6.} Let $L_1(x)$ and $\hr_1(x)$ ($x\in \ll$) be as in Lemma 10.7. 
Fix  an arbitrary $x\in \ll$ and set $x_j = \varphi^{-j}(x)$ for any integer $p \geq 0$.
In what follows we use the maps $F^p_y$ ($y\in \ll$, $p \geq 1$) from  Lemma 10.7, and also
the notation from the proof of Lemma 10.7.

\ms

(a) Given $u \in E^u_1(x, \hr_1(x))$ and $p \geq 0$, set $u_p = F^p_x(u) \in E^u_1(x; \hr(x))$.  
To show that the sequence  $\{ u_p\}$ is Cauchy, consider any $q > p$ and set 
$$v = \hf_x^{-p}(u) \in E^u_1(x_p, \hr_1(x_p)) .$$
By (3.12)  we have
\be
\|v\|'_{x_p} \leq \frac{\|u\|'_x}{\mu_1^p} .
\ee
Set 
$$v_{q-p} = d \hf_{x_q}^{(q-p)}(0)\cdot (\hf^{-(q-p)}_{x_p}(v)) .$$
From Lemma 10.7 we know that   $\|v_{q-p} - v\| \leq L_1(x_p)\, \|v\|^{1+\beta} , $ i.e.
$$\|d \hf_{x_q}^{(q-p)}(0)\cdot (\hf^{-(q-p)}_{x_p}(v)) - v\| \leq L_1(x_p)\, \|v\|^{1+\beta} .$$
Applying $d\hf_{x_p}^p(0)$ to the above and using (3.14), we get
\begin{eqnarray*}
\|u_q - u_p\| 
& =    &    \|d\hf^q_{x_q}(0)\cdot (\hf^{-q}_x(u)) - d\hf^p_{x_p}(0)\cdot v\|  \leq \Gamma(x)\, \|d\hf^p_{x_q}(0)\cdot  (v_{q-p} - v) \|'_x\\
& \leq & L_1(x_p) \Gamma(x) (\nu_1)^p \, (\|v\|'_{x_p})^{1+\beta} \\
& \leq & \Gamma(x) L_1(x)  (\nu_1 e^{\hep}/\mu_1^{1+\beta})^p \, (\|u\|'_{x})^{1+\beta}  \leq \Gamma(x)^{2+\beta}\, L_1(x) \gamma_2^p\, \|u\|^{1+\beta} ,
\end{eqnarray*}
where as above, $\gamma_2 = e^{2\hep}\nu_1/\mu_1^{1+\beta}  < 1$.
Thus, $\{ u_p\}$ is Cauchy, so there exists $F_x(u) = \lim_{p\to\infty} u_p$. 
Moreover, letting $q\to \infty$ in the above gives
$$\|F_x(u) - u_p\| \leq  \Gamma(x)^{2+\beta}\, L_1(x) \gamma_2^p\, \|u\|^{1+\beta} 
\quad, \quad u\in E^u(x; \hr_1(x))\:, \: p\geq 0  .$$

\medskip

(b) Given $u, v\in E^u_1(x; \hr_1(x))$ and $p \geq 0$, it follows from Lemma 10.7(c) that
$$\|(u_p - v_p) - (u-v)\| \leq L_1(x)\, [ \|u-v\|^{1+\beta} + \|v\|^\beta \, \|u-v\|] .$$
Letting $p \to \infty$, proves the desired estimate.

\ms

(c) Let $v\in E^u_1(x_1, \hr_1(x_1))$ be such that $\|\hf_{x_1}(v)\| \leq \hr_1(x)$.
It is enough to show that 
\be
d\hf_{x_1}(0) \circ F_{x_1} (v) = F_x \circ \hf_{x_1} (v) .
\ee
Set $u = \hf_{x_1}(v)$. For any integer $p \geq 0$ we have $v_p = d\hf^p_{x_{p+1}}(0) \cdot \hf^{-p}_{x_1}(v)$, so
$$d\hf_{x_1}(0) (v_p) =  d\hf^{p+1}_{x_{p+1}}(0) \cdot \hf^{-p}_{x_1}(v) 
= d\hf^{p+1}_{x_{p+1}}(0) \cdot \hf^{-(p+1)}_{x}(u) = u_{p+1} .$$
Letting $p \to \infty$ gives $d\hf_{x_1}(0) (F_{x_1}(v)) = F_x(u)$, which is exactly (10.16).
\endofproof

\bs

\def\tBo{\widetilde{B}^{u,1}}

\index{$\tB^{u,1}_p(z,\ep)$}

In what follows we assume that $\tr_1(x)$ ($x\in \ll$) is a {\bf fixed $\frac{5\hep_1}{\beta}$-slowly varying radius function}
with the properties described in Proposition 10.6.

For $z\in \ll$, $\ep \in (0, \hr_1(z)]$ and an integer $p \geq 0$ set
$$\tB^{u,1}_p(z,\ep) = F_z(\hB^{u,1}_p(z,\ep)) \subset E^u_1(z; \hr(z)) .$$
Then, using  Proposition 10.6(c) we get
\be
d\hf^{-1}_x(0) (\tBo_{p+1}(x, \delta )) \subset \tBo_p({\varphi^{-1}(x)}, \delta)  \quad , \quad x\in \ll\;, \; p \geq 1 .
\ee
Indeed, if $\eta \in \tBo_{p+1}(x, \delta)$, then $\eta = F_x(v)$ for some 
$v \in \hBo_{p+1}(x,\delta)$, and then clearly
$w = \hf^{-1}_x(v) \in \hBo_p(x, \delta)$. Setting $y = \varphi^{-1}(x)$, by Proposition 10.6(c), 
$$\eta = F_x(v) = F_x (\hf_y(w)) = d\hf_y(0) \cdot (F_y(w)) ,$$ 
so $d\hf_x^{-1}(0)\cdot \eta  = F_y(w) \in \tBo_p(y,\delta)$.
Moreover, locally near $0$ we have an equality in (10.17), i.e. if $\delta' \in (0,\delta)$ is sufficiently small, then 
$$d\hf^{-1}_x(0) ( \tBo_{p+1}(x,\delta)) \supset \tBo_p(\varphi^{-1}(x),\delta') .$$

To prove part (a) of Lemma 10.3 we have to establish the following lemma which is similar to Lemma 4.4 
in \cite{St5} (see also the Appendix in \cite{St5}).

\bs

\noindent
{\bf Lemma 10.8.} {\it Let $x\in \ll$ and let  $0 < \delta \leq  \ep \leq \hr_1(x)$. Then 
$$\ell \left( \tB^{u,1}_p( \varphi^{-p}(x),\ep)  \right) 
\leq  8 n_1  \frac{\ep}{\delta}   \, \ell \left(\tB^{u,1}_p ( \varphi^{-p}(x) , \delta) \right) $$
for any  integer $p \geq 0$, where $n_1 = \dim (E^u_1(x))$.}

\ms

\noindent
{\it Proof of Lemma 10.8.} 
Choose an orthonormal basis $e_1 ,e_2 , \ldots,e_{n_1}$ in $E^{u}_1 (x)$ and set $u_i = \frac{\delta}{4} e_i$.

Consider an arbitrary integer  $p \geq 1$ and set $z = \varphi^{-p}(x)$. Given
$v \in  \tBo_p( z , \ep)$, we have $v = F_z(w)$ for some $w \in \hBo_p(z,\ep)$. 
Then $\|\hf^p_z(w)\| \leq \ep$. 
Now it follows from Proposition 10.6 that 
$$\|d\hf_z^p(0)\cdot v\| = \|d\hf^p_z(0)\cdot F_z(w)\| = \| F_x(\hf_z^p(w))\|\leq 2\|\hf^p_z(w)\| \leq 2\ep .$$
So, $u = d\hf_z^p(0)\cdot v \in  E^{u}_1(x, 2\ep)$.
We have $u = \sum_{s=1}^{n_1} c_s\, u_s$ for some real numbers $c_s$ and
$u = \sum_{s=1}^{n_1} \frac{\delta c_s}{4} e_s$,
so $\sqrt{\sum_{s=1}^{n_1} c_s^2} = \frac{4}{\delta} \|u\|$ 
and therefore   $|c_s| \leq 8 \ep\, \frac{1}{\delta} $ for all $ s = 1, \ldots, n_1$. 

By (10.17), $v_j = d\hf_x^{-p}(0) \cdot u_j \in \tBo_p( z , \delta)$. 
Indeed,  since $\|u_j\| \leq \frac{\delta}{4}$, we have $u_j = F_x(u'_j)$  for some $u'_j \in E^u_1(x, \delta/2)$. 
Set $v'_j = \hf^{-p}_x(u'_j)$; then $\| \hf_z^p(v'_j)\|\leq \frac{\delta}{2}$, so
$v'_j \in \hBo_p(x, \delta/2)$ and therefore $v_j = F_z(v'_j) \in \tBo_p(z, \delta)$. Using Proposition 10.6(b), we get
$$d\hf_x^{-p}(0) \cdot u_j = d\hf_x^{-p}(0) \cdot F_x(u'_j) = F_z(\hf^{-p}_x(u'_j)) = F_z(v'_j) = v_j  .$$

It now follows that
\begin{eqnarray*}
\|v\| 
& =     &  \|d \hf^{-p}_x (0) \cdot u\| =  \left\| \sum_{s=1}^{n_1} c_s\,  d\hf^{-p}_x(0) \cdot u_s\right\|\\
& \leq  & n_1 \, 8 \ep \, \frac{1}{\delta} \, \max_{1\leq s\leq n_1} \|v_s\|   \leq  8 n_1  \frac{\ep}{\delta} \, \ell(\tBo_p(z,\delta)) .
 \end{eqnarray*}
Therefore 
$\di \ell \left(\tBo_p( z ,\ep)  \right)    \leq  8 n_1 \frac{\ep}{\delta} \, \ell \left(\tBo_p ( z , \delta) \right) ,$
which proves the lemma.
\endofproof

\bs

Lemma 10.3(b) is a consequence of the following.

\bs

\noindent
{\bf Lemma 10.9.} {\it Let $x\in \ll$ and let  $0 <  \ep \leq \hr_1(x)$ and $\rho \in (0,1)$. 
Then for any $\delta$ with
$0 < \delta \leq  \frac{\rho \, \ep}{8 n_1} $
 we have
$$\di \ell \left( \tB^{u,1}_p( \varphi^{-p}(x),\delta)  \right)  \leq \rho  \, \ell \left(\tB^{u,1}_p ( \varphi^{-p}(x) , \ep) \right) $$
for any  integer $p \geq 0$.}

\ms

\noindent
{\it Proof of Lemma 10.9.} 
We will essentially repeat the argument in the proof of Lemma 10.8.

Fix $x\in \ll$,  $0 <  \ep \leq \hr_1(x)$ and $\rho \in (0,1)$.  Set  $u_i = \frac{\ep}{4} e_i$ and let
$0 <  \delta  \leq  \frac{\rho \, \ep}{8 n_1}$. Then given an  integer $p \geq 1$, set $z = \varphi^{-p}(x)$. Given
$v \in  \tBo_p( z , \delta)$, as before we have $u = d\hf_z^p(0)\cdot v \in  E^{u}_1(x, 2\delta)$.
Then $u = \sum_{s=1}^{n_1} c_s\, u_s$ for some real numbers $c_s$ and $u = \sum_{s=1}^{n_1} \frac{\ep}{4} e_s$,
so again we get $|c_s| \leq 8 \delta\, \frac{1}{\ep} $ for all $ s = 1, \ldots, n_1$.  As in the proof of Lemma 10.8 we obtain 
$$v_j = d\hf_x^{-p}(0) \cdot u_j \in \tBo_p( z , \ep) ,$$ 
therefore
\begin{eqnarray*}
\|v\| 
 =      \|d \hf^{-p}_x (0) \cdot u\| =  \left\| \sum_{s=1}^{n_1} c_s\,  d\hf^{-p}_x(0) \cdot u_s\right\|
 \leq  n_1 \, 8 \delta \, \frac{1}{\ep} \, \max_{1\leq s\leq n_1} \|v_s\| \leq  \rho \, \ell(\tBo_p(z,\delta)) .
 \end{eqnarray*}
 Thus, $ \di \ell(\tB^u_p(z,\delta)) \leq \rho \, \ell(\tB^u_p(z,\ep))$.
\endofproof

\def\htau{\hat{\tau}}
\def\tcc{\tilde{\cc}}

\subsection{Consequences for cylinders in Markov partitions}

Here we prove Lemma 4.4 using arguments similar to these in Sect. 4 in \cite{St5}.  We use the notation from Sect. 4.

Let $c_0$, $\gamma$ and $\gamma_1$ be the constants from (2.1). Recall that by the choice of $\hep_1$ we have
$\hep = \frac{6\hep_1}{\beta} .$
Assume $\hep > 0$ is chosen small enough so that $e^{\hep}/\gamma < 1$.

Let $\hr_1(x)$ ($x\in \ll$)  be an $\hep$-slowly varying radius function  with the properties described in Theorem 10.4 and Corollary 10.5. 
Then there exists  a constant $\hr_0 > 0$ such that $\hr_1(x) \geq \hr_0$ for all $x\in P_0$ and $x \in \tP_0$. 
{\bf Fix $\hr_0$ with this property.}

Let $B_1 > 0$ be a {\it Lipschitz constant} for the {\it projection along the flow}
$$\psi : \cup_{i=1}^{k_0} \phi_{[-\ep,\ep]}(D_i) \longrightarrow \cup_{i=1}^{k_0} D_i ,$$
i.e. for all $i = 1,\ldots,k_0$ and all $x\in \phi_{[-\ep,\ep]}(D_i) $ we have $\psi(x) = \pr_{D_i}(x)$.  
\index{$B_1$}

 Next,  fix an integer $d_0 \geq 1$ such that 
\be
\frac{2 \tk \Gamma^3_0 e^{2\hep} r_1}{\hr_0} < (\mu_1 e^{\hep})^{d_0} 
\quad , \quad  \frac{1}{c_0 (\gamma e^{-\hep})^{d_0}} < \frac{\hr_0}{2} .
\ee
Set
\be
r'_0 = \frac{\hr_0 e^{- (d_0 + 1) \hep}}{B_1} .
\ee

\ms

\noindent
{\it Proof of Lemma} 4.4. 
First note the following. Let $z\in \tR_j$ be such that $\tpp^{d_0+1}(z) \in \tP_0$. Then
$z\in C_V[\ii']$ for some $\ii' = [i_0, \ldots,i_{d_0+1}]$ with $i_0 = j$, where
$V = W^u_{\tR}(z)$.  Set $\ii = [i_0, \ldots,i_{d_0}]$. We claim that
\be
C_V[\ii] \subset B_V(z, r'_0) \quad \mbox{\rm and } \quad \hr_1(z) \geq r'_0 .
\ee
Indeed, by (2.1) and (10.18), 
$$\diam(C_V[\ii]) \leq \frac{1}{c_0\gamma^{d_0+1}} < r'_0/2 .$$
On the other hand, $\hr_1(x)$ is an  $\hep$-slowly varying radius function and $y = \tpp^{d_0+1}(z) \in \tP_0$ 
and the definition of $\tP_0$ show that $\hr_1(y) \geq \hr_0 $. Also recall that $0 < \tau(x) \leq 1$ for all 
$x\in \tR$ by the choice of the Markov family.  Now using (10.19), we get
$$\hr_1(z) \geq \hr_1(y) e^{-\tau_{d_0+1}(z) \hep} \geq \hr_0 e^{-(d_0+1)\hep } = r_0' >   2\, \diam(C_V[\ii]) .$$
This proves (10.20).

\medskip

(a)  Assume that $m > d_0$, and let $\ii = [i_0,i_1, \ldots,i_m]$ and 
$\ii' = [i_0,i_1, \ldots,i_m, i_{m+1}]$ be admissible sequences. Let  $\cc = C[\ii]$ and $\cc' = C[\ii']$
be the corresponding cylinders in $\tR$. Assume that  there exists $z \in \cc' \cap P_0$
with $\pp^{m+1}(z) \in P_0$.

Fix such a point $z \in \cc'$; then  $y = \tpp^{m+1}(z)  \in \tP_0$ and $\tpp^j(z) \in \tR_{i_j}$ 
for all $j = 0,1,\ldots,m+1$. Set $\tpp^{m-d_0}(z) = x$, $V = W^u_{\tR}(x)$, Since $\tpp^{d_0+1}(x) = y\in \tP_0$, 
we have $\hr_1(y) \geq \hr_0$, so $\hr_1(x) \geq r_0$.

Consider the cylinders
$$\tcc' = C_{V}[i_{m-d_0}, i_{m-d_0+1}, \ldots, i_m, i_{m+1}] \subset 
\tcc = C_{V}[i_{m-d_0}, i_{m-d_0+1}, \ldots, i_m]  \subset V .$$
Since $\tpp^{d_0+1}(x) = y$, using (10.20) we get $\tcc \subset B_{V}(x, r'_0)$.
On the other hand it is easy to see using (2.1) that 
$\tcc' \supset B_{V}(x,c_0  r_1/\gamma_1^{d_0+1})$.

We will now use Corollary 10.5(a) with $x$ and $z$ as above, $T = \tau_{m-d_0}(z) > 0$ and
$$0 < \delta = \delta_3 = \frac{c_0 r_1}{B_1 \gamma_1^{d_0+1}} < \ep = r'_0 .$$
This, combined with (10.20), gives
$$\diam(B^u_T(z,\delta_3)) \geq \frac{\delta_3}{B_1 L_1 \Gamma_0^3 R_0 r'_0}\, \diam(B^u_T(z, r'_0)) .$$
However, using the above information about $\tcc$ and $\tcc'$, as in the proof of 
Proposition 3.3 in \cite{St2},  one easily observes that $\cc' \supset B^u_T(z,\delta_3)$ 
and $\cc \subset B^u_T(z, B r'_0)$. Thus,
$$\diam (\cc') \geq \frac{\delta_3}{B_1 L_1 \Gamma_0^3 R_0 r'_0}\, \diam(\cc) .$$

This proves part (a) for $m > p_0$. Since there are only finitely many cylinders of length $\leq p_0$,
it follows immediately that there exists 
$$\rho_1 \in (0, \frac{\delta_3}{B_1 L_1 \Gamma_0^3 R_0 r'_0}\,]$$
which satisfies the requirements of part (a).

\ms

(b)  Let $\rho'\in (0,1)$. 
It follows from Corollary 10.5(b) that for $z\in \tR \cap \ll$ with $\Phi_T(z) \in \tP_0$ for some $T > 0$  we have
\be
\diam( B^u_T(z, B\delta)) \leq \rho' \, \diam( B^u_T(z, r_1/B)) ,
\ee
provided 
$$ 0 < \delta \leq \delta_4 = \frac{\rho' r_1}{B^2 L_1 \Gamma^3_0 R_0} .$$
Fix an integer $q' \geq 1$ so large that
\be
\frac{1}{c_0 \gamma^{q'}} \leq \delta_4 . 
\ee

Consider the cylinders
$$\cc = C[\ii] = C[i_0, \ldots,i_m] \supset \cc' = C[\ii'] = C[i_0,i_{1}, \ldots, i_{m+1}, \ldots, i_{m+q'}] $$
in some unstable leaf in $\tR$. 
Assume that there exists $z \in \cc' \cap P_0$ with $\pp^{m+q'}(z) \in P_0$. 
Set $T = \ttau_{m+q'}(z)$. Then (10.22) holds.
Set $x = \tpp^{m}(z)$, $V = W^u_{\tR}(x)$; then  $\tpp^{q'}(x) \in \tP_0$, so $\hr_1(x) \geq r_0$. 
Consider also the cylinder 
$$\tcc = C_{V}[i_{m}, i_{m+1}, \ldots, i_m, i_{m+q'}] \subset V .$$
Since $x\in \tcc$ and $\diam (\tcc) \leq \frac{1}{c_0 \gamma^{q'}}$, (10.22) 
implies $\tcc \subset B_{V}(x, \delta_4)$.

Next, we have $C[\ii'] \subset B^u_T(z, B_1 \delta_4)$. Indeed, if $u \in C[\ii']$, then  $\tpp^m(u) \in \tcc$, so 
$$d(\tpp^m(z), \tpp^m(u)) < \frac{1}{c_0 \gamma^{q'}} < \delta_4 .$$ 
Thus, $C[\ii'] \subset B^u_T(z, B_1\delta_4)$ and therefore by (10.22),
$$\diam (C[\ii']) \leq \diam (B^u_T(z, B_1 \delta_4) ) \leq \rho' \, \diam( B^u_T(z, r_1/B_1)) .$$
On the other hand,  $W^u_{r_1}(x) \subset \tR_{i_m}$ (see Sect. 4.2 for the choice of $r_1 > 0$), so
$B^u_T(z, B_1 r_1) \subset C[\ii]$, and therefore $\diam (C[\ii']) \leq \rho'\, \diam(C[\ii])$.

\ms

(c)  Let again $m > d_0$,  let $\ii = [i_0,i_1, \ldots,i_m]$ be an admissible sequence, let $\cc = \cc_W[\ii]$
be the corresponding cylinder in an unstable leaf $W$ in $\tR$.  Let $z\in \cc \cap \tP_0$ and let $\tpp^m(z) = z' \in \tP_0$.
Set 
$$z''= \tpp^{m-d_0}(z) \quad , \quad V = W^u_{\tR}(z'') .$$
If $z' = \phi_T(z)$ and $z'' = \phi_t(z)$; 
then $\phi_{T-t}(z'') = z'$, so $T-t = \ttau_{d_0}(z'') < d_0$.
Thus, 
$$\hr(z'') \geq \hr(z') e^{-d_0 \hep} \geq \hr_0 e^{-d_0\hep} > r'_0 .$$
 As in part (a), for the cylinder
$$\tcc = C_{V}[i_{m-d_0}, i_{m - d_0+1}, \ldots, i_m]$$ 
in $V$, we have 
$$z'' \in B_V(z'', c_0\hr_0/ \gamma_1^{d_0}) \subset \tcc = \tpp^{m-d_0}(\cc) \subset B_V (z'', r'_0) .$$
Setting $\ep' = c_0\hr_0/ \gamma_1^{d_0} < \ep = B_1 r'_0$, it follows from Corollary 10.5(c) that 
for $0 < \delta_5 = \frac{\ep'}{100 n_1}$ there exists $x \in B^u_t(z,\ep')$ such that for every $y \in W^u_\delta (z)$ with
$d(\phi_t(y),\phi_t(x)) \leq  \delta_5$ we have 
\be
d(z,y) \geq \frac{\ep'}{L_1 \ep \Gamma_0^3R_0 } \, \diam(B^u_t (z, B_1 r'_0)) 
\geq \frac{c_0 \hr_0}{L_1 B_1 r'_0\Gamma_0^3R_0 \gamma_1^{d_0} } \, \diam(\cc),
\ee
since $\cc \subset B^u_t(z, B_1 r'_0)$.

Take the integer $q_0 \geq 1$ so large that 
$$\frac{1}{c_0 \gamma^{d_0+q_0}} < \frac{\delta_5}{B_1} = \frac{\ep'}{100 B_1 n_1} ,$$ 
where $n_1 = \dim(E^u_1)$.

Let  $x\in \cc$ and let
$$\cc' = C[\ii'] = C[i_0,i_{1}, \ldots, i_{m+1}, \ldots, i_{m+q_0}]$$ 
be the sub-cylinder of $\cc$ of co-length $q_0$ containing $x$.
Then for the cylinder
$$\tcc' = C_{V}[i_{m-d_0}, i_{m-d_0+1}, \ldots, i_{m}, i_{m+1}, \ldots, i_m, i_{m+q_0}] \subset V $$
we have $\tpp^{m-d_0}(x) \in \tcc'$ and 
$$\diam(\tcc') < \frac{1}{c_0 \gamma^{d_0+q_0}} < \frac{\delta_5}{B_1} .$$
Since for any $y \in \cc'$ we have $\tpp^{m-d_0} (y) \in \tcc'$, it follows that 
$$d(\tpp^{m-d_0}(x), \tpp^{m-d_0}(y)) < \delta_5/B_1 ,$$
and therefore $d(\phi_t(x), \phi_t(y)) < \delta_5$. Thus, $y$ satisfies (10.23). This proves the assertion with
$\rho_1 =  \frac{c_0\hr_0}{L_1 B_1 r'_0 \gamma_1^{d_0} \Gamma_0^3 R_0} $.

\ms

(d) This follows from Theorem 10.4(a).
\endofproof

\newpage

\section{Appendix: Proofs of some technical lemmas}
\setcounter{equation}{0}

In this chapter we prove Lemmas 5.2 and 6.10.

\bs

\noindent
{\it Proof of Lemma} 5.2. Let  $u, u' \in \hU_i $ for some $i = 1, \ldots,k_0$ 
and  let $m \geq 1$ be an integer.  Given $v \in \hU$ with $\sigma^m(v) = u$, let $C[\ii] = C[i_0, \ldots,i_m]$
be the  cylinder of length $m$ containing $v$.  Set 
$$\hC[\ii] = C[\ii]\cap \hU .$$
Since the sequence $\ii = [i_0, \ldots,i_m]$ is 
admissible, the Markov property implies $i_m = i$ and $\sigma^m(\hC[\ii]) = \hU_i$. Moreover,  
$\sigma^m : \hC[\ii] \longrightarrow \hU_i $ 
is a homeomorphism, so there  exists a unique $v' = v'(v)\in \hC[\ii]$ such that $\sigma^m(v') = u'$.  Clearly,
$$\dte (\sigma^j(v),\sigma^j(v'(v))) = \theta^{m-j}\, \dte (u,u') \quad, \quad j = 0,1, \ldots, m-1 .$$
Consequently, using (5.3),
\begin{eqnarray*}
|\fa_m(v) - \fa_m(v')| 
& \leq & \sum_{j=0}^{m-1} |\fa(\sigma^j(v)) - \fa (\sigma^j(v'))|\\
& \leq &
\sum_{j=0}^{m-1} \Lip_\theta(\fa) \,\theta^{m-j}\, \dte (u,u')\nonumber\\
& \leq & \frac{ \theta \, T_0}{1-\theta}\, \dte (u,u')\;. \nonumber
\end{eqnarray*}

Also notice that  $\dte(v,v'(v)) \leq \theta^m \dte (u,u')$.

Recall the definition of the Ruelle operator $\ma = L_{\fa}$ from Sect. 5.1.

We have
\begin{eqnarray*}
&         & |\lab^m h(u) - \lab^m h(u')|
 =       \di \left| \sum_{\sigma^m v = u} \left( e^{\fa_m(v) - \i b \tau_m(v)}\, h(v) -
  e^{\fa_m(v') - \i b \tau_m(v')}\, h(v')\right)\right|\\
& \leq & \di \left| \sum_{\sigma^m v = u} e^{\fa_m(v) - \i b \tau_m(v)}\, [h(v) - h(v')]\right|\\
&       &          \di  +  \left| \sum_{\sigma^m v = u} \left( e^{\fa_m(v) - \i b \tau_m(v)} -
  e^{\fa_m(v') - \i b \tau_m(v')}\right)\, h(v')\right|\\
& \leq & \di \sum_{\sigma^m v = u} e^{\fa_m(v) }\, | h(v) - h(v')|\\
&        &  +  
\sum_{\sigma^m v = u} \left| e^{\fa_m(v) - \i b \tau_m(v) -
  \fa_m(v') + \i b \tau_m(v')} - 1\right|\, e^{\fa_m(v') }\, |h(v')| .
\end{eqnarray*}
From the assumptions in the lemma, $|h(v) - h(v')| \leq B\, H(v')\, \dte (v,v')$. Also,
using elementary inequalities one checks that $|e^{x+\i y} - 1| \leq e^{|x|} (|x|+ |y|)$  for real  $x$ and $y$. 
By Lemma 5.1 and (5.3), 
$$|\tau_m(v) - \tau_m(v')| \leq T_0\, \dte(u,u') .$$
Assuming 
$$A_0 \geq e^{ \frac{\theta T_0}{1-\theta}}\,\max\left\{  1 , \frac{2\,\theta\, T_0}{1-\theta}\right\}$$ 
and  $|b| \geq 1$, the above and (5.3) give
$$\left| e^{\fa_m(v) - \i b \tau_m(v) -\fa_m(v') + \i b \tau_m(v')} - 1\right| \leq \frac{2\,\theta\, T_0\, e^{ \frac{\theta T_0}{1-\theta}}}{1-\theta}\,  
|b|  \, \dte(u,u') \leq A_0\,|b|  \, \dte (u,u') .$$

Next,
\begin{eqnarray*}
\sum_{\sigma^m v = u} e^{\fa_m(v) }\, | h(v) - h(v')|
& \leq & \sum_{\sigma^m v = u} e^{\fa_m(v) - \fa_m(v')}\, e^{\fa_m(v')}\, B\, H(v')\, \dte (v,v')\\
& \leq &  B \theta^m \, e^{\frac{\theta T_0}{1-\theta }}\,\dte (u,u')\, (\ma^mH)(u') .
\end{eqnarray*}
Therefore
\begin{eqnarray*}
|\lab^m h(u) - \lab^m h(u')|
& \leq & \di A_0\left[ B\,\theta^m \, (\ma^m H)(u') + |b|\, (\ma^m |h| )(u')\right]\, \dte (u,u') .
\end{eqnarray*}
This proves the assertion.
\endofproof

\bs

\noindent
{\it Proof of Lemma} 6.10. (a) Let $u, u'\in \dd'_j$ for some $j \leq j_0$. 
Let $\dd_j \subset \cc_m$, $m \leq m_0$.
Then for  $v = \vli(u)\in  \xijl$ and $v' = \vli(u')\in \xijl$, 
we have $\ell(v,v') \geq N$  and $\sigma^{N}v, \sigma^N v' \in \cc_m$.
This and $H \in \kk_E$ imply
$$|\ln H(v) - \ln H(v')| \leq \frac{|H(v) - H(v')|}{\min\{H(v), H(v')\}} \leq 
E\, \dd(v,v')  = E\,\theta^q \dd(u,u') \leq E \theta^q  < \ln 2\; ,$$
Hence $|\ln H(x') - \ln H(x'')|\leq \ln 2$, so $\frac{1}{2} \leq \frac{H(\vl_i(u'))}{H(\vl_i(x''))} \leq 2$.

\ms

(b) 
Consider the case when for some $v\in \xijl$ we have $|h(v)|\geq \frac{3}{4}H(v)$.
Fix $v$ with this property and consider an arbitrary $v'\in \xijl$. It follows from (ii) in Sect. 6.4 that
$$ |h(v') -h(v)| \leq E|b|\,\theta_2^{p+N} H(v)\, \diam (\Psi(\cc_m)) \leq E|b|\,\theta_2^N H(v)\, \frac{\ep_1}{|b|}
= E \ep_1 \theta_2^N H(v) .$$
Using $2H(v) \geq H(v')$ which follows from (a), one obtains
\begin{eqnarray*}
|h(v')| 
& \geq & |h(v)| - E \ep_1 \theta_2^N H(v) \geq (3/4  - E \ep_1 \theta_2^N) H(v) \geq \frac{1}{4} \, H(v'),
\end{eqnarray*}
since $E\epsilon_1 \theta_2^N \leq 1/4$ by (5.5). 
Thus, in this case the second alternative in (b) holds for all $v \in \xijl$.

In the same way one shows that if $|h(v)| \leq \frac{1}{4}H(v)$ for some
$v\in \xijl$, then the first alternative in (b) holds. \endofproof

\printindex


\end{document}